\documentclass[12pt]{article}
\usepackage{latexsym}
\usepackage{amsfonts}
\usepackage{amssymb}
\usepackage{amscd}
\usepackage{amsbsy}
\usepackage{amsmath}
\usepackage[dvips]{graphicx}
\usepackage[all]{xy}
\usepackage{epsfig}    
\input epsf.tex
\parskip=\medskipamount                 
\thicklines                     
\newcommand{\bimn}[7]{\bibitem{#1}#2,
{\em #3},
{ #4}\hspace{0.25em}{\bf
#5}\hspace{0.25em}(#6)\hspace{0.25em}{#7}.}

\def\inbar{\vrule height1.5ex width.4pt depth0pt}
\def\IC{\relax\,\hbox{$\inbar\kern-.3em{\rm C}$}}
\def\IN{\relax{\rm I\kern-.18em N}}
\def\IQ{\relax\,\hbox{$\inbar\kern-.3em{\rm Q}$}}
\def\IR{\relax{\rm I\kern-.18em R}}
\def\ZZ{\relax{\sf Z\kern-.4em Z}}
\def\a{\alpha} \def\b{\beta}   \def\e{\epsilon} 
 \def\l{\lambda} 
   
\def\cA{{\cal A}} \def\cB{{\cal B}} \def\cC{{\cal C}} \def\cD{{\cal D}}
\def\cE{{\cal E}}
  \def\cH{{\cal H}} 
 \def\cK{{\cal K}} \def\cL{{\cal L}} 
 \def\cO{{\cal O}}  \def\cQ{{\cal Q}}
  \def\cT{{\cal T}} 

\newtheorem{theorem}{Theorem}[section]

\newtheorem{corollary}[theorem]{Corollary}
\newtheorem{conjecture}[theorem]{Conjecture}
\newtheorem{lemma}[theorem]{Lemma}
\newtheorem{definition}[theorem]{Definition}
\newtheorem{remark}[theorem]{Remark}
\newtheorem{example}[theorem]{Example}

\newtheorem{crl}{Corollary of Conjecture~4.1}

\catcode`\@=11

\newif\if@fewtab\@fewtabtrue

\catcode`\@=11

\newif\if@fewtab\@fewtabtrue

{\count255=\time\divide\count255 by 60
\xdef\hourmin{\number\count255}
\multiply\count255 by-60\advance\count255 by\time
\xdef\hourmin{\hourmin:\ifnum\count255<10 0\fi\the\count255}}
\def\ps@draft{\let\@mkboth\@gobbletwo
    \def\@oddhead{}
    \def\@oddfoot
      {\hbox to 7 cm{\footnotesize {\em Draft of \jobname:} \draftdate
       \hfil}\hskip -7cm\hfil\rm\thepage \hfil}
    \def\@evenhead{}\let\@evenfoot\@oddfoot}

\def\ceqno{\global\@fewtabfalse
    \ifcase\@eqcnt \def\@tempa{& & &}\or \def\@tempa{& &}
      \or \def\@tempa{&}
      \or\def\@tempa{}\fi\@tempa
{\rm(\theequation)}}

\def\aeqno#1{\global\@fewtabfalse
    \ifcase\@eqcnt \def\@tempa{& & &}\or \def\@tempa{& &}
      \or \def\@tempa{&}
      \or\def\@tempa{}\fi\@tempa
{\rm(\theequation,#1)}}

\def\label#1{\ifnum\draftcontrol=1
 \global\def\draftnote{$\scriptstyle #1$}\fi
 \@bsphack\if@filesw {\let\thepage\relax
   \def\protect{\noexpand\noexpand\noexpand}%
\xdef\@gtempa{\write\@auxout{\string
      \newlabel{#1}{{\@currentlabel}{\thepage}}}}}\@gtempa
   \if@nobreak \ifvmode\nobreak\fi\fi\fi
  \@esphack}

\def\alabel#1#2{\label{#1}\global\@fewtabfalse
    \ifcase\@eqcnt \def\@tempa{& & &}\or \def\@tempa{& &}
      \or \def\@tempa{&}
      \or\def\@tempa{}\fi\@tempa
{\hbox to 3cm{\phantom{\rm(\theequation,#2)}
\draftnote \hfil}\hskip -3cm {\rm(\theequation,#2)}}}

\def\clabel#1{\label{#1}\global\@fewtabfalse
    \ifcase\@eqcnt \def\@tempa{& & &}\or \def\@tempa{& &}
      \or \def\@tempa{&}
      \or\def\@tempa{}\fi\@tempa
{\hbox to 3cm{\phantom{\rm(\theequation)}
\draftnote \hfil}\hskip -3cm{\rm(\theequation)}}}

\def\eqnarray{\def\draftnote{{}}\global\@fewtabtrue
\stepcounter{equation}\let\@currentlabel=\theequation
\global\@eqnswtrue
\global\@eqcnt\z@\tabskip\@centering\let\\=\@eqncr
$$\halign to \displaywidth\bgroup\@eqnsel\hskip\@centering\@eqcnt\z@
  $\displaystyle\tabskip\z@{##}$&\global\@eqcnt\@ne
  \hskip 1\arraycolsep \hfil$\displaystyle{##}$\hfil
  &\global\@eqcnt\tw@ \hskip 1\arraycolsep
$\displaystyle\tabskip\z@{##}$
\hfil  \tabskip\@centering&\global\@eqcnt\thr@@\llap{##}\tabskip\z@
\cr}

\def\endeqnarray{\@@eqncr\egroup
      \global\advance\c@equation\m@ne$$\global\@ignoretrue}

\def\@eqnnum{\hbox to 3cm{\phantom{\rm(\theequation)} \draftnote
                         \hfil}\hskip -3cm {\rm(\theequation)}}

\def\@@eqncr{\let\@tempa\relax
    \ifcase\@eqcnt \def\@tempa{& & &}\or \def\@tempa{& &}
      \or \def\@tempa{&}
      \or\def\@tempa{}
\fi\@tempa
\if@eqnsw
\if@fewtab\@eqnnum\fi
\stepcounter{equation}\fi\global
\@eqnswtrue\global\@eqcnt\z@\global\@fewtabtrue\cr}

\def\draftcite#1{\ifnum\draftcontrol=1#1\else{}\fi}

\def\@lbibitem[#1]#2{\item{}\hskip -3cm \hbox to 2cm
{\hfil$\scriptstyle\draftcite{#2}$}\hskip
1cm[\@biblabel{#1}]\if@filesw
     {\def\protect##1{\string ##1\space}\immediate
      \write\@auxout{\string\bibcite{#2}{#1}}}\fi\ignorespaces}

\def\@bibitem#1{\item\hskip -3cm \hbox to 2cm
{\hfil $\scriptstyle\draftcite{#1}$}\hskip 1cm
\if@filesw \immediate\write\@auxout
       {\string\bibcite{#1}{\the\value{\@listctr}}}\fi\ignorespaces}

\def\nsection#1{\section{#1}\setcounter{equation}{0}}

\def\draftdate{\number\month/\number\day/\number\year\ \ \ \hourmin }

\global\def\draftcontrol{0}
\catcode`\@=12

\def\theequation{{\thesection.\arabic{equation}}}

\def\qq{\begin{eqnarray}}
\def\qqq{\end{eqnarray}}
\def\rx#1{~(\ref{#1})}
\def\rxw#1{(\ref{#1})}
\def\ex#1{eq.\hspace*{-3pt}\rx{#1}}
\def\eex#1{eqs.\hspace*{-3pt}\rx{#1}}
\def\cx#1{~\cite{#1}}
\def\rw#1{~\ref{#1}}

\def\fg#1{Fig.~\ref{#1}}
\def\tb#1{Table~\ref{#1}}

\newlength{\shiftwidth}
\def\shift#1{&&\hbox to \shiftwidth{\hfill $\displaystyle#1$}}
\newlength{\sshiftwidth}
\def\sshift#1{\lefteqn{\hbox to
\sshiftwidth{\hfill$\displaystyle#1$}}}

\def\qbezier{\bezier{120}}
\def\DottedCircle{
\bezier{4}(0.966,-0.259)(1.04,0)(0.966,0.259)
\bezier{4}(0.966,0.259)(0.897,0.518)(0.707,0.707)
\bezier{4}(0.707,0.707)(0.518,0.897)(0.259,0.966)
\bezier{4}(0.259,0.966)(0,1.04)(-0.259,0.966)
\bezier{4}(-0.259,0.966)(-0.518,0.897)(-0.707,0.707)
\bezier{4}(-0.707,0.707)(-0.897,0.518)(-0.966,0.259)
\bezier{4}(-0.966,0.259)(-1.04,0)(-0.966,-0.259)
\bezier{4}(-0.966,-0.259)(-0.897,-0.518)(-0.707,-0.707)
\bezier{4}(-0.707,-0.707)(-0.518,-0.897)(-0.259,-0.966)
\bezier{4}(-0.259,-0.966)(0,-1.04)(0.259,-0.966)
\bezier{4}(0.259,-0.966)(0.518,-0.897)(0.707,-0.707)
\bezier{4}(0.707,-0.707)(0.897,-0.518)(0.966,-0.259)
}
\def\Endpoint[#1]{
\ifcase#1
\put(1,0){\circle*{0.15}}
\or\put(0.866,0.5){\circle*{0.15}}
\or\put(0.5,0.866){\circle*{0.15}}
\or\put(0,1){\circle*{0.15}}
\or\put(-0.5,0.866){\circle*{0.15}}
\or\put(-0.866,0.5){\circle*{0.15}}
\or\put(-1,0){\circle*{0.15}}
\or\put(-0.866,-0.5){\circle*{0.15}}
\or\put(-0.5,-0.866){\circle*{0.15}}
\or\put(0,-1){\circle*{0.15}}
\or\put(0.5,-0.866){\circle*{0.15}}
\or\put(0.866,-0.5){\circle*{0.15}}
\fi}
\def\Arc[#1]{
\thicklines         
\ifcase#1
\bezier{25}(0.966,-0.259)(1.04,0)(0.966,0.259)
\or
\bezier{25}(0.966,0.259)(0.897,0.518)(0.707,0.707)
\or
\bezier{25}(0.707,0.707)(0.518,0.897)(0.259,0.966)
\or
\bezier{25}(0.259,0.966)(0,1.04)(-0.259,0.966)
\or
\bezier{25}(-0.259,0.966)(-0.518,0.897)(-0.707,0.707)
\or
\bezier{25}(-0.707,0.707)(-0.897,0.518)(-0.966,0.259)
\or
\bezier{25}(-0.966,0.259)(-1.04,0)(-0.966,-0.259)
\or
\bezier{25}(-0.966,-0.259)(-0.897,-0.518)(-0.707,-0.707)
\or
\bezier{25}(-0.707,-0.707)(-0.518,-0.897)(-0.259,-0.966)
\or
\bezier{25}(-0.259,-0.966)(0,-1.04)(0.259,-0.966)
\or
\bezier{25}(0.259,-0.966)(0.518,-0.897)(0.707,-0.707)
\or
\bezier{25}(0.707,-0.707)(0.897,-0.518)(0.966,-0.259)
\fi}
\def\DottedArc[#1]{
\ifcase#1
\bezier{4}(0.966,-0.259)(1.04,0)(0.966,0.259)
\or
\bezier{4}(0.966,0.259)(0.897,0.518)(0.707,0.707)
\or
\bezier{4}(0.707,0.707)(0.518,0.897)(0.259,0.966)
\or
\bezier{4}(0.259,0.966)(0,1.04)(-0.259,0.966)
\or
\bezier{4}(-0.259,0.966)(-0.518,0.897)(-0.707,0.707)
\or
\bezier{4}(-0.707,0.707)(-0.897,0.518)(-0.966,0.259)
\or
\bezier{4}(-0.966,0.259)(-1.04,0)(-0.966,-0.259)
\or
\bezier{4}(-0.966,-0.259)(-0.897,-0.518)(-0.707,-0.707)
\or
\bezier{4}(-0.707,-0.707)(-0.518,-0.897)(-0.259,-0.966)
\or
\bezier{4}(-0.259,-0.966)(0,-1.04)(0.259,-0.966)
\or
\bezier{4}(0.259,-0.966)(0.518,-0.897)(0.707,-0.707)
\or
\bezier{4}(0.707,-0.707)(0.897,-0.518)(0.966,-0.259)
\fi}
\def\Chord[#1,#2]{
\thinlines
\ifnum#1>#2\Chord[#2,#1]
\else\ifnum#1<#2
\ifcase#1
\ifcase#2
\or\qbezier(1,0)(0.516,0.138)(0.866,0.5)
\or\qbezier(1,0)(0.45,0.26)(0.5,0.866)
\or\qbezier(1,0)(0.327,0.327)(0,1)
\or\qbezier(1,0)(0.179,0.311)(-0.5,0.866)
\or\qbezier(1,0)(0.0536,0.2)(-0.866,0.5)
\or\put(1, 0){\line(-2, 0){2}}
\or\qbezier(1,0)(0.0536,-0.2)(-0.866,-0.5)
\or\qbezier(1,0)(0.179,-0.311)(-0.5,-0.866)
\or\qbezier(1,0)(0.327,-0.327)(0,-1)
\or\qbezier(1,0)(0.45,-0.26)(0.5,-0.866)
\or\qbezier(1,0)(0.516,-0.138)(0.866,-0.5)
\fi
\or\ifcase#2\or
\or\qbezier(0.866,0.5)(0.378,0.378)(0.5,0.866)
\or\qbezier(0.866,0.5)(0.26,0.45)(0,1)
\or\qbezier(0.866,0.5)(0.12,0.446)(-0.5,0.866)
\or\qbezier(0.866,0.5)(0,0.359)(-0.866,0.5)
\or\qbezier(0.866,0.5)(-0.0536,0.2)(-1,0)
\or\put(0.866, 0.5){\line(-5, -3){1.73}}
\or\qbezier(0.866,0.5)(0.146,-0.146)(-0.5,-0.866)
\or\qbezier(0.866,0.5)(0.311,-0.179)(0,-1)
\or\qbezier(0.866,0.5)(0.446,-0.12)(0.5,-0.866)
\or\qbezier(0.866,0.5)(0.52,0)(0.866,-0.5)
\fi
\or\ifcase#2\or\or
\or\qbezier(0.5,0.866)(0.138,0.516)(0,1)
\or\qbezier(0.5,0.866)(0,0.52)(-0.5,0.866)
\or\qbezier(0.5,0.866)(-0.12,0.446)(-0.866,0.5)
\or\qbezier(0.5,0.866)(-0.179,0.311)(-1,0)
\or\qbezier(0.5,0.866)(-0.146,0.146)(-0.866,-0.5)
\or\put(0.5, 0.866){\line(-3, -5){1}}
\or\qbezier(0.5,0.866)(0.2,-0.0536)(0,-1)
\or\qbezier(0.5,0.866)(0.359,0)(0.5,-0.866)
\or\qbezier(0.5,0.866)(0.446,0.12)(0.866,-0.5)
\fi
\or\ifcase#2\or\or\or
\or\qbezier(0,1.)(-0.138,0.516)(-0.5,0.866)
\or\qbezier(0,1.)(-0.26,0.45)(-0.866,0.5)
\or\qbezier(0,1.)(-0.327,0.327)(-1,0)
\or\qbezier(0,1.)(-0.311,0.179)(-0.866,-0.5)
\or\qbezier(0,1.)(-0.2,0.0536)(-0.5,-0.866)
\or\put(0, 1){\line(0, -2){2}}
\or\qbezier(0,1.)(0.2,0.0536)(0.5,-0.866)
\or\qbezier(0,1.)(0.311,0.179)(0.866,-0.5)
\fi
\or\ifcase#2\or\or\or\or
\or\qbezier(-0.5,0.866)(-0.378,0.378)(-0.866,0.5)
\or\qbezier(-0.5,0.866)(-0.45,0.26)(-1,0)
\or\qbezier(-0.5,0.866)(-0.446,0.12)(-0.866,-0.5)
\or\qbezier(-0.5,0.866)(-0.359,0)(-0.5,-0.866)
\or\qbezier(-0.5,0.866)(-0.2,-0.0536)(0,-1)
\or\put(-0.5, 0.866){\line(3, -5){1}}
\or\qbezier(-0.5,0.866)(0.146,0.146)(0.866,-0.5)
\fi
\or\ifcase#2\or\or\or\or\or
\or\qbezier(-0.866,0.5)(-0.516,0.138)(-1,0)
\or\qbezier(-0.866,0.5)(-0.52,0)(-0.866,-0.5)
\or\qbezier(-0.866,0.5)(-0.446,-0.12)(-0.5,-0.866)
\or\qbezier(-0.866,0.5)(-0.311,-0.179)(0,-1)
\or\qbezier(-0.866,0.5)(-0.146,-0.146)(0.5,-0.866)
\or\put(-0.866, 0.5){\line(5, -3){1.73}}
\fi
\or\ifcase#2\or\or\or\or\or\or
\or\qbezier(-1,0)(-0.516,-0.138)(-0.866,-0.5)
\or\qbezier(-1,0)(-0.45,-0.26)(-0.5,-0.866)
\or\qbezier(-1,0)(-0.327,-0.327)(0,-1)
\or\qbezier(-1,0)(-0.179,-0.311)(0.5,-0.866)
\or\qbezier(-1,0)(-0.0536,-0.2)(0.866,-0.5)
\fi
\or\ifcase#2\or\or\or\or\or\or\or
\or\qbezier(-0.866,-0.5)(-0.378,-0.378)(-0.5,-0.866)
\or\qbezier(-0.866,-0.5)(-0.26,-0.45)(0,-1)
\or\qbezier(-0.866,-0.5)(-0.12,-0.446)(0.5,-0.866)
\or\qbezier(-0.866,-0.5)(0,-0.359)(0.866,-0.5)
\fi
\or\ifcase#2\or\or\or\or\or\or\or\or
\or\qbezier(-0.5,-0.866)(-0.138,-0.516)(0,-1)
\or\qbezier(-0.5,-0.866)(0,-0.52)(0.5,-0.866)
\or\qbezier(-0.5,-0.866)(0.12,-0.446)(0.866,-0.5)
\fi
\or\ifcase#2\or\or\or\or\or\or\or\or\or
\or\qbezier(0,-1.)(0.138,-0.516)(0.5,-0.866)
\or\qbezier(0,-1.)(0.26,-0.45)(0.866,-0.5)
\fi
\or\ifcase#2\or\or\or\or\or\or\or\or\or\or
\or\qbezier(0.5,-0.866)(0.378,-0.378)(0.866,-0.5)
\fi\fi\fi\fi}
\def\FullChord[#1,#2]{
\Endpoint[#1]
\Endpoint[#2]
\Arc[#1]
\Arc[#2]
\Chord[#1,#2]
}
\def\EndChord[#1,#2]{
\Endpoint[#1]
\Endpoint[#2]
\Chord[#1,#2]
}
\def\Picture#1{
\begin{picture}(2,1)(-1,-0.167)
#1
\end{picture}
}
\def\DottedChordDiagram[#1,#2]{
\Picture{\DottedCircle \FullChord[#1,#2]}
}
\def\mapright#1{\smash{ \mathop{\longrightarrow}\limits^{#1}}}
\def\mapdown#1{\Big\downarrow\rlap
   {$\vcenter{\hbox{$\scriptstyle#1$}}$}}
\def\mapne#1{\nearrow\rlap
   {$\vcenter{\hbox{$\scriptstyle#1$}}$}}
\def\mapse#1{\searrow\rlap
   {$\vcenter{\hbox{$\scriptstyle#1$}}$}}

\def\ZZ{ \mathbb{Z} }
\def\IQ{ \mathbb{Q} }
\def\IC{ \mathbb{C} }
\def\IR{ \mathbb{R} }

\def\u#1{ \underline{#1} }

\def\ux{ {\u{x}} }
\def\uy{ {\u{y}} }
\def\ua{ {\u{a}} }
\def\ual{ {\u{\a}} }

\def\um{ {\u{m}} }
\def\ut{ {\u{t}} }

\def\ut{ {\u{t}} }
\def\uo{ {\u{o}} }

\def\vx{ \vec{x} }
\def\vy{ \vec{y} }
\def\va{ \vec{a} }
\def\vb{ \vec{b} }

\def\vr{ \vec{r} }
\def\vl{ \vec{\lambda} }
\def\vrho{ \vec{\rho} }
\def\vua{ \vec{\ua} }

\def\val{ \vec{\a} }
\def\vbet{ \vec{\b} }

\def\bx{ \bar{x} }
\def\by{ \bar{y} }

\def\tb{ \tilde{b} }

\def\mfg{ \mathfrak{g} }
\def\mfh{ \mathfrak{h} }
\def\mfr{ \mathfrak{r} }
\def\mfv{ \mathfrak{v} }
\def\mfvs{ \mfv^* }

\def\ie{{\it i.e.\ }}
\def\eg{{\it e.g.\ }}
\def\cf{{\it cf.\ }}
\def\rhs{{\it r.h.s.\ }}
\def\lhs{{\it l.h.s.\ }}

\def\bl{BL\ }

\def\End{\mathop{{\rm End}}\nolimits}

\def\Tr{\mathop{{\rm Tr}}\nolimits}
\def\STr{\mathop{{\rm STr}}\nolimits}
\def\Vol{\mathop{{\rm Vol}}\nolimits}

\def\span{\mathop{{\rm span}}\nolimits}

\def\deg{ \mathop{{\rm deg}}\nolimits }
\def\p{^{\prime}}

\def\sdim{\mathop{{\rm sdim}}\nolimits}

\def\lk{\mathop{{\rm lk}}\nolimits}

\def\End{\mathop{{\rm End}}\nolimits}

\def\pr#1#2{ \noindent{\em Proof of #1~\ref{#2}.} }
\def\proof{ \noindent{\em Proof.} }

\def\qed{ \hfill $\Box$ }
\def\const{ {\mbox{const}} }

\def\lrbc#1{ \left( #1 \right) }
\def\lrbs#1{ \left[ #1 \right] }

\def\atv#1#2{ \left. #1\right|_{#2} }
\def\atvb#1#2{  #1\big|_{#2} }

\def\ldtc{ ,\ldots, }

\def\Ad{ {\rm Ad} }
\def\Adv#1{ \Ad_{#1} }

\def\ad{ {\rm ad} }
\def\adv#1{ \ad_{#1} }

\def\eqsym{ \mathop{\;=\;}\limits }
\def\dfs{ {\rm def} }
\def\eqdef{ \eqsym_{\dfs} }

\def\widevec#1{ \mathop{#1}\limits^{\longrightarrow} }

\def\smvi#1{ \sum_{#1}^\infty }
\def\snzi{ \smvi{n=0} }
\def\snoi{ \smvi{n=1} }
\def\snti{ \smvi{n=2} }
\def\smoi{ \smvi{m=1} }
\def\smzi{ \smvi{m=0} }
\def\smti{ \smvi{m=2} }
\def\smnzi{ \smvi{m,n=0} }
\def\smthi{ \smvi{m=3} }

\def\svldp{ \sum_{\vl\in\Delp} }

\def\prvldp{ \prod_{\vl\in\Delp} }

\def\Delp{ \Delta_+ }

\def\xmap#1#2#3{ #1\colon #2\to #3 }
\def\ximap#1#2#3{ #1\colon #2\hookrightarrow #3 }
\def\xmapt#1{ \mathop{{\longmapsto}}^{#1} }

\def\urcc{$U(1)$-RCC }
\def\unknot{ {\rm unknot} }

\def\hlf{ {1\over 2} }

\def\Invr#1#2{ #1^{#2} }
\def\Invrs#1#2{ \Invr{(#1)}{#2} }
\def\Invrss#1#2{ \Invr{\lrbcs{#1}}{#2} }

\def\limv#1#2{ #1 \longrightarrow #2 }
\def\limone#1{ \limv{#1}{1} }
\def\limzero#1{ \limv{#1}{0} }

\def\qp#1{ q^{#1} }
\def\qal{ \qp{\a} }
\def\qali#1{ \qp{\a_{#1}} }

\def\pdv#1#2{ #1^{#2} - #1^{-#2} }
\def\pdhv#1#2{ \pdv{#1}{#2/2} }
\def\pdhq{ \pdv{q}{1} }
\def\pdhqa{ \pdv{q}{\a} }

\def\pdht{ \pdhv{t}{1} }

\def\ssSt#1{ #1\subset S^3 }
\def\kst{ \ssSt{\cK} }
\def\lst{ \ssSt{\cL} }

\def\lst{ \cL\subset S^3 }
\def\stml{ S^3\setminus \cL }

\def\Jv#1#2#3{ J_{#1}(#2;#3) }
\def\Jvq#1#2{ \Jv{#1}{#2}{q} }
\def\JaK{ \Jvq{\a}{\cK} }

\def\AP#1{ \Delta_{\rm A}^{#1} }
\def\APv#1#2{ \AP{}(#1;#2) }
\def\APKv#1{ \APv{\cK}{#1} }
\def\APKt{ \APKv{t} }
\def\APKqa{ \APKv{\qal} }

\def\APpv#1#2#3{ \AP{#1}(#2;#3) }
\def\APKpv#1#2{ \APpv{#1}{\cK}{#2} }
\def\APKpqa#1{ \APKpv{#1}{\qal} }

\def\AF#1{ \nabla_{\rm A}^{#1} }
\def\AFv#1#2{ \AF{}(#1;#2) }
\def\AFLv#1{ \AFv{\cL}{#1} }
\def\AFLtN{ \AFLv{\toN} }
\def\AFLut{ \AFLv{\ut} }
\def\AFLeua{ \AFLv{e^{\ua}} }

\def\tAF#1{ \tilde{\nabla}_{\rm A}^{#1} }
\def\tAFLv#1{ \tAF{}(\cL;#1) }
\def\tAFLua{ \tAFLv{\ua} }
\def\tAFLhue{ \tAFLv{\uhed} }

\def\AFpv#1#2#3{ \AF{#1}{(#2;#3)} }
\def\AFLpv#1#2{ \AFpv{#1}{\cL}{#2} }

\def\AFLptN#1{ \AFLpv{#1}{\toN} }

\def\xP{ P }
\def\Pv#1#2#3{ \xP_{#1}(#2;#3) }
\def\PnKv#1{ \Pv{n}{\cK}{#1} }
\def\PnKt{ \PnKv{t} }
\def\PnKqa{ \PnKv{\qal} }

\def\PnLv#1{ \Pv{n}{\cL}{#1} }
\def\PnLtN{ \PnLv{\toN} }

\def\Zti{ \ZZ[t,t^{-1}] }

\def\ZK{ Z }
\def\ZKiv#1#2{ \ZK_{#1}(#2) }
\def\ZKAv#1{ \ZKiv{\cA}{#1} }
\def\ZKAL{ \ZKAv{\cL} }
\def\ZKALLp{ \ZKAv{\cL\cup\cLp} }
\def\ZKALM{ \ZKAv{\cL,M} }
\def\ZKALp{ \ZKAv{\cLp} }
\def\ZKALpv#1{ \ZKAv{\cLp;#1} }
\def\ZKBv#1{ \ZKiv{\cB}{#1} }
\def\ZKBCv#1{ \ZKiv{\cBC}{#1} }
\def\ZKBL{ \ZKBv{\cL} }
\def\ZKBLM{ \ZKBv{\cL,M} }

\def\ZKBLpe{ \ZKBv{\cLp;\DLx,\ubfx} }
\def\ZKBLx{ \ZKBv{\cL;\ubfx} }
\def\ZKBLxt{ \ZKBv{\cL;\ubfx;t} }
\def\ZKBCLv#1{ \ZKBCv{\cL;#1} }
\def\ZKBCLx{ \ZKBCLv{\ubfx} }

\def\ZKOBv#1{ \ZK_{\cB}^{\Wh}(#1) }
\def\ZKOBL{ \ZKOBv{\cL} }
\def\ZKOBLxt{ \ZKOBv{\cL;\ux;t} }
\def\ZKOBCv#1{ \ZK_{\cBC}^{\Wh}(#1) }
\def\ZKOBCLv#1{ \ZKOBCv{\cL;#1} }
\def\ZKOBCLx{ \ZKOBCLv{\ubfx} }

\def\ZKDCv#1{ \ZKiv{\cDC}{#1} }
\def\ZKDCLv#1{ \ZKDCv{\cL;#1} }
\def\ZKDCLx{ \ZKDCLv{\ubfx} }

\def\ZKriv#1#2{ \ZK_{#1}\irom(#2) }
\def\ZKDCrv#1{ \ZKriv{\cDC}{#1} }
\def\ZKDCrLv#1{ \ZKDCrv{\cL;#1} }
\def\ZKDCrLx{ \ZKDCrLv{\uxz} }
\def\ZKDCrLxk{ \atv{\ZKDCrLx}{x_k=0} }

\def\ZKOriv#1#2{ \ZK_{#1}\iromW(#2) }
\def\ZKODCrv#1{ \ZKOriv{\cDC}{#1} }
\def\ZKODCrLv#1{ \ZKODCrv{\cL;#1} }
\def\ZKODCrLx{ \ZKODCrLv{\uxz} }

\def\xuo{ {\rm r} }
\def\xuoc{ {\rm c} }
\def\dxuor{ _\xuo }
\def\dxuoc{ _\xuoc }
\def\uo{ ^{\xuo} }

\def\ZKuriv#1#2{ \ZK_{#1}\uo(#2) }
\def\ZKDCurv#1{ \ZKuriv{\cQDC}{#1} }

\def\ZKDCurLa{ \ZKDCurv{\cL} }
\def\ZKDCurKa{ \ZKDCurv{\cK} }
\def\ZKBK{ \ZKiv{\cB}{\cK} }
\def\ZKBKx{ \ZKiv{\cB}{\cK;\bfx} }
\def\ZKBKa{ \ZKiv{\cB}{\cK;a} }

\def\ZKOuriv#1#2{ \ZK_{#1}^{\Wh,\xuo}(#2) }
\def\ZKODCurv#1{ \ZKOuriv{\cQDC}{#1} }

\def\ZKODCurLa{ \ZKODCurv{\cL} }

\def\WKuriv#1#2{ \WK_{#1}\uo(#2) }
\def\WKDCurLa{ \WKuriv{\cQDC}{\cL} }

\def\WKriv#1#2{ \WK_{#1}\irom(#2) }
\def\WKDCrv#1{ \WKriv{\cDC}{#1} }
\def\WKDCrLv#1{ \WKDCrv{\cL;#1} }
\def\WKDCrLx{ \WKDCrLv{\uxz} }

\def\WKOriv#1#2{ \WK_{#1}\iromW(#2) }
\def\WKODCrv#1{ \WKOriv{\cDC}{#1} }
\def\WKODCrLv#1{ \WKODCrv{\cL;#1} }
\def\WKODCrLx{ \WKODCrLv{\ux} }

\def\WKOuriv#1#2{ \WK\irorW_{#1}(#2) }
\def\WKODCurLa{ \WKOuriv{\cQDC}{\cL} }

\def\WKBCv#1{ \WKiv{}{\cBC}{#1} }
\def\WKBCLv#1{ \WKBCv{\cL;#1} }
\def\WKBCLx{ \WKBCLv{\ubfx} }

\def\WKFDriv#1#2{ \WK_{#1}\iromFD(#2) }
\def\WKFDDCrv#1{ \WKFDriv{\cDC}{#1} }
\def\WKFDDCrLv#1{ \WKFDDCrv{\cL;#1} }
\def\WKFDDCrLx{ \WKFDDCrLv{\uxz} }

\def\ZKODCv#1{ \ZK_{\cDC}^{\Wh}(#1) }
\def\ZKODCLv#1{ \ZKODCv{\cL;#1} }
\def\ZKODCLx{ \ZKODCLv{\ubfx} }

\def\WK{ W }
\def\WKv#1{ \WK(#1) }
\def\WKL{ \WKv{\cL} }
\def\WKLX{ \WKv{\cL;\ubfx} }
\def\WKLstr{ \WK\istr(\cL) }
\def\WKLx{ \WKv{\cL;\ux} }
\def\WKLte{ \WKv{\cL;\ubfx;t} }
\def\WKLpe{ \WKv{\cLp;\DLx;\ubfx} }
\def\WKOv#1{ \WK^{\Wh}(#1) }
\def\WKOL{ \WKOv{\cL} }
\def\WKOLX{ \WKOv{\cL;\ubfx} }
\def\WKODCLx{ \WKOv{\cL;\ubfx} }
\def\WKOLstr{ \WK^{\Wh}\istr(\cL) }
\def\WKiv#1#2#3{ \WK^{#1}_{#2}(#3) }
\def\WKOiv#1#2{ \WKiv{\Wh}{#1}{#2} }
\def\WKOLv#1{ \WKOiv{#1}{\cL} }
\def\WKOLmn{ \WKOLv{m,n} }

\def\xX{ \mathfrak{L} }
\def\xN{ L }
\def\xNp{ L\p }
\def\xD{ D }
\def\xDo{ \xD_1 }
\def\xDt{ \xD_2 }
\def\xDth{ \xD_3 }
\def\xDot{ \xDo\amalg\xDt }
\def\txD{ \tilde{\xD} }
\def\xDir{ \xD_{i,\xuo} }
\def\xDic{ \xD_{i,\xuoc} }
\def\txDir{ \txD_{i,\xuo} }
\def\txDic{ \txD_{i,\xuoc} }

\def\xCD{ \cC_D }
\def\xchi{ \hat{\chi} }
\def\echi{ \chi }
\def\xF{ F }
\def\xfl{ \flat }
\def\xfle{ \,\xfl\, }
\def\xG{ \mathfrak{G} }
\def\xT{ \mathfrak{T} }

\def\xk{ k }
\def\xvr{ \vx }

\def\bfa{ \mathbf{a} }
\def\bfn{ \mathbf{n} }
\def\bfD{ \mathbf{D} }
\def\bfY{ \mathbf{Y} }
\def\vbfY{ \vec{\bfY} }

\def\Jc{ \check{J} }
\def\Jcr{ \Jc^{\rm (r)} }
\def\Jcrv#1#2#3{ \Jcr(#1;#2;#3) }
\def\JcrLth{ \Jcrv{\cL}{\toN}{h} }
\def\JcrLqaN{ \Jcrv{\cL}{\qaoN}{q-1} }

\def\JaNL{ \Jvq{\aloN}{\cL} }
\def\JvaNL{ \Jvq{\uval}{\cL} }

\def\xXv#1{ \xX_{\rm #1} }
\def\xXt{ \xXv{t} }
\def\xXr{ \xXv{r} }
\def\xXc{ \xXv{c} }

\def\toN{ t_1\ldtc t_{\xN} }
\def\aloN{ \a_1\ldtc \a_{\xN} }
\def\qaoN{ \qali{1}\ldtc \qali{\xN} }
\def\valoN{ \val_1\ldtc \val_{\xN} }
\def\hvaloN{ \hb\val_1\ldtc \hb\val_{\xN} }
\def\vbetoN{ \vbet_1\ldtc \vbet_{\xN} }
\def\vboN{ \vbm{1}\ldtc\vbm{\xN} }
\def\vaoN{ \vam{1}\ldtc\vam{\xN} }
\def\dvbetoN{ d\vbet_1\cdots d\vbet_{\xN} }

\def\xoN{ x_1\ldtc x_{\xN} }

\def\ut{ \underline{t} }
\def\ual{ \underline{\a} }
\def\qua{ \qp{\ual} }
\def\uval{ \underline{\val} }
\def\huval{ \hb\uval }
\def\uvbet{ \underline{\vbet} }
\def\uvb{ \underline{\vb} }
\def\uva{ \underline{\va} }
\def\duvbet{ d\uvbet }

\def\ux{ \underline{x} }

\def\uvbrk{ \uvb\rv{k} }
\def\duvbrk{ d\uvbrk }

\def\aloN{ \ual }
\def\qaoN{ \qua }
\def\valoN{ \uval }
\def\hvaloN{ \huval }
\def\vbetoN{ \uvbet }
\def\vboN{ \uvb }
\def\vaoN{ \uva }
\def\dvbetoN{ \duvbet }

\def\xoN{ \ux }

\def\vbetm#1{ \vbet_{#1} }
\def\vbetj{ \vbetm{j} }
\def\vbm#1{ \vb_{#1} }
\def\vam#1{ \va_{#1} }
\def\vbj{ \vbm{j} }

\def\ZZtNh{ \ZZ[t_1^{\pm 1/2}\ldtc t_{\xN}^{\pm 1/2},1/2] }
\def\QIQuua{ \Qf{\IQuua} }
\def\IQuua{ \IQ[[\ua]] }
\def\evQIQuua{ \evn{\QIQuua} }

\def\evn#1{ \Big[ #1 \Big]^{\rm even} }

\def\lk{ \mathop{{\rm lk}}\nolimits }
\def\lkf#1#2{ \lk(#1,#2) }
\def\lkLf#1{ \lkf{\cL}{#1} }

\def\lkLuva{ \lkLf{\vua} }

\def\lkn#1#2{ l_{#1 #2} }
\def\lkij{ \lkn{i}{j} }

\def\snoijxNk{ \sum_{1\leq i\neq j\leq \xN\atop i\neq k} }
\def\sjoxNk{ \sum_{1\leq j\leq \xN\atop j\neq k} }

\def\smnoi{ \smvi{m,n=1} }
\def\sjoxN{ \sum_{j=1}^{\xN} }
\def\sioxN{ \sum_{i=1}^{\xN} }
\def\sijoxN{ \sum_{i,j=1}^{\xN} }
\def\sijoxNp{ \sum_{i,j=1}^{\xNp} }
\def\sijpxN{ \sum_{i,j=\xNp+1}^{\xN} }
\def\sijoxX{ \sum_{i,j=1}^{\orX} }

\def\smnd{ \sum_{ m\geq 0 \atop n\geq -1} }
\def\smnd{ \smvi{m,n=0} }
\def\smnd{ \smnzi }

\def\sDDL{ \sum_{\xD\in\cgDL} }
\def\sDDLn{ \sum_{\xD\in\cgDL\atop \dgth(\xD)=n } }

\def\pjoxN{ \prod_{j=1}^{\xN} }
\def\pioxN{ \prod_{i=1}^{\xN} }
\def\pldp{ \prod_{\vl\in\Delta_+} }

\def\arN{ \uparrow_{1}\cdots\uparrow_{\xN} }

\def\cAxy{ \cA(\uparrow_{X},\uparrow_{Y}) }
\def\cAz{ \cA(\uparrow_{Z} ) }
\def\cBstv#1{ B_{\{#1\}} }
\def\cBxy{ \cBstv{X,Y} }
\def\cBz{ \cBstv{Z} }

\def\cAL{ \cA(\arN) }
\def\cALs{ \cA(\uparrow_{\cL}) }
\def\cBL{ \cB_{\cL} }

\def\xupr#1{ \uparrow_{#1} }

\def\cALp{ \cA(\xupr{\DLxo},\xupr{X_1}\ldtc\xupr{X_{\xN}}) }
\def\cBLp{ \cBstv{ \DLxo,X_1\ldtc X_{\xN} } }

\def\hb{ \hbar }

\def\pnvaL{ p_n(\cL;\valoN) }

\def\Vm#1{ V_{#1} }
\def\Vval{ \Vm{\val} }
\def\Vmva#1{ \Vm{\val_{#1}} }

\def\Vvalk{ \Vmva{k} }
\def\VmvaoN{ \Vmva{1}\ldtc \Vmva{\xN} }

\def\TVmvaoX{ \Vmva{1}\otimes\cdots\otimes\Vmva{\orX} }

\def\TVuva{ \Vm{\uval} }
\def\TVua{ \Vm{\ua} }

\def\Tg{ \cT_{\mfg} }
\def\Tuoo{ \cT_{\uoo} }
\def\Tgh{ \cT_{\mfg}^{1,\hb} }
\def\Tgth{ \cT_{\mfg}^{2,\hb} }
\def\tTg{ \tilde{\cT}_{\mfg} }

\def\orX{ |\xX| }
\def\orXp{ |\Xp| }

\def\UgL{ \Invrs{(\Ug)^{\otimes \xN}}{\xG} }
\def\UgX{ \Invrs{(\Ug)^{\otimes \orX}}{\xG} }
\def\UgLh{ \UgL\,[[\hb]] }
\def\Ug{ \Urs\mfg }
\def\Sg{ \Srs\mfg }

\def\SgL{ \Invrs{(\Sg)^{\otimes \xN}}{\xG} }
\def\SgLo{ (\Sg)^{\otimes (\xN-1)} }
\def\SgX{ \Invrs{(\Srs g)^{\otimes \orX}}{\xG} }

\def\oSgXT{ (\Sh)^{\otimes |\xXc|}\otimes \Invrs{(\Sra)^{\otimes |\xXr|}\otimes
(\Sg)^{\otimes|\xXt|} }{\xT} }
\def\QfShXc{ \Qf{\ShXc} }
\def\tSgXt{ \Qf{\ShXc} \otimes (\Sh)^{\otimes |\xXc\setminus\Xp|}
\otimes \Invrs{(\Sra)^{\otimes |\xXr|}\otimes
(\Sg)^{\otimes|\xXt|} }{\xT} }
\def\tSgXts{ (\Sh)^{\otimes |\xXc\setminus\Xp|}
\otimes \Invrs{(\Sra)^{\otimes |\xXr|}\otimes
(\Sg)^{\otimes|\xXt|} }{\xT} }

\def\Sh{ S\mfh }
\def\Sra{ S\mfr }

\def\ShXc{ (\Sh)^{\otimes |\Xp|} }

\def\hm{ h }

\def\ICh{ \IC[[\hb]] }

\def\mfgs{ \mfg^{\ast} }
\def\mfgsC{ \mfg^{\ast}_{\IC} }
\def\mfhs{ \mfh^{\ast} }

\def\mfgC{ \mfg_{\IC} }
\def\mfrC{ \mfr_{\IC} }

\def\Va{ V_{\val} }

\def\TrVua{ \Tr_{\TVuva} }
\def\TrVa{ \Tr_{\Va} }

\def\dgh#1{ \deg_1 #1 }
\def\dgth#1{ \deg_2 #1 }

\def\bmg{ \b_{\mfg} }

\def\duflo{ D(j^{1/2}_{\mfg}) }
\def\dflL{ D_{\cL}(j^{1/2}_{\mfg}) }
\def\Wh{ \Omega }

\def\Whh{ \hat{\Wh} }
\def\Whhin#1{ \Whh_{#1} }
\def\Whhj{ \Whhin{j} }

\def\Whhjv#1{ \Whhin{j}{#1} }

\def\Whin#1{ \Wh_{#1} }
\def\Whj{ \Whin{j} }

\def\WhhL{ \Whhin{\cL} }

\def\btn{ b_{2n} }
\def\omtn{ \omega_{2n} }

\def\corbv#1{ \cO_{#1} }
\def\corbal{ \corbv{\val} }
\def\corbalv#1{ \corbv{\val_{#1}} }
\def\icorbal{ \int_{\corbal} }

\def\corbuval{ \corbv{\uval} }
\def\iOuval{ \int_{\corbuval} }
\def\corbuvalrk{ \corbv{\uvalrk} }
\def\iOuvalrk{ \int_{\corbuvalrk} }
\def\corbhuvalrk{ \corbv{\hb\uvalrk} }
\def\iOhuvalrk{ \int_{\corbhuvalrk} }

\def\corbuvark{ \corbv{\uvark} }
\def\iOuvark{ \int_{\corbuvark} }

\def\corba{ \corbv{\va} }

\def\rv#1{_{(#1)}}

\def\uvalrk{ \uval\rv{\xk} }
\def\uvark{ \uva\rv{\xk} }

\def\dgf#1{ d_{\mfg}(#1) }
\def\dgva{ \dgf{\va} }
\def\dgvaj{ \dgf{\va_j} }
\def\dgval{ \dgf{\val} }
\def\dgvak{ \dgf{\vak} }
\def\dguva{ \dgf{\uva} }
\def\dgrho{ \dgf{\vrho} }
\def\dguvarho{ {\dguva\over\dgrho} }

\def\vak{ \va_{\xk} }
\def\vbk{ \vb_{\xk} }

\def\del{ \partial }
\def\parv#1{ \del_{#1} }
\def\parxv#1{ \parv{x_{#1}} }
\def\parxj{ \parxv{j} }

\def\Lv#1#2#3{ L_{#1}(#2;#3) }
\def\LLv#1#2{ \Lv{#1}{\cL}{#2} }
\def\LLbetv#1{ \LLv{#1}{\vbetoN} }
\def\LLbetmn{ \LLbetv{m,n} }
\def\LLbetmz{ \LLbetv{m,0} }
\def\LLbetmo{ \LLbetv{m,1} }
\def\LLbetzz{ \LLbetv{0,0} }
\def\LLbettz{ \LLbetv{2,0} }
\def\LLamo{ \LLv{m,1}{\uva} }
\def\LLamos{ \LLv{m,1}{\ua} }

\def\LLbv#1{ \LLv{#1}{\vboN} }
\def\LLbmn{ \LLbv{m,n} }
\def\LLbmz{ \LLbv{m,0} }
\def\LLbmo{ \LLbv{m,1} }

\def\LLmzuva{ \LLv{m,0}{\uva} }
\def\LLmzuvb{ \LLv{m,0}{\uvb} }

\def\Lr{ L^{\rm R} }
\def\Lrv#1#2#3#4{ \Lr_{#1}(#2;#4;#3) }
\def\LLrav#1#2{ \Lrv{#1}{\cL}{#2}{\uva} }
\def\LLrarkv#1{ \LLrav{#1}{\uxvrk} }
\def\LLrarkmn{ \LLrarkv{m,n} }
\def\LLrarkmz{ \LLrarkv{m,0} }
\def\LLrarkzz{ \LLrarkv{0,0} }
\def\LLrarkon{ \LLrarkv{1,n} }
\def\LLrarkoz{ \LLrarkv{1,0} }
\def\LLrarktz{ \LLrarkv{2,0} }
\def\LLrarkzo{ \Lr_{0,1}(\cL;\uva;\uxvrk) }
\def\LLrarkzos{ \Lr_{0,1}(\cL;\ua) }

\def\LLradrkv#1{ \LLrav{#1}{\del_{\uxvrk}} }
\def\LLradrkmz{ \LLradrkv{m,0} }
\def\LLradrkmn{ \LLradrkv{m,n} }

\def\avbav#1#2{ \left. #1 \right|_{\vb_{\rm{#2}} = \va_{#2}} }
\def\LLbmnk{ \avbav{\LLbmn}{k} }
\def\LLbmzk{ \avbav{\LLbmz}{k} }

\def\Jfv#1#2#3#4{ J(#1;#2;#3,#4) }
\def\Jfvk#1#2#3{ \Jfv{#1}{#2}{#3}{\xk} }
\def\Jfhv#1#2{ \Jfvk{#1}{#2}{\hb} }
\def\JfLv#1{ \Jfhv{\cL}{#1} }
\def\JfvaNL{ \JfLv{\vaoN} }

\def\Ir{ I^{({\rm r})} }
\def\Ir{ I\uo_{\mfg} }
\def\Irfv#1#2#3#4{ \Ir(#1;#2;#3,#4) }
\def\Irfvk#1#2#3{ \Irfv{#1}{#2}{#3}{\xk} }
\def\Irfhv#1#2{ \Irfvk{#1}{#2}{\hb} }
\def\IrfLv#1{ \Irfhv{\cL}{#1} }
\def\IrfvaNL{ \IrfLv{\vaoN} }
\def\IrsufaNL{ I\uo_{su(2)}(\cL;\ua;\hb,k) }

\def\Irv#1#2#3{ \Ir(#1;#2;#3) }
\def\Irhv#1#2{ \Irv{#1}{#2}{\hb} }
\def\IrhLv#1{ \Irhv{\cL}{#1} }
\def\IrhLuva{ \IrhLv{\uva} }

\def\istl#1{ \int\limits_{[#1]} }
\def\ist#1{ \int_{[#1]} }
\def\ilvrkz{ \istl{\uxvrk=0} }
\def\ivrkz{ \ist{\uxvrk=0} }

\def\uxvr{ \u{\xvr} }
\def\uxvrk{ \uxvr\rv{\xk} }
\def\duxvrk{ d\uxvrk }

\def\uxrvk{ \ux\rv{k} }
\def\duxrvk{ d\uxrvk }

\def\RiL{ \exp (\Tgh(\WKOL)) }

\def\ohb{ {1\over \hb} }

\def\nrtg{ \nrt{\mfg} }
\def\nrtg{ \nDp }

\def\hbpLg{ \hb^{-\xN\nDp} }

\def\mohbpLg{ (-2\pi\hb)^{(1-\xN)\nDp} }
\def\hbpLgo{ \hb^{ - (\xN - 1)\, \nrtg} }
\def\hbmnDp{ \hb^{-\nDp} }

\def\scp#1#2{ (#1,#2) }

\def\GmT{ \xG/\xT }

\def\xvrj{ \xvr_j }
\def\xvri{ \xvr_i }

\def\vaj{ \va_j }
\def\vai{ \va_i }
\def\vajv#1{ \va_{j_{#1}} }

\def\msr#1#2{ \msj_{#1}(#2) }
\def\msra#1{ \msr{\va}{#1} }
\def\msravr{ \msra{\xvr} }
\def\msraz{ \msra{\vzr} }

\def\msj{ \mu }
\def\omsj{ \omega }
\def\bomsj{ \boldsymbol{\omega} }
\def\umsj{ \boldsymbol{\msj} }
\def\msjv#1{ \msj(#1) }
\def\omsjv#1{ \omsj(#1) }
\def\bomsjv#1{ \bomsj(#1) }
\def\umsjv#1{ \umsj(#1) }
\def\msgr{ \msjv{\crr} }
\def\msgx{ \msjv{x} }
\def\msgvx{ \msjv{\vx} }
\def\omsgx{ \omsjv{x} }
\def\bomsgx{ \bomsjv{\ux} }
\def\umsgx{ \umsjv{\ux} }

\def\xXi{ \Xi }

\def\vzr{ \vec{0} }

\def\AdgT{ \Adv{g\xT} }
\def\Ader{ \Adv{\er} }
\def\Adey{ \Adv{\ey} }

\def\er{ e^{\xvr} }
\def\ey{ e^{\vy} }
\def\eFLyr{ e^{\FLyr} }

\def\eyv#1{ e^{\vy_{#1}} }
\def\eyo{ \eyv{1} }
\def\eyt{ \eyv{2} }

\def\FLyto{ \FLv{\vy_2}(\vy_1) }
\def\eFLyto{ e^{\FLyto} }

\def\FL{ F^{\rm L} }
\def\FL{ {\rm F} }
\def\FLv#1{ \FL_{#1} }
\def\FLy{ \FLv{\vy} }
\def\FLyr{ \FLy(\xvr) }
\def\FLyv#1{ \FLv{\vy_{#1}} }
\def\FLyo{ \FLyv{1} }
\def\FLyt{ \FLyv{2} }

\def\FLgy{ \FLv{\cry} }
\def\FLgx{ \FLv{x} }
\def\FLgmx{ \FLv{-x} }

\def\tF{ \tilde{\FL} }

\def\tFLv#1{ \tF_{#1} }
\def\tFLga{ \tFLv{\cra} }
\def\tFLgy{ \tFLv{\cry} }
\def\tFLgta{ \tFLv{t\cra} }

\def\tFLyx{ \tFLv{\vy}(\xvr) }

\def\cry{ y }

\def\FLdyot{ \FLv{\FLyto} }

\def\FLp{ {\FL}^{\rm L,\prime} }
\def\FLp{ {\FL}^{\prime} }
\def\FLpv#1{ \FLp_{#1} }

\def\drFLv#1#2{ \FLpv{#1}(#2) }
\def\drFLb#1{ \FL_{#1,\star} }
\def\drFLb#1{ \FLpv{#1} }
\def\drFLyvv#1{ \drFLv{\vy}{#1} }

\def\drFLyr{ \drFLyvv{\xvr} }
\def\drFLyv#1{ \drFLb{\vy_{#1}} }
\def\drFLyo{ \drFLyv{1} }
\def\drFLyt{ \drFLyv{2} }
\def\drFLdyot{ \drFLb{\FLyto} }
\def\drFLvrz{ \drFLv{\xvr}{\vzr} }
\def\drFLmvrr{ \drFLv{-\xvr}{\xvr} }

\def\vFhyto{ \tFLv{\vy_2}(\vy_1) }

\def\AdvFhyto{ \Adv{\exp\big(\vFhyto\big)} }

\def\dtmfr{ \det_{\mfr} }

\def\dxvr{ \Delta\xvr }

\def\drr#1{ {I-e^{-\adv{#1}}\over \adv{#1}} }
\def\drrrw{ \drr{\xvr} }
\def\drrr{ \drrrw\;}
\def\drrrd{ \drrr\dxvr }
\def\drrrds{ \lrbc{\drrrd\,} }
\def\prjr{ P_{\mfr} }
\def\prjh{ P_{\mfh} }

\def\adxvr{ \adv{\xvr} }

\def\Idop{ I }

\def\TgD{ \Tg D }
\def\TgDv#1{ \Tg D(#1) }
\def\vbom{ \vb_1\ldtc\vb_m }
\def\vaom{ \va_1\ldtc\va_m }

\def\vatm{ \va_2\ldtc\va_m }
\def\TgDbm{ \TgDv{\vbom} }
\def\TgDxbm{ \TgDv{\vb,\vatm} }
\def\TgDbmx{ \TgDv{ [\vb,\va_1],\va_2\ldtc\va_m} }

\def\lstbiv#1#2#3{ #1_{#2}\; (1\leq #2\leq #3) }
\def\lstbivj#1#2{ \lstbiv{#1}{j}{#2} }
\def\lstbivjxN#1{ \lstbivj{#1}{\xN} }

\def\ojxNi{ 1\leq j\leq \xN }
\def\oixNi{ 1\leq i\leq \xN }
\def\oijxNi{ 1\leq i,j\leq \xN }

\def\dlxj{ \del_{x_j} }
\def\dlxi{ \del_{x_i} }
\def\dlyi{ \del_{y_i} }
\def\dluy{ \del_{\uy} }
\def\dlux{ \del_{\ux} }
\def\dlvux{ \del_{\vux} }
\def\dluz{ \del_{\uz} }
\def\dlz{ \del_z }

\def\ndv#1{ \nabla_{#1} }
\def\ndxi{ \ndv{\xi} }
\def\ndxia{ \ndv{\xia} }
\def\ndmada{ \ndv{\mad_a} }

\def\dvad{ {\rm dvd} }
\def\dvadv#1{ \dvad_{#1}\, }

\def\xmltv#1#2{ {\rm m}_{#1}^{#2} }
\def\mxyz{ \xmltv{X,Y}{Z} }
\def\myxz{ \xmltv{Y,X}{Z} }
\def\mcxyz{ \xmltv{[X,Y]}{Z} }

\def\DLt{ \Delta t }
\def\tDLt{ t+\DLt }
\def\DLx{ \Delta X }
\def\DLxo{ \DLx_1 }
\def\dlxo{ \del_{X_1} }
\def\cLp{ \cL\p }

\def\mfrxN{ \mfr^{\oplus(\xN-1)} }

\def\xQv#1#2#3{ Q_{#1}(#2;#3) }
\def\gQv#1#2#3{ \xQv{\mfg,#1}{#2}{#3} }
\def\gQLv#1#2{ \gQv{#1}{\cL}{#2} }
\def\gQLuav#1{ \gQLv{#1}{\uva} }
\def\gQLuak{ \gQLuav{k} }
\def\gQsuLuak{ \xQv{su(2),k}{\cL}{\ua} }

\def\xQiv#1#2#3{ Q^{-1}_{#1}(#2;#3) }
\def\gQiv#1#2#3{ \xQiv{\mfg,#1}{#2}{#3} }
\def\gQLiv#1#2{ \gQiv{#1}{\cL}{#2} }
\def\gQLuaiv#1{ \gQLiv{#1}{\uva} }
\def\gQLuaik{ \gQLuaiv{k} }

\def\xQLxav#1{ \xQv{#1}{\cL}{\ux\rv{#1}} }
\def\xQLxak{ \xQLxav{k} }

\def\xQiv#1#2#3{ Q_{#1}^{-1}(#2;#3) }
\def\xQiLxav#1{ \xQiv{#1}{\cL}{\ux\rv{#1}} }
\def\xQiLxak{ \xQiLxav{k} }

\def\uxrk{ \ux\rv{k} }
\def\uxrkz{ \uxrk }

\def\FD{ \mathop{\rm FD}\nolimits }
\def\FDv#1#2{ \FD(#1;#2) }
\def\FDLv#1{ \FDv{\cL}{#1} }
\def\FDLua{ \FDLv{\uva} }
\def\FDLgua{ \FD(\cL) }

\def\FDLn{ \FD_{\rm nrw} }
\def\FDLngua{ \FDLn(\cL) }

\def\pv#1#2#3{ p_{#1}(#2;#3) }
\def\pLuva#1{ \pv{#1}{\cL}{\uva} }
\def\pQnLuva{ \pLuva{{\rm Q},n} }
\def\pnQnLuva{ \pLuva{n} }

\def\FN#1#2{ { F}_{\nabla,#1}^{#2} }
\def\FNgv#1#2#3{ \FN{\mfg}{#3}(#1,#2) }
\def\FNgLua#1{ \FNgv{\cL}{\uva}{#1} }

\def\FNdg{ \dguva\FNgLua{} }

\def\Qf#1{ \cQ(#1) }

\def\xCC{ {\rm CC} }

\def\cBC{ \cB^{\rm C} }
\def\cBCL{ \cBC_{\cL} }
\def\cBCX{ \cBC_{\xX} }
\def\cBCXD{ \cBCX(\xD) }
\def\tcBC{ \tilde{\cB}^{\rm C} }
\def\tcBCX{ \tilde{\cB}^{\rm C}_{\xX} }
\def\tcBCXD{ \tcBCX(\xD) }
\def\tcBCXIHX{ \tilde{\cB}^{\rm C}_{\xX,\IHX} }
\def\tcBCXIHXv#1{ \tilde{\cB}^{\rm C}_{\xX,\IHX,#1} }
\def\tcBCXIHXi{ \tcBCXIHXv{i} }
\def\tcBCXIHXo{ \tcBCXIHXv{1} }
\def\tcBCXIHXoD{ \tcBCXIHXo(\xD) }
\def\tcBCXIHXt{ \tcBCXIHXv{2} }
\def\tcBCXIHXtD{ \tcBCXIHXt(\xD) }
\def\tcBCXIHXth{ \tcBCXIHXv{3} }
\def\tcBCXIHXf{ \tcBCXIHXv{4} }
\def\tcBCXIHXz{ \tcBCXIHXv{0} }
\def\tcBCXCC{ \tilde{\cB}^{\rm C}_{\xX,\xCC} }
\def\tcBCXCCo{ \tilde{\cB}^{\rm C}_{\xX,\xCC_1} }
\def\tcBCXCCt{ \tilde{\cB}^{\rm C}_{\xX,\xCC_2} }
\def\tcBCXoC{ \tcBCX }
\def\tcBX{ \tilde{\cB}_{\xX} }
\def\tcBXIHX{ \tilde{\cB}_{\xX,\IHX} }

\def\ttcBCXv#1{ \tilde{\tilde{\cB}}^{\rm C}_{\xX}(#1) }
\def\ttcBCXD{ \ttcBCXv{\xD} }

\def\IHX{ {\rm IHX} }
\def\CC{ {\rm CC} }

\def\cAX{ \cA(\uparrow_{\xX}) }

\def\cBv#1{ \cB_{#1} }
\def\cBX{ \cBv{\xX} }
\def\cBL{ \cBv{\cL} }

\def\nev#1#2{ n_{#1}(#2) }
\def\neDv#1{ \nev{#1}{D} }
\def\neDo{ \neDv{1} }
\def\neDt{ \neDv{3} }
\def\neDe{ \neDv{\rm e} }

\def\nehDo{ \neDv{1,\mfh} }
\def\nerDo{ \neDv{1,\mfr} }

\def\zhi{ h^{-1} }
\def\zf{ f }

\def\fmfg{ f_{\mfg} }
\def\hmfg{ h_{\mfg} }

\def\hiv#1{ h^{-1}_{#1} }
\def\himfrv#1{ \hiv{\mfr_{#1}} }
\def\himfg{ h^{-1}_{\mfg} }
\def\himfr{ h^{-1}_{\mfr} }
\def\himfrl{ h^{-1}_{\mfrl} }
\def\himfrc{ h^{-1}_{\mfr,\xc} }
\def\himfh{ h^{-1}_{\mfh} }
\def\himfhij{ h^{-1}_{\mfh;i,j} }
\def\hi{ h^{-1} }
\def\hiD{ \hi_D }

\def\wdgtg{ \bigwedge^3\mfgs }
\def\wdgtgC{ \bigwedge^3\mfgsC }
\def\wdgtgx{ \big(\wdgtg\big)^{\otimes \neDt} }
\def\stg{ \Srs^2\mfg }
\def\stgC{ \Srs^2\mfgC }
\def\stgx{ \big(\stg\big)^{\otimes \neDe} }
\def\stgs{ \Srs^2\mfgs }

\def\str{ \Srs^2\mfr }
\def\sth{ \Srs^2\mfh }
\def\tmfgx{ \mfg^{\otimes\neDo} }
\def\tmfhx{ \mfh^{\otimes\nehDo} }
\def\tmfrx{ \mfr^{\otimes\nerDo} }
\def\tmfc{ \tmfhx\otimes\tmfrx }

\def\evv#1{ \mathop{\rm ev}\nolimits_{#1} }
\def\evvx#1#2{ \evv{#1;#2} }
\def\evkx#1{ \evvx{\xk}{#1} }
\def\evkal{ \evkx{\valk} }
\def\eval{ \evkx{\val} }

\def\valk{ \val_{\xk} }

\def\finc{ f }
\def\ftr{ \finc_{\rm trv} }
\def\ftri{ \ftr^{-1} }
\def\tftr{ \tilde{\finc}_{\rm trv} }
\def\tftri{ \tftr^{-1} }
\def\fct{ \finc_{\rm Crtn} }
\def\fcti{ \fct^{-1} }
\def\tfct{ \tilde{\finc}_{\rm Crtn} }
\def\tfcti{ \tfct^{-1} }
\def\fcto{ \fct^{(1)} }
\def\fctt{ \fct^{(2)} }
\def\tfct{ \tilde{\finc}_{\rm Crtn} }
\def\tfcto{ \tfct^{(1)} }
\def\tfctt{ \tfct^{(2)} }

\def\Hov#1{ H^1(#1) }
\def\Hobv#1{ \Hov{#1,\del #1} }
\def\HobZv#1{ \Hov{#1,\del #1;\ZZ} }
\def\HobZtD{ \HobZv{\txD} }
\def\ZHobZt{ \ZZ(\HobZtD) }
\def\Hhobv#1{ H_1(#1,\del #1) }
\def\HhobDp{ \Hhobv{\xD\p} }
\def\HobD{ \Hobv{\xD} }
\def\HobtD{ \Hobv{\txD} }
\def\HobDp{ \Hobv{\xD\p} }
\def\HobDo{ \Hobv{\xD_1} }
\def\HobDt{ \Hobv{\xD_2} }

\def\HobtxDir{ \Hobv{\txDir} }
\def\HobtxDic{ \Hobv{\txDic} }
\def\sG#1{ G_{#1} }
\def\GD{ \sG{D} }
\def\GDp{ \sG{D\p} }
\def\GDo{ \sG{\xDo} }
\def\GDt{ \sG{\xDt} }
\def\GDot{ \sG{\xDot} }
\def\PGv#1{ P_{#1} }
\def\PGD{ \PGv{\GD} }
\def\PGDp{ \PGv{\GDp} }
\def\PGDot{ \PGv{\GDot} }
\def\PGDv#1{ \PGv{\sG{#1}} }

\def\SHobv#1{ S^{\ast}H^1(#1,\del #1) }
\def\SHobD{ \SHobv{D} }

\def\HCobv#1{ H^1_{\rm C}(#1) }
\def\HCobD{ \HCobv{D} }

\def\HCobDp{ \HCobv{D\p} }
\def\HCobDj{ \HCobv{D_j} }

\def\HCfobv#1{ \tilde{H}^1(#1) }
\def\HCfobD{ \HCfobv{D} }
\def\HCfobDp{ \HCfobv{D\p} }
\def\HCfobDd{ \Hobv{\txD} }
\def\HCfobDo{ \HCfobv{\xDo} }
\def\HCfobDt{ \HCfobv{\xDt} }
\def\HCfobDot{ \HCfobv{\xDot} }

\def\HCfeobv#1{ \tilde{H}_\rxp^1(#1) }
\def\rxp{{\rm exp} }
\def\HCfeobD{ \HCfeobv{\xD} }

\def\SHCfobv#1{ \Srs^{\ast}\HCfobv{#1} }
\def\SHCfobXpv#1{ (\SHCfobv{#1})^{\otimes\orXp} }
\def\SHCfobD{ \SHCfobv{D} }
\def\SHCfobDp{ \SHCfobv{D\p} }
\def\SHCfobDpj{ \SHCfobv{\Dpj} }
\def\SHCfobDj{ \SHCfobv{D_j} }
\def\SHCfobDtX{ (\SHCfobD)^{\otimes|\xX|} }
\def\SHCfobDtXp{ (\SHCfobD)^{\otimes|\Xp|} }
\def\SHCfobDjtXp{ (\SHCfobDj)^{\otimes|\Xp|} }
\def\SHCfobDptXp{ (\SHCfobDp)^{\otimes|\Xp|} }
\def\SHCfobDpjtXp{ (\SHCfobDpj)^{\otimes|\Xp|} }

\def\SHCfobvd#1{ S^{\ast}\Hobv{\tilde{#1}} }
\def\SHCfobDd{ \SHCfobvd{D} }
\def\SHCfobDtXpd{ (\SHCfobDd)^{\otimes|\Xp|} }

\def\QeSHCfobDtX{ \cQev{\SHCfobDtXp} }

\def\cQev#1{ \cQ_{\rm e}(#1) }

\def\ed{ e }
\def\edz{ e_0 }
\def\hed{ \hat{\ed} }
\def\uhed{ \hat{ \mathbf{\ed} } }
\def\fed{ f_{\ed} }

\def\edx{ \ed_{(x)} }
\def\eda{ \ed_a }
\def\heda{ \hed_a }
\def\hedza{ \hed_{0,a} }
\def\hedx{ \hed_{(x)} }

\def\mex{ m_{\edx} }
\def\mea{ m(\eda) }

\def\fir{ f_{i,\xuo} }
\def\fic{ f_{i,\xuoc} }

\def\preDr{ \prod_{\ed\in\bEtD} }

\def\hedxXp{ (\heda)_{a\in\Xp} }
\def\hedzxXp{ (\hedza)_{a\in\Xp} }

\def\pexXp{ p_{\ed}(\hedxXp)}
\def\pezxXp{ p_{\edz}(\hedzxXp)}

\def\CCv#1{{\rm CC}_#1}
\def\CCo{\CCv{1}}
\def\CCt{\CCv{2}}

\def\mF{ f_{\cH} }
\def\tmF{ \tilde{f}_{\cH} }
\def\tmFi{ \tmF^{-1} }
\def\tmFsi{ \tilde{f}_{\cH,s}^{-1} }
\def\mFi{ \mF^{-1} }
\def\mFq{ f_{\cQ} }
\def\mFiq{ f_{\cQ\cH} }
\def\mFiqi{ \mFiq^{-1} }

\def\zA{ A }

\def\blt{ \circ }
\def\blx{ \bullet }
\def\blxb{ \blx\,\bl }

\def\cHvv#1#2{ \cH_{#1}(#2) }
\def\cXpHv#1{ \cHvv{\Xp}{#1} }
\def\cXpHD{ \cXpHv{\xD} }

\def\cXpHb{ \cXpHv{\blt} }
\def\cXpHx{ \cXpHv{\blxb} }
\def\cXpHc{ \cXpHv{\rcrc} }
\def\cXpHl{ \cXpHv{\lcrc} }

\def\cHevv#1#2{ \cH^\rxp_{#1}(#2) }
\def\cXepHv#1{ \cHevv{\Xp}{#1} }
\def\cXepHD{ \cXepHv{\xD} }
\def\cLepHv#1{ \cHevv{\cL}{#1} }
\def\cLepHD{ \cLepHv{\xD} }

\def\cQXepHD{ \cQ\cXepHD }

\def\cHv#1{ \cH(#1) }
\def\cHD{ \cHv{\xD} }
\def\cHDp{ \cHv{\xD\p} }
\def\cHDo{ \cHv{\xDo} }
\def\cHDt{ \cHv{\xDt} }
\def\cHDot{ \cHv{\xDot} }
\def\cQHDp{ \cQ\cHv{D\p} }

\def\cHDo{ \cHv{D_1} }
\def\cHDt{ \cHv{D_2} }
\def\cQHDo{ \cQHv{D_1} }
\def\cQHDt{ \cQHv{D_2} }
\def\cHDr{ \cHv{\Dr} }
\def\cHDc{ \cHv{\Dc} }
\def\cHDrc{ \cHDr\oplus\cHDc }

\def\cQHv#1{ \cQ\cH(#1) }
\def\cQHD{ \cQHv{D} }
\def\cQHev#1#2{ \cQ\cH_{#1}(#2) }
\def\cQHeev#1{ \cQHev{\ed}{#1} }
\def\cQHDe{ \cQHeev{D} }
\def\cQHDje{ \cQHeev{D_j} }

\def\cQHDev#1{ \cQHeev{D_{#1}} }
\def\cQHDr{ \cQHv{\Dr} }
\def\cQHDc{ \cQHv{\Dc} }
\def\cQHDrc{ \cQHDr\oplus\cQHDc }

\def\cQHDXp{ \cQ\cXpHD }

\def\dslv#1{ #1_{\rm rt} }
\def\dslD{ \dslv{D} }

\def\gDv#1{ \mathbf{D}_{#1} }
\def\gDX{ \gDv{\xX} }
\def\gDXp{ \gDv{\xXp} }
\def\gDL{ \gDv{\cL} }

\def\cgDv#1{ \mathbf{D}\p_{#1} }

\def\cgDXp{ \cgDv{\xXp} }
\def\cgDL{ \cgDv{\cL} }

\def\tcDCv#1{ \tilde{\cD}^{\rm C}_{#1} }
\def\tcDCX{ \tcDCv{\xX} }
\def\tcDCXp{ \tcDCv{\xXp} }

\def\tcDCXv#1{ \tcDCv{\xX,{\rm #1}} }

\def\tcDCXIHX{ \tcDCXv{\IHX} }
\def\tcQDCXIHX{ \cQ\tcDCXIHX }

\def\tcQDCX{ \cQ\tcDCX }
\def\tcQDCXp{ \cQ\tcDCXp }

\def\cQD{ \cQ\cD }
\def\cDC{ \cD^{\rm C} }
\def\cQDC{ \cQ\cD^{\rm C} }
\def\cDCv#1{ \cDC_{#1} }
\def\cDCX{ \cDCv{\xX} }
\def\cDCXp{ \cDCv{\xXp} }
\def\cDCL{ \cDCv{\cL} }
\def\cQDCL{ \cQDCv{\cL} }

\def\cDCe{ \cE^{\rm C} }
\def\cDCev#1{ \cDCe_{#1} }
\def\cDCeXp{ \cDCev{\xXp} }
\def\cDCeX{ \cDCev{\xX} }
\def\cDCeL{ \cDCev{\cL} }

\def\cQDCeXp{ \cQ\cDCeXp }
\def\cQDCeX{ \cQ\cDCeX }
\def\cQDCeL{ \cQ\cDCeL }

\def\tcDCe{ \tilde{\cE}^{\rm C} }

\def\tcQDCe{ \cQ\tcDCe }
\def\tcQDCev#1{ \tcQDCe_{#1} }
\def\tcQDCeX{ \tcQDCev{\xX} }
\def\tcQDCeXp{ \tcQDCev{\xXp} }
\def\tcQDCeXpIHX{ \tcQDCev{\xX,\IHX} }

\def\cQDCv#1{ \cQ\cD^{\rm C}_{#1} }
\def\cQDCX{ \cQDCv{\xX} }
\def\cQDCXp{ \cQDCv{\xXp} }

\def\boDgDXp{\bigoplus_{D\in\gDXp} }
\def\boDgDXpp{\bigoplus_{D\in\cgDXp} }

\def\Dr{ D_{\rm r} }
\def\Dc{ D_{\rm c} }

\def\Hr{ H_{\rm (r)} }
\def\Hc{ H_{\rm (c)} }

\def\mtrc#1{ \mtr{#1}_{\rm c} }
\def\mtrr#1{ \mtr{#1}_{\rm r} }
\def\mtri#1{ \mtr{#1}^{-1} }

\def\mtrv#1#2{ \mtr{#1}_{#2} }
\def\mtrk#1{ \mtrv{#1}{k} }

\def\mtriv#1#2{ \mtri{#1}_{#2} }
\def\mtrik#1{ \mtriv{#1}{k} }

\def\VC{ V_{\rm C} }
\def\VCIHX{ V_{{\rm C},\IHX} }

\def\IHX{ {\rm IHX} }

\def\gll{ \flat }
\def\gllv#1#2#3{ \gll_{#1,#2;#3} }
\def\gllsv#1#2#3{ \;\;\gll_{#1,#2;#3}\;\; }
\def\gllxym{ \gllv{x}{y}{m} }
\def\gllxyLR{ \gllv{x}{y}{\iLR} }
\def\gllxymj{ \gllv{x_1}{x_2}{m|j} }

\def\fglj{ f_{\gll,j} }

\def\iLR{ {\rm LR} }

\def\PI{PI}

\def\lgt{ \log_\times }
\def\xpt{ \exp_\times }

\def\lrbcs#1{ \Big( #1 \Big) }
\def\lrbss#1{ \Big[ #1 \Big] }

\def\glyxLR{ \gllsv{\uy}{\ux}{\iLR} }
\def\glyyLR{ \gllsv{\uy}{\uy}{\iLR} }
\def\glxyLR{ \gllsv{\ux}{\uy}{\iLR} }
\def\glxxLR{ \gllsv{\ux}{\ux}{\iLR} }
\def\glxxLRo{ \gllsv{x}{x}{\iLR} }
\def\glxxkLR{ \gllsv{\uxrvk}{\uxrvk}{\iLR} }
\def\glxxkLRb{ \gllv{\uxrvk}{\uxrvk}{\iLR} }
\def\gluyo{ \gllsv{\uy}{\uy}{1} }
\def\gluzo{ \gllsv{\uz}{\uz}{1} }

\def\glmdy{ \gllsv{d\uy}{\dluy}{1} }
\def\glxxo{ \gllsv{\ux}{\ux}{1} }

\def\Dpj{ D^\prime_j }

\def\bS{ \mathbf{S} }
\def\bSv#1{ \bS_{#1} }
\def\bSD{ \bSv{D} }

\def\bE{ \mathbf{E} }
\def\bV{ \mathbf{V} }
\def\bEv#1{ \bE_{#1} }
\def\bErv#1{ \bEv{{\rm r},#1} }
\def\bEDr{ \bErv{\xD} }
\def\bED{ \bEv{\xD} }
\def\bEtD{ \bEnv{\txD} }

\def\bEnv#1{ \bE(#1) }
\def\bED{ \bEnv{\xD} }

\def\Env#1{ E(#1) }
\def\ED{ \Env{\xD} }
\def\EtD{ \Env{\txD} }

\def\SEtD{ \Srs^*\EtD }
\def\SEtDX{ (\SEtD)^{\otimes\orXp} }

\def\SEtDG{ \Invrs{\SEtDX}{\GD} }

\def\Evnv#1{ E_{\rm v}(#1) }
\def\EvD{ \Evnv{\xD} }
\def\EvtD{ \Evnv{\txD} }
\def\IvtD{ I_{\rm v}(\txD) }

\def\bVv#1#2{ \bV_{#1}(#2) }
\def\bVmv#1{ \bVv{m}{#1} }
\def\bVmD{ \bVmv{\xD} }

\def\Dg{ \Delta_{\mfg} }

\def\Dpj{ D^{\prime}_j }

\def\xc{ c }

\def\hc{ \hat{c} }

\def\otorXp{ {\otimes\orXp} }

\def\ShoXp{ (\Sh)^\otorXp }
\def\QShoXp{ \Qf{\ShoXp} }

\def\mfrl{ \mfr_{\vl} }
\def\mfrml{ \mfr_{-\vl} }

\def\Xp{ \xXc\p }
\def\xXp{ \xX,\Xp }

\def\prc#1#2{ (#1,#2) }
\def\prcox{ \prc{\omega}{\xi} }

\def\uF{ F }
\def\uG{ G }
\def\uH{ H }
\def\uK{ K }

\def\uz{ \u{z} }

\def\lji{ l_{ji} }
\def\lij{ l_{ij} }
\def\lsij{ \stl_{ij} }

\def\liij{ l^{-1}_{ij} }

\def\lijua{ \lij(\ua) }

\def\lin{ _{\rm lin} }
\def\str{ _{\rm str} }

\def\uFaff{ \uF\str  }
\def\uFlin{ \uF\lin  }

\def\ltiv#1#2{ #1^{(#2)} }
\def\ltiH#1{ \ltiv{\uH}{#1} }
\def\ltiK#1{ \ltiv{\uK}{#1} }

\def\mtr#1{ || #1 || }
\def\mtrrv#1#2{ \mtr{#1}\rv{#2} }
\def\rstrt{ \smile }
\def\tstrt{ \smile }

\def\rstv#1#2{ \xymatrix{ #1 \ar@{--}[r]& #2} }
\def\rstxot{ \rstv{x_1}{x_2} }
\def\rstx{ \rstv{x}{x} }

\def\tstv#1#2{ \xymatrix{ #1 \ar@{-}[r]& #2} }

\def\cstv#1#2{  \xymatrix{#1 \ar@{.}[r] & #2} }

\def\rstyxij{ \rstv{y_i}{x_j} }
\def\cstyxij{ \cstv{y_i}{x_j} }
\def\tstyxij{ \tstv{y_i}{x_j} }
\def\tstyxj{ \tstv{y_j}{x_j} }
\def\tstuxy{ \tstv{\uy}{\ux} }

\def\Id{ {\rm Id} }

\def\Diff{ {\rm Diff} }
\def\DiffxN{ \Diff_{\xN} }
\def\DiffxNn{ \Diff_{\xN}^{\rm nrw} }

\def\Vect{ {\rm Vect} }
\def\VectxN{ \Vect_{\xN} }

\def\Ldv#1{ L_{#1} }
\def\Ldxi{ \Ldv{\xi} }
\def\LdXi{ \Ldv{\Xi} }
\def\Ldxia{ \Ldv{\xia} }
\def\Ldmada{ \Ldv{\mad_a} }
\def\Lmadz{ \Ldv{\madz} }

\def\dvrg{ {\rm div} }
\def\dvrgv#1{ \dvrg_{#1}\, }
\def\dvxxi{ \dvrgv{\ux} \xi }
\def\dvxXi{ \dvrgv{\ux} \Xi }

\def\dvxxia{ \dvrgv{x} \xia }
\def\dvxmada{ \dvrgv{x} (\mad_a) }

\def\txi{ \tilde{\xi} }

\def\attz#1{ \atv{#1}{t=0} }

\def\SP{SP}
\def\SPN{SPND}
\def\nrw{narrow}

\def\pri#1#2{ {}^{#1}_{#2} }
\def\pryx{ \pri{y}{x} }
\def\prxy{ \pri{x}{y} }
\def\prxx{ \pri{x}{x} }
\def\prxxM{ \prxx M }
\def\przyM{ \pri{z}{y}M }
\def\pryxM{ \pri{y}{x}M }
\def\yxFp{ \pryx F\p }
\def\xyFip{ \prxy (F^{-1})\p }
\def\pMv#1#2{ {}^{#1}_{#2} M }
\def\yxM{ \pMv{y}{x} }
\def\xxM{ \pMv{x}{x} }
\def\xxMp{ \xxM\p }
\def\qfv#1{ (#1)_{\rm qf} }
\def\Mqf{ \qfv{M} }
\def\xxaM{ M }
\def\Mij{ M^i_j }
\def\Mv#1#2#3#4{ {}^{#1}_{#2} M^{#3}_{#4} }
\def\Fpv#1#2#3#4{ {}^{#1}_{#2} (F\p)^{#3}_{#4} }
\def\yxMij{ \Mv{y}{x}{i}{j} }
\def\xxMij{ \Mv{x}{x}{i}{j} }
\def\yxFpij{ \Fpv{y}{x}{i}{j} }

\def\istr{ _{\rm str} }
\def\ostr{ _{\rm 1,str} }
\def\tstr{ _{\rm 2,str} }

\def\xxaMstr{ \xxaM\istr }

\def\xxaMtstr{ \xxaM\tstr }
\def\xxaMostr{ \xxaM\ostr }

\def\gxG{ G }
\def\gxQ{ Q }
\def\gxP{ P }

\def\gxGp{ \gxG\p }
\def\gxQp{ \gxQ\p }
\def\gxPp{ \gxP\p }

\def\gxQi{ \gxQ^{-1} }
\def\hgxQi{ \hgxQ^{-1} }

\def\hgxQ{ \hat{\gxQ} }

\def\gxQLxa{ \gxQ(\cL;\uxz) }

\def\hQkua{ \hgxQ_k(\cL) }

\def\Qstkss{ \gxQ_{su(2),k}(\cL;\ua) }

\def\FGi{ \int_{\rm FG} }

\def\iexp{ ^{\rm (exp)} }
\def\Fiexp{ F\iexp }

\def\eQix{ e^{-\hlf\gxQi(\ux)} }
\def\eQixo{ e^{-\hlf\gxQi(x)} }

\def\eQi{ e^{-\hlf\gxQi} }
\def\eQx{ e^{\gxQ(\ux)} }
\def\eQxt{ e^{\gxQ(\ux;t)} }
\def\eQx{ e^{ \hlf\,\gxQ(\ux) } }
\def\eQ{ e^{\hlf\gxQ} }

\def\sFst{ F\drst }
\def\drst{^{\ast}}

\def\crr{ x }
\def\cra{ a }

\def\xia{ \xi_\cra }

\def\Qpr{$\IQ$-parameter}

\def\dlt{ \del_t }

\def\BCH{ {\rm E} }
\def\Ev#1#2{ \BCH(#1,#2) }
\def\Esv#1{ \BCH(#1) }
\def\Esvb#1{ \BCH\Big(#1\Big) }
\def\BCHr{ \BCH\xdr }
\def\BCHc{ \BCH\xdc }
\def\Erv#1#2{ \BCHr(#1,#2) }
\def\Ecv#1#2{ \BCHc(#1,#2) }

\def\SPL{ {\rm S} }
\def\Sv#1#2{ \SPL(#1,#2) }
\def\Sr{ \SPL\xdr }
\def\Sc{ \SPL\xdc }
\def\Srv#1#2{ \Sr(#1,#2) }
\def\Scv#1#2{ \Sc(#1,#2) }

\def\Fr{ F\xdr }
\def\Fc{ F\xdc }

\def\xdr{ _{\rm r} }
\def\xdc{ _{\rm c} }

\def\tBCH{ \tilde{\BCH} }
\def\tEv#1#2{ \tBCH(#1,#2) }

\def\bfx{ X }
\def\bfy{ Y }
\def\bfz{ Z }
\def\bfu{ U }
\def\bfw{ W }
\def\bfm{ \mathbf{m} }
\def\bfn{ \mathbf{n} }

\def\ubfx{ \mathbf{\bfx} }

\def\bfun{ \mathbf{1} }

\def\mad{ \vec{\ad} }
\def\madz{ \mad_z }
\def\mady{ \mad_y }

\def\mmad{ \widehat{\ad} }

\def\Mrx{ M_{{\rm r},x} }
\def\Mrmx{ M_{{\rm r},-x} }

\def\tmad{ \widetilde{\ad} }
\def\tmad{ \ad }
\def\tmadz{ \tmadv{z} }
\def\tmadv#1{ \tmad_{#1} }

\def\mAdv#1{ \Ad_{#1} }
\def\mAdev#1{ \mAdv{\exp(#1)} }
\def\mAdesv#1{ \mAdv{\exp\big(#1\big)} }
\def\mAdz{ \mAdev{z} }
\def\mAdtz{ \mAdev{tz} }
\def\mAda{ \mAdev{a} }
\def\mAdta{ \mAdev{ta} }

\def\xpLie{ \exp_{\rm Lie} }

\def\DRl{ {\rm R} }
\def\DRv#1{ \DRl_{#1} }
\def\DRax{ \DRa(x) }
\def\DRa{ \DRv{a} }
\def\DRav#1{ \DRl_{a_{#1}} }

\def\vDRl{ \DRl }
\def\vDRv#1{ \vDRl_{#1} }
\def\vDRa{ \vDRv{\va} }
\def\vDRav#1{ \vDRl_{\va_{#1}} }
\def\vDRax{ \vDRa(\vx) }

\def\tG{ \tilde{G} }

\def\zchA{ \check{\zA} }
\def\Gfz{ \xXi_\bfz }
\def\tGfz{ \tGf_\bfz }
\def\tGf{ \tilde{\xXi} }

\def\iromx{ ^{\rm R} }
\def\irom{ ^{\cQ{\rm R}} }
\def\irorW{ ^{\Omega,{\rm r} } }
\def\iromFD{ ^{\rm s,FD} }
\def\iromW{ ^{\Omega,\cQ{\rm R}} }

\def\Gr{ G\irom }

\def\Grua{ \Gr }
\def\frua{ f\irom }

\def\Gruax{ G\iromx }

\def\uxz{ \ux }

\def\dxtxN {dx_3\cdots dx_\xN}
\def\jxtxN { \msjv{x_3}\cdots \msjv{x_\xN} }

\def\cLSt{ \cL\in S^3 }
\def\cLM{ \cL\in M }

\def\spx#1{ x^{(#1)} }
\def\spxi{ \spx{i} }

\def\xl{ l }
\def\xlij{ \xl_{ij} }
\def\xlijt{ \xl_{ij|2} }
\def\xxl{ \ell }
\def\xxlij{ \xxl_{ij} }

\def\nDp{ |\Delta_+| }

\def\mmrcr{ (-2\pi)^{-\nDp} }
\def\mfrvl{ \mfr_{\vl} }
\def\mfrvlj{ \mfr_{\vl,j} }
\def\TgMvl{ (\Tg M)_{\vl} }

\def\pvlDg{ \prod_{\vl\in \Dg } }
\def\svlDg{ \sum_{\vl\in \Dg} }

\def\vdfm#1{ \vec{#1} }
\def\vdfmw#1{ \widevec{#1} }
\def\vF{ \vdfm{F} }

\def\xiFv#1{ \xi_{{\rm F},#1} }
\def\xiFy{ \xiFv{y} }
\def\bxiFy{ \boldsymbol{\xi}_{{\rm F},y} }

\def\detp{ \det\nolimits\p }

\def\stv#1{ {#1}^{\ast} }
\def\stl{ \stv{l} }

\def\bcup{ \cup }
\def\bcap{ \cap }

\def\pja{ \pjoxN a_j }
\def\pjba{ \lrbc{ \pja } }
\def\pjbai{ \pjba^{-1} }

\def\Fso{ F_{o} }
\def\Fsuo{ F_{\yuo} }

\def\dgsv#1{ \deg_{#1} }
\def\dgsj{ \dgsv{(j)} }
\def\xPsv#1{ P_{#1} }
\def\xPsj{ \xPsv{j} }

\def\ndgap{ \AFLut\not\equiv 0 }

\def\xrov#1{ _{\{#1\}} }
\def\xroj{ \xrov{j} }
\def\cLroj{ \cL\xroj }

\def\bcL{ \bar{\cL} }

\def\mpcd#1#2#3{ \begin{CD} {#1} @>{#2}>> {#3} \end{CD} }

\def\xFDL{ F_{\xD}(\cL) }
\def\xtFDL{ \tilde{F}_{\xD}(\cL) }
\def\xtFODL{ \tilde{F}^{\Omega}_{\xD}(\cL) }
\def\xFODL{ F^\Wh_{\xD}(\cL) }
\def\xFzLua{ F_{0}(\cL;\ua) }
\def\xFWzLue{ F^\Wh_{0}(\cL;\uhed) }
\def\xFWzLav{ F^\Wh_{0}(\cL;\scp{\vua}{\vl}) }
\def\hxFzLua{ \hlf\,\xFzLua }
\def\dtplijua{ \detp\mtr{\lijua} }
\def\dtplijue{ \detp\mtr{\lij(\uhed)} }

\def\xFWzLua{ F^\Wh_{0}(\cL;\ua) }
\def\hxFWzLua{ \hlf\,\xFWzLua }

\def\xslija{ \sum_{1\leq j\leq \xN \atop j\neq i} \xxlij a_j }
\def\hxslija{ \hlf \xslija }

\def\lknb{ \Wh_{\nabla} }
\def\lknv#1#2{ \lknb(#1;#2) }
\def\lknLa{ \lknv{\cL}{\ua} }

\def\saE{ E }
\def\saF{ F }

\def\saPp{ \Psi^+ }
\def\saPm{ \Psi^- }
\def\saPpm{ \Psi^{\pm} }

\def\sap{ \hat{\psi}
 }
\def\sapp{ \sap_+ }
\def\sapm{ \sap_- }

\def\savp{ \psi }
\def\savpp{ \savp_+ }
\def\savpm{ \savp_- }

\def\sapr{ \hat{\psi}^\prime
 }
\def\saprp{ \sapr_+ }
\def\saprm{ \sapr_- }

\def\savpv#1{ \savp_{#1} }
\def\savppj{ \savpv{+,j} }

\def\savppi{ \savpv{+,i} }
\def\savpmi{ \savpv{-,i} }
\def\savppk{ \savpv{+,k} }
\def\savpmk{ \savpv{-,k} }

\def\usapm{ \upsi_- }
\def\usapp{ \upsi_+ }

\def\shlfv#1{ {\sinh(#1/2) \over (#1/2) } }
\def\shlfph{ \shlfv{\phi} }

\def\xfB{ f_{\rm B} }

\def\STrVa{ \STr_{V_a} }

\def\uoo{ u(1|1) }
\def\Uoo{ U(1|1) }

\def\frm{ {\rm F} }
\def\bsn{ {\rm B} }

\def\xDF{ \xD_\frm }

\def\mfgb{ \mfg_\bsn }
\def\mfgf{ \mfg_\frm }

\def\himfgb{ h_{\mfg,\bsn}^{-1} }
\def\himfgf{ h_{\mfg,\frm}^{-1} }

\def\Vbsn{ V_\bsn }
\def\Vfrm{ V_\frm }

\def\TgDF{ \Tg(\xD,\xDF) }

\def\nor{ {n_{\rm cycl}} }
\def\norl{ {n_{\rm lin}} }
\def\neDb{ n_{\bsn} }
\def\neDf{ n_{\frm} }

\def\gra{ a }

\def\lpijua{ l\p_{ij}(\ua) }

\def\vgrav#1{ |\gra,#1\rangle }
\def\vgrap{ \vgrav{+\hlf} }
\def\vgram{ \vgrav{-\hlf} }
\def\vgrapm{ \vgrav{\pm\hlf} }

\def\STrvuak{ \STr_{\TVuva;k} }
\def\STrv#1{ \STr_{#1} }
\def\STrsvuak{ \STrv{\uval;k} }
\def\STrvua{ \STrv{\TVuva} }
\def\STrua{ \STrv{\TVua} }
\def\STrsuak{ \STrv{\ua;k} }
\def\STrV{ \STr_{V} }
\def\STrVav#1{ \STr_{V_{a_{#1}}} }
\def\STrVai{ \STrVav{i} }
\def\STrVak{ \STrVav{k} }

\def\pO{ O }

\def\pOC{ \pO }
\def\pOCk{ \pO_{k} }
\def\pOk{ \pO_k }

\def\eaphe{ e^{{a\over\phi}\,\etm\etp} }

\def\gxf{ \phi }

\def\dgxfv#1{ \del_{\gxf_{#1}} }
\def\dgxfj{ \dgxfv{j} }
\def\etp{ \eta^+ }

\def\detpv#1{ \del_{\etp_{#1}} }
\def\detpj{ \detpv{j} }
\def\etm{ \eta^- }

\def\detmv#1{ \del_{\etm_{#1}} }
\def\detmi{ \detmv{i} }

\def\lijQ{ \lij\in\left\{
\begin{array}{cl}
\QIQuua, &\mbox{if $1\leq i,j\leq \xNp$}
\\
\IQ, &\mbox{if $\xNp\leq i,j\leq \xN$}
\end{array}
\right.
}

\def\Vuoa{ V_a }

\def\Suoo{ \Srs\,\uoo }
\def\SuootN{ (\Suoo)^{\otimes \xN} }
\def\Uuoo{ \Urs\,\uoo }

\def\Srs{ {\rm S} }
\def\Urs{ {\rm U} }

\def\xgm{ \gamma }
\def\xgmg{ \xgm_{\mfg} }

\def\Dcst{ \mathfrak{D}_{\rm const} }
\def\Dcstg{ \Dcst(\mfg) }

\def\vY{ \vec{Y} }

\def\buoo{ \b_{\uoo} }

\def\rT{ {\rm T} }

\def\vd{ v }
\def\inav#1{ a_{#1} }
\def\inave{ \inav{\vd,\ed} }

\def\spanv#1{ \span\lrbc{#1} }

\def\bopgD{ \bigoplus_{\xD\in\gDXp} }

\def\bl{ b }
\def\Xpfr{$\Xp$-frame}

\def\xxs{ s }

\def\xfo{ f_1 }
\def\xft{ f_2 }

\def\IQuva{ \IQ[[\uva]] }
\def\IQuvark{ \IQuva[\uxvrk] }
\def\IQQuvahb{ \Qf{\IQ[[\uva]]}[[\hb]] }

\def\ual{ \boldsymbol{\a} }
\def\ueta{ \boldsymbol{\eta} }
\def\upsi{ \boldsymbol{\psi} }
\def\uval{ \vec{\ual} }
\def\uvbet{ \vec{\ubet} }

\def\uvb{ \vec{\mathbf{b}} }
\def\ua{ \mathbf{a} }
\def\yuo{ \mathbf{o} }
\def\uva{ \vec{\ua} }
\def\ux{ \boldsymbol{x} }
\def\uz{ \boldsymbol{z} }
\def\ut{ \boldsymbol{t} }
\def\vux{ \vec{\ux} }
\def\uy{ \boldsymbol{y} }
\def\vuy{ \vec{\uy} }
\def\uxvr{ \vec{\mathbf{x}} }

\def\bD{ \mathbf{D} }

\def\zA{ T }
\def\vb{ \vec{\bfx} }
\def\uvb{ \vec{\ubfx} }
\def\vbet{ \vec{\chi} }
\def\uvbet{ \vec{\boldsymbol{\chi}} }

\def\orcirc{
\thinlines
\unitlength=1mm
\begin{picture}(4,3)
\put(1.5,1.2){\circle{4}}
\put(3.5,2){\vector(0,1){0}}
\end{picture}
}

\def\rcrc{\begin{picture}(450,226)(-105,60)\epsfig{file=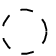}\end{picture}}
\def\lcrc{\begin{picture}(450,226)(-105,60)\epsfig{file=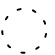}\end{picture}}
\def\tcrc{\begin{picture}(450,226)(-105,60)\epsfig{file=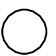}\end{picture}}

\begin{document}

\setlength{\unitlength}{3947sp}

\begin{titlepage}
\vfill
\begin{center}

{\large \bf A universal \urcc invariant of links and rationality conjecture.}\\

\bigskip

\bigskip
\centerline{L. Rozansky\footnote{
This work was supported by NSF Grants DMS-0196235 and DMS-0196131}
}



\centerline{\em Department of Mathematics, University of North Carolina}
\centerline{\em CB \#3250, Phillips Hall}
\centerline{\em Chapel Hill, NC 27599}
\centerline{{\em E-mail address:} {\tt rozansky@math.unc.edu}}

\vfill
{\bf Abstract}

\end{center}
\begin{quotation}

We define a graph algebra version of the stationary phase integration over
the coadjoint orbits in the Reshetikhin formula for the colored
Jones-HOMFLY polynomial. As a result, we obtain a `universal' \urcc
invariant of links in rational homology spheres, which determines the
\urcc invariants based on simple Lie algebras. We formulate a rationality
conjecture about the structure of this invariant.

\end{quotation} \vfill \end{titlepage}

\pagebreak

\tableofcontents

\nsection{Introduction}
\label{s1}
\hyphenation{Re-she-ti-khin}
\hyphenation{Tu-ra-ev}
\hyphenation{sub-stan-ti-al}
\hyphenation{in-va-ri-ant}

\subsection{Motivation}
\label{ss1.1}




A classical topology interpretation of `quantum' invariants of knots, links and
3-manifolds still remains an open question. A discovery of the
Alexander polynomial inside the colored Jones polynomial as well as a
discovery of the Casson-Walker invariant inside the
Witten-Reshetikhin-Turaev invariant of rational homology spheres suggests
that quantum invariants are somehow packed up with classical invariants
(both known and unknown). Thus the appropriate unpacking of quantum
invariants may lead to a discovery of new classical invariants which may,
in turn, provide a solution to long-standing puzzles of 3-dimensional
topology.

The main method of unpacking the quantum invariants is splitting them into
invariants of finite type. For rational homology spheres, this procedure seems
to correspond to the isolation of the contribution of the trivial connection to
the Witten-Reshetikhin-Turaev invariant in the semi-classical limit of
$\limone{q}$. It turns out that the best way of implementing the same procedure
for the Jones polynomial of knots and links is
\emph{not} to study the trivial
connection contribution in the limit $\limone{q}$,
$\a_{1}\ldtc\a_{\xN}=\const$ (which would
amount to simply expanding the colored Jones polynomial in powers of $q-1$),
but rather to isolate
a contribution of the $U(1)$-reducible flat connections
in
the knot or link complement
(we abbreviate this contribution as
$U(1)$-RCC) to the colored Jones polynomial in the limit
\qq
\limone{q},\qquad
\qali{1}\ldtc\qali{\xN}=\const,
\label{1.1}
\qqq
where $\a_1\ldtc\a_N$ are colors (dimensions of $SU(2)$ modules) assigned to
the components of an $\xN$-component link.
This approach yields the following formula\cx{Ro0}
for the colored Jones polynomial of
a knot $\kst$:
%
\qq
\JaK = {[\a]\over \APKqa}
\lrbc{ 1 + \snoi  {\PnKqa \over \APKpqa{2n} }\, (q-1)^n
}.
\label{1.2}
\qqq
Here $[\a]= (\pdhqa) / (\pdhq)$, $\APKt$ is the Alexander polynomial of $\cK$
normalized by the condition $\APv{\unknot}{t}=1$ and
\qq
\PnKt \in \Zti
\label{1.3}
\qqq
are polynomial invariants of $\cK$ which, due to
the manner of their appearance in \ex{1.2},
may have a classical topology interpretation similar to that of the Alexander
polynomial and Casson-Walker invariant. The formula\rx{1.2} itself should be
interpreted as an equation between the expansion of both sides in powers of
$q-1$ and $\a$ (known as the Melvin-Morton expansion of the colored Jones
polynomial).

Since the expansion of $\JaK$ in powers of $q-1$ can be deduced from the
Kontsevich integral of $\cK$, then equation\rx{1.2} implies that
Kontsevich integral has a special `rational' structure which matches the
structure of the \rhs of \ex{1.2}. The corresponding conjecture about
Kontsevich integral was formulated in\cx{RoC} and A.~Kricker
proved it in\cx{Kr} (a stronger version was proved by S.~Garoufalidis
and A.~Kricker in\cx{GaKr}). The conjecture states that the coefficients at
individual graphs in the Kontsevich integral of $\cK$ can be assembled into
`universal' polynomials which determine the polynomials $\PnKt$ of expansions
of the type\rx{1.2} for colored Jones-HOMFLY polynomials based on any simple
Lie algebra.

The case of links is less simple than that of knots. On one
hand, S.~Garoufalidis and A.~Kricker\cx{GaKr} proved the rationality of
Kontsevich integral for boundary links. On the other hand,
the $R$-matrix
calculation, which led to \ex{1.2} for knots,
can be repeated for the opposite class of `sufficiently connected' links
(that is, links whose Alexander polynomial is not identically equal to 0).
This $R$-matrix calculation yields a \urcc
invariant of an oriented link $\cL$, which has the form
%
\qq
\JcrLth = {1 \over h}\,
{1\over \AFLtN}
\lrbc{1 +
\snoi {\PnLtN \over \AFLptN{2n} }\,h^n
},
\label{1.4}
\qqq
where $\AFLtN$ is the Alexander `polynomial' of $\cL$ in the normalization
$$\AFv{\unknot}{t} = {1 \over (\pdht)},$$
which is appropriate for links, while
\qq
\PnLtN\in\ZZtNh
\qqq
are polynomial invariants of $\cL$. Path integral arguments suggest that
the expression
\qq
q^{ {1\over 4}
\sijoxN
\lkn{i}{j} \a_i \a_j }
\JcrLqaN
\label{1.5}
\qqq
represents a contribution of one of $2^{\xN-1}$ $U(1)$-reducible flat
connections in $\stml$ with appropriate monodromies to the colored Jones
polynomial $\JaNL$ in the semi-classical limit\rx{1.1}. In
contrast to the knot case, the
expression\rx{1.5} does not determine the whole Jones polynomial through
a simple relation similar
to\rx{1.2}\footnote{A knot relation\rx{1.2} is rather an
exception which, from the path integral point of view, is due to the fact
that a $U(1)$-reducible connection is the only flat connection in the knot
complement when the monodromy along the meridian of
the knot is sufficiently small.}. Also, deriving the series\rx{1.4} from
the Kontsevich integral of $\cL$ requires a few extra steps outlined
in\cx{Ro2}:
%
\qq
\mbox{link $\cL$}
\mapright{1}
\begin{array}{c}
\mbox{Kontsevich} \\
\mbox{intergral} \\
\mbox{of $\cL$ in $\cB$}
\end{array}
\mapright{2}
\begin{array}{c}
\mbox{Wheeled} \\
\mbox{Kontsevich} \\
\mbox{intergral} \\
\mbox{of $\cL$ in $\cB$}
\end{array}
\mapright{3}
\begin{array}{c}
\mbox{Reshetikhin} \\
\mbox{integrand} \\
\mbox{for $\JaNL$}
\end{array}
\mapright{4}
\begin{array}{c}
\mbox{\urcc} \\
\mbox{invariant}\\
\mbox{of $\cL$}
\end{array}
\label{1.6}
\qqq
Here $\cB$ is the space of (1,3)-valent graphs with cyclic ordering at
3-valent vertices, and 1-valent vertices labeled by (that is, assigned
to) components of $\cL$. Map~1 is a standard combination of
Kontsevich integral (which maps $\cL$ into the space $\cA$) and a PBW
symmetrization map (which identifies $\cA$ with $\cB$). Map 2 is the
wheeling map defined in\cx{Wh} (it maps $\cB$ into $\cB$). Map 3 is a link
version of the application of the `wheeled' $su(2)$ weight system. Instead
of simply attaching the shifted highest weights of $su(2)$ modules to the
1-valent vertices of the graphs (as we did in the case of knots), we also
have to perform a Kirillov-type integration over their (co-)adjoint orbits
in order to recover the colored Jones polynomial. In the
limit\rx{1.1}, the stationary phase approximation may be applied
to this integral. Map 4 picks up a contribution of a particular stationary
phase point corresponding to all integrands belonging to the same Cartan
subalgebra.

The goal of this paper is to remove the Lie algeba from the maps 3,4. In
other words, we will describe the stationary phase integration in
Reshetikhin integral purely in diagrammatic terms. This integration will
become a map from the gaussian expressions in $\cB$ into a new space
$\cQD$. The image of $\cL$ in $\cQD$ will be a `universal' \urcc invariant
in the following sense: for any simple Lie algebra $\mfg$ and for any
assignment of its modules to components of $\cL$,
the $\mfg$-based \urcc invariant can be obtained by applying an appropriate
weight system to the image of $\cL$ in $\cQD$, so that instead
of the sequence\rx{1.6} we may use a new one
\qq
\mbox{link $\cL$}
\mapright{1}
\begin{array}{c}
\mbox{Kontsevich} \\
\mbox{intergral} \\
\mbox{of $\cL$ in $\cB$}
\end{array}
\mapright{2}
\begin{array}{c}
\mbox{Wheeled} \\
\mbox{Kontsevich} \\
\mbox{intergral} \\
\mbox{of $\cL$ in $\cB$}
\end{array}
\mapright{5}
\begin{array}{c}
\mbox{Universal} \\
\mbox{\urcc} \\
\mbox{invariant} \\
\mbox{of $\cL$ in $\cQD$}
\end{array}
\mapright{6}
\begin{array}{c}
\mbox{\urcc invariant} \\
\mbox{related to} \\
\mbox{the $\mfg$-Jones}\\
\mbox{polynomial of $\cL$}
\end{array}
\label{1.7}
\qqq
where map 5 is the diagrammatic stationary phase integration over the
(co-)adjoint orbits
and
map 6 is the application of $\mfg$ weight system to an element of $\cQD$.

\subsection{Results}
\label{ss1.2}
Here is the outline
of the paper. In Section\rw{xs2} we define the Lie algebra
based \urcc invariant by performing the stationary phase integration in a
familiar context of integration over the coadjoint orbits. In
Sections\rw{xs3}--\rw{xs5} we translate this calculation into the pure language
of graphs and then in Section\rw{xs6} we prove the topological invaraince of
graph $\urcc$ invariant.

We begin Section\rw{xs2} by recalling the basic facts about a
Kontsevich integral of a link in $S^3$ and in a rational homology sphere. We
look at the relation between Kontsevich integral and the Melvin-Morton
expansion of the Jones-HOMFLY polynomial. We modify it by using Kirillov
integrals instead of traces in representations and thus come to the Reshetikhin
formula for the Melvin-Morton expansion. By rescaling the integration
variables, we make the Reshetikhin integral suitable for the stationary phase
approximation and define the $\mfg$-based \urcc invariant as a contribution of
a particular stationary phase point (when all the variables belong to the same
Cartan subalgebra).

In Section\rw{xs3} we introduce a new graph algebra in which the edges are
split into their root and Cartan parts. Then we encode Cartan legs of a
graph as elements of that graph's cohomology. Thus Cartan legs are converted
into a symmetric algebra of the graph cohomology space. We extend that
symmetric algebra by allowing
division by polynomials depending on individual edges.
Finally, we formulate basic graph operations of disconnected union and leg
gluing in terms of the cohomology.

In Section\rw{xs4} we define the basic elements of differential geometry (such
as functions, tensors, vector fields and matrix fields) on our graph algebras.
We also define a determinant of a matrix field and a formal gaussian integral
and establish their properties.

In Section\rw{xs5} we concentrate on the differential geometry of the graph
algebra version of a coadjoint orbit. We find an invariant measure on the orbit
and study the properties of related gaussian integrals.

In Section\rw{xs6} we apply all this knowledge to define a universal \urcc
invariant of an oriented link as an integral over the `graph' coadjoint orbits.
We describe the structure of this invariant and conjecture its rationality.

The main results of the paper are Theorems \rw{t6.6} and\rw{t6.7}, which
establish the topological invariance of the \urcc invariant, Theorem\rw{t6.9}
which relates the universal \urcc invariant to its Lie algebra based
analog, and
the rationality conjecture\rw{6c.1}, which suggests that the $\urcc$ invariant
is a source of interesting polynomial invariants of links, possibly
with a  nice
topological interpretation. As a by-product, we establish Theorem\rw{t6.11}
which expresses the Alexander polynomial of a link in terms of a tree and
1-loop parts of its Kontsevich integral. We give an alternative proof of this
theorem in Appendix\rw{A1}. Finally, while considering the simplest
example of the graph coadjoint orbit integral, we conjecture a graph
version of the Duistermaat-Heckmann theorem (Conjecture\rw{cdh}).

\subsection{Notations}
\label{ss1.3}

Throughout the paper we will work with (1,3)-valent graphs. We assume that
every graph is endowed with a cyclic order of edges at 3-valent vertices. If
that order changes at one vertex, then we assume that the whole graph acquires
a minus sign. We refer to the egdes incident to 1-valent vertices as
\emph{legs} and we call the edges incident to two 3-valent vertices
\emph{internal}. A graph which consists of a single edge, is called a
\emph{strut}.

An edge of a graph is called a \emph{bridge}, if its removal increases the
number of connected components of the graph. A sequence of edges is called a
\emph{path}, if the beginning of the next edge coincides with the end of the
previous edge. A path without self-intersections is called a \emph{chain}. A
closed chain is called a \emph{cycle}. A \emph{haircomb} is a graph consisting
of an open chain and legs attached to it. A \emph{wheel} is a graph consisting
of a cycle and legs attached to it.

We call an element of a graph algebra \emph{narrow} if it is a
linear combination of connected graphs.

We use multi-index notations: $\ux=x_1\ldtc x_{\xN}$, where in topological
applications $\xN$ is usually a number of components of a link $\cL$.
$\ux\rv{i}$ denotes the same sequence $x_1\ldtc x_{\xN}$ but with $x_i$
removed.
We also use the following notations for exponentials and scalar products
\qq
e^{\ux} = e^{x_1}\ldtc e^{x_{\xN}},\qquad
\scp{\ux}{y} = \scp{x_1}{y}\ldtc\scp{x_{\xN}}{y}.
\label{1.8x}
\qqq

$\xG$ is a simple Lie group, $\mfg$ is its Lie algebra. $\xT\subset\xG$ is a
maximal torus and $\mfh\subset\mfg$ is a corresponding Cartan subalgebra.
$\mfgs$ and $\mfhs$ are dual spaces of $\mfg$ and $\mfh$. Lie algebra splits
into a sum $\mfg = \mfh \oplus \mfr$, where $\mfr$ is the span of the root
spaces.
We assume that $\mfg$
has an invariant scalar product (normalized in such a way that the length of
short roots is $\sqrt{2}$), so we do not always distinguish $\mfg$ from
$\mfgs$. We use vector notations $\vx$ for their elements.
$\Dg\subset\mfhs$ denotes the set of all roots of $\mfg$,
$\Delp\subset\mfhs$ denotes the set of positive roots of $\mfg$, and we
denote $\vrho = \svldp\vl$. Sometimes we will assume that
$\Dg,\Delp\subset\mfh$, root being transferred from $\mfhs$ to $\mfh$ with the
help of the scalar product. $\Va$ denotes a module of $\mfg$ with highest
weight $\val-\vrho$.

$\cL$ denotes an $\xN$-component link either in $S^3$ or in a rational
homology sphere $M$. The numbers $\xxlij$ denote the linking numbers
between its components.

For a ring $R$, $\Qf{R}$ denotes the field of its fractions.

\nsection{A Lie algebra based \urcc invariant}
\label{xs2}
\subsection{Kontsevich integral and colored Jones-HOMFLY polynomial}
\label{2xs.1}

Let us recall some basic facts about Kontsevich integral of links (we refer
the reader to\cx{A2} and references therein). We will work with links in
rational homology spheres, but we start with a link in $S^3$.
Let $\cL$ be an
$\xN$-component oriented \emph{dotted Morse link}.
In other words, we assume that $\cL$ is
imbedded in $\IR^3$, and a dot is marked on each component. Then Kontsevich
integral $\ZKAL$ of $\cL$ takes value in the space $\cAL$. This is a
(factored over the STU relations) space of (1,3)-valent graphs,
whose 1-valent vertices are placed on $\xN$ oriented segments, which
correspond to the link components. The paper\cx{A2} uses the notation
$\cAX$, where $\xX$ is a finite set of labels,
in our case $X$ obviously denotes
the set of link components, so
\qq
\xX=\{1\ldtc \xN\};
\label{1.11z}
\qqq
we
will abuse notations by
using $\cL$ directly instead of $\xX$, thus writing $\cALs$ when
we deal with Kontsevich integral of $\cL$.

Note that $\ZKAL$ is not an invariant of the
link $\cL$, because it depends on the positions of the dots on the link
components. This dependence could be eliminated, had we defined Kontsevich
integral as taking values in the space $\cA(\orcirc_1\ldtc\orcirc_{\xN})$,
in which the 1-valent vertices of (1,3)-valent graphs are placed on
oriented circles. However, we have to keep the dots, since we are about to
pass from $\cALs$ to $\cBL$.

The space $\tcBX$ is a space of (1,3)-valent graphs, whose
1-valent vertices (or, equivalently, legs)
are assigned to the elements of $\xX$. The vertices which
carry the same label (that is, element of $\xX$),
are considered equivalent and
may be permuted when an isomorphism between two graphs is considered. The
space $\cBX$ is a quotient of $\tcBX$ over its subspace
$\tcBXIHX$, which is a span of the IHX relation.


The space $\cBX$ is a bi-graded
algebra. The multiplication of two graphs is
defined as their disjoint union. The
unit of this multiplication is the empty graph, and we denote it as
$\bfun$. The two gradings that are preserved by the multiplication of $\cBX$
are
\qq
\dgh{\xD} & = & \# (\mbox{edges}) - \# (\mbox{3-vertices})
= \hlf\; \#(\mbox{all vertices})
= \#(\mbox{1-vertices}) + \echi(\xD),
\label{1.11}
\\
\dgth{\xD} & = & \#(\mbox{edges}) - \#(\mbox{all vertices}) = \echi(\xD) =
\dgh{\xD} - \#(\mbox{1-vertices}),
\label{1.11x}
\qqq
where $\echi(D) = \#(\mbox{edges}) - \#(\mbox{all vertices})$ is the
(\emph{opposite})
Euler characteristic of $D$. The grading $\dgh{}$ is also compatible with the
multiplication in the algebra $\cALs$.

A transition from
$\cAX$
to $\cBX$ is produced by the Poincare-Birkhoff-Witt
symmetrization of the order in which 1-valent vertices appear on the
oriented segments. The inverse map is denoted as
\qq
\xmap{\xchi}{\cBX}{\cAX}.
\label{1.11y}
\qqq
%

We define $\ZKBL$ as the
Kontsevich integral of $\cL$ with value in $\cBL$, which in line with\rx{1.11z}
is, of course, a space
of (1,3)-valent graphs, whose 1-valent vertices are labeled by components of
$\cL$ or simply by the numbers $1\ldtc \xN$.
%
In other words,
\qq
\ZKAL = \xchi (\ZKBL).
\label{1.7*}
\qqq

If $\cL$ is in a rational homology sphere $M$, then we present it as a
surgery on the framed $\cLp$ part of a combined link $\cL\cup\cLp\in S^3$.
Then $\ZKALM$ is defined with the help of the LMO map\cx{LMO} or its \AA
rhus version\cx{A2}. Namely, we first symmetrize all legs of $\ZKALLp$
which are attached to the components of $\cLp$, and then we apply the LMO
or \AA rhus map to those legs. The result is $\ZKALM\in\cALs$. Then we
pass from it to $\ZKBLM\in\cBL$ in exactly the same way as we did it for
$\cL\in S^3$. Therefore from now on we consider $\cL$ to be a link in a
rational homology sphere $M$, and since $M$ will always be the same, we
drop it from our notations.


By using the fact that
$\ZKBL$ (by definition) has coefficient 1 at the empty graph, we can
define the logarithm of $\ZKBL$
\qq
\WKL = \log \ZKBL
\label{1.8}
\qqq
by the formula $\log (1+x) = \snoi (-1)^{n-1} x^n/n$.
We call an element of a graph algebra
\emph{\nrw}, if it can be presented as a linear
combination of connected graphs.
An important
property of Kontsevich integral is that
$\WKL$ is \nrw.

As we have mentioned, $\WKL$ depends on a presentation of $\cL$ as a
dotted
Morse link. However, the strut part of $\WKL$ is topologically invariant
and, in fact, is well-known to be
\qq
\WKLstr = \hlf \sijoxN \xxlij\;\; \tstv{i}{j},
\label{1.8*}
\qqq
where $\xxlij$ are the linking numbers of $\cL$.

Let $\xG$ be a compact Lie group,
$\mfg$ be its Lie algebra and $\mfh\subset \mfg$ its Cartan
subalgebra. We assume that $\mfg$ is equipped with a positive-definite Killing
form, so we will not make a consistent
effort to distinguish between $\mfg$,
$\mfh$ and their dual spaces $\mfgs$, $\mfhs$.
Let $\valoN\in\mfhs$ be the highest
weights (shifted by a half-sum of positive roots) of $\xG$-modules
$\VmvaoN$
assigned to the components of $\cL$. We denote by $\JvaNL$ the colored
Jones polynomial of $\cL$ which corresponds to this data. It is well-known
that it has an expansion in powers of $\hb = \log q$
\qq
\JvaNL = \lrbc{\pjoxN \dim V_{\val_j} } \snzi \pnvaL \hb^n,
\label{1.9}
\qqq
%
where the link invariants $\pnvaL$ are Weyl-invariant polynomials of
$\uval$. This expansion can be derived from $\ZKBL$. The standard
method is to step back to $\ZKAL$ through \ex{1.7*}
and then apply the $\mfg$ weight system.
The weight system is applied in two steps:
\qq
\begin{CD}
\cALs @>{\Tgh}>> \UgLh @>{\TrVua}>> \ICh
\end{CD}
\label{1.10}
\qqq
The map $\Tgh$ is a combination of two maps. First, a graph
$\xD\in\cALs$ is
multiplied by $\hb^{\dgh{\xD}}$.
Then the map $\cALs\mapright{\Tg}\UgL$ is applied. Let us recall the
definition\cx{BN}
of two similar maps (which we denote by the same symbol)
\qq
\cAX\mapright{\Tg}\UgX,\qquad \cBX\mapright{\Tg}\SgX,
\label{1.11*}
\qqq
where $\orX$ denotes the number of elements in $\xX$.
Let $\fmfg\in \wdgtg$ be the structure constant tensor of $\mfg$ and
let $\hmfg\in \stgs$ be the
positive definite Killing metric on $\mfg$ normalized in such
a way that short roots have length $\sqrt{2}$. We denote by
$\himfg\in \stg$ its inverse matrix. Now let $D$ be a (1,3)-valent
graph with $\neDo$ 1-valent vertices, $\neDt$ 3-valent vertices and $\neDe$
edges. Assume for a moment that all 1-valent vertices of $D$ have
distinct labels.
Associate a space $\wdgtg$ with each 3-vertex of $D$ and a
space $\Srs^2\mfg$ with each edge of $D$. Consider their tensor products
$\wdgtgx$ and $\stgx$. There is a natural `index contraction' map
\qq
\begin{CD}
\wdgtgx\otimes\stgx @>{\xCD}>> \mfg^{\otimes\neDo},
\end{CD}
\label{1.11*1}
\qqq
which pairs up the spaces $\mfg$ and $\mfgs$ from incident egdes and 3-valent
vertices, so that the spaces $\mfg$ of edges which are incident to 1-valent
vertices remain uncontracted. Then we define $\tTg(D)$ as
\qq
\fmfg^{\otimes\neDt}\otimes(\himfg)^{\otimes\neDe}
\mathop{{\longmapsto}}^{\xCD} \tTg(D)\in\tmfgx.
\label{1.11*2}
\qqq
Now, as we know, if the graph $D$ is an element of $\cAX$, then the STU
relations project $\tTg(D)$ naturally to $\Tg(D)\in\UgX$, whereas if $D$ is
an element of $\cBX$, then the equivalences between 1-vertices project $\tTg(D)$
to $\Tg(D)\in\SgX$.


The map $\TrVua$ of\rx{1.10} evaluates
the traces of the elements of $\Ug$ in $\mfg$-modules $\VmvaoN$ assigned
to the link components. Finally,
\qq
\JvaNL = \TrVua\Tgh (\ZKAL).
\label{1.12}
\qqq
\begin{remark}
\label{r2.1}
\rm
Since there is a natural inclusion $\mfg\hookrightarrow\mfgC$, then we may
consider $\fmfg$ and $\himfg$ as elements of $\wdgtgC$ and $\stgC$.
Obviously, both maps\rx{1.10} commute with complexification. The complexified
version of\rx{1.10} will be more convenient later when we split $\mfg$ into a
Cartan subalgebra $\mfh$ and the space of roots
$\mfr$, because the adjoint action
of $\mfh$ on $\mfr$ could be diagonalized over $\IC$.
\end{remark}

\subsection{Wheeling, Kirillov integral and Reshetikhin formula}
\label{2xs.2}

An alternative way of deriving the expansion\rx{1.9} from Kontsevich
integral is based on Kirillov's integral formula for the trace map
$\TrVua$. Let $\Sg\mapright{\bmg}\Ug$ be the Poincare-Birkhoff-Witt
symmetrization isomorphism and let $\xmap{\duflo}{\Sg}{\Sg}$ be the Duflo
isomorphism (see\cx{Wh} and references therein). For an element
$\val\in\mfhs$ let $\corbal\subset \mfgs$
be its coadjoint orbit equipped
with the integration measure derived from the Kirillov-Kostant
symplectic form. Then the following square is commutative
\qq
%
\begin{CD}
\Sg  @>{\bmg}>>  \Ug
\\
@VV{\duflo}V   @VV{\TrVa}V
\\
\Sg @>{\icorbal}>> \IC
\end{CD}
\label{1.13}
\qqq
Here the map $\xmap{\icorbal}{\Sg}{\IC}$ means the following: the space
$\Sg$ is naturally isomorphic to the space
of polynomial functions
on $\mfgs$ and these polynomials can be integrated over $\corbal$.

According to\cx{Wh}, Duflo isomorphism can be performed at the graph level
as a wheeling map $\Whh$. Here is the definition of $\Whhj$ acting on
the legs of the graphs of $\cBL$ assigned to the link component $\cL_j$.
We define the `modified Bernoulli numbers' $\btn$ by the expansion
\qq
\hlf \log \lrbc{
\sinh(x/2) \over (x/2) } = \snoi \btn\,x^{2n}.
\label{1.14}
\qqq
Then we define $\Wh\in\cB$ as the exponential
\qq
\Wh = \exp\lrbc{ -\snoi \btn \omtn},
\label{1.14*}
\qqq
where $\omtn$ is a circle with $2n$
radial edges attached to it\footnote{In contrast to\cx{Wh}, we placed an extra
minus sign in the exponent of the formula for $\Wh$. Thus our `wheeling' in the
inverse of the wheeling of\cx{Wh}.}. Following the notations of\cx{A2}, we label
the 1-valent vertices of the graphs of $\cBL$ with link components by
assigning to them formal variables $\xoN$. We also assign the dual variable
$\parxj$ to all the 1-valent vertices of $\Wh$ thus producing
$\Whj = \Wh(\parxj)$. Now for any $\xF\in\cBL$
\qq
\Whhjv{\xF} = \Whj \xfle \xF,
\label{1.15}
\qqq
where the operation $\xfl$ was defined in\cx{A2} as the `application of
derivatives' $\parxj$ of $\Whj$ to the variables $x_j$ sitting at some
1-valent vertices of $\xF$. This application of the derivatives means that
the corresponding edges of $\Whj$ and $\xF$ are joined into single edges:
\qq
\begin{picture}(4500,500)
\thinlines
\unitlength=1mm
\put(0,0){\line(1,0){20}}
\put(10,0){\circle*{1.5}}
\put(10,0){\line(0,1){10}}
\put(18,-5){\mbox{$\parxj$} }
\put(25,0){\mbox{$\xfl$} }
\put(30,0){\line(1,0){20}}
\put(40,0){\circle*{1.5}}
\put(40,0){\line(0,1){10}}
\put(30,-5){\mbox{$x_j$} }
\put(55,0){\mbox{$=$} }
\put(60,0){\line(1,0){30}}
\put(70,0){\circle*{1.5}}
\put(80,0){\circle*{1.5}}
\put(70,0){\line(0,1){10}}
\put(80,0){\line(0,1){10}}
\end{picture}
\label{1.16}
\qqq

Denote by $\WhhL$ the composition of the wheeling maps $\Whhj$ for all the
link components of $\cL$.
This map allows us to extend a commutative square\rx{1.13} to a
commutative diagram (\cf the monster diagram of\cx{Wh})
\qq
\xymatrix{
& \cALs \ar@/^/[dr]^{\Tg} &
\\
\cBL\ar@/^/[ur]^{\xchi} \ar[r]^-{\Tg} \ar[d]^{\WhhL}
& \SgL \ar[r]^{\bmg} \ar[d]^{\dflL} & \UgL \ar[d]^{\TrVua}
\\
\cBL \ar[r]^-{\Tg} & \SgL \ar[r]^-{\iOuval} & \IC
}
\label{1.17}
\qqq
Since the map\rx{1.11y}
is the graph algebra counterpart
of the PBW symmetrization map $\bmg$, then the upper triange of\rx{1.17} is
commutative.
The maps $\xmap{\Tg}{\cBL}{\SgL}$ convert (1,3)-valent graphs into
the elements of $\SgL$ by placing the structure constants of $\mfg$ at
3-valent vertices and contracting indices along the internal edges,
while $\dflL$ is the Dulfo isomorphism applied to all components of the
tensor product $\SgL$. We also used a shortcut notation
$\corbuval$ for $\corbalv{1}\times\cdots\times\corbalv{\xN}$.
%
%

The $\xG$-invariance of the elements of $\SgL$ means that we can omit one of
the integrals in $\iOuval$, replacing it by the volume of the coadjoint orbit.
Since $\Vol\corbal = \dim \Vval$, then in the context of the diagram\rx{1.17}
\qq
\iOuval = \dim \Vvalk \evkal \iOuvalrk, \qquad 1\leq k\leq \xN,
\label{1.17*1}
\qqq
%
where $\xmap{\eval}{\SgL}{\SgLo}$ evaluates the elements of the $\xk$-th
component of $\SgL$ (considered as polynomials on $\mfhs$) on $\val$.

It is easy to see that the wheeling map $\Whh$ does not change the
degree\rx{1.11} of graphs, therefore we can replace $\Tg$ by $\Tgh$ and
add $[[\hb]]$ to the spaces $\SgL$, $\UgL$ and $\IC$ in the
diagram\rx{1.17} without spoiling its commutativity.

A combination of equations\rx{1.7*} and\rx{1.12} demonstrates that
\qq
\JvaNL = \TrVua\circ\Tgh \circ \xchi (\ZKBL).
\label{1.18}
\qqq
At the same time, the $[[h]]$ version of the diagram\rx{1.17} shows that
\qq
\TrVua \circ  \Tgh \circ \xchi =
\iOuval
\Tgh \circ \WhhL.
\label{1.19}
\qqq
Thus we obtain the alternative version of \ex{1.12}
\qq
\JvaNL =
\iOuval
\Tgh \circ \WhhL (\ZKBL).
\label{1.20}
\qqq

Let us look more closely at the wheeled Kontsevich integral
\qq
\ZKOBL = \WhhL (\ZKBL).
\label{1.20*}
\qqq
According to\cx{A2},
the operation $\xfl$ between two exponentials of connected graphs is
itself an exponential of connected graphs, whose structure is described in
Exercise~2.4 of\cx{A2}
(one just has to expand both exponentials, apply $\xfl$ and then
pick only connected graphs -- they will make up the exponent of the resulting
exponential). Therefore
\qq
\ZKOBL = \exp (\WKOL)
\label{1.21}
\qqq
and $\WKOL\in\cBL$ is \nrw.
The
strut part of $\WKOL$ is the same as that of $\WKL$:
\qq
\WKOLstr = \hlf \sijoxN \xxlij\;\;
\tstv{i}{j},
\label{1.21*}
\qqq

The maps
$\Tg$ and $\Tgh$ are algebra homomorphisms, so
\qq
\Tgh (\ZKOBL) = \RiL,
\label{1.22}
\qqq
and we proved the following
\begin{theorem}[Reshetikhin's formula]
The expansion of the colored Jones-HOMFLY polynomial of a link $\lst$ in powers
of $\hb=\log q$ can be presented as Kirillov's integral over the coadjoint
orbits of the (shifted) highest weights of $\xG$-modules assigned to the link
components
\qq
\JvaNL = \iOuval \RiL.
\label{1.23}
\qqq
\end{theorem}
%
Let us write \ex{1.23} more explicitly. First, we split $\WKOL$ into a sum
\qq
\WKOL = \smnd \WKOLmn,
\label{1.24}
\qqq
where $\WKOLmn$ are linear combinations of graphs $D$ such that
$\#(\mbox{1-vertices}) = m$ and $\echi(D)=n-1$ so that according to \ex{1.11},
\qq
\dgh{\WKOLmn} = m+n-1.
\label{1.25}
\qqq
Then $\Tg(\WKOLmn)$ is an element of $\SgL$ of total degree $m$ and
$\Tgh(\WKOLmn)$ is proportional to $\hb^{m+n-1}$. Since $\Sg$ is
canonically isomorphic to the algebra of polynomials on $\mfgs$, we can
present $\Tg(\WKOLmn)$ as a $\xG$-invariant
homogeneous polynomial $\LLbetmn$ of degree $m$,
where $\vbetoN\in\mfgs$.  Then in view of \ex{1.25},
\qq
\Tgh \WKOL = \smnd \LLbetmn\,\hb^{m+n-1},
\label{1.25*}
\qqq
and
we can rewrite \ex{1.23} as
\qq
\JvaNL = \iOuval \dvbetoN \exp\lrbc{
\smnd \LLbetmn\,\hb^{m+n-1}
}.
\label{1.26}
\qqq

Since all graphs with a single leg are equal to zero in $\cBL$, we
conclude that
\qq
\LLbetmo = 0.
\label{1.27}
\qqq
Also, since a tree graph has at least 2 legs, $\LLbetzz = 0$.
We deduce from \ex{1.21*} that
\qq
\LLbettz =  \hlf \sijoxN \xxlij\; \scp{\vbet_i}{\vbet_j},
\label{1.27*1}
\qqq
Finally, since $\cBL$ does not contain a $\echi(D)=0$ graph without legs,
then $\LLbetmo=0$.
Since $\LLbetmz = \LLbetmo = 0$ for $m\leq 1$, then
the exponent of the \rhs of \ex{1.26} contains only
strictly positive powers of $\hb$. If we expand the exponential in Taylor
series and then integrate the coefficients at each power of $\hb$, then we will
reproduce the Melvin-Morton expansion\rx{1.9}.

\subsection{Invariant measure on a coadjoint orbit}
\label{2xs.3}

Let us perform a substitution\rx{1.17*1} and
change the integration variables in \ex{1.26} from $\vbetj$ to
$\vbj=\hb\vbetj$, $1\leq j\leq \xN$. Since the polynomials $\LLbetmn$ are
homogeneous of degree $m$, then \ex{1.26} becomes
\qq
\JvaNL = \hbpLgo \dim \Vvalk \iOhuvalrk \duvbrk \exp \lrbc{
\smnd \LLbmn\,\hb^{n-1},
}
\label{1.28}
\qqq
where $\nrtg=(1/2)\dim \corbal$ is the number of positive roots of $\mfg$.
Now for the variables $\vaoN\in\mfhs$ let us consider a formal expression
\qq
\JfvaNL = \hbpLg\, {\dgvak\over\dgrho} \iOuvark \duvbrk \exp  \lrbc{
\ohb \smnd \LLbmnk\,\hb^n
},
\label{1.29}
\qqq
%
where
\qq
\dgva = \prvldp \scp{\va}{\vl},
\label{1.29*}
\qqq
so that $\dgval/\dgrho=\dim\Vval$.
Expression\rx{1.29}
does not make sense in itself, because the coefficients at
the powers of
$\hb$ in the exponent are formal power series in $\vboN$, and we can not make
sense of the integral. We can only make sense of the expression
$\JfLv{\hvaloN}$, because in view of \ex{1.28} it can be defined as a power
series in $\hb$ and
\qq
\JfLv{\hvaloN} = \JvaNL.
\label{1.30}
\qqq
%

Despite the lack of proper definition of the whole integral\rx{1.29}, its
integrand looks suitable for the stationary phase approximation in
the limit of
\qq
\limzero{\hb},\qquad(\limone{q}),\qquad
\vam{j}\;(=\hb\val_j)\;=\const,\qquad 1\leq j\leq \xN.
\label{1.31}
\qqq
Indeed, we may assume that a formal parameter $\hb$ is purely imaginary;
actually, this is the case in the Quantum Field Theory approach to quantum
invariants of 3d topology.
The advantage of the stationary phase approximation for us is that within its
calculations the `honest' integrals are essentially replaced by simple formal
combinatorial manipulations with the integrand. This combinatorics is
called `Feynman rules' in Quantum Field Theory.

The first step in the stationary phase approximation is to identify the
stationary points of $\smti\LLbmzk$
on $\corbuvark$. Generally, this may be
tricky, since $\smti\LLbmz$ is only a formal power series: one might try to
find the stationary points of its lowest terms and then correct them
perturbatively by taking into account the higher order terms.
Luckily, there exists a point
\qq
\uvb = \uva 
\label{1.32}
\qqq
(that is, when all $\uvb$ belong to the same Cartan subalgebra) which is
manifestly stationary (we will check this a bit later). We are going to
calculate (or, more precisely, define) its contribution to the
integral\rx{1.29} by following the formal rules of the stationary phase
approximation. We will call the resulting power series in $\hb$ `the
$\mfg$-based \urcc invariant'.

Let us find a convenient parametrization of $\corbuvark$ around the
point\rx{1.32}. A map
\qq
g\xT \mapsto \AdgT\,\va
\label{1.33}
\qqq
establishes an isomorphism between a coadjoint orbit $\corba$ and the
quotient $\GmT$. A composition of maps
\qq
\begin{CD}
\mfr @>{\exp}>> \xG @>>> \GmT,
\end{CD}
\label{1.34}
\qqq
identifies the vicinity of the origin in
the space of roots $\mfr\subset\mfg$ with the
vicinity of the identity at $\GmT$. Thus a composition of the
maps\rx{1.33} and\rx{1.34} makes $\mfr$ a coordinate space for a vicinity
of $\va$ in $\corba$:
\qq
\vb = \vDRax,\qquad\mbox{where}\qquad
\vDRax =
\Ader\, \va
= \smzi {\adxvr^m\,\va\over m!}.
\label{1.35}
\qqq
%
%
Let us find the measure $\msravr$ on $\mfr$ which corresponds to the
Kirillov-Kostant measure on $\corba$. It is easy to find the measure at the
origin:
\qq
\msraz = \mmrcr\;\dgva.
\label{1.36}
\qqq
(the minus sign compensates for the fact that the metric on $\mfg$ is
negative-definite).

We are going to use the invariance of the
Kirillov-Kostant measure under the adjoint action of  the elements $\ey$
($\vy\in\mfr$) on $\corba$ in order to find $\msravr$ at a general point
$\xvr$.
Let us rewrite this action in terms of $\mfr$ coordinates. We define a
(non-linear) map $\xmap{\FLy}{\mfr}{\mfr}$ by
%
a commutative diagram
\qq
\begin{CD}
\mfr @>\FLy>> \mfr
\\
@VV{\vDRa}V @VV{\vDRa}V
\\
\corba @>{\Adey}>> \corba
\end{CD}
\label{1.37}
\qqq
where the map
$\xmap{\vDRa}{\mfr}{\corba}$
is defined by \ex{1.35}. This diagram
implies a defining relation for $\FLyr$ and its companion $\tFLyx$
\qq
\ey\er = \eFLyr\,
e^{\tFLyx},
\qquad
\tFLyx\in\mfh,\;\;\vy,\xvr,\FLyr\in\mfr.
\label{1.38}
\qqq
Since the integration measure $\msravr$ has to be $\FLy$-invariant, it must be
of the `left-invariant' form
\qq
\msravr =  \msraz\big/ \det \drFLvrz,
\label{1.39}
\qqq
where $\drFLyr = \del \FLyr/\del\xvr$ is the derivative of the map $\FLyr$.
Since $T_{\xvr}\mfr\cong\mfr$, we assume that $\xmap{\drFLyr}{\mfr}{\mfr}$.
%
\begin{theorem}
\label{t1.1}
The measure\rx{1.39} is $\FL$-invariant.
\end{theorem}

\noindent{\em Quick proof of Theorem\rw{t1.1}.}
The $\FL$-invariant measure exists: it is a pull-back of the Kirillov-Kostant
measure by the map $\vDRa$
of\rx{1.35}. Also it is unique (up to a constant factor), since
it must satisfy \ex{1.39}. Therefore, the measure\rx{1.39} is
$\FL$-invariant.\qed

This proof is based on the uniqueness of the
conjugation-invariant measure on $\corba$. Since our
ultimate goal is to strip away Lie algebras from our calculations, we will
present another proof, which is essentially based on the combinatorics of
commutators.

\noindent{\em Combinatorial proof of Theorem\rw{t1.1}.}
Since $\exp(\mfr)\subset\xG$ does not form a subgroup, then the invariance
proof is a bit more subtle than the standard proof of the left-invariance of
measures like\rx{1.39} on the whole group $\xG$. We have to learn how to
compose $\FLyo$ and $\FLyt$. A relation
\qq
\eyt \eyo = \eFLyto
e^{ \tFLv{\vy_2}(\vy_1) }
\label{1.40}
\qqq
implies that
\qq
\FLyt \circ \FLyo = \FLdyot\circ
\AdvFhyto.
\label{1.41}
\qqq
Since $\xmap{
\AdvFhyto
}{\mfr}{\mfr}$ is a
linear map and $\dtmfr\AdvFhyto = 1$, we conclude that
\qq
\det\drFLyt\,\det\drFLyo =  \det\drFLdyot
\label{1.42}
\qqq
and in particular
\qq
\det\drFLyt(\vy_1)\,\det\drFLyo(\vzr) = \det\drFLdyot(\vzr).
\label{1.43}
\qqq
The latter relation means that the measure\rx{1.39} is indeed
adjoint-invariant.\qed

Now let us calculate the determinant in the measure\rx{1.39}. For a stationary
phase calculation it is sufficient to express it as a power series in $\xvr$.
Actually, it is easier to calculate the determinant of the inverse operator
$\drFLmvrr$. In order to find the action of that operator on $\dxvr\in\mfr$, we
perform the following approximate calculation up to $\cO(\dxvr^2)$
\qq
e^{-\xvr}\, e^{\xvr+\dxvr}  & \approx & \exp \drrrds
\nonumber
\\
& \approx & \exp\lrbc{ \prjr\,\drrrd\, }\;
\exp\lrbc{ \prjh\,\drrrd\, },
\label{1.44}
\qqq
where $I$ is the identity operator and by definition
\qq
\drrr = \snzi {(-\adxvr)^n \over (n+1)!}
= I + \snoi {(-\adxvr)^n \over (n+1)!},
\label{1.45}
\qqq
while $\prjr$, $\prjh$ are the orthogonal projectors of $\mfg$ onto $\mfr$ and
$\mfh$. Equation\rx{1.44} means that
\qq
\drFLmvrr = \prjr\,\drrr.
\label{1.46}
\qqq
in the sense that the \rhs of this formula maps $\mfr\subset\mfg$ into $\mfr$.
We see from \eex{1.45} and\rx{1.46} that as an operator on $\mfr$
\qq
\drFLmvrr = \Idop + \prjr\snoi {(-\adxvr)^n \over (n+1)!},
\label{1.47}
\qqq
so we can calculate its determinant as a power series in $\xvr$ through the
formula
\qq
\det (I + A) = \exp\lrbcs{\Tr \log (I+A)}
= \exp \lrbc{ -\snoi {(-1)^{n}\over n}
\,\Tr A^n
}
\label{1.48}
\qqq
Thus combining \eex{1.36},\rx{1.39},\rx{1.47} and\rx{1.48} we find that
\qq
\msravr = \mmrcr\;\dgva \exp\lrbs{- \snoi {(-1)^{n}\over n}\,
\Tr_{\mfg}\lrbc{
\smoi{\prjr\, (-\adxvr)^m \over (m+1)! }
}^n\,
}.
\label{1.49}
\qqq
Note that all calculations with operators in this equation are performed in
$\mfg$, while the projector $\prjr$ essentially reduces $\Tr_{\mfg}$ to
$\Tr_{\mfr}$.

\subsection{$\mfg$-based \urcc invariant as a stationary phase integral}
\label{2xs.4}

Having presented the integration measure for the substitution\rx{1.35} as an
exponential of the power series in $\xvr$, we proceed to substitute
\qq
\vbj=
\vDRav{j}(\vx_j)
\qquad (1\leq j\leq \xN,\; j\neq \xk), \qquad \vbk=\vak
\label{1.50}
\qqq
\vspace{-0.0005in}
\nopagebreak in the exponent
of\rx{1.29}, expand it in powers of $\xvr$ and combine this expansion
with the
exponent of \ex{1.49}. This results in a formula
\qq
\lefteqn{ \IrfvaNL }
\nonumber\\
& = &
\hbpLg\, {\dgvak\over\dgrho}
\nonumber
\\
&&\times
\ilvrkz \duxvrk
\atv{\lrbs{
\exp  \lrbc{
\smnd 
\LLv{m,n}{\vDRav{1}(\vx_1)\ldtc\vDRav{\xN}(\vx_{\xN})}
\,\hb^{n-1}
}   \pjoxN \msr{\va_j}{\vx_j}
}
}{\vx_{k}=0}
\nonumber\\
& = &
\hbmnDp\,
\mohbpLg\, \dguvarho
\nonumber
\\
& &\qquad\times
\ilvrkz \duxvrk \exp\lrbc{\ohb \smnzi \LLrarkmn \,
\hb^n
}.
\label{1.51}
\qqq
Here
%
$\dguva = \pjoxN \dgvaj$,
%
$\ivrkz$ denotes the stationary phase contribution of the point
$\uxvrk=\vzr$ to the integral, while $\LLrarkmn$ are homogeneous polynomials of
$\uxvrk$ of degree $m$ whose coefficients are formal power series of $\uva$:
$\LLrarkmn\in\IQuvark$.

Let us check the structure of the terms $\LLrarkmz$, $m\leq 2$ which are
crucial for the stationary phase approximation. We will need two simple lemmas.
Let $D$ be a connected graph with $m$ legs from $\cBX$, all of its legs
having distinct labels, and let $\TgD$ be
its image in $\Invrs{(\Ug)^{\otimes m}}{\xG}$.
Since $\TgD$ is of degree $m$, then we may consider it to
be an $m$-linear $\xG$-invariant function on $\mfgs$,
which we denote as $\TgDbm$.
\begin{lemma}
\label{l1.1}
If $D$ is a tree graphs, $\vatm\in\mfh$ and $\vb\in\mfg$, then $\TgDxbm=0$,
unless $D$ is a strut.
\end{lemma}
\proof
If $D$ is not a strut, then out of any $m-1$ of its $m$ 1-valent vertices one
can always find a pair of 1-valent vertices connected to the same 3-valent
vertex. This means that each term in the expression of $\TgDxbm$ contains a
commutator of the elements $\vaj$ attached to those 1-valent vertices. Then
$\TgDxbm=0$, because all elements of $\mfh$ commute.\qed
\begin{lemma}
\label{l1.2}
If $\vaom\in\mfh$ and $\vb\in\mfg$, then
\qq
\TgDbmx = 0.
\label{1.52}
\qqq
\end{lemma}
\proof
Since $\TgDbm$ is invariant under adjoint action of $\xG$ on its agruments,
then
\qq
\TgDbmx = -\sum_{j=2}^m \TgDv{\vb,\va_2\ldtc[\va_j,\va_1]\ldtc\va_m}.
\label{1.53}
\qqq
The \rhs of this equation is zero, since all $[\va_j,\va_1]=\vzr$.\qed

Obviously, $\LLrarkzz = \smti\LLmzuva$ (we dropped $m=0,1$
in the sum, because
for them $\LLmzuvb=0$). The polynomials $\LLmzuvb$ come from $\Tg$
acting on the tree graphs of $\WKOL$, which coincide with the tree graphs
of $\WKL$, since wheeling does not produce new tree graphs.
Therefore, according to
Lemma\rw{l1.1}, the only contribution to $\LLrarkzz$ comes from the strut
part of $\WKL$, and in view of \ex{1.8*},
\qq
\LLrarkzz =\hlf
\sijoxN\xxlij\, \scp{\va_i}{\va_j} \eqdef \lkLuva.
\label{1.55}
\qqq

The terms $\LLrarkon$ (and in particular $\LLrarkoz$) come from all graphs
of $\WKOL$. The formula
$\vDRax = \va + [\xvr,\va] + \cO(\xvr^2)$ suggests
that the contribution of a graph $D$ from $\WKOL$ is calculated by placing
$[\xvr_j,\va_j]$ at one of its 1-valent vertices and $\uva$ at all other
1-valent vertices. Then Lemma\rw{l1.2} says that
\qq
\LLrarkon = 0\qquad\mbox{for all $n$}.
\label{1.56}
\qqq
In particular, this means that there are no terms of order $\hb^{-1}$
in the exponent of\rx{1.51} which are linear in $\uxvrk$.
Thus we proved the following
\begin{lemma}
\label{l1.3}
The point\rx{1.32} is a stationary point of the exponent of\rx{1.29}.
\end{lemma}

Now we turn to $\LLrarktz$. Since
$\vb = \va + [\xvr,\va] + (1/2)[\xvr,[\xvr,\va]] + \cO(\xvr^3)$, the terms of
$\LLrarktz$ come from tree graphs of $\WKL$ in two ways: either by placing
the double commutator
$(1/2)[\xvrj,[\xvrj,\vaj]]$ at one of their 1-valent vertices and $\uva$ at the
others, or by placing $[\xvrj,\vaj]$ at two 1-valent vertices and $\uva$ at the
others. According to Lemma\rw{l1.1}, only the strut graphs may contribute in
the first way, so the contribution of the first way is
\qq
{1\over 4}\snoijxNk\lkij\, \scp{[\xvri,[\xvri,\va_i]}{\va_j}
=
{1\over 4}\snoijxNk\lkij\, \scp{ \adv{\va_i}\adv{\vaj} \xvri}{\xvri}
\label{1.57}
\qqq
(we neglected the first
way contribution of strut graphs $\tstv{j}{j}$
since it is canceled by their
second way contribution).
It is easy to see that the only graphs contributing to $\LLrarktz$ a term
proportional to $\xvri\xvrj$ are the ones which have the `haircomb' shape
of \fg{f1.1} up to
the cyclic order at 3-valent vertices (note that the graph of \fg{f1.1}
contributes also to the terms $\xvr_{j_1}\xvr_j$, $\xvr_i\xvr_{j_m}$ and
$\xvr_{j_1}\xvr_{j_m}$).
\begin{figure}[hbt]
\unitlength=1mm
\begin{center}
\begin{picture}(45,15)(7,0)
\put(0,-2){ \makebox{$\cL_i$} }
\put(5,0){ \line(1,0){45} }
\put(10,0){ \circle*{1.5} }
\put(10,0){ \line(0,1){5} }
\put(8,7){ \makebox{$\cL_{j_1}$} }
\put(17,0){ \circle*{1.5} }
\put(17,0){ \line(0,1){5} }
\put(15,7){ \makebox{$\cL_{j_2}$} }
\put(24,0){ \circle*{1.5} }
\put(24,0){ \line(0,1){5} }
\put(31,0){ \circle*{1.5} }
\put(31,0){ \line(0,1){5} }
\put(38,0){ \circle*{1.5} }
\put(38,0){ \line(0,1){5} }
\put(45,0){ \circle*{1.5} }
\put(45,0){ \line(0,1){5} }
\put(43,7){ \makebox{$\cL_{j_m}$} }
\put(52,-2){ \makebox{$\cL_j$} }
\end{picture}
\end{center}
\caption{A haircomb graph.}
\label{f1.1}
\end{figure}
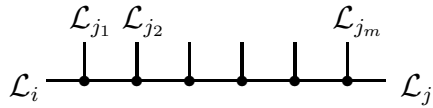
The contribution of the graph of \fg{f1.1} is
\qq
\scp{
[\vajv{m},[\ldots[\vajv{2},[\vajv{1},[\xvri,\vai]\ldots]
}{[\xvr_j,\vaj]} =
-\scp{
\adv{\vaj}\adv{\vajv{m}}\cdots\adv{\vajv{1}}\adv{\vai}\xvri}{\xvrj}
\label{1.58}
\qqq
%
Thus adding up expression\rx{1.57} and the contributions\rx{1.58} we obtain
$\LLrarktz$ as a quadratic form of $\uxvrk$.
%

Now we are almost ready to write a well-defined expression for the stationary
phase integral\rx{1.51}.
With a slight abuse of notation, let us think of $\uxvrk$ as a single vector in
a direct sum of $\xN-1$ spaces $\mfr$:
\qq
\uxvrk = \sjoxNk \xvrj \in \mfrxN.
\label{1.59}
\qqq
Then we can write
\qq
\LLrarktz =
\hlf\sum_{1\leq i,j \leq N\atop i,j\neq k} \lrbcs{ \vx_i,
(\gQLuak)_{ij}\,\vx_j},
\label{1.60}
\qqq
where $\gQLuak$
is a linear operator acting on $\mfrxN$, whose coefficients
$(\gQLuak)_{ij}$
are
formal power series in $\vua$. Let us assume that $\gQLuak$
is non-degenerate,
or in other words, that
\qq
\det \gQLuak \not\equiv 0.
\label{1.61}
\qqq
We denote by $\gQLuaik$
the inverse operator. Then the formula for the stationary
phase integral\rx{1.51} is
\qq
\IrfvaNL = 
\hbmnDp\;
\dguvarho  \;\;
e^{\ohb\lkLuva}\,
\lrbcs{ \det\gQLuak }^{-1/2}
e^{\LLrarkzo}
\,
\FDLua,
\label{1.62}
\qqq
where $\lkLuva$ was defined by \ex{1.55},
%
%
while the most interesting factor $\FDLua$ is defined by the formula
\qq
\lefteqn{
\FDLua =
}
\label{1.64}
\\
&&\hspace*{-37pt}
=
\exp \ohb \lrbc{\smthi \LLradrkmz +
\smnoi
\LLradrkmn \,\hb^n
}
\exp\Big(\!\! -{\hb\over 2}\,
\scp{\uxvrk}{\gQLuaik\, \uxvrk}
\Big)
\Big|_{\uxvrk=
0}
\nonumber
\qqq
and is known in the context of Quantum Field Theory as a sum of Feynman
diagrams
(hence the name FD).
It is easy to see that $\FDLua$ is a well-defined element of
$\IQQuvahb$, where for a ring $R$, $\Qf{R}$ denotes its
field of fractions. The coefficients of $\LLrarkmn$ and $\gQLuak$ belong
to $\IQuva$, so denominators in $\FDLua$ are due exclusively to the
denominators of $\gQLuaik$
which are, of course, equal to the determinant
$\det\gQLuak$.
Then an easy combinatorics (which we will review in
subsection\rw{6xs.2})
demonstrates that $\FDLua$ can be put in a form
\qq
\FDLua = \exp\lrbc{
\snoi   {\pQnLuva \over \Big(\det\gQLuak\Big)^{3n} }\;\hb^n
},\qquad \pQnLuva\in\IQuva.
\label{1.65}
\qqq
\begin{remark}
\rm
Since the sum in the exponent of \ex{1.49} starts at $n=1$, then it does
not contain the terms of 0-th order in $\vx$. Therefore, the term
$\LLrarkzo$, which is singled out in \ex{1.62}, comes exclusively from the
polynomials $\LLbmo$
\qq
\LLrarkzo = \smzi \LLamo
\label{1.65z}
\qqq
(in fact, one can show that $\LLv{0,1}{\uva}=\LLv{1,1}{\uva}=0$, so the
sum in \ex{1.65z} begins at $m=2$).
\end{remark}

In\cx{Ro2} we studied $\IrfvaNL$ when $\mfg=su(2)$. Since $su(2)$ has a
1-dimensional Cartan subalgebra, then in that case it is convenient to
introduce the variables $\ua=\scp{\uva}{\vl}$,
where $\vl$ is the single positive root of $su(2)$. We proved in\cx{Ro2} that
\qq
\IrsufaNL =\exp\lrbc{
{1\over 4\hb}\sijoxN \xxlij\,a_i a_j
} \Jcrv{\cL}{e^{\ua}}{e^{\hb}-1},
\label{1.65y}
\qqq
where $\JcrLth$ is defined by \ex{1.4}. By comparing \eex{1.62} and\rx{1.4}
through the relation\rx{1.65y} and taking into account that in the case of
$su(2)$ $\dguva=\ua$ and $\dgrho=1$, we find that
\qq
\pjbai\lrbcs{ \gQsuLuak}^{1/2}\, e^{-\LLrarkzos} = \AFLeua,
%
\label{1.65y1}
\qqq
or, in view of \ex{1.65z},
\qq
\pjbai\lrbcs{ \gQsuLuak}^{1/2}\;
\exp\lrbc{ - \smzi \LLamos } = \AFLeua.
\label{1.65y2}
\qqq
%

%

Let us introduce a notation
\qq
\FNgLua{} = \pldp
\AFLv{ e^{\scp{\uva}{\vl}} }.
\label{1.67}
\qqq
In subsection\rw{6xs.2} we will use the properties of the universal \urcc
invariant in order to prove the following
\begin{theorem}
\label{t1.2}
If $\AFLut{\not\equiv} 0$, then
the $\mfg$-based \urcc invariant
defined by \eex{1.62} and\rx{1.64} is well-defined and
does
not depend on a presentation of $\cL$ as a dotted Morse link and on a
choice of $\xk$ ($1\leq \xk\leq \xN$) thus being a topological invariant
of an oriented link $\cL$. It can be presented in a form
\qq
\IrhLuva  & = &
{\hbmnDp\over\dgrho}\;\;
{e^{\ohb\lkLuva}\over \FNgLua{} }
\;\;\exp\lrbc{\snoi {\pnQnLuva \over
\lrbs{\FNdg}^{3n} }\;\hb^n
},
\label{1.65y3}
\qqq
where the series
\qq
\pnQnLuva\in\IQuva
\label{1.69}
\qqq
are invariants of oriented links.
\end{theorem}

In subsection\rw{6xs.3} we will formulate the rationality conjecture for
links and derive from it the following corollary which sharpens the
formula\rx{1.65y3}.
\begin{corollary}[corollary of the rationality conjecture]
\label{3c.1}
The $\mfg$-based \urcc invariant can be presented in a form
\qq
\IrhLuva  & = &
{\hbmnDp\over\dgrho}\;\;
{e^{\ohb\lkLuva}\over \FNgLua{} }
\;\;\exp\lrbc{\snoi {
P_n\lrbc{\cL;(e^{\scp{\va_i}{\vl_j}})_{1\leq i\leq \xN\atop 1\leq j\leq r}
}
\over
\FNgLua{3n} }\;\hb^n
},
\label{6.65}
\qqq
where $\vl_1\ldtc\vl_r$ are simple roots of $\mfg$ and
\qq
P_n\lrbcs{\cL;(t_{ij})_{1\leq i\leq \xN\atop 1\leq j\leq r} } \in
\IQ\lrbs{  (t^{\pm 1}_{ij})_{1\leq i\leq \xN\atop 1\leq j\leq r}  }.
\label{6.66}
\qqq
\end{corollary}

\nsection{Graph algebras $\cBCX$, $\cDCX$ and $\cQDCX$}
\label{xs3}
\subsection{Root and Cartan edges of algebra $\cBCX$}
\label{3xs.1}

Now our main goal is to rewrite the calculations of the previous section solely
in terms of (1,3)-valent graphs while avoiding any use of Lie algebras.
However, it seems that the albebra $\cBL$ is not
sufficient for this purpose, because the integration in \ex{1.51} does not
go over
the whole Lie algebra $\mfg$ but only over its space of roots $\mfr$.
Therefore we will construct a bigger algebra $\cBCL$ which reflects a
distinction between the root space and Cartan subalgebra.

%
%
\begin{figure}[htb]
\begin{center}
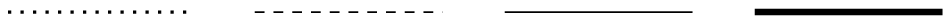
\end{center}
\caption{Types of edges in graphs of $\cBC$}
\label{f4.1}
\end{figure}

We begin by defining a bigger space $\tcBCX$. Similarly to $\cBX$,
it consists of
the formal linear combinations of (1,3)-valent graphs with cyclic order fixed
at 3-valent vertices, while the graphs which differ by the cyclic orientation
at one vertex, are considered opposites of each other (the AS relation).
As usual,
the 1-valent vertices of the graphs are labeled by the elements of $\xX$.
However, in contrast to $\cBX$, we allow two different basic
types of edges: Cartan
egdes (which are depicted by dotted lines) and root edges (which are depicted
by dashed lines). For convenience, we define two `auxiliary'
types of edges which are expressed in terms of the basic types.
First, a `total' edge (which we depict by a standard thin line) is a `sum' of
the root edge and the Cartan edge. In other words, a graph with a total edge is
a sum of two graphs in which the total edge is replaced by the root edge and by
the Cartan edge. Second, an `any' edge can be
of any type, we depict it by
a thick solid line. All types of edges are depicted in
\fg{f4.1}.

\begin{figure}[htb]
\begin{center}
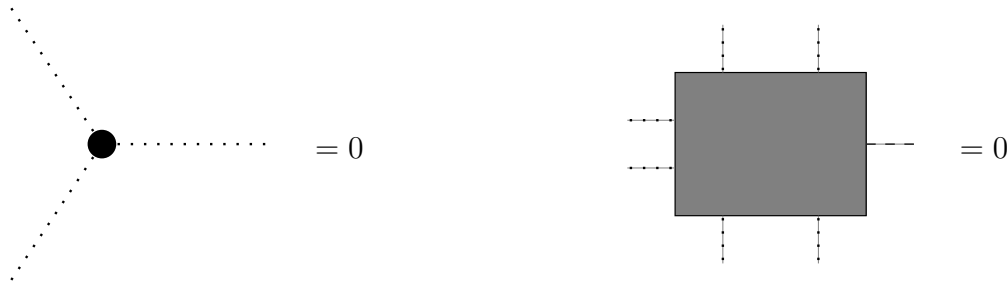
\end{center}
\caption{Cartan-commutator relations $\CCo$ and $\CCt$}
\label{f4.3}
\end{figure}

Next, we exlude the `Cartan-commutator' or $\CC$ graphs. A graph is called
$\CCo$ if it contains a 3-vertex which is incident to 3 Cartan edges. We call a
subgraph $D\p$ of $D$ \emph{proper}, if all vertices of $D\p$ are either
1-valent or 3-valent. A graph is called $\CCt$, if it contains a proper
subgraph, all of whose legs are Cartan except for one leg, which is root. A
graph is called $\CC$, if it is either $\CCo$ or $\CCt$, and we exclude all
these graphs from $\tcBCX$.

Eliminating $\CC$ graphs from the set of all graphs
can be equivalently described by setting these graphs to zero through factoring
over their span. Namely, let $\tcBCXCC$ be the span of $\CC$ relations, then
$\tcBCXoC$ is the span of all possible graphs with Cartan and
root edges and $\tcBCX$ is defined as its quotient over $\tcBCXCC$.
This point of view is useful, for example, when a $\CC$ graph
appears in an IHX relation.

For future calculations it is convenient to supplement the set of graphs
defining $\tcBCX$ with two more
objects, which are not actually graphs (although we also call them graphs):
a root circle $\rcrc$ and a Cartan circle $\lcrc$. They will be
needed in particular in \rw{4xs.5} in order to define determinants of
matrices in graph calculus.

%
The set of labels $\xX$ splits into three subsets:
$\xX = \xXr\bcup \xXc\bcup \xXt$ of root, Cartan and total labels.
If a root label is placed on a Cartan leg, or a Cartan
label is placed on a root
leg, then the graph is considered to be zero. A total label can be placed on
both types of legs.


Thus we defined the space $\tcBCX$. We define a familiar subspace $\tcBCXIHX$ as
a span of linear combination of the graphs of \fg{f4.2}. As we just
mentioned, if a $\CC$ graph appears there, then it is considered to be equal to
zero.
Now we define $\cBCX$ as a quotient space
\qq
\cBCX = \tcBCX / \tcBCXIHX.
\label{4.1}
\qqq
\begin{figure}[htb]
\begin{center}
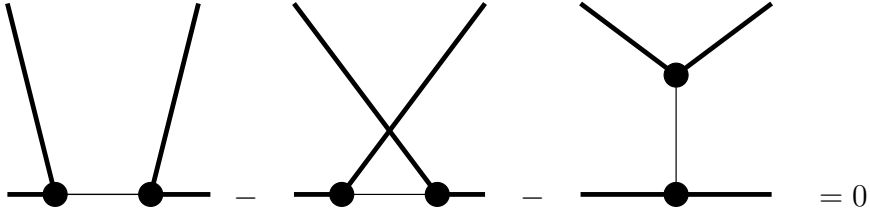
\end{center}
\caption{The IHX relation}
\label{f4.2}
\end{figure}
%

\begin{lemma}[Cartan commutativity]
\label{lcc}
Any graph having a 3-vertex incident to two Cartan edges, is $\CC$ and
therefore is zero in $\cBCX$.
\end{lemma}
\proof
If the third edge is Cartan, then the graph is $\CCo$, and if the third edge is
root, then the graph is $\CCt$.\qed

A $\mfg$-based weight system for $\cBX$ is a linear map
\qq
\xmap{\Tg}{\cBCX}{\oSgXT},
\label{4.3*}
\qqq
which is similar to the second map of\rx{1.11*}. In order to define it, we
modify\rx{1.11*2}.
Projectors $\prjr$ and $\prjh$ project $\stg$ naturally  onto $\str$ and
$\sth$. Thus we construct an element $\hiD\in\stgx$ by taking a tensor
product of $\himfr = \prjr(\himfg)$ for every root edge and $\himfh =
\prjh(\himfg)$ for every Cartan edge and then we change the
definition\rx{1.11*2}:
\qq
\fmfg^{\otimes\neDt}\otimes\hiD
\mathop{{\longmapsto}}^{\xCD} \tTg(D)\in \tmfc
\subset\tmfgx,
\label{4.3}
\qqq
where $\nehDo$ and $\nerDo$ are the numbers of Cartan and root legs of $\xD$.
The symmetrization over 1-valent
vertices which have the same label, projects $\tTg(D)$ to
$\Tg(D)\in\oSgXT$.

For special circle
graphs we define $\Tg(\rcrc) = \dim\mfr$ and $\Tg(\lcrc)=\dim\mfh$.

\subsection{Injections}
\label{3xs.1*}
In this subsection
we adopt the definition of $\tcBCX$, in which it
includes the $\CC$ graphs, so that now
\qq
\cBCX =  \tcBCX/(\tcBCXIHX + \tcBCXCC),  \qquad
\tcBCXCC = \tcBCXCCo + \tcBCXCCt.
\label{4.2y}
\qqq
where $\tcBCXCCo$ and $\tcBCXCCt$ are the spans of $\CCo$ and $\CCt$ graphs.
Also,
as usual,
$\tcBX$ is the space of all (1,3)-valent graphs before it was factored over
its IHX relation subspace $\tcBXIHX$.

Let us consider two injective linear maps
\qq
&\xmap{\tftr}{\tcBX}{\tcBCX} , \qquad \xX=\xXt,
\label{4.2t}
\\
&\xmap{\tfct}{\tcBX}{\tcBCX}, \qquad \xX=\xXc,
\label{4.2xt}
\qqq
which map each graph $\xD$ of $\tcBX$ into the same graph of $\tcBCX$, such that
all edges of $\tftr(\xD)$ are total, while all internal edges of $\tfct(\xD)$
are total and all legs are Cartan. Since
$\tftr(\tcBXIHX),\tfct(\tcBXIHX)\subset\tcBCXIHX$, then these maps extend to
\qq
&\xmap{\ftr}{\cBX}{\cBCX} , \qquad \xX=\xXt,
\label{4.2}
\\
&\xmap{\fct}{\cBX}{\cBCX}, \qquad \xX=\xXc.
\label{4.2x}
\qqq
The weight system\rx{4.3*} obviously commutes with $\ftr$.

\begin{theorem}
\label{t4.1}
The map $\ftr$ is an injection. If $\xX$ consists of a single element,
then the map $\fct$ is also an injection.
\end{theorem}

In order to prove this Theorem, we need a couple of lemmas about the structure
of $\cBX$.

\begin{lemma}
\label{brleg}
If a graph $\xD$ of $\tcBX$ contains a proper non-strut subgraph with a single
leg, then this grpah belongs to $\tcBXIHX$ and therefore is zero in
$\cBX=\tcBX/\tcBXIHX$.
\end{lemma}
\proof
We follow the $\cBX$ transformations of \fg{f4.3yx1}.
\begin{figure}[htb]
\begin{center}
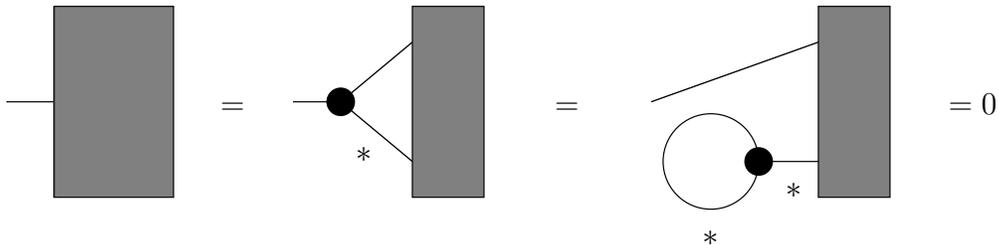
\end{center}
\caption{The IHX slide of the asterisk-marked edge onto itself}
\label{f4.3yx1}
\end{figure}
Since the subgraph is non-strut, then we can `pull out' a 3-valent vertex from
it (first equation). The IHX relations in $\cBCX$ allow us to slide the
asterisk-marked edge all the way through the grey box onto itself (second
equation). Finally, any graph which has a proper subgraph consisting of a
circle with a single leg is equal to zero in view of the AS relation (third
equation).\qed

Recall that an edge of a graph is called \emph{a bridge},  if it connects two
otherwise disconnected components (or, in other words, if its removal increases
the number of connected components of the graph). Note that a leg is always a
bridge.

\begin{lemma}
\label{3-bridge}
If $\xX$ consists of a single element, then
a graph of $\tcBX$, in which a 3-valent vertex is incident to 3 bridges,
belongs to $\tcBXIHX$ and therefore
is zero in $\cBX=\tcBX/\tcBXIHX$.
\end{lemma}
\proof
The proof relies on IHX-slides depicted in \fg{f4.3y1}.
\begin{figure}[htb]
\begin{center}
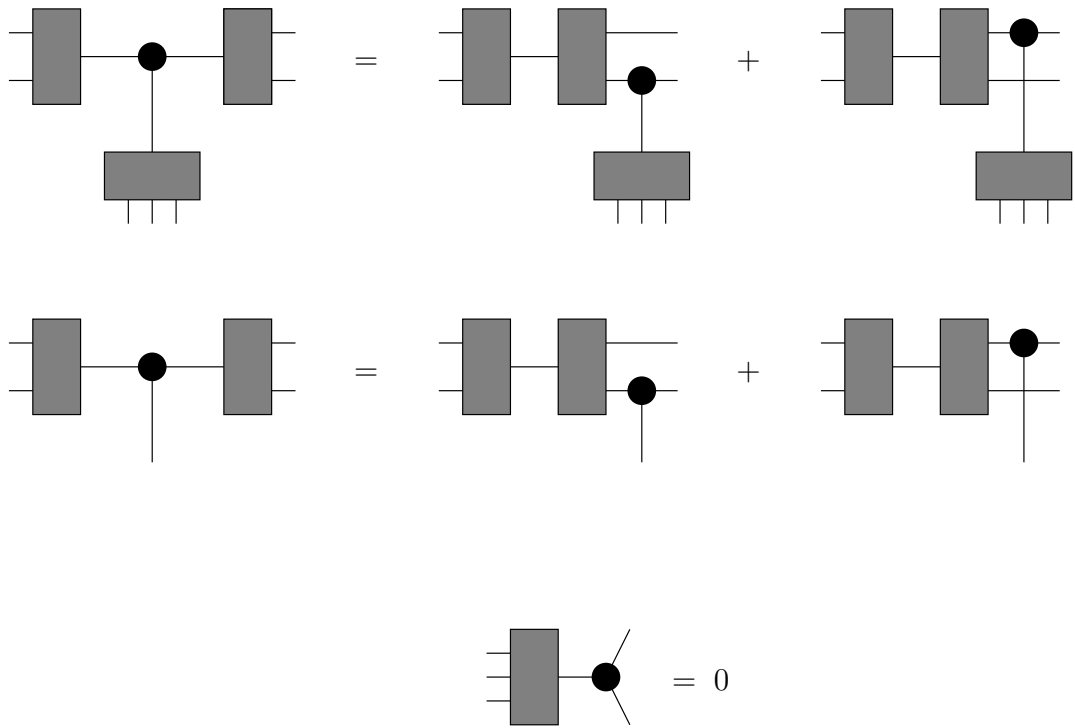
\end{center}
\caption{The IHX-slide of a bridge and the IHX slide of a leg}
\label{f4.3y1}
\end{figure}
First, we IHX-slide one of the bridges onto the legs. The summands in the \rhs
of the first equation of \fg{f4.3y1} have a 3-vertex, which is incident to a
leg and two bridges. Now we IHX-slide that leg onto the legs, as in the second
equation of \fg{f4.3y1}. Each graph in the \rhs of the second equation has a
3-vertex, which is incident to two legs. It is easy to see that, due to the
anti-symmetry of cyclic orientations at 3-vertices, a graph which has a
3-vertex incident to two legs of the same label, is equal to zero.\qed

\pr{Theorem}{t4.1}\footnote{We are very thankful to D.~Thurston, who
considerably streamlined the proof of this theorem.}
Our strategy for proving the injectivity of the maps $\ftr$, $\fct$ is to
define the left-inverses of $\tftr$ and $\tfct$:
\qq
\xmap{\tftri,\tfcti}{\tcBCX}{\tcBX},\qquad \tftri\circ\tftr = \tfcti\circ\tfct
= I.
\label{4.2*x1}
\qqq
Then we will show that
\qq
\tftri(\tcBCXIHX),\tftri(\tcBCXCC)\subset\tcBXIHX,\qquad
\tfcti(\tcBCXIHX),\tfcti(\tcBCXCC)\subset\tcBXIHX.
\label{4.2*x2}
\qqq
This implies that the maps $\tftri$ and $\tfcti$ can be extended to the maps
$\xmap{\ftri,\fcti}{\cBCX}{\cBX}$, which in view of \ex{4.2*x1} are
left-inverses of $\ftr$ and $\fct$:
\qq
\ftri\circ\ftr=\fcti\circ\fct = I.
\label{4.2*x3}
\qqq
The latter equation means that $\ftr$ and $\fct$ are injective and thus proves
the theorem.

Thus it remains to define the maps $\tftri$ and $\tfcti$ which satisfy
\eex{4.2*x1}, and prove the inclusions\rx{4.2*x2}.

In order to define the maps $\tftri$ and $\tfcti$ we pick a particular basis in
$\tcBCX$: it consists of graphs each of whose edges is either root or Cartan.
By definition, $\tftri$ maps a basis graph of $\tcBCX$ into the same graph of
$\tcBX$, if all of its edges are total, and maps it to zero otherwise. Then,
obviously, $\tftri\circ\tftr=I$. The space $\tcBCXIHX$ is a span of relations
of \fg{f4.2}, in which all participating graphs are basis. Depending on the
nature of their common edges, $\tftri$ either maps the graphs in each triplet
to zero, or to the corresponding graphs of $\tcBX$, which form the $\tcBXIHX$
triplets. Therefore, $\tftri(\tcBCXIHX)\subset\tcBXIHX$.
Since
$\tftri$ maps all graphs, which have at least one Cartan edge, to zero, then
$\tftri(\tcBCXCCo)=0$. Almost all $\CCt$ graphs also have Cartan edges and
hence are also mapped to zero. The only exception are those $\CCt$ graphs,
whose $\CCt$ subgraph has only one leg (which is root). But $\tftri$-image of
such graph is proportional to a graph of $\tcBX$ which, according to
Lemma\rw{brleg}, belongs to $\tcBXIHX$.

By definition, $\tfcti$ maps a basis graph of $\tcBCX$ into the same graph of
$\tcBX$ if all of its non-bridge edges are total, and maps it to zero
otherwise. Obviously, $\tfcti\circ\tfct=I$. Also
$\tfcti(\tcBCXIHX)\subset\tcBXIHX$ for the same reason as in the case of
$\tftri$. The $\tfcti$ image of a $\CCo$ graph can be non-zero only if all three
Cartan edges of its $\CCo$ subgraph are bridges, but according to
Lemma\rw{3-bridge}, the image of such graphs belongs to $\tcBXIHX$, so
$\tfcti(\tcBCXCCo)\subset\tcBXIHX$. Similarly, the $\tfcti$ image of a $\CCt$
graph can be non-zero only if all Cartan legs of its $\CCt$ subgraph are
bridges of the graph, but then the root leg of that subgraph must also be a
bridge. However, the $\tfcti$ image of a graph with a root bridge is zero
(indeed, we can present it as a difference of two graphs, in which this bridge
is total and Cartan, but their $\tfcti$ images are the same, since the action
of $\tfcti$ on basis elements does not depend on the nature of their bridge
edges).\qed

\subsection{Symmetric algebras of cohomologies of graphs}
\label{3xs.2}

For our purposes it will be more convenient to use a different
presentation for the space $\cBCX$. Let $\Xp\subset\xXc$ be a subset of
Cartan labels.
We are going to convert Cartan legs
with labels of $\Xp$ into the elements of cohomology of (1,3)-valent graphs
in exactly the same way as we did it with all legs of the graphs of $\cB$
in\cx{RoC}. In other words, we will define a new algebra $\cDCX$ and
then show that it is canonically isomorphic to $\cBCX$. The only
difference between our case and that of\cx{RoC} is that the space $\cD$ is
based there on 3-valent graphs and ours will be based on (1,3)-valent
graphs, while the leg commutativity lemma~3.3 of\cx{RoC} is replaced by
Cartan leg commutativity lemma\rw{lcc}. Thus all proofs remain exactly the
same, but we will repeat them here for convenience.

Let $\xD$ be a general graph. Thinking of it as a $CW$-complex, its
boundary $\del\xD$ being the set of 1-valent vertices, we may consider a
rational relative cohomology space $\HobD$. An oriented edge $\ed$ of $\xD$
represents an integral element
$\hed\in\HobD$, a pairing between $\hed$ and a cycle of $\xD$ being the
coefficient at $\ed$ in a presentation of the cycle as a linear
combination of edges. Moreover, $\HobD$ can be presented as a quotient of
a linear space, whose basis elements are oriented edges of $\xD$. Let us
fix an orientation of the edges of $\xD$ and let $\bED$ be the set of
these oriented edges. $\ED$ denotes the corresponding linear space, whose
basis is formed by the elements of $\bED$
(we also assume that edges with opposite orientation represent the
opposite elements of $\ED$).  For a vertex $\vd$
and an oriented edge $\ed$ we define an incidence number $\inave$ to be 0
if $\ed$ is not incident to $\vd$ or if it is a loop attached to $\vd$,
$\inave=1$ if $\ed$ goes into $\vd$ and $\inave=-1$ if $\ed$ goes out of
$\vd$. Then
\qq
\HobD = \ED/ \EvD,
\label{grco}
\qqq
where
\qq
\EvD =
\bigoplus_{m\geq 2}
\spanv{\sum_{\ed\in\bED}\inave\,\ed\,| \vd\in\bVmD},
\label{grco1}
\qqq
and $\bVmD$ denotes the set of $m$-valent vertices of $\xD$. We excluded
the 1-valent vertices from the \rhs of \ex{grco}, because they represent
the boundary of $\xD$ and we are interested in the relative cohomology.

Now let $\xD$ be a (1,3)-valent graph of $\cBCX$ which has no $\Xp$
Cartan legs and
let $\HCobD\subset\HobD$ denote a linear span of Cartan edges of $\xD$.
We will use a quotient space
%
\qq
\HCfobD=\HobD/\HCobD.
\label{4.3*1}
\qqq
It has an alternative description. Let $\txD$ be the graph (or,
rather,
the cell complex) constructed from $\xD$ by removing Cartan edges.
The $\CCt$ condition guarantees
that no 3-valent vertices of $\xD$ are incident to
two Cartan edges. Therefore all 1-valent vertices of $\txD$ come from
those of $\xD$ and as a result in view of \ex{grco}
\qq
\HCfobD = \Hobv{\txD}.
\label{4.3*2}
\qqq

\begin{remark}
\rm
\label{r4.1*1}
By using the cohomology spaces that we have just introduced, we can rewrite
condition $\CCt$ of \fg{f4.2}. Namely, a graph $D$ contains a proper
$\CCt$ subgraph
iff it has a root edge $\ed$
such that $\hed\in\HCobD$.
\end{remark}

Let $\GD$ be the symmetry group of $D$: $\GD$ maps edges to edges of the same type,
3-vertices to 3-vertices and 1-vertices to 1-vertices of the same label.
The elements of $\GD$ may change orientation
at 3-valent vertices. $\GD$ acts on $\HobD$,
and $\HCobD$ is invariant under that
action, so the action of $\GD$ on $\HCfobD$ is well-defined. This action can be
extended to the algebra
$\SHCfobDtXp$, which is the $\orXp$-th tensor power of the symmetric algebra
of $\HCfobD$. We modify this action of $\GD$ by
multiplying the action of $g\in\GD$ by a sign factor $(-1)^{|g|}$, where $|g|$
is the number of 3-vertices of $\xD$ whose cyclic order was changed by $g$.
Now we define the space $\cXpHD$ associated to a graph $\xD$ as
the $\GD$-invariant part of $\SHCfobDtXp$
%
\qq
\cXpHD = \Invrss{\SHCfobDtXp}{\GD},
\label{4.5}
\qqq
if $\xD$ does not contain proper $\CC$ subgraphs,
and $\cXpHD$ is by definition
0-dimensional otherwise. In the former case there is also a natural
symmetrization projector
\qq
\mpcd{\SHCfobDtXp}{\PGD}{\cXpHD}.
\label{4.5m}
\qqq

We define the space $\cH_{\Xp}$ also for four
special `graphs'. The space $\cXpHD$ for the first two of them is defined
in accordance with \ex{4.5}.
The first special graph is a root circle  $\rcrc$, so
$\HCfobv{\rcrc}=\Hov{\rcrc}$ and the symmetry group is $\sG{\rcrc}=\{1,-1\}$,
its element $-1$ flipping the circle. Therefore $\cXpHc=\SPL^{\rm even}\,
\IQ^{\orXp}$.
The second special graph is a Cartan circle $\lcrc$, so $\dim\HCfobv{\lcrc}=0$
and by definition
$\cXpHl=\IQ$.
The third special graph is an empty unlabeled dot $\blt$, it
consists of a single vertex and has no edges. We define
$\cXpHb = \Srs^2(\IQ^{|\Xp|})$. The fourth special graph is a filled
labeled dot
$\blxb$, it consists of a single vertex labeled by elements of
$\xXc\setminus\Xp$ and we define $\cXpHx=\IQ^{\orXp}$.
If $\xD$ is a disjoint union of $m_1$ graphs $\rcrc$, $m_2$ graphs
$\lcrc$, $m_3$ graphs $\blt$, $m_4(\bl)$ graphs $\blxb$,
$\bl\in\xXc\setminus\Xp$
and a graph $\xD\p$ of $\cBCX$, which has no
$\Xp$ Cartan legs, then we define
\begin{multline}
\cXpHD = \Srs^{m_1}\cXpHc \otimes \Srs^{m_2}\cXpHl
\otimes \Srs^{m_3}\cXpHb
\\
\otimes
\lrbc{ \bigotimes_{\bl\in\xXc\setminus\Xp} \Srs^{m_4(x)}\cXpHx}
\otimes \HCfobv{\xD\p}.
\end{multline}
Let $\gDXp$ be the set of all such graphs $\xD$. Then we define the space
\qq
\tcDCXp = \boDgDXp \cXpHD.
\label{4.7}
\qqq

\begin{remark}
\label{r4.1*}
\rm
We will abbreviate the notations $\cXpHD$, $\gDXp$ and $\tcDCXp$ down to $\cHD$,
$\gDX$ and $\tcDCX$ when it is clear that a particular
set $\Xp$ has been chosen.
\end{remark}

$\tcDCX$ has a commutative algebra structure. In order to define it, we
observe that if $\xDot$ is a disjoint union of the graphs $\xDo$ and
$\xDt$, then
$\HCfobDot=\HCfobDo\oplus\HCfobDt$
and as a result
$\SHCfobv{\xDot} = \SHCfobv{\xDo}\otimes\SHCfobv{\xDt}$ and
$\SHCfobXpv{\xDot} = \SHCfobXpv{\xDo}\otimes\SHCfobXpv{\xDt}$. Therefore
$\GDot$ acts on $\SHCfobXpv{\xDo}\otimes\SHCfobXpv{\xDt}$ and since
$\GDo\times\GDt\subset\GDot$, then in view of the definition\rx{4.5},
$\cHDot\subset\cHDo\otimes\cHDt$ and the symmetrization projector of the
type\rx{4.5m} establishes a projection
\qq
\mpcd{\cHDo \otimes\cHDt}{\PGDot}{ \cHDot}.
\label{4.7m}
\qqq
The latter map defines the multiplication structure of $\tcDCX$.

The spaces $\cXpHD$ inherit $|\Xp|$ independent gradings $\dgsj$
($1\leq j\leq |\Xp|$) from the individual symmetric algebras $\SHCfobD$. We
define $|\Xp|$ parity operators $\xPsj$ acting on an element $x\in\cXpHD$ with
definite degree $\dgsj$ as
\qq
\xPsj(x) = (-1)^{\dgsj(x)}\,x.
\label{4.7*}
\qqq
Also, we define a grading $\dgth$ of $\tcDCXp$ by assigning
\qq
\dgth \cXpHD = \echi(\xD)
\label{1.11xx}
\qqq
(\cf \ex{1.11x}). All gradings $\dgsj$ as well as $\dgth$ are respected by the
multiplication\rx{4.7m}.

Suppose that $D\in\gDXp$ and $D\p$ is obtained from $D$ by replacing one of its
root edges $\ed$
with a Cartan edge (if $\ed$ is a leg, then we assign to it a color from
$\xX\setminus\Xp$).
Obviously, $\txD\p\subset\txD$,
%
%
so there is a natural map
%
\qq
\begin{CD}
\HCfobD @>\fed>> \HCfobDp,
\end{CD}
\label{4.7*1}
\qqq
which can be extended to an algebra homomorphism
$\SHCfobD\longrightarrow\SHCfobDp$
and further to
\qq
\begin{CD}
\cHD @>\fed>> \cHDp,
\end{CD}
\label{4.8}
\qqq
which
is a combination of the map $\xmap{}{\cHD}{\SHCfobDtXp}$ and a symmetrization
projection $\xmap{\PGDp}{\SHCfobDtXp}{\cHDp}$.  If $D\p$ is $\CC$, then by
definition $\fed=0$.

The map\rx{4.8} allows us to define the space $\cHD$ for a graph $\xD$
which, in
addition to root and Cartan edges, may also have total
edges. In order to simplify notations, suppose that $\xD$ has exactly
one total edge
$\ed$. Let $\Dr$ and $\Dc$ be the graphs constructed from $\xD$
by replacing the total edge by a root edge and by a Cartan edge.
We define $\cHD$ as a subspace of $\cHDrc\subset\tcDCX$:
\qq
\cHD =
\left\{
\begin{array}{ll}
\begin{CD}
\cHDc @>0\oplus I>> \cHDrc, &
\qquad\mbox{if $\Dr$ contains a $\CC$ subgraph,}
\end{CD}
\\
\begin{CD}
\cHDr@>I\oplus \fed>>\cHDrc &\qquad\mbox{otherwise},
\end{CD}
\end{array}
\right.
\label{4.9}
\qqq
where $I$ is the identity map and $\fed$ is the map\rx{4.8}. Thus $\cHD$ is
naturally injected into $\tcDCX$.


Now we define an ideal $\tcDCXIHX\subset\tcDCX$.
Consider a (1,3,4)-valent
graph $\xD$
which has a single 4-valent vertex. We assume that a cyclic order is fixed
at its 3- and 4-valent vertices, that 1-valent vertices are colored by the
elements of $\xX$, that its edges are either root or Cartan and that it has no
$\Xp$ Cartan legs. We can `resolve' a 4-valent vertex into two 3-valent vertices
connected by a total edge $\ed$ in three different ways depicted in \fg{f4.2}.
We denote the corresponding graphs of $\gDX$ as $\xD_i$ ($i=1,2,3$). The graphs
$\xDir$ and $\xDic$ are constructed from $\xD_i$ by declaring $\ed$ to be
either root or Cartan. For any $i$, $\txD$ can be constructed from
$\txDir$ by contracting the new root edge and also from $\txDic$ by gluing
together two 2-vertices, which were 3-vertices of $\xDic$ incident to the
Cartan edge. Therefore, there are natural maps
\qq
&\mpcd{\HobtD}{\fir}{\HobtxDir},
\nonumber\\
&\mpcd{\HobtD}{\fic}{\HobtxDic},
\label{4.10*}
\qqq
which can be extended to symmetric
algebra homomorphisms composed with symmetry projectors
%
\qq
&\mpcd{\SHCfobv{\xDir}}{\PGDv{\xDir}}{\cHv{\xDir}}
\nonumber
\\
&\mpcd{\SHCfobv{\xDic}}{\PGDv{\xDic}}{\cHv{\xDic}}.
\label{4.10*1}
\qqq
Since, according to the definition\rx{4.9}, the spaces $\cHv{\xD_i}$ are either
$\cHv{\xDir}$ or $\cHv{\xDic}$ (depending on whether $\xDir$ contains a $\CC$
subgraph), the maps\rx{4.10*1} define maps
\qq
\begin{CD}
\SHCfobDtXp @>{f_i}>> \cHv{D_i},\qquad i=1,2,3.
\label{4.10*1x}
\end{CD}
\qqq
Note that if $\xDir$ does not contain a $\CC$ subgraph, then
$\fic=\fed\circ\fir$ and therefore a composition of maps
$$
\begin{CD}
\SHCfobDtXp @>{f_i}>> \cHv{D_i} \hookrightarrow \cDCX
\end{CD}
$$
is equal to the sum $\fir+\fic$ of maps\rx{4.10*1}.

The ideal $\tcDCXIHX$ is defined as a
span of all elements $f_1(x)-f_2(x)-f_3(x)$ for all graphs $D$ and all elements
$x\in\SHCfobDtX$. Finally, we define the algebra $\cDCX$ as a quotient
\qq
\cDCX = \tcDCX/
\tcDCXIHX.
\label{4.10}
\qqq

\subsection{Converting Cartan legs into cohomology}
\label{3xs.2x}

Now we will show that $\cDCX$ is canonically isomorphic to the algebra
$\cBCX$. First of all, we define a map from the graphs of $\tcBCX$ to the
graphs of $\gDXp$, which removes the $\Xp$ legs. A Cartan strut with two
$\Xp$ labels maps to $\blt$. A Cartan strut with one $\Xp$ label and one
$\xXc\setminus\Xp$ label $\bl$
maps to $\blxb$. In any other graph, $\Xp$ legs are
attached to root edges. We remove these legs and `dissolve' their incident
vertices, thus obtaining the graphs of $\gDXp$. We call them \emph{\Xpfr s} of
the graphs of $\cBCX$. A removal of $\Xp$ legs preserves the homotopy
class of the $CW$-complex and hence its Euler characteristic, so if $\xD\p$ is
a \Xpfr\ of $\xD$, then, according to definitions\rx{1.11x} and\rx{1.11xx},
$\dgth{\xD}=\dgth{\xD\p}$.

Let $\ttcBCXD\subset\tcBCX$ denote the span of all
graphs whose \Xpfr\ is $\xD\in\gDXp$. Then
$\tcBCX = \bopgD \ttcBCXD$.

Let $\tcBCXIHXi\subset\tcBCXIHX$ ($0\leq i\leq 4$)
denote a span of IHX relations of \fg{f4.2}, in which $i$ connecting edges
of the IHX graphs are $\Xp$ legs. If at least 3 connecting edges are
Cartan legs, then the corresponding IHX graphs are $\CC$, so
$\tcBCXIHXth=\tcBCXIHXf=0$ and $\tcBCXIHX = \bigoplus_{i=0}^{i=2}
\tcBCXIHXi$. All graphs in 2 and 3 $\Xp$-legged $\tcBCX$ IHX triplets
reduce to the same $\gDXp$ graph after the removal of $\Xp$ legs.
Therefore, if we define
\qq
\tcBCXD & = & \ttcBCXD/\tcBCXIHXtD,
\label{4.10x1}
\qqq
where $\tcBCXIHXtD = \ttcBCXD\cap\tcBCXIHXt$, and if we also define
\qq
\cBCXD & = & \ttcBCXD/(\ttcBCXD\cap(\tcBCXIHXo + \tcBCXIHXt))
\nonumber
\\
& = & \tcBCXD/\tcBCXIHXoD,
\label{4.10x2}
\qqq
where $\tcBCXIHXoD = (\tcBCXIHXo/(\tcBCXIHXo\cap\tcBCXIHXt))\cap\cBCXD$,
then
\qq
\tcBCX/\tcBCXIHXt & = & \bopgD \tcBCX(\xD),
\label{4.10x3}
\\
\tcBCX/(\tcBCXIHXo+\tcBCXIHXt) & = & \bopgD \cBCX(\xD)
\label{4.10x4}
\qqq
and therefore
\qq
\cBCX = \lrbc{\bopgD \cBCX(\xD)}/\lrbcs{\tcBCXIHXz/
((\tcBCXIHXo+\tcBCXIHXt)\cap\tcBCXIHXz)}.
\label{4.10x5}
\qqq

In each triplet of IHX graphs with two connecting edges being $\Xp$ legs,
one graph has these legs attached to the same 3-vertex, and so it is
$\CC$. The remaining two graphs differ only in order in which the two
$\Xp$ legs are attached to the (root) edge. Therefore the graphs which
differ only in the order in which their $\Xp$ legs are attached to root
edges, are equal in $\tcBCX/\tcBCXIHXt$.

Next, we establish the isomorphisms $\xmap{\mF}{\cHD}{\cBCXD}$
for all $\xD\in\gDXp$. First, we define it for $\blt$ and $\blxb$. According to
our definition, $\cXpHx=\IQ^{\orXp}$. We choose the basis $e_a$, $a\in\Xp$ in
this space and then define
$$e_a \xmapt{\mF} \cstv{a}{\bl}.$$ Similarly, we choose the basis $e_{a_1}
e_{a_2}$, $a_1,a_2\in\Xp$ of $\cXpHb = \Srs^2(\IQ^{|\Xp|})$ and define
$$e_{a_1} e_{a_2} \xmapt{\mF}\cstv{a_1}{a_2}.$$ For the rest of the graphs
$\xD\in\gDXp$ we define $\mF$ with the help of the following
commutative diagram:
\qq
\xymatrix{
%
%
%
\\
&\SEtDX \ar[dr]^{\tmF} \ar[dl]_{\PGD}
\ar[dd]|\hole^(.7){\xfo}
\\
\SEtDG 
\ar[rr]^(.7){\tmF}
\ar[dd]^{\xfo} & & \tcBCXD   \ar[dd]^{\xft} 
\ar@{-->}@/_{2pc}/[ul]_-{\tmFsi}
\\
&\SHCfobDtXp \ar[dr]^{\mF} \ar[dl]_{\PGD}&
\\
\cHD \ar[rr]^{\mF} & &\cBCXD  
}
\label{1.17x}
\qqq
Recall that $\EtD$ is a linear space, whose basis vectors are oriented root
edges of $\xD$. We define the map
\qq
\xymatrix@1{\SEtDX\ar[r]^-{\tmF} & \tcBCXD}
\label{isom0}
\qqq
by its
action on the monomials of $\SEtDX$: a monomial
$\prod_{\ed\in\bEtD} \prod_{a\in\Xp}\lrbc{\eda}^{\mea}$
is
mapped to a graph constructed from $\xD$ by attaching $\mea$ legs of color
$a\in\Xp$ to the right side of each oriented root edge $\ed$. This map is
well-defined, because the order of attaching $\Xp$ legs to a root edge does not
matter after we take a quotient over $\tcBCXIHXt$. The map
\qq
\xymatrix@1{\SEtDG\ar[r]^-{\tmF} & \tcBCXD}
\label{isom1}
\qqq
is the restriction of the former
map to the subspace $\SEtDG\subset\SEtDX$.

\begin{lemma}
\label{iso1}
The map\rx{isom1} is an isomorphism.
\end{lemma}
\proof
Suppose that $\xD\p$ is a graph of $\tcBCXD$. Generally, there may be many ways
of identifying its \Xpfr\ $\xD$ with a `strandard copy' of $\xD$, but they all
differ by a composition with the elements of $\GD$. Let $\xxs$ be a way of
identifying the \Xpfr s of every graph of $\tcBCXD$ with the standard copy of
$\xD$. Then we can construct the left-inverse
$\xymatrix@1{\tcBCXD\ar[r]^-{\tmFsi} & \SEtDX}$ of $\tmF$. The composition
$\tmFi=\PGD\circ\tmFsi$ does not depend on the choice of $\xxs$ and serves as
the left-inverse of\rx{isom1}. Since the
images of $\tmFsi$ for all possible identifications $\xxs$ span the whole
$\SEtDX$, then $\tmFi$ is surjective, hence\rx{isom1}
is an isomorphism.\qed

Let us check that
the spaces of the lower triangle of the diagram\rx{1.17x} are the quotients of
the spaces of the upper triangle. According to the definition\rx{4.10x2},
$\cBCXD = \tcBCXD/\tcBCXIHXoD$.
According to\rx{grco},
$\HCfobD=\EtD/\EvtD$, so $\SHCfobDtXp = \SEtDX/\IvtD$, where
$\IvtD\subset\SEtDX$ is the ideal generated by $\EvtD$. Since $\EvtD$ is
invariant under the action of $\GD$, then
\qq
\cHD = \Invrss{\SHCfobDtXp}{\GD}  & = &
\SEtDG/\lrbcs{\SEtDG\cap\IvtD}.
\nonumber
\\
& = & \SEtDG/\Invrs{\IvtD}{\GD}.
\qqq
This follows from a simple
\begin{lemma}
Suppose that a finite group $G$ acts on a linear space $V$ and its subspace
$W\subset V$ is invariant under this action. Then $\Invrs{V/W}{G} =
V^G/(V^G\cap W)\subset V/W$.
\end{lemma}
\proof
Consider a commutative diagram
\qq
\xymatrix{
V \ar[d]^f\ar[r]^{P_G} & V^G \ar[d]^f &\hspace{-1.8in}\subset V
\\
V/W \ar[r]^-{P_G} & V^G/(V^G\cap W) &\subset V/W
}
\qqq
where $f$ maps $V$ naturally onto its quotient $V/W$ and $P_G$,
as usual, is a $G$-symmetrization projector.
Since all elements of $V^G/(V^G\cap W)$ are $G$-invariant, then the lower $P_G$
is surjective,
hence $V^G/(V^G\cap W) = (V/W)^G$.\qed

Thus, the vertical maps $\xfo,\xft$ of the diagram\rx{1.17x} are surjections of
linear spaces onto their quotients. It is easy to see, that $\tmF$ maps the
$\EvtD$-generated ideal $\IvtD$ into the subspace
$\tcBCXD\cap\tcBCXIHXo\subset\tcBCXD$. Therefore $\tmF$ of\rx{isom1}
descends to the map
%
$\xymatrix@1{\SHCfobDtXp\ar[r]^-{\mF} & \tcBCXD}.$
%
\begin{lemma}
\label{iso2}
The restriction of $\mF$
\qq
\xymatrix@1{\cHD\ar[r]^-{\mF} & \cBCXD}
\label{4.3y1}
\qqq
is an isomorphism.
\end{lemma}
\proof
In view of Lemma\rx{iso1}, it is sufficient to show that the inverse map
$\tmFi$ descends to a map
$\xymatrix@1{\cBCXD\ar[r]^-{\mFi} & \cHD}.$ The latter statement would follow
from the fact that $\tmFi(\ker \xft)\subset \ker \xfo$. In order to see this,
observe, that $\ker\xft = \tcBCXIHXoD$. Consider a triplet of IHX graphs
from $\tcBCXIHXoD$. Suppose that an identification assignment $\xxs$ identifies
their \Xpfr s consistently with the standard copy of $\xD$. Then $\tmFsi$ maps
their IHX linear combination into $\IvtD$. Since $\tmFi = \PGD\circ\tmFsi$ does
not depend on the choice of $\xxs$, then all $\tcBCXIHXoD$ triplets are mapped
into $\Invrs{\IvtD}{\GD}=\ker \xfo$.\qed

Thus, if we apply a sum of individual maps $\mF$ to the direct sum
$\boDgDXp \cHD$, then in view of relations\rx{4.7} and\rx{4.10x4} we establish
an isomorphism
\qq
\xymatrix@1{\tcDCX\ar[r]^-{\mF}& \tcBCX/(\tcBCXIHXo+\tcBCXIHXt).}
\label{4.3y2}
\qqq

We leave it for the reader to verify the following
\begin{lemma}
\label{iso3}
The map\rx{4.3y2} establishes the isomorphism between the IHX subspaces\\
$\tcDCXIHX\subset\tcDCX$ and
$\tcBCXIHXz/((\tcBCXIHXo+\tcBCXIHXt)\cap\tcBCXIHXz)\subset\tcBCX/(\tcBCXIHXo+\tcBCXIHXt)$.
\end{lemma}
\proof
Recall the definition of $\tcDCXIHX$. Let $\xD$ be a (1,3)-valent graph without
$\Xp$ legs and with a single 4-vertex. Consider a monomial
\qq
\prod_{\ed\in\bEtD}
\prod_{a\in\Xp}\lrbc{\heda}^{\mea}\in\SHCfobDtXp
\label{mon1}
\qqq
and let $\xD\p$ be a graph
representing its $\mF$ image (in other words, $\xD\p$ is a result of attaching
$\Xp$ legs to the root edges of $\xD$ according to the powers $\mea$). Let us
compare the IHX elements constructed from the monomial and from $\xD\p$.

The root edges of $\xD$ are naturally identified with the root edges of the
graphs $\xDir$ and $\xDic$. The naturality of the maps\rx{4.10*} guarantees
that if a root edge $\ed$ of $\xD$ corresponds to a root edge $\ed\p$ of
$\xDir$ (or $\xDic$), then $\fir(\hed) = \hed\p$ (or $\fic(\hed)=\hed\p$).
This means that if we apply the maps\rx{4.10*1x} to the monomial\rx{mon1} and
then apply $\mF$ to the resulting monomials, the result will be the graphs
constructed by IHX resolving the 4-vertex of $\xD\p$. This proves
that
\qq
\mF(\tcDCXIHX)=\tcBCXIHXz/((\tcBCXIHXo+\tcBCXIHXt)\cap\tcBCXIHXz).
\label{inc1}
\qqq
\qed
%
%
%

Now, in view of the definition\rx{4.10} and relation\rx{4.10x5} we come to the
following
\begin{theorem}
\label{t4.2}
The map $\tmF$ of\rx{isom0} descends to the isomorphism
\qq
\begin{CD}
\cDCX @>{\mF}>> \cBCX,
\end{CD}
\label{4.11}
\qqq
which preserves the grading $\dgth{}$.
\end{theorem}

\subsection{Edge-related denominators}
\label{3xs.3}

For the purpose of describing the universal \urcc invariant we have to consider
a modified algebra $\cQDCXp$
(which we abbreviate down to $\cQDCX$). This algebra is similar to
$\cDCX$, except that in constructing the analogs of the spaces $\cHD$ we allow
the `edge related' denominators. As a result, whenever we define the analogs of
operations in $\cDCX$, we have to make sure that we do not produce zeroes in
these denominators.

We begin by defining the analogs of the spaces $\cHD$, which we call
$\cQHD$. Let $D\in\gDXp$ be a regular graph. If it
contains a $\CC$ subgraph, then $\cQHD$ is by
definition zero-dimensional. Suppose now that $\xD$ does not contain $\CC$
subgraphs. According to Remark\rw{r4.1*1}, this means that for any root edge
$\ed$ of $\xD$, $\hed\not\in\HCobD$. Therefore, as an element of the
quotient space $\HCobD$, $\hed\neq 0$, hence any polynomial of $\hedx$ is
non-zero as an element of $\SHCfobDtXp$. Thus we can define $\QeSHCfobDtX$
as an extension of the algebra $\SHCfobDtX$, which is a linear span of
fractions $f\over g$, where $f\in\SHCfobDtX$ and
%
%
\qq
g=\preDr \pexXp,
\label{4.11*}
\qqq
where $\pexXp\in\SHCfobDtXp$
are formal power series of
$\hedxXp$ (that is, each formal power
series depends on a particular root edge of all possible $\Xp$ labels).
We denote the $\GD$-invariant part of this algebra as
\qq
\cQHDXp  = \Invr{\QeSHCfobDtX}{\GD}
\label{4.12}
\qqq

Let us consider a couple of useful examples of these spaces.
First, let $D$ be a root strut with different colors
at 1-vertices: $D = \rstxot$. The single edge of $D$ is a natural
basis element of $\HobD$, $\HCobD$ is trivial and $\GD$ is also trivial.
Hence there is a canonical isomorphism
\qq
\cQHv{\rstxot} = \QIQuua,
\label{4.13**}
\qqq
where $\bfa = (a_1\ldtc a_{|\Xp|})$
and $\cQ(R)$ denotes a field of quotients of a ring $R$.
Consider a map $D\rightarrow D$ which
flips the root strut.
This map is not a symmetry of the labeled graph, but it acts on
$\HobD$ by changing the sign of the basis element. We extend its action to the
whole space $\cQHD$. We denote the image of $l\in\cQHD$ as $\stl$. If we
represent $l$ as an element of $\QIQuua$, then obviously
\qq
\stl(\ua) = l(-\ua).
\label{4.13**x}
\qqq

Our second example is a root strut with the same colors at 1-vertices:
$D= \rstx$. The difference with the first example is that the graph
symmetry group is non-trivial: $\GD=\{1,-1\}$ where the $-1$ element
reverses the edge. Therefore, $\cQHD$ is just the (simultaneously) even
part of the algebra\rx{4.13**}
\qq
\cQHv{\rstx} = \evQIQuua.
\label{4.13**1}
\qqq

It remains to define $\cQHDXp$ for special `graphs'. We define
\qq
\cQ\HCfobv{\lcrc} = \HCfobv{\lcrc},\qquad
\cQ\cXpHb = \cXpHb,\qquad
\cQ\cXpHx = \cXpHx.
\label{4.13**x1}
\qqq
Since the cohomology and symmetry of the root circle $\rcrc$ are exactly
the same as those of
$\rstx$, then we define
\qq
\cQHv{\rcrc} = \evQIQuua.
\label{4.13**2}
\qqq
%


\begin{remark}
\rm
\label{r4.sp1}
There is an obvious inclusion and projection which relate the
algebras\rx{4.13**}, \rx{4.13**1} and\rx{4.13**2}
\qq
\evQIQuua\hookrightarrow \QIQuua,\qquad \QIQuua\longrightarrow \evQIQuua.
\label{4.13**3}
\qqq
\end{remark}
\begin{remark}
\label{r4.sp2}
\rm
Actually, we assume that the algebras $\cQHv{\rcrc}$ and $\cQHv{\lcrc}$ are
extended by the
logarightms of $\evQIQuua$ and $\IQ$.
\end{remark}

Now we define the analog of $\tcDCXp$
\qq
\tcQDCXp = \boDgDXp \cQHDXp.
\label{4.13}
\qqq
As usual, we will tend to drop $\Xp$ from these notations.

Similarly to $\tcDCX$, $\tcQDCX$ has a
structure of a graded commutative algebra, its multiplication and grading
$\dgth$ being defined by the obvious analogs of\rx{4.7m} and\rx{1.11xx}.
Although the spaces $\cQHD$ do not have
gradings $\dgsj$ (because of denominators), we can
still define the action of parity operators $\xPsj$ of \ex{4.7*} on $\tcQDCX$
by the
formula
\qq
\xPsj(x/y) = \xPsj(x)/\xPsj(y),
\qquad x,y\in\cHD.
\label{4.13**3x}
\qqq
%


The definition of
$\tcQDCXIHX$ involves a total edge and therefore requires some care.
Let $\ed$ be a root edge of a graph $D$. Let $D\p$ denote a graph constructed
from $D$ by replacing $\ed$ with a Cartan edge.
We call an element $x\in\cQHD$
\emph{non-singular at a root edge $\ed$ of $D$} if
either $D\p$ contains a $\CCt$ subgraph or
$x$ can be presented as a sum of fractions, whose denominators do not
belong to the kernel of the unsymmetrized version
of the map $\fed$ of\rx{4.8}:
\qq
\begin{CD}
\SHCfobDtXp   @>\fed>> \SHCfobDptXp.
\end{CD}
\label{4.13*1}
\qqq
\begin{lemma}
\label{l4.1*}
If for a root edge $\edz$ of $\xD$, an element $x\in\cQHD$ has a presentation as
a sum of fractions such that
\qq
\pezxXp\Big|_{\hedza=0\;\;{\rm for\;\; all}\;\; a\in\Xp} \neq 0
\label{4.13*}
\qqq
in all their denominators\rx{4.11*}, then $x$ is non-singular at $\edz$.
\end{lemma}
\proof
A product\rx{4.11*} belongs to the kernel of\rx{4.13*1}, if at least one of the
polynomials $\pexXp$ belongs to it, and that happens only if the corresponding
edge $\hed$ belongs to the kernel of the map\rx{4.7*1}. That would mean
that $\hed\in\HCobDp$, where $D\p$  is the graph constructed from $D$ by
declaring $\edz$ to be Cartan. But then, according to Remark\rw{r4.1*1},
$D\p$ contains a $\CCt$ subgraph, so $x$ is still non-singular at $\edz$.
\qed

All elements of $\cQHD$, which are non-singular at $\ed$, form a
subalgebra $\cQHDe\subset\cQHD$.
The map\rx{4.8} extends to $\cQHDe$
\qq
\begin{CD}
\cQHDe @>\fed>> \cQHDp,
\end{CD}
\label{4.14}
\qqq
if we define $\fed$ to be indentically zero if $D\p$ contains a $\CCt$ subgraph
(this is the only possible definition because in this case the space $\cQHDp$
is trivial).

The map\rx{4.14}
allows us to define a space $\cQHD$ for a graph $D$ with a total
edge $\ed$ as a subspace of $\cQHDrc$:
\qq
\cQHD =
\left\{
\begin{array}{ll}
\begin{CD}
\cQHDc @>0\oplus I>> \cQHDrc, &
\qquad\mbox{if $\Dr$ contains a $\CC$ subgraph}
\end{CD}
\\
\begin{CD}
\cQHeDr@>I\oplus \fed>>\cQHDrc &\qquad\mbox{otherwise},
\end{CD}
\end{array}
\right.
\label{4.14x}
\qqq
where $\Dr$ and $\Dc$ are again the graphs constructed
from $D$ by declaring $\ed$ to be root or Cartan.

Now we define $\tcQDCXIHX$ by copying the definition of $\tcDCXIHX$. Again,
let $D$ be a graph with a 4-valent vertex and $D_1$, $D_2$ and $D_3$ be
(1,3)-valent
graphs of \fg{f4.2} with a single total edge $\ed$ constructed by resolving the
4-valent vertex of $D$. This time we extend the linear maps\rx{4.10*} to the
maps
\qq
\mpcd{\QeSHCfobDtX}{\fic}{\cQHv{\xDic}}&\qquad&\mbox{if $\xDir$ contains
a $\CC$ subgraph}
\nonumber
\\
\mpcd{\QeSHCfobDtX}{\fir}{\cQHeev{\xDir}}
&\qquad&\mbox{otherwise.}
\label{4.10*1y}
\qqq
\begin{lemma}
The homomorphisms\rx{4.10*1y} are well-defined.
\end{lemma}
\proof
If $\xDir$ contains a $\CC$ subgraph, then $\hed\in\HCobv{\xDir}$ and this
means that the map $\fic$ of\rx{4.10*} has no kernel. Therefore, denominators
of the elements of $\QeSHCfobDtX$ do not map into zero, and the map $\fic$
of\rx{4.10*1y} is well-defined.

Lemma\rw{l4.1*} guarantees that the image of the map $\fir$
is
contained in $\cQHeev{\xDir}\subset\cQHv{\xDir}$, because the edge $e$ is not
present in $\xD$ and therefore does not contribute polynomials to the
denominators of the elements of $\QeSHCfobDtX$.\qed

The definition\rx{4.14x} allows us to combine the homomorphisms\rx{4.10*1y}
into the homomorphisms
\qq
\begin{CD}
\QeSHCfobDtX @>{f_i}>> \cQHv{D_i},\qquad i=1,2,3.
\label{4.14*}
\end{CD}
\qqq
%
Then the ideal  $\tcQDCXIHX\subset\tcQDCX$ is a span of all elements
$f_1(x)-f_2(x)-f_3(x)$ for all graphs $\xD$ with a single 4-valent vertex, and we
define the new algebra
\qq
\cQDCX = \tcQDCX/
\tcQDCXIHX.
\label{4.15}
\qqq

Naturally, $\cHD\subset\cQHD$ for any $\xD\in\gDXp$, so $\tcDCX$ is a
subalgebra of $\tcQDCX$. It is easy to see that $\tcDCXIHX\subset\tcQDCXIHX$,
therefore there is a natural map
\qq
\mpcd{\cDCX}{\mFq}{\cQDCX},
\label{4.15*}
\qqq
and, in fact, we have a sequence of natural maps
\qq
\begin{CD}
\cBX @>{\ftr}>> \cBCX @>{\mFiq}>> \cQDCX,
\end{CD}
\qquad\mbox{where $\mFiq = \mFq\circ\mFi$.}
\label{4.15*1}
\qqq
%
\begin{conjecture}
$\tcDCXIHX = \tcQDCXIHX\cap\tcDCX$, and hence the map $\mFq$ of\rx{4.15*} is
an injection.
\end{conjecture}



\subsection{Gluing of legs}
\label{3xs.4}
The spaces $\cDCX$ and $\cQDCX$ have algebra structure, but for
our future purposes we will need another useful operation on graph
spaces: the gluing of legs. The definition of gluing $\gll$ is very
straightforward on $\cBX$ (see \eg\cx{A2})
and on $\cBCX$. First, we define it as a unary
operation (that is, linear map)
$\mpcd{\cBCX}{\gllxym}{\cBCX}$ by its action on graphs.
Suppose that we have a graph $D$ with one 1-valent vertex
labeled $x$ and another labeled $y$. Then
$\gllv{x}{y}{1}(\xD)$ is a graph constructed by joining and `dissolving'
these vertices. In other words, we join two legs into a single edge. More
precisely, if both legs are of the same type, then we glue them, if one
leg is root and the other Cartan, then $\gllv{x}{y}{1}(\xD)=0$ by
definition. If one
leg is either root or Cartan and the other is total, then since a total edge is
a sum of a root edge and a Cartan edge, then we must define
$\gllv{x}{y}{1}{\xD}$ by gluing them and declaring
the resulting edge correspondingly either root or Cartan.

For a positive integer $m$, if $D$ has at least $m$ legs of each label $x$
and $y$, then $\gllxym(\xD)$ is a sum of graphs constructed by gluing
$m$ legs $x$ with $m$ legs $y$ in all the possible ways, and if the number
of legs is insufficient, then $\gllxym(\xD)=0$. If a graph $\xD$ is
$\CC$, then so are all the graphs of $\gllxym(\xD)$. Therefore
$\gllxym$ is defined as an
operation on $\tcBCX$. If $\xD$ is a graph with a single 4-valent vertex, which
generates and IHX triplet $\xDo-\xDt-\xDth$ of \fg{f4.2}, then the graphs of
$\gllxym(\xD)$ generate $\gllxym(\xDo-\xDt-\xDth)$. Hence,
$\gllxym(\tcBCXIHX)\subset\tcBCXIHX$ and $\gllxym$ is well-defined as an
operation on $\cBCX=\tcBCX/\tcBCXIHX$.

The binary
operation $D_1\gllxym D_2$ is a composition of two operations:
first, taking a product (that is, disconnected union) of $D_1$ and $D_2$
and then applying the unary $\gllxym$ to it.

\begin{remark}
\label{r4.1}
\rm
Gluing changes the grading $\dgth{}$ in a simple way:
\qq
\dgth{\gllxym(D)} = \dgth{D} + m.
\label{4.16*}
\qqq
\end{remark}

Now we have to define $\gll$ on $\cQDCXp$ in such a way that its
restriction to $\cDCXp$ coincides with $\gll$ on $\cBCX$ transferred by
the isomorphism\rx{4.11}. We will define the gluing of only the legs which
are either root or Cartan of labels $\xX\setminus\Xp$. So let $D$ be a
graph of $\gDXp$ which has at least $m$ legs labeled $x$ and at least $m$
legs of the same type which are labeled by $y$.  Let $\Dpj$ $j=1,\ldots$
be the graphs constructed by gluing $m$ legs of labels $x$ and $y$
pairwise in all possible ways, $j$ indexing particular ways of gluing.
The map
\qq
\begin{CD}
\cQHD @>\gllv{x_1}{x_2}{m}>> \bigoplus_{j}\cQHv{\Dpj}.
\label{4.17}
\end{CD}
\qqq
is a sum of individual maps
\qq
\begin{CD}
\cQHD @>\gllxymj>> \cQHv{\Dpj}.
\label{4.18}
\end{CD}
\qqq
These maps are defined in the following way. If a graph $\Dpj$ has a
proper $\CC$ subgraph, then the space $\cQHv{\Dpj}$ is
zero-dimensional and we define the image of\rx{4.18} as zero. If $\Dpj$ does
not contain proper CC subgraphs, then we consider
a natural map
\qq
\mpcd{\Hhobv{\Dpj}}{\fglj}{\Hhobv{D}},
\label{4.19}
\qqq
which simply cuts the cycles at gluing points. Then we extend the dual map
\qq
\mpcd{\HobD}{\fglj}{\Hobv{\Dpj}}
\label{4.20}
\qqq
to the algebra homomorphism
$\mpcd{\SHCfobDtXp}{\gllxymj}{\SHCfobDjtXp}$ and then to\rx{4.18}.
We have to check that\rx{4.18} is
defined correctly. More precisely, we have to verify that this map does
not produce zeroes in the denominators of $x\in\cQHD$ in case when $D_j$
does not contain proper $\CC$ subgraphs.
In other words, we
have to verify that for any denominator\rx{4.11*},
\qq
\gllxymj  \lrbc{
\preDr \pexXp
}\neq 0.
\label{4.20*}
\qqq
Zero may appear only if
$\gllxymj(\hed) = 0$
for some root edge $\ed$ of $D$, but that would mean that the
corresponding root edge of $D_j$ is in $\HCobDj$. Then, according to
Remark\rw{r4.1*1}, the graph $D_j$ has a $\CCt$ subgraph and therefore
the map\rx{4.18} is defined as identical zero irrespective of
denominators.
Thus we defined a unary operation $\mpcd{\tcDCX}{\gllxym}{\tcDCX}$.

Gluing can also be defined naturally
for graphs containing total internal edges, because
it does not introduce zeroes into denominators. This is
gruaranteed by the following
\begin{lemma}
\label{l4.1*5}
For a root edge $\ed$ of $D$, let $\ed\p$ denote the corresponding root
edge of $D_j$. Then $\gllxymj(\cQHDe)\subset\cQHDje$.
\end{lemma}
\proof
Let us denote by $D\p$ and $\Dpj$ the graphs obtained from $D$ and $D_j$
by declaring $\ed$ and $\ed\p$ to be Cartan. If $\Dpj$ contains a
proper $\CCt$
subgraph, then all of $\cQHD$ is non-singular by definition.
Therefore we assume that neither $D_j$ nor $\Dpj$ contain a proper $\CCt$
subgraph.
If $\Dpj$ does not contain a proper $\CCt$ subgraph,
then neither does $D\p$, because proper
$\CCt$ subgraphs survive gluing.
The maps\rx{4.7*1} and\rx{4.20} commute, hence their symmetric algebra extensions\rx{4.13*1}
and\rx{4.18} also commute:
\qq
\begin{CD}
\SHCfobDtXp @>\gllxymj>> \SHCfobDjtXp \\
@VV{\fed}V          @VV{\fed}V \\
\SHCfobDptXp   @>\gllxymj>> \SHCfobDpjtXp
\end{CD}
\qqq
Since $D\p$ and $\Dpj$ do not contain $\CCt$ subgraphs, then no root edges
belong to the kernels of the horizontal maps. Therefore since the
denominators of $x$ do not belong to the kernel of the left map, then
their image by the upper map can not belong to the kernel of the right
map. Therefore $\gllxymj(x)$ is non-singular at $\ed\p$.\qed

\begin{lemma}
A gluing $\gllxym$ maps the IHX ideal $\tcDCXIHX$ into itself:
\qq
\gllxym(\tcDCXIHX) \subset \tcDCXIHX.
\label{gllihx}
\qqq
\end{lemma}
\proof
Let $\xD$ be a graph with a single 4-vertex, and it
produces 3 IHX related graphs $\xD_i$, $i=1,2,3$.
There is a natural identification between the legs of all 4 graphs
$\xD,\xDo,\xDt,\xDth$, so we can index consistently the ways in which the $x$
and $y$ legs are glued for each of these graphs. Let us fix a
particular gluing $j$. Obviously, the graph $\xD\p=\gllxymj(\xD)$
produces the IHX triplet of graphs $\xD\p_i=\gllxymj(\xD_i)$, $i=1,2,3$.

Let $w\in\QeSHCfobDtX$. Since maps\rx{4.10*} commute with the maps\rx{4.20},
then the maps\rx{4.14*} and\rx{4.18} also commute. Therefore, if
$x\in\QeSHCfobDtX$, then $f_i(\gllxymj(w)) = \gllxymj(f_i(w))$, so
\qq
\gllxymj(f_1(w)-f_2(w)-f_3(w)) = f_1(z)-f_2(z)-f_3(z),
\label{gllp}
\qqq
where $z=\gllxym(w)$. Equation\rx{gllp} means that a glued IHX triplet is
itself an IHX triplet, which implies\rx{gllihx}.\qed

The latter lemma indicates that gluing $\gllxym$ is well-defined on the
quotient space $\cDCX=\tcDCX/\tcDCXIHX$.

\begin{remark}
\rm
\label{rcomm}
The maps\rx{4.15*1} commute with
leg gluing.
\end{remark}

\subsection{Weight system}
\label{3xs.5}

Finally, let us describe an algebra homomorphism
\qq
\begin{CD}
\cQDCX @>
{\Tg}
>> \tSgXt,
\end{CD}
\label{4.21}
\qqq
which coincides with the homomorphism $\Tg$ of\rx{4.3*} on
$\cDCX\subset\cQDCX$. Let us orient the edges of a graph $D\in\gDXp$. We
denote by $\bEDr$ the set of root edges of $D$ and we denote by
$\Dg\subset\mfh$ the set of roots of $\mfg$. A map
$\xmap{\xc}{\bEDr}{\Dg}$ is called a root assignment. An assignment is
called consistent if the sum of incoming roots is equal to the sum of
outgoing roots at every 3-valent vertex of $D$. Let $\bSD$ be the set of
all consistent root assignments for $D$. For a root $\vl\in\Dg$ let
$\mfrl\subset\mfrC$ be the corresponding root space. Obviously,
$\mfrC = \bigoplus_{\vl\in\Dg}\mfrl$. We denote the projection of
$\himfr\in\mfrC^{\otimes 2}$ onto
$\mfrl\otimes\mfrml$ as $\himfrl$ and for a consistent assignement
$\xc\in\bSD$ we construct a tensor
$\himfrc=(\bigotimes_{\ed\in\bEDr}\himfrv{\xc(\ed)})\otimes \himfh$.

Think of $D$ as a cell complex. Let $D\p$ denote a complex that is
obtained from $D$ by removing all Cartan edges. A consistent assignment
$c$ defines an element $\hc\in\HhobDp\otimes\mfh$. According to the
definition\rx{4.3*1}, $\HCfobD = \HobDp$, hence an element
$x\in\cQHD\subset\QeSHCfobDtX$ can be evaluated on $\hc^\otorXp$ to
produce an element $x(\hc)\in\QShoXp$. Then we can define a map
\qq
\begin{CD}
\cQHD @>\tTg>> \QShoXp\otimes\tmfc
\end{CD}
\label{4.22}
\qqq
%
%
by the formula
\qq
\tTg(x)   =\sum_{c\in\bSD} x(\hc)\otimes
\xCD(\fmfg^{\otimes\neDt}\otimes\himfrc),
\label{4.23}
\qqq
where $\xCD$ is the contraction map\rx{1.11*1}. The symmetrization over 1-valent
vertices of $D$ which have the same label, projects $\tTg(x)$ to
$\Tg(x)\in\tSgXt$.
If we restrict $\Tg$ to
$\cHD\subset\cQHD$, then it takes the values in
\qq
\ShoXp\otimes\tSgXts\cong\oSgXT
\nonumber
\qqq
and coincides with its previous
definition which followed\rx{4.3}.

Let us describe explicitly the action of $\Tg$ on spaces associated to special
`graphs'. An element $x\in\cXpHb = S^2(\IQ^{|\Xp|})$ can be presented as a
symmetric matrix $\mtr{\lij}$, $1\leq i,j\leq |\Xp|$. Then
$\Tg(x) = \sum_{i,j=1}^{|\Xp|}\lij\,\himfhij\in
S^2\lrbcs{\bigoplus_{j=1}^{|\Xp|}\mfh_j}$,
where $\himfhij\in \mfh_i\otimes\mfh_j$ is a copy of $\himfh$.
For $x\in\cXpHl=\IQ$ we define $\Tg(x)=x\in\IQ$. Finally,
if we think of $x\in\cQHv{\rcrc} = \evQIQuua$ as a `function' of $\ua$, then
\qq
\Tg(x) = \sum_{\vl\in\Dg\subset \mfh} x( \vl_1\ldtc\vl_{|\Xp|} ) \in \QShoXp,
\label{4.23*}
\qqq
where $\vl_j$ is a copy of $\vl$ in the $j$-th factor of $\ShoXp$. Since the
function $x(\ua)$ is even, then the last expression can be written as a sum
over only positive roots
\qq
\Tg(x) = 2\sum_{\vl\in\Delp}  x( \vl_1\ldtc\vl_{|\Xp|} ).
\label{4.23*1}
\qqq

\nsection{Calculus and differential geometry of graph algebras}
\label{xs4}
\subsection{Functions and general tensor fields}
\label{4xs.1}

We have to introduce a few basic calculus definitions and theorems for
algebras $\cBX$, $\cBCX$, $\cDCXp$ and $\cQDCXp$ in order to define a
universal graph version of the stationary phase integral of
subsection\rw{2xs.4}
and prove its invariance under the shifts of marks on the components of a
marked Morse link.
Our general strategy is to translate as
much differential geometry of $\mfgs$ and of coadjoint orbits $\corbal$
as possible
into the language of graph algebras. We are going to define a graph algebra
version of
functions, vector fields and matrix fields as well as related operations such
as a Lie derivative, a determinant and a (stationary phase) integral. An
application of a weight system should translate graph objects into their Lie
coalgebra relatives in such a way that the operations on graphs become the
usual operations of differential geometry of $\mfgs$ and $\corbal$.
Actually, the definition of the stationary phase (\ie gaussian)
integral was already given in\cx{A2}. Fubini's theorem and the integration
by parts formula were also proved there. In\cx{Ron} we hope to introduce a
graph version of differential forms with the help of grassman
variables and to use them to prove the Duistermaat-Heckmann theorem for
coadjoint orbits.

All our objects are defined as graphs with a particular labeling of legs.
Most of our operations are defined either as relabeling or
as gluing.
Such definitions work equally well for all graph
algebras $\cBX$, $\cBCX$, $\cDCXp$ and $\cQDCXp$.
We just assume that all edges in the
definitions are of `any' type and whenever there is gluing, we invoke an
appropriate definition.
 However, the definition of
an inverse map, a determinant and a gaussian integral require inverting the
coefficients at strut graphs. As a result, these objects are
not well-defined on
$\cDCXp$, but are well-defined on other three graph algebras.
We will work with the space
$\cQDCXp$, because it includes $\cDCXp$, while $\cDCXp$ coincides literally
with $\cBCX$ if we take $\Xp=\emptyset$. Whenever necessary, we will explain
how to modify a definition for the case of $\cBX$.

It will be obvious that our definitions and operations commute with the
maps\rx{4.15*1}. Only determinants and gaussian integrals present a slight
problem. Their definition involves an explicit separation of a `strut' part of
a linear combination of graphs. However, the definition of struts is
not respected by
the isomorphism $\mFiq$: the latter `shaves off' some legs
and thus creates new struts out of some haircomb graphs (\fg{f1.1}) of $\cBCX$.
We will address this issue separately.

All our objects are graphs and their linear combinations. When we say
`graph' in a definition, we actually mean a (possibly infinite) linear
combination of graphs. The distinction between the objects is made by the
way in which we label their legs. We also allow a labeling of a leg by a linear
combination of labels. In this case, following\cx{A2}, we assume that the
result is a linear combination of graphs in which this leg is labeled by
individual labels. In particular, this means that if a leg of $D$ is labeled
with $0$, then the whole graph is equal to $0$.

The legs of our graphs will be labeled with variables $x,y,\ldots$ as well as
with `differentials' $dx,dy,\ldots$ and with derivatives
$\del_x,\del_y,\ldots$. The variables, as labels, may be either root, or
Cartan, or total. A total variable $x$ is essentially a sum of a root variable
$x\dxuor$ and a Cartan variable $x\dxuoc$. We just have to remember that the
spaces $\mfr$ and $\mfh$ related to $x\dxuor$ and $x\dxuoc$ through the
application of the weight system $\Tg$, are parts of the same space $\mfg$.

Let $\ux = (x_1\ldtc x_{\xN})$
be our `coordinate' labels. They may be of any type. Other variable colors (such as
$y,z,\ldots$) may be interpreted as parameters. A \emph{function} $f(\ux)$
is a graph whose legs are labeled by the variables. Let $\mfv_j$
($1\leq j\leq\xN$) be the spaces $\mfg$, $\mfr$ or $\mfh$ depending on whether
$x_j$ is a total, root or Cartan variable. A weight system\rx{4.21}
maps $f$ into an element
\qq
\Tg f\in \QfShXc\otimes\bigotimes_{j=1}^{\xN}S^*\mfv_j,
\label{5.x1}
\qqq
which can be also
interpreted as a formal power series
$\Tg f(\vux;\vua)\in\Qf{\IQ[[\vua]]}[[\vux]]$, where $\vx_j\in\mfvs$
($1\leq j\leq \xN$) and $\va_j\in\mfhs$ ($1\leq j\leq |\Xp|$ ). Thus if we
forget about a distinction between a formal power series and a function, then we
may think of $\Tg F$ as a function of variables $\vux$ which also depends on
Cartan parameters $\vua$. Obviously, $\Tg$ converts a product of graphs in
$\cQDCXp$ into a product of functions.

For positive integer numbers $m_j, n_j$ ($(1\leq j\leq \xN$)
we define
a \emph{tensor field} $\zA(\ux)$ of ranks $\bfm,\bfn$ in $\ux$ as a graph which
has exactly one leg labeled by each of the colors $d\spxi_j$ ($1\leq i\leq
m_j$, $1\leq j\leq\xN$), $\del_{\spxi_j}$, ($1\leq i \leq n_j$, $1\leq j\leq
\xN$), all other legs being labeled by $\ux$. When applying the weight system
$\Tg$ to $\zA$, we use extra metric tensors $h$ in order to convert the spaces
$\mfv$ of $\del_x$ legs into $\mfvs$. Thus,
$\Tg \zA\in \Qf{\IQ[[\vua]]}[[\vux]]\otimes \bigotimes_{j=1}^{\xN}
(\mfv^{\otimes m_j}\otimes(\mfvs)^{\otimes n_j})$ is indeed a tensor field
(formal power series) of variables $\vux$. $\Tg$ converts a product of graphs
in $\cQDCXp$ into a pointwise tensor product of tensor fields.

Next, we define a \emph{contraction of indices} $d\spx{i_1}_j$ and
$\del_{\spx{i_2}_j}$ of $A$ simply as gluing of the corresponding legs
\qq
\zA\longmapsto   \gllv{ d\spx{i_1}_j }{ \del_{\spx{i_2}_j} }{1} (\zA).
\label{5.1*}
\qqq
\begin{lemma}
\label{l5.1*}
The weight system $\Tg$ turns this operation into a
standard definition of contracting upper and lower indices of a tensor
$\mfv_j\otimes\mfvs_j\longrightarrow \IQ$.
\end{lemma}
\proof
This is obvious, since, according to our convention, a weight system puts
a tensor $h^{-1}$ only on the $dx$ leg, effectively omitting it on the
$\del_x$ leg, while a new internal edge produced by gluing carries the
same tensor $h^{-1}$.\qed

The action of the maps\rx{4.15*1} on a tensor field $\zA$ is prescribed by
their definitions. In particular, to define the action of $\mFi$, we suppose
that the labels in the set $\Xp$ are some Cartan variables $\ua=(a_1\ldtc
a_{\xN})$. Assume that the graphs of a tensor $\zA(\ux,\ua)$ of $\cBCX$ do
not have labels $\del_{\ua}$ or $d\ua$. Then the action of $\mFi$ on $\zA$
`shaves off' the legs $\ua$, converting them into the elements of the
spaces\rx{4.12}. Thus, $\ua$ are no longer treated as coordinates, and
we can write
\qq
\zA(\ux,\ua)\xmapt{\mFi}\zA(\ux)\in\cDCX,
\label{5.1u}
\qqq
keeping the same name $\zA$ for the image of the tensor field in $\cDCXp$.

A \emph{1-form} $\omega(\ux)$ is a tensor field all legs of which are colored
by coordinates
except for one
leg which is colored by one of the differentials
$\lstbivjxN{dx}$.


A \emph{partial differential} $\dlxj \zA$ of a tensor field
$A$ is a sum of graphs constructed from $\zA$ by
changing a color at one of the legs of $\zA$ from $x_j$ to $dx_j$ in all
the possible ways. Also we define a `total differential'
\qq
\nabla \zA = \sjoxN \dlxj \zA.
\label{5.1*1}
\qqq
For a function $f$ we define a
\emph{differential} $df=\nabla f$ which is a 1-form.


\subsection{Vector fields}
\label{4xs.2}

A \emph{vector field} $\xi(\ux)$ is a tensor
field all of whose legs are colored by
the variables and parameters except for one leg which is colored by a
derivative $\dlxj$ ($\ojxNi$). We define a pairing $\prcox$ between a
1-form $\omega$ and a vector field $\xi$
as a contraction between the $dx$ and $\del_x$ indices in the tensor product
$\omega\,\xi$
%
\qq
\prcox = \sjoxN \omega \gllsv{dx_j}{\dlxj}{1} \xi.
\label{5.1}
\qqq
In other words, if $\omega$ and $\xi$ are both presented by a single
graph, then $\prcox$ is constructed by gluing the differential leg of
$\omega$ with the derivative leg of $\xi$ if their variables match, and
$\prcox=0$ otherwise. If $\omega$ and $\xi$ are sums of graphs, then
$\prcox$ is defined by bilinearity.

A vector field $\xi$ defines a derivation $\ndxi$ acting on any object $\zA$
by the formula
\qq
\ndxi \zA = \sjoxN \xi \gllsv{\dlxj}{x_j}{1} \zA.
\label{5.2}
\qqq
Obviously, for a function $f$
\qq
\ndxi f = \prc{df}{\xi}.
\label{5.3}
\qqq

A \emph{divergence} of a vector field $\xi(\ux)$ (relative to the coordinates
$\ux$) is a result of gluing the $\dlux$ leg of $\xi$ to one of its $\ux$ legs
in all the possible ways:
\qq
\dvxxi =  \gllv{\dlux}{\ux}{1} (\xi).
\label{5.3*}
\qqq
$\Tg$ transforms this divergence into a standard divergence of a vector field.

\subsection{Diffeomorphisms}
\label{4xs.3}

Next, we define differentiable maps and
diffeomorphisms which we will also call substitutions.
For two sets of coordinate labels $\ux$, $\uy$,
a \emph{map} $\uy = F(\ux)$ is defined as a sum of graphs all legs of which are
labeled by $\ux$ and parameters except for one leg which is labeled by
(one of the coordinates) $\uy$. In applying $\Tg$ to $F$ we treat the $\ux$
legs similar to $\del_x$ legs in tensor fields: we apply an extra metric tensor
$\hm$ to the $\uy$ legs so that
$\Tg F_j\in \Qf{\IQ[[\vua]]}[[\vux]]\otimes\mfvs_j$, where $F_j$ is the
sum of graphs of $\uy=F(\ux)$ which have the  $y_j$ leg.
Thus a set of functions $\Tg F_j$ defines a
local smooth map (or, more precisely, a formal power series) from $\vux$ to
$\vuy$.

A definition of a composition of maps is
obvious, but we will write it in a formal way which is a bit messy.
First, for an element $\zA\in \cQDCXp$ we define a gluing operation
\qq
\gllxyLR (\zA) = \smzi\atv{\gllxym(\zA)}{x=y=0}.
\label{5.4y1}
\qqq
In other words, $\gllxyLR$ glues all $x$ legs to all $y$ legs and
it
yields 0 if the numbers of $x$ and $y$ legs of $\zA$ do not match. We also define
a multiple gluing $\glxyLR$ as a composition of individual gluings
$\gllv{x_j}{y_j}{\iLR}$, $1\leq j\leq \xN$.
Then for a map $\uy=F(\ux)$
we define an element
\qq
\Fiexp = \sjoxN \snzi {(F_j)^n\over n!}.
\label{5.4x}
\qqq
Now a
composition of maps
$\uy=F(\ux)$ and $\uz=G(\uy)$ is a map $\uz=G\circ F\,(\ux)$ defined by
the formula
\qq
G\circ F = G \glyyLR \Fiexp,
\label{5.4x*}
\qqq
which says that all $\uy$ legs of the graphs of $G$ have to be glued to
the matching $\uy$ legs of $F$ in all possible ways. There exists an identity
map
\qq
\Id = \tstuxy,
\label{5.4x1}
\qqq
where we introduced an abbreviated multi-index notation
\qq
\tstuxy = \sjoxN \tstyxj.
\label{5.4x1n}
\qqq
Obviously, $F\circ \Id=\Id\circ F=F$. It is easy to verify that
$\Tg$ converts the composition\rx{5.4x*} into a composition of local
maps.

We will use an obvious shortcut notation $\uz=G(F(\ux))$ for a composition
$\uz= G\circ F (\ux)$.
Moreover, we will use the same `substitution' notation for
any tensor field $\zA(\uy)$ in which $\uy$ is replaced by $\ux$ through
$\uy=F(\ux)$:
\qq
\zA(F(\ux)) = \zA  \glyyLR \Fiexp
\label{5.4xx1}
\qqq
%

A map $\uy=F(\ux)$ is called \emph{linear} if it consists only of struts, one
vertex of which is colored by $\ux$ and the other by $\uy$.
Let us assume for a moment that all our coordinates are
either root or Cartan (that is, if we have a total coordinate, then we split it
into a sum of a root and cartan coordinates) and that the first $L\p$
coordinates ($0\leq L\p\leq L$) are root, while the rest are Cartan.
Since $\cQHv{\cstv{y}{x}} = \IQ\subset \QIQuua$ and
$\cQHv{\rstv{y}{x}}=\QIQuua$ (\cf \eex{4.13}), then a
linear map has a form
\qq
F = \sijoxNp \lij \rstyxij + \sijpxN \lij \cstyxij, \qquad \lij\in\QIQuua,
\label{5.4}
\qqq
and we define
\qq
\det F = \det \mtr{\lji(\bfa)},
\label{5.5}
\qqq
where $\mtr{\lji(\bfa)}$ denotes an $\xN\times\xN$ matrix with entries
$\lji(\bfa)$.
If $\det \uF\not\equiv 0$, then a linear map $F$ can be inverted
%
\qq
%
F^{-1}
= \sijoxNp \liij \rstyxij + \sijpxN \liij \cstyxij,
\qquad\uF^{-1}\circ \uF = \Id,
\label{5.6}
\qqq
where $\mtr{\liij(\bfa)} = \mtr{\lij(\bfa)}^{-1}$. Note that the block-diagonal
structure of $\mtr{\lij(\bfa)}$ implies that $\lij\in\IQ$ for $i,j>\xNp$.

\begin{remark}
\rm
In case of the space $\cBX$ a linear map has a form
\qq
F=\sijoxN \lij\tstv{x_i}{x_j},\qquad \lij\in\IQ
\qqq
and we define $\det F=\det\mtr{\lij}$.
\end{remark}

Now we again assume that coordinates may be of any type.
We call a map $\uF$ \emph{strut} if it contains only strut graphs. Those strut
graphs are of two types: the linear ones $\tstyxij$ and the `constant'
ones $\tstv{y_i}{w_j}$, where $w_j$ are parameters on which $\uF$ depends. For
any map $\uF$, let us denote by $\uFaff$ the strut part of $\uF$ and by
$\uFlin$ the linear part of $\uF$. We call a map $F$ \emph{non-degenerate}
iff $\det \uFlin\not\equiv 0$.
The following is obvious
\begin{lemma}
\label{l5.1}
A strut map $\uF$ is invertible iff it is non-degenerate.
\end{lemma}

Now we can prove a general statement
\begin{theorem}
\label{t5.1}
A map $\uF$ is invertible iff it is non-degenerate.
\end{theorem}
\proof
Since the operation of taking a linear part of a map commutes with composition,
then $\uF^{-1}\circ\uF=\Id$ implies $\uFlin^{-1}\circ\uFlin=\Id$. Thus if $\uF$
is invertible, then $\uFlin$ is invertible and hence non-degenerate.

Now suppose that $\uFlin$ is non-degenerate. Then, according to Lemma\rw{l5.1},
$\uFaff$ is invertable.
It is easy to see that the composite map $\uG=(\uFaff)^{-1}\circ\uF$ does not
have constant struts and $\uG\lin=\Id$. In other words, $\uG = \Id + \uG\p$, and
the map $\uG\p$ contains graphs with at least two edges. Therefore $\uG$ can be
inverted perturbatively: if $\ltiH{n}$ is such that
$(\Id + \ltiH{n})\circ(\Id + \uG\p) = 1 + \ltiK{n+1}$, where $\ltiK{n+1}$
contains only graphs with at least $n+1$ edges, then we can choose
$\ltiH{n+1}=\ltiH{n} + \ltiK{n+1}$. Then
\qq
\uG^{-1} = \Id + \snti \ltiH{n},
\label{5.7}
\qqq
and $\uF^{-1} = \uG^{-1}\circ\uFaff^{-1}$.\qed

\begin{remark}
\rm
Although the map $\mFiq$ of\rx{4.15*1} related the identity maps of $\cBCX$ and
$\cQDCX$ and commutes with the composition of maps, yet the
definitions of inverse maps are a bit different: some maps, which are
degenerate in $\cBCX$, may have non-degenerate $\mFiq$ images in $\cQDCX$.
Indeed, if $\uF$ is a map in $\cBCX$ and $\uF\p=\mFiq(\uF)\in\cQDCX$, then
$\uFlin\p$ has a form\rx{5.4}
\qq
\uFlin\p = \sijoxNp \lij(\ua) \rstyxij + \sijpxN \lij \cstyxij,
\label{5.7y}
\qqq
while
\qq
\uFlin = \sijoxNp \atv{\lij(\ua)}{\ua=0} \rstyxij + \sijpxN \lij \cstyxij.
\label{5.7y1}
\qqq
Thus, generally speaking, $\det\uF\p\neq\det\uF$, rather
$\atv{\det\uF\p}{\ua=0}=\det\uF$, and the reason for this is that the strut
part of $\uF\p$ includes the haircomb graphs of $\uF$ with legs from $\Xp$
(these legs are shaved off by $\mFiq$ and converted into powers of $\ua$).
Therefore it may happen that $\det\uF=0$, whereas $\det\uF\p\not\equiv 0$ and
so $\uF\p$ is non-degenerate. However, in this case its inverse diffeomorphism
$(\uF\p)^{-1}$ is not a $\mFiq$ image of any diffeomorphism of $\cBCX$.
\end{remark}

We call invertible maps \emph{diffeomorphisms} or \emph{substitutions}.
Theorem\rw{t1.1} shows that diffeomorphisms form a group which we call
$\DiffxN$. $\Tg$ generates a homomorphism of this group into a group of
local diffeomorphisms of $\bigoplus_{j=1}^{\xN}\mfvs_j$ depending on
parameters $\bfa$.
We will define the \emph{contragradient}
action of $\DiffxN$ on some tensor fields.
This action is converted by $\Tg$
into a natural contragradient action of the group of local diffeomorphisms
$\bigoplus_{j=1}^{\xN}\mfvs_j$.

First of all,
$\DiffxN$ acts naturally on functions:
a substitution $\uy=\uF(\ux)$ converts a
function $f(\uy)$ it into a function
$\sFst f(\ux) = f(F(\ux))$
of $\ux$
by gluing the $\uy$ legs of $\uF$ to all $\uy$ legs of the graphs of $f$:
\qq
\sFst f(\ux) = f(F(\ux)) = f \glyyLR \Fiexp.
\label{5.7*}
\qqq
%

Formally, we can turn a map $\uy=F(\ux)$ into a vector field by replacing the
labels $\uy$ with $\dlux$. We denote this vector field as $\vF(\ux)$.
Conversely,
we can turn a vector field $\xi(\ux)$ into a map $\uy=\txi(\ux)$ by replacing
the labels $\dlux$ with $\uy$.

The Lie algebra of $\DiffxN$ can be identified with the space of vector
fields $\VectxN$. If a diffeomorphism $\uy=F(\ux;t)$ depends on a \Qpr\ $t$ and
$F(\ux;0)=\Id$, then we identify the tangent vector to
the curve $\uy=F(\ux;t)$ in $\DiffxN$ at $t=0$ with
the vector field $\vdfmw{\dlt F}\big|_{t=0}$.
We leave it for the reader to check that the Lie
algebra commutator induced on $\VectxN$ coincides with the usual vector field
Lie bracket defined as
\qq
[\xi_1,\xi_2] = \ndv{\xi_1} \xi_2 - \ndv{\xi_2} \xi_1.
\label{5.8}
\qqq

Since $\VectxN$ is a Lie algebra of $\DiffxN$, then the action of vector fields
on various objects (such as functions, 1-forms and vector fields) is determined
by the action of diffeomorphisms. The corresponding action of a vector field
$\xi$ is called a \emph{Lie derivative} and we denote it as $\Ldxi$. More
precisely, suppose that $\DiffxN$ acts on a certain type of objects
(say, tensor fields). Denote an action of $F\in\DiffxN$ on $\zA$ as $\sFst \zA$.
For a vector field $\xi(x)$ we construct a 1-parametric family of
diffeomorphisms $F(x;t) = \Id\,+\,t\txi$, so that if $\DiffxN$ acts on a tensor
field $\zA$, then we define $\Ldxi \zA = \atvb{\dlt\sFst(t) \zA}{t=0}$.
Obviously, $\Ldxi f = \ndxi f$ for a function $f$ and $\Ldxi \eta =
[\xi,\eta]$ for a vector field $\eta$. Also for any 1-parametric family of
diffeomorphisms $y=F(x;t)$,
\qq
\dlt \sFst(t) \zA = \Ldv{ \widevec{ (\dlt F(t))\circ F^{-1}(t)} }
(\sFst(t) \zA).
\label{5.8x}
\qqq

\begin{example}
\rm
\label{e5.1}
Consider a map $y=\tmadz(x)$ defined on \fg{f5.1}, $x,y$ being total variables.
\begin{figure}[htb]
\begin{center}
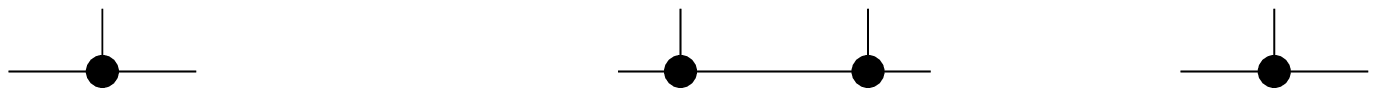
\end{center}
\caption{A map $\tmadz$, a map $\mAdz$ and a vector field $\madz$}
\label{f5.1}
\end{figure}
Its composition exponent defines another map $y = \mAdz(x)$ as
\qq
\mAdz= \exp_\circ (\tmadz) = \snzi {(\tmadz)^{\circ n}\over n!}.
\label{5.8*}
\qqq
The maps $\mAdtz$, where $t$ is a \Qpr, form a 1-parametric subgroup of
$\Diff_x$, which is generated by the vector field $\madz(x)$ of \fg{f5.1}.
This means that $\mAdtz = \xpLie(t\,\madz)$,
where $\xpLie$ means a standard
exponential map from a Lie algebra to a Lie group. More specifically, the
action of $\mAdtz$ on any object is equal to the operator product
exponential $\exp(t\Lmadz)$.

The IHX relation says that
\qq
\tmadv{\tmadv{y}(x)}  = \tmadv{y}\circ\tmadv{x} - \tmadv{x}\circ\tmadv{y}.
\label{5.8*1}
\qqq
Then it is easy to see that
\qq
\tmadv{\mAdev{y} x} = \mAdev{y} \circ\tmadv{x}\circ\mAdev{-y},
\label{5.8*2}
\qqq
and with the help of \ex{5.8*} we derive a useful formula
\qq
\mAdesv{\mAdev{y} x}  = \mAdev{y}\circ\mAdev{x}\circ\mAdev{-y}.
\label{5.8*3}
\qqq

\end{example}

The following lemma follows easily from IHX and $\CC$ relations.
\begin{lemma}
\label{l5.4}
If a function $f(x)$ depends only on a single coordinate $x$ (that is, all its
legs are labeled by $x$), then
\qq
\ndv{\mad_y} f(x) = 0,
\label{5.8**2}
\qqq
if one of the two conditions holds: either if all edges of the graphs of $f$ are
total (that is $f$ belongs to the image of the map\rx{4.2} and $x$ is total),
or if $y$ is Cartan. \end{lemma}

\begin{remark}
\rm
If $f_1(x)$ and $f_2(x)$ are two functions which may
depend on other parameters, then
\qq
\ndv{\mad_y}\lrbcs{f_1(x) \gllsv{x}{x}{m} f_2(x)} =
\lrbcs{ \ndv{\mad_y} f_1(x) } \gllsv{x}{x}{m} f_2(x) +
f_1(x) \gllsv{x}{x}{m}
\lrbcs{ \ndv{\mad_y} f_2(x) }.
\label{5.8**1}
\qqq
\end{remark}


\subsection{Derivatives and matrix fields}
\label{4xs.4}

Let us define a derivative of a diffeomorphism. A derivative is
a `matrix field'. A \emph{matrix field} $\yxM(\ux)$ is a sum of graphs all of
whose legs are
labeled by coordinates $\ux$ and parameters except for two legs, one
carrying a color $dx_j$ ($\ojxNi$) and the other carrying a color $\dlyi$
($\oixNi$). We will use an abbreviated notation $M(\ux)$ for $\xxM(\ux)$. A
derivative of a map $y=F(x)$ is a matrix field $\yxFp(x)$ whose graphs are
constructed from the graphs of $F$ by replacing the colors $y_i$ with $\dlyi$
and a color $x_j$ at one leg by $dx_j$ in all possible ways:
\qq
\yxFp =
(\tstv{\dluy}{\uy})
\gluyo F \gllsv{\ux}{\ux}{1}
(\tstv{\ux}{d\ux}).
\label{5.9}
\qqq
Obviously, $\Tg$ converts \ex{5.9} into a standard derivative of a
diffeomorphism.
A vector field $\sFst \xi(\ux)$ produced by the contragradient
(adjoint) action of a
diffeomorphism $\uy=\uF(\ux)$ on a vector field $\xi(\uy)$ can be expressed as
\qq
\sFst \xi(x) = \lrbc{
\xyFip
\glmdy
\xi(y)
}
\glyyLR
F\iexp,
\label{5.10}
\qqq
where $\ux=\uF^{-1}(\uy)$ is the diffeomorphism inverse to $\uy=\uF(\ux)$.

The matrix fields can be multiplied in a matrix way. We denote this product by
the symbol $\times$ in order to distinguish it from the multiplication in the
graph algebra. The formula for the product is transparent: for example,
\qq
\przyM_1\times \pryxM_2 = \przyM_1 \glmdy \pryxM_2.
\label{5.11}
\qqq
The identity matrix field is obviously
\qq
\prxx I =
\tstv{\dlux}{d\ux}.
\label{5.11*}
\qqq
\begin{theorem}[chain rule]
\label{t5.2}
For two maps $\uy=\uF(\ux)$
and $\uz=\uH(\uy)$, the derivative of their composition is
given by the formula
\qq
(\uH\circ \uF)\p(\ux) = \uH\p (F(\ux))\times \uF(\ux).
\label{5.12}
\qqq
\end{theorem}
\proof
Let us change the order of gluings in the formula for $(\uH\circ \uF)\p$
\qq
(\uH\circ \uF)\p
& = &
(\tstv{\dluz}{\uz})
\gluzo
(G \glyyLR \Fiexp)
\glxxo (\tstv{\ux}{d\ux})
\nonumber\\
& = &
(\tstv{\dluz}{\uz})
\gluzo
\Big[ G \glyyLR \Big( \Fiexp
\glxxo (\tstv{\ux}{d\ux})
\Big) \Big].
\label{5.13}
\qqq
The gluing $\Fiexp \glxxo(\tstv{\ux}{d\ux})$
replaces a label $\ux$ at one of the
legs of one of the graphs of $\Fiexp$ with the corresponding label $d\ux$. Let
us relabel the single $\uy$ leg of that same graph with the label $\dluy$ and
then relabel the $\uy$ leg of $G$ which was glued to it by $\glyxLR$ with the
label $d\uy$. Then we come directly to the \rhs of \ex{5.12}. \qed

The chain rule has an obvious corollary:
\begin{corollary}[derivative of the inverse function]
If a map $\uy=F(\ux)$ is nondegenerate, then the derivative of the inverse map
$\ux = F^{-1}(\uy)$ is equal to the matrix inverse of the derivative of $F$
\qq
(F^{-1})\p(\uy) = \lrbcs{ F\p(F^{-1}(\uy)) }^{\times (-1)}.
\label{5.13**}
\qqq
\end{corollary}

A \emph{trace} of a matrix field $\xxM$ is a function defined in an
obvious way
\qq
\Tr \xxM = \gllv{\dlux}{d\ux}{1} (\xxM).
\label{5.14}
\qqq
It is easy to see that
\qq
\Tr (\xxM_1\times \xxM_2) = \Tr (\xxM_2\times \xxM_1),\quad
\Tr\,[\xxM_1,\xxM_2] = 0,
\label{5.15}
\qqq
where the commutator is defined relative to the matrix product $\times$.

We define a \emph{transposed matrix field} $\xxM^\rT$ by switching the $d\ux$
and $\dlux$ labels at the legs of the graphs of $\xxM$. Obviously,
\qq
(\xxM_1\times\xxM_2)^\rT = \xxM_2^\rT\times\xxM_1^\rT,\qquad
\Tr \xxM^\rT = \Tr \xxM.
\label{5.15x}
\qqq

Let $\xi(\ux)$ be a vector field. Consider a map
$\uy = \txi(\ux)$ constructed from
$\xi$ by replacing its $\dlux$ colors with $\uy$. Then 
\qq
\Tr \txi\p = \dvxxi.
\label{5.15*}
\qqq

\subsection{Determinant and its properties}
\label{4xs.5}

A \emph{determinant} of a matrix field is a function which we will define
in stages, while making sure that it is multiplicative
\qq
\det\, (\prxxM_1\times\prxxM_2) = (\det \prxxM_1)\,(\det\prxxM_2).
\label{5.16}
\qqq
First, we define a determinant of a \emph{strut matrix field}
\qq
\xxM = \sijoxN \lji\;
\tstv{\dlxj}{dx_i},
\label{5.17}
\qqq
Similarly to \ex{5.4}, we assume that coordinates $x_j$, $1\leq j\leq \xNp$ are
root, while coordinates $x_j$, $\xNp+1\leq j \leq \xN$ are Cartan, so that
\qq
\hspace*{-25pt}
\xxM = \sijoxNp \lij \rstv{\dlxi}{dx_j} + \sijpxN \lij \cstv{\dlxi}{dx_j},
\qquad
\lijQ
\label{5.17*}
\qqq
Then we define
%
%
$\det\xxM$ as
\qq
\det\xxM =
\exp\lrbs{\;
\hlf\,\lrbcs{\log(\det\mtrr{\lij}) + \log(\det\mtrr{\lsij}) }\rcrc
+ \log(\det \mtrc{\lij})\lcrc},
\label{5.18}
\qqq
if $\det\mtr{\lij}\not\equiv 0$, and $\det\xxM=0$ if $\det\mtr{\lij}\equiv 0$.
Here $\mtrr{\lij}$ and $\mtrc{\lij}$ are root and Cartan blocks of the
block-diagonal matrix $\mtr{\lij}$.
This definition makes sense, 
because $\cQHv{\rcrc}$ and $\cQHv{\rstv{}{}}$ are related as algebras (see
Remark\rw{r4.sp1}), so we can calculate the products of $\lji$ in
$\QIQuua$ and then project them into $\evQIQuua$, while Remark\rw{r4.sp2}
warned the reader that we would have to include logarightms in the spaces
$\cQHv{\rcrc}$ and $\cQHv{\lcrc}$. The multiplicativity of the
determinants\rx{5.18} follows easily from the multiplicativity of the matrices
$\mtrr{\lij}$ and $\mtrc{\lij}$.

\begin{remark}
\rm
If we work in $\cBX$, then for a strut matrix field\rx{5.17} with $\lij\in\IQ$
we define
\qq
\det\xxM = \exp\lrbcs{ \log(\det\mtr{\lij})\tcrc }.
\qqq
\end{remark}

We call a strut matrix field non-degenerate if its determinant is non-zero.
A non-degenerate strut matrix
field can be inverted with the help of a formula
\qq
(\xxM)^{-1} =
\sijoxN \liij\;\tstv{\dlxi}{dx_j}.
\label{5.19}
\qqq

Next, we define a determinant of a matrix field which is close to identity, by
using a well-known formula which holds for ordinary matrices:
\qq
\det M = \exp (\Tr \log M).
\label{5.19*}
\qqq
In order to abbreviate our notations we will temporarily drop the indices
$\prxx$ from the matrix field notation.

A matrix field $\xxaM$ is called a
\emph{perturbative identity} (or simply \PI) if it is of the form
\qq
\xxaM = I + \xxaM\p,
\label{5.20}
\qqq
where the matrix field $\xxaM\p$ contains only the graphs with at least two
edges. For such a matrix we can define the logarithm
\qq
\lgt(\xxaM) = \snoi (-1)^{n+1}\, {(\xxaM\p)^{\times n}\over n}
\label{5.21}
\qqq
and then define the determinant
\qq
\det\xxaM = \exp \Big(\Tr \lgt \xxaM   \Big).
\label{5.22}
\qqq
\begin{lemma}
\label{l5.2}
The determinant of \PI\ matrix fields defined by \ex{5.22} is multiplicative.
\end{lemma}
\proof
Let $\xxaM_i=
I + \xxaM\p_i$ ($i=1,2$) be two \PI\ matrices.
According to the BCH formula,
\qq
\lefteqn{
\xpt \lrbcs{ \lgt(\xxaM_1) }\; \xpt \lrbcs{ \lgt(\xxaM_2) }
}
\label{5.23}
\\
&
\hspace{1in}
= &
\xpt\lrbcs{ \log(\xxaM_1) + \log(\xxaM_2) + \mbox{(sum of commutators)}
}
\nonumber
\qqq
Taking the $\times$-logarithms of both sides, we find that
\qq
\lgt(\xxaM_1\times\xxaM_2) = \lgt(\xxaM_1) + \lgt(\xxaM_2) +
\mbox{(sum of commutators)}.
\label{5.24}
\qqq
Then in view of \ex{5.15},
\qq
\Tr \lgt(\xxaM_1\times\xxaM_2)  = \Tr \lgt \xxaM_1 + \Tr \lgt \xxaM_2.
\label{5.25}
\qqq
The multiplicativity\rx{5.16} follows from this formula and the
definition\rx{5.22}.\qed
\begin{lemma}
\label{l5.3}
Determinant\rx{5.22} is
conjugation-invariant: if $\xxaM_1$ is \PI\ and $\xxaM_2$ is a
strut matrix field, then
\qq
\det\, (\xxaM_2\times \xxaM_1\times \xxaM_2^{-1}) = \det \xxaM_1
\label{5.26}
\qqq
\end{lemma}
\proof
First, note that the matrix field
$\xxaM_2\times\xxaM_1\times\xxaM_2^{-1}$ is \PI, so
\ex{5.26} makes sense. Then it is easy to see from the definition\rx{5.21} of
$\lgt$ that
\qq
\lgt (\xxaM_2\times\xxaM_1\times\xxaM_2^{-1} ) = \xxaM_2\,
(\lgt\xxaM_1)\, \xxaM_2^{-1}.
\label{5.27}
\qqq
Now it follows from \ex{5.15} that
\qq
\Tr \lgt (\xxaM_2\times\xxaM_1\times\xxaM_2^{-1} ) = \Tr \lgt \xxaM_1.
\label{5.28}
\qqq
This equation together with \ex{5.22} implies \ex{5.26}.\qed

For a matrix field $\xxaM$, let $\xxaMstr$ denote its strut part. We call $\xxaM$
non-degenerate if $\xxaMstr$ is non-degenerate, that is, if
$\det\xxaMstr\not\equiv 0$. If $\xxaM$ is non-degenerate, then it is easy to see
that $\xxaMstr^{-1}\times\xxaM$ is \PI. Then we define
\qq
\det \xxaM =
\left\{
\begin{array}{ll}
\det\xxaMstr \;\det (\xxaMstr^{-1}\times\xxaM)
, & \mbox{if $\det \xxaMstr\not\equiv 0$,}
\\
0, & \mbox{if $\det \xxaMstr\equiv 0$,}
\end{array}
\right.
\label{5.29}
\qqq

Let us establish the properties of the graph determinant.
\begin{theorem}
For a matrix field $M$, $\det( \Tg M) = \Tg (\det M)$.
\end{theorem}
\proof
Let us first prove this for a strut matrix field\rx{5.17}. Let us assume for
simplicity that all coordinates are root (a general case is easy to consider,
since struts do not mix root and cartan coordinates, so that a general strut
matrix field has a block-diagonal form\rx{5.17*}). A weight system $\Tg$
applied to a
root strut matrix field
\qq
M = \sijoxN
\lji\;
\rstv{\dlxj}{dx_i}
\label{5.35*}
\qqq
takes value in $\QShoXp\otimes\End\lrbcs{\bigoplus_{j=1}^{\xN}\mfr_j}$, where
$\mfr_j$ are $\xN$ copies of the space $\mfr$. In fact, the weight system of a
matrix
field\rx{5.35*} is diagonal with respect to the root spaces $\mfrvl$, so
\qq
\Tg M = \svlDg \TgMvl,\qquad
\TgMvl = \mtr{\lij(\vl)} \in \QShoXp\otimes
\End\lrbc{
\bigoplus_{j=1}^{\xN} \mfrvlj }.
\label{5.35*1}
\qqq
%
As a result,
\qq
\det(\Tg M) = \pvlDg \det\TgMvl = \pvlDg \det\mtr{\lij(\vl)}.
\label{5.35*2}
\qqq
A combination of the definition of determinant\rx{5.18} and the
formula\rx{4.23*} for the application of $\Tg$ to $\cQHv{\rcrc}$ yields the same
expression in view of \ex{5.19} applied to $\mtr{\lij(\vl)}$. Thus we proved
the theorem for strut matrix fields.

A general definition of a determinant is based on formulas\rx{5.22}
and\rx{5.29}. The application of weight system $\Tg$ converts \ex{5.22} into
\ex{5.19} and converts \ex{5.29} into a multiplicativity of the ordinary
determinant. Therefore $\Tg$ `commutes' with these defining equations, and
hence it commutes with graph determinant.\qed

\begin{theorem}
\label{t5.3}
The determinant\rx{5.29} is multiplicative.
\end{theorem}
\proof Let $\xxaM_1$ and $\xxaM_2$ be non-degenerate matrix fields. Then
\qq
\det\lrbcs{\xxaM_1\times\xxaM_2} = \det\lrbcs{\xxaMostr\times\xxaMtstr}\;
\det\lrbcs{(\xxaMostr\times\xxaMtstr)^{-1}\times\xxaM_1\times\xxaM_2}.
\label{5.30}
\qqq
Determinant of strut matrix fields is multiplicative, so
\qq
\det\lrbcs{\xxaMostr\times\xxaMtstr} = (\det \xxaMostr)\; (\det \xxaMtstr).
\label{5.31}
\qqq
As for the second factor in the \rhs of \ex{5.30},
\qq
\lefteqn{
\det\lrbcs{(\xxaMostr\times\xxaMtstr)^{-1}\times\xxaM_1\times\xxaM_2}
}  \hspace{1in} &
\nonumber\\
& = &
\det\lrbcs{\xxaMtstr^{-1}\times\xxaMostr^{-1}\times\xxaM_1\times\xxaM_2}
\nonumber \\
& = &
\det \Big[
\lrbcs{\xxaMtstr^{-1}\times(\xxaMostr^{-1}\times\xxaM_1)\times
\xxaMtstr }
\times(\xxaMtstr^{-1}\times\xxaM_2) \Big]
\nonumber \\
& = &
\det \lrbcs{\xxaMtstr^{-1}\times(\xxaMostr^{-1}\times\xxaM_1)\times
\xxaMtstr }
\det (\xxaMtstr^{-1}\times\xxaM_2)
\nonumber \\
& = &\det (\xxaMostr^{-1}\times\xxaM_1)  \;
\det (\xxaMtstr^{-1}\times\xxaM_2).
\label{5.32}
\qqq
Here the third equality is due to Lemma\rw{l5.2} and the fourth equality is due
to Lemma\rw{l5.3}. The multiplicativity\rx{5.16} follows from this equation and
the definition\rx{5.29}.       \qed

Since an exponential is non-zero, then it follows from \ex{5.29} that a matrix
field $\xxaMstr$ is non-degenerate iff $\det\xxaMstr\not\equiv 0$.
\begin{theorem}
\label{t5.4}
A matrix field $\xxaM$ can be inverted iff it is non-degenerate.
\end{theorem}
\proof
The easiest way to prove this theorem is to deduce it from Theorem\rw{t5.1}.
Indeed, one just has to declare the coordinate colors $\ux$ to be parameters
and call $\dlux$ and $d\ux$ new outgoing and incoming coordinates. \qed

An explicit
formula for the inverse matrix field is
\qq
\xxaM^{-1} =
\xxaMstr^{-1}\times\snzi (-1)^n (\xxaMstr^{-1}\times\xxaM)^{\times n}.
\label{5.33}
\qqq

\begin{remark}
\rm
A determinant of a non-degenerate matrix field $M$ is invertible as an element
of the graph algebra. Indeed, \eex{5.29},\rx{5.22} and\rx{5.18} define it as a
product of two functions, both of which are exponentials, and an
inverse of an exponential function is constructed by adding a minus sign to the
exponent.
\end{remark}

\begin{theorem}[derivative of a determinant]
\label{t5.5}
If a matrix field $M(t)$ depends on a \Qpr\ $t$ and $M(t)$ is non-degenerate,
then
\qq
\dlt\det M = \det M \Tr (M^{-1}\, \dlt M).
\label{5.34}
\qqq
\end{theorem}
\proof
Let us evaluate $\det M(t+dt)$ to the linear order in $dt$:
\qq
\det M(t+dt) & = &\det\lrbss{M(t) + \lrbcs{\dlt M(t)}dt}
= \det M(t) \det \lrbss{I + \lrbcs{M^{-1}(t)\,\dlt M(t)}dt }
\nonumber
\\
& = & \lrbcs{\det M(t)}\lrbss{ 1 + \Tr\lrbcs{ M^{-1}(t)\, \dlt M(t)}dt },
\label{5.35}
\qqq
where the last equation follows from \ex{5.22}. Equation\rx{5.34} follows from
\ex{5.35}.\qed

A derivative of a diffeomorphism $\uy=F(\ux)$ is a matrix field $\yxFp(\ux)$.
We define its determinant by first replacing the colors $\dluy$ by the
corresponding colors $\dlux$ and then applying the definition\rx{5.29}.

It is obvious that the definition\rx{5.29} of a determinant commutes with the
maps $\ftr$ and $\mFq$ of\rx{4.15*1}. However the commutativity with $\mFiq$ is
not immediately obvious and requires a special consideration.

\begin{theorem}
\label{tdet}
Let $\xxM(\ux;\ua)$ be a matrix field in $\cBCX$, which depends on Cartan
parameters $\ua$. Suppose that $\ua\in\Xp$ and denote
$\xxMp(\ux) = \mFiq(\xxM(\ux;\ua))$ (\cf\rx{5.1u}, we used notation $\xxMp$ for
the image of $\xxM$ in $\cDCXp$ instead of $\xxM$ in order to avoid any mix up).
Suppose that $\xxM$ is non-degenerate. Then $\xxMp$ is also non-degenerate and
\qq
\det\xxMp(\ux) = \mFiq(\det\xxM(\ux;\ua)).
\label{5.35z1}
\qqq
\end{theorem}
\proof
Let us drop ${}^{x}_{x}$ from our notations of matrix fields. It is easy to see
that the strut graphs of $M\p$ come from the haircomb graphs (\fg{1.1}) of $M$
whose vertical legs are colored by $\ua$ and horizontal legs are colored by
$\del_{\ux}$ and $d\ux$. Therefore,
if following
\eex{5.17} and\rx{5.17*} we denote
\qq
(M\p(\ux))\str = \sijoxN \lij(\ua) \tstv{\del_{x_i}}{dx_j},\qquad
\lij(\ua)\in\IQ[[\ua]],
\label{5.35z2}
\qqq
then
\qq
(M(\ux;\ua))\str = \sijoxN \atv{\lij(\ua)}{\ua=0}
\tstv{\del_{x_i}}{dx_j}.
\label{5.35z3}
\qqq
As a result,
\qq
\det M\str
= \atv{\det M\p\str
}{\ua=0},
\label{5.35z3x}
\qqq
and
non-degeneracy of $M(\ux;\ua)$ indeed implies the non-degeneracy of $M\p(\ux)$.

Since the determinant is multiplicative in both $\cBCX$ and $\cDCXp$, we can
prove \ex{5.35z1} separately for $M\p\str$ and for $(M\p\str)^{-1} M\p$. Since
\qq
\det\mtr{\lij(\ua)} = \det\mtr{\lij(\ua)}\Big|_{\ua=0}
\exp\lrbcs{\Tr \log(\mtr{\lij(\ua)}^{-1}\Big|_{\ua=0}\mtr{\lij(\ua)} },
\label{5.35z4}
\qqq
then
\qq
\mFiqi(\det M\p\str(\ux)) = \det \mFiqi(M\p\str(\ux;\ua)).
\label{5.35z5}
\qqq
Also since the strut parts of both $(M\p\str(\ux))^{-1} M\p(\ux)$ and
$\mFiqi(\,(M\p\str(\ux))^{-1} M\p(\ux)\,)$ are identity matrices, then their
determinants are calculated just by the formula\rx{5.19*} and hence these
determinants are also related by $\mFiq$.\qed

\begin{remark}
\label{rq1}
\rm
Equation\rx{5.35z3x} indicates that it may happen that
$\det M(\ux;\ua)\equiv 0$, while $\det M\p(\ux)\not\equiv 0$. This means that
the definition of a determinant in $\cQDCXp$ is finer than in $\cBCX$, and the
reason for this is that a notion of a strut in $\cQDCXp$ is wider than that in
$\cBCX$, since the former includes the haircomb graphs of $\cBCX$ with legs
$\Xp$ legs.
\end{remark}

\subsection{Integration measure}
\label{4xs.6}
Now we can define
an integration measure. As a graph with labeled legs, an
\emph{integration measure} is the same as a function. The only difference is in
the action of substitutions. A diffeomorphism
$\uy=F(\ux)$ acts on an integration measure $G(\ux)$ as
\qq
G(\uy) \xmapt{F} \sFst G(\ux) = { G(F(\ux)) \det F\p(\ux) }.
\label{5.36}
\qqq
Chain rule and the multiplicativity property of the determinant indicate
that this is a genuine (contragradient) group action.
It is easy to verify that a graph algebra
product $f(\ux)\, G(\ux)$ of a function $f(\ux)$ and an integration
measure $G(\ux)$ is an integration measure.
\begin{theorem}[Lie derivative of integration measure]
\label{t5.6}
The action of a vector field $\xi(\ux)$ on an integration measure $G(\ux)$ is
described by a formula
\qq
\Ldxi\,G = \ndxi G + G\,\dvxxi.
\label{5.37}
\qqq
\end{theorem}
\proof
We define a family of diffeomorphisms $\uy=F(\ux,t)$ by a formula
$F(\ux,t) = \Id + t\,\txi(\ux)$, where $\uy=\txi(\ux)$ is the vector field $\xi$
converted into a map by replacing the $\dlux$ colors with $\uy$. Then by
definition
\qq
\Ldxi\,G = \dlt \attz{ \lrbcs{ {G(F(\ux,t)) \det F\p(\ux,t) } } }.
\label{5.38}
\qqq
As we already know,
\qq
\dlt \attz{ \lrbcs{ G(F(\ux,t)) } } = \ndxi G.
\label{5.39}
\qqq
At the same time, $F\p = I + t\, \nabla\xi$,
so $\det F\p(x,0) = 1$ and according
to \eex{5.34} and\rx{5.15*},
\qq
\attz{\dlt\det F\p} = \Tr\, (\nabla\xi)  = \dvxxi.
\label{5.40}
\qqq
The formula\rx{5.37} follows easily from \eex{5.38},\rx{5.39} and\rx{5.40}.\qed

%


\subsection{A formal gaussian integral}
\label{4xs.7}

Now we define the integral. The only integral that we can define is the
gaussian or the stationary phase one, because, as it turns out, its definition
can be given purely in terms of combinatorics of the integrand without any
references to Riemann sums. We will use the names `gaussian' and `stationary
phase' as synonims, although there is a slight difference between them:
the stationary phase integral allows infinite formal power series as a
preexponential factor.

A gaussian
integral for a graph algebra $\cB$ has been already defined in\cx{A2}.
Its definition can be transferred \emph{verbatim} to any other graph algebra,
and we are going to do it. The only difference is that the paper\cx{A2}
neglected the 1-loop determinant, because it played a rather trivial role in
those calculations. However, we must reinstate it, since it participates in the
proof of topological invariance of \urcc invariant and also contributes to the
Alexander polynomial.

We say that a function,
a tensor field or an integration measure $\gxP$ is $\ux$-substantial if
there exists a positive number $N$ such that $\gxP$ can be presented as a
(possibly infinite) linear combination of elements of spaces $\cQHD$ such that
each graph $D$ has at most $N$ struts $\tstv{x_i}{x_j}$, $1\leq i,j\leq \xN$,
among its connected components. In short, $\gxP$ has a polynomial dependence
on these struts. 
A function, a tensor field or an integration measure $\gxG$ is called
\emph{stationary phase} (\SP), if it has a form
\qq
\gxG(\ux) = \eQx \gxP(\ux),
\label{5.42}
\qqq
where
\qq
\gxQ(\ux) = \sijoxN \lij\; \tstv{x_i}{x_j},    
\label{5.43}
\qqq
while $\gxP(\ux)$
is $\ux$-substantial (we put an extra $1/2$ in the exponent for
future convenience). We 
choose the coefficients $\lij$ so that they satisfy a condition
\qq
\lsij = \lji.
\label{5.43*1}
\qqq
%
We say that $\gxG$ is \emph{stationary phase
non-degenerate} (\SPN) if the matrix
$\mtr{\lij}$ is non-degenerate.
%
We associate to $\gxQ(\ux)$ an `inverse quadratic form'
\qq
\gxQi(\ux) = \sijoxN \liij\; \tstv{x_i}{x_j}
\label{5.43*}
\qqq
( $\mtr{\liij} = \mtr{\lij}^{-1}$) and
a matrix field
\qq
\hgxQ = \sijoxN \lji\;\tstv{\dlxi}{dx_j}.
\label{5.44}
\qqq

If $\gxG(\ux)$
of \ex{5.42} is \SPN, then following\cx{A2} we define its \emph{formal
gaussian} (\ie stationary phase) integral by the formula
\qq
\FGi d\ux\;
\eQx \gxP(\ux)
= \lrbcs{ \det\hgxQ }^{-1/2}
\lrbc{
\gxP(\ux) \glxxLR \eQix
}
,
\label{5.45}
\qqq
%
The \SPN\ condition on $\gxG(\ux)$ guarantees that the gluing in
\ex{5.45} is well-defined, that is, each particular graph in the second
brackets can be constructed by gluing the legs of $\gxP(\ux)$ and
$\eQix$ in only finitely many ways.
The difference between our definition and that of\cx{A2} is in the
`1-loop' determinant factor $\lrbcs{ \det\hgxQ }^{-1/2}$ which was dropped
in\cx{A2}.

\begin{remark}
\rm
Condition\rx{5.43*1} implies that $\det\mtr{\lsij} = \det\mtr{\lij}$, hence
$\det\hgxQ$ can be written simply as
\qq
\det\hgxQ = \exp\lrbcs{ \log (\det\mtrr{\lij})\rcrc +
\log (\det\mtrc{\lij})\lcrc }
\label{5.45*y1}
\qqq
(\cf \ex{5.18}).
\end{remark}

\subsection{Properties of the integral}
\label{4xs.8}

Now let us establish the basic properties of integrals for the
definition\rx{5.45}. Most of the work in this direction has been already done
in\cx{A2}.

First of all, let us use the weight system $\Tg$ in order to relate the graph
formula\rx{5.45} to the usual calculus definition of the stationary phase
integral. Assume for simplicity that $\gxG(\ux)$ and $\gxP(\ux)$ are both either
functions or integration measures. Then, according to\rx{5.x1} we may think of
$\Tg\gxG$ and $\Tg \gxP$ as formal power series of $\vux$ (let us ignore their
dependence on parameters $\vua$).
Usually in calculus the preexponential factor $\Tg \gxP$ would be
presented as a power series in $\hb$, so that the whole integral would be
well-defined as such a power series.
First, let us not commit ourselves to any
particular way of inserting $\hb$ into $\gxP(\ux)$ or
$\gxQ(\ux)$
and consider
the
case of a `polynomial' $\Tg \gxP$.
\begin{lemma}
Suppose that there exists a \emph{finite} set of graphs
$\bD_{\gxP}$, such that
\qq
\gxP(\ux)\in\bigoplus_{D\in\bD_{\gxP}}\cQHD.
\label{5.45*x6}
\qqq
Then $\Tg \gxP$ is a polynomial of $\vux$ and
\qq
\Tg\lrbc{ \FGi d\ux\; \gxG(\ux) } = (-2\pi)^{-\hlf\sjoxN\dim\mfvs_j}\int d\vux\;
\Tg \gxG (\vux).
\label{5.45*x}
\qqq
\end{lemma}
\proof
Obviously, $\Tg \gxP$ is a polynomial of $\ux$ of degree which is equal to the
maximum number of $\ux$ legs in the graphs of $\bD_{\gxP}$. $\Tg \gxQ$ is a
quadratic form
\qq
\Tg \gxQ(\vux) = \sijoxN \lij \;\scp{\vx_i}{\vx_j},
\label{5.45*x1}
\qqq
so if we assume that this quadratic form is negative-definite, then the
\rhs of \ex{5.45} is a well-defined gaussian integral with a polynomial
prefactor. This integral is calculated with the help of an explicit formula
\qq
(-2\pi)^{-\hlf\sjoxN\dim\mfvs}\int d\vux\;
e^{\hlf \Tg\gxQ(\vux_j)}\, \Tg \gxP(\vux) =  \lrbcs{\det \mtr{\lij}}^{-1/2}
\atv{
\Tg \gxP(\dlvux)\;  e^{-\hlf\Tg\gxQi(\vux)}
}{\vux=0}.
\label{5.45*x2}
\qqq
It is easy to see that the \rhs of this formula is equal to $\Tg$ applied to the
\rhs of \ex{5.45}. \qed

\begin{remark}
\rm
The prefactor $(-2\pi)^{-\hlf\sjoxN\dim\mfvs}$ in equation\rx{5.45*x2} indicates
that \ex{5.45} defines $\FGi {d\ux\over (-2\pi)^{1/2}}$ rather than
$\FGi d\ux$. Alternatively, we could remove that prefactor by replacing
$\lrbcs{ \det\hgxQ }^{-1/2}$ in the \rhs of \ex{5.45} with
$\lrbcs{ \det(-2\pi\hgxQ) }^{-1/2}$.
\end{remark}

In order to bring in $\hb$, we introduce a weight system
$\Tgth$, which acts in the same way as $\Tgh$, except that it multiplies the
graphs by $\hb^{\dgth{D}}$ rather than by $\hb^{\dgh{D}}$.
\begin{theorem}
\label{t5.11}
For a \SPN\ function\rx{5.42}
\qq
\Tgth
\lrbc{ \FGi d\ux\; \gxG(\ux) } = (-2\pi\hb)^{-\hlf\sjoxN\dim\mfvs_j}\int d\vux\;
\Tgth \gxG (\vux).
\label{5.45*x3}
\qqq
\end{theorem}
\proof
Let $\gxP_N$ denote the part of $\gxP$ which contains only the contributions of
the graphs $D$ with $\dgth{D}\leq N$. Since there are only finitely many graphs
with such property for a given $N$, then $\gxP_N$ satisfies the
property\rx{5.45*x6} and \ex{5.45*x} holds for
$\gxG_N = \eQx \gxP_N$.
Now \ex{5.45*x3} for $\gxG_N$ follows from \ex{5.45*x} and
from
a relation
\qq
\Tgth
\lrbc{
\gxP_N(\ux) \glxxLR \eQix
}
=
\atv{
\Tgth \gxP_N(\dlvux)\;  e^{-{\hb\over2 }\Tg\gxQi(\vux)}
}{\vux=0},
\qqq
which is established by an easy counting of the powers of $\hb$
(note that gluing a strut to two legs
of a graph $D$ reduces $\dgth{D}$ by 1). Equation\rx{5.45*x3} for $\gxG$
follows by taking a limit $N\longrightarrow\infty$, which is well-defined on
both sides of this equation.\qed


\begin{theorem}[parity invariance]
If $\gxG(\ux)$ is \SPN, then
\qq
\FGi d\ux\; \gxG(-\ux) = \FGi d\ux\; \gxG(\ux)
\label{5.45*}
\qqq
\end{theorem}
\proof This theorem was proved in\cx{A2}. We can use this proof verbatim,
since changing the sign of $\ux$ does not change the strut part and
therefore does not change the 1-loop determinant.\qed

\begin{theorem}[Fubini's theorem]
\label{t5.7}
Consider an \SP\ integrand $\gxG(\ux,\uy)$ which depends on two sets of
coordinates $\ux$ and $\uy$. Suppose that both the full quadratic form
$\gxQ(\ux,\uy)$ and its
$\tstv{\ux}{\ux}$ part $\gxQ_x(\ux)$ are non-degenerate.
Then $\FGi d\ux \;\gxG(\ux,\uy)$ is a \SPN\ integrand for $\uy$
and
\qq
\FGi d\ux d\uy\; \gxG(\ux,\uy) =
\FGi d\uy \lrbc{ \FGi d\ux\;\gxG(\ux,\uy) }.
\label{5.46}
\qqq
\end{theorem}
\proof
The theorem was proved in\cx{A2} for the case when 1-loop determinants
were neglected. Accounting for the determinants is easy. If we split the matrix
$\mtr{\lij}$ of $\gxQ(\ux,\uy)$ into $\ux$ and $\uy$ blocks as
\qq
\mtr{\lij} =  \left(
\begin{array}{c|c}
A & B \\
\hline
B^{\rm T} & C
\end{array}
\right),
\label{5.46*}
\qqq
then it is easy to see that the corresponding $y$-matrix of the exponent of
$\FGi d\ux\;\gxG(\ux,\uy)$
is $C - B^{\rm T} A^{-1} B$ and the determinant part of
Fubini's formula follows from the well-known identity
$\det \mtr{\lij} = \det A\;\det(C- B^{\rm T} A^{-1} B)$.
\qed
\begin{remark}
\label{r5.3}
\rm
The condition that the  $\tstv{\ux}{\ux}$ part $\gxQ_x(\ux)$ of
$\gxQ(\ux,\uy)$
should be non-degenerate, does not restrict the applicability
of Fubini's theorem too much. If needed, one can consider a
deformed quadratic form $\gxQ(\ux,\uy;t)$ which depends on a parameter $t$ in
such a way that $\gxQ(\ux,\uy;0)=\gxQ(\ux,\uy)$ and $\gxQ_x(\ux;t)$ is
non-degenerate for general values of $t$. Then one can apply Fubini's
theorem for a general value of $t$, obtain all necessary results and at
the very end set $t=0$.
\end{remark}
If $\gxG(\ux)$ is an \SPN\ integration measure and $\xi$ is a
$\ux$-substantial vector field,
then $\ndxi \gxG$,
$\gxG\,\dvxxi$ and
$\Ldxi\,\gxG$ are all \SPN\ (the latter also being an
integration measure).
\begin{theorem}[integration by parts]
\label{t5.8}
If $\gxG(\ux)$ is an \SPN\ integration measure
and $\xi(\ux)$ is a $\ux$-substantial
vector field, then
\qq
\FGi d\ux \; \Ldxi \gxG = 0.
\label{5.47}
\qqq
\end{theorem}
\proof
The claim of this theorem obviously does not depend on the inclusion of
the 1-loop determinant factor in \ex{5.45}, so the proof of\cx{A2} applies
without modifications.\qed

This theorem has a simple corollary
\begin{corollary}
\label{c5.1}
Let $\Xi(\ux)$ be an \SPN\ vector field and $\mu(\ux)$ be an $\ux$-substantial
integration measure. If $\LdXi \mu=0$, then
\qq
\FGi d\ux\; (\dvxXi)\,\mu = 0.
\label{5.47*}
\qqq
\end{corollary}

If an \SPN\ integrand
$\gxG(\ux;t)$ depends on a \Qpr\ $t$, then $\dlt \gxG(\ux;t)$ is
also an \SPN\ integrand.
\begin{theorem}[derivative of an integral over a parameter]
\label{t5.9}
If an \SPN\ integrand\\ $G(\ux;t)$ depends on a \Qpr\ $t$, then
\qq
\dlt \lrbc{ \FGi d\ux \; \gxG(\ux;t) } =
\FGi d\ux\; \dlt \gxG(\ux;t).
\label{5.48}
\qqq
\end{theorem}
\proof
Let $\gxG(\ux;t)= \eQxt \gxP(\ux;t)$.
Since
\qq
\dlt \gxG(\ux;t) = \eQxt \lrbcs{
\hlf\;\gxP(\ux;t)\,\dlt\gxQ(\ux;t) + \dlt\gxP(\ux;t)},
\label{5.48yy1}
\qqq
then proving \ex{5.48} amounts to showing the following:
\qq
\lefteqn{
\dlt\lrbs{ \lrbcs{ \det\hgxQ }^{-1/2} \lrbc{ P \glxxLR \eQi } }
}
\nonumber
\\
& = & \hlf  \lrbcs{ \det\hgxQ }^{-1/2}
\Bigg[
\lrbcs{
(\gxP\,\dlt\gxQ)\glxxLR \eQi
}
\;+\;
2
\lrbcs{
\dlt\gxP \glxxLR \eQi
}
\Bigg].
\label{5.48yy2}
\qqq
We calculate explicitly the \lhs of this equation
\qq
\lefteqn{
\dlt\lrbs{ \lrbcs{ \det\hgxQ }^{-1/2} \lrbc{ P \glxxLR \eQi } }
}
\label{5.48yy3}
\\
&&
\hspace{-0.5in}=
\hlf \lrbcs{ \det\hgxQ }^{-1/2} \Bigg[ -\Tr (\hgxQi \dlt \hgxQ)
\lrbc{ P \glxxLR \eQi }
\;-\;
\gxP \glxxLR \lrbcs{ (\dlt \gxQi)\, \eQi }
\nonumber
\\
&&\qquad +\;  2
\lrbcs{\dlt\gxP \glxxLR \eQi}
\Bigg],
\nonumber
\qqq
where $\dlt\gxQi$ can be determined with the help of the formula
\qq
\widehat{\dlt \gxQi} = - \hgxQi\times\dlt\hgxQ\times\hgxQi.
\label{5.48yy4}
\qqq
The last terms in the \rhs of \eex{5.48yy2} and\rx{5.48yy3} coincide, so in
order to prove \ex{5.48yy2} we have to show that
\qq
\lefteqn{
(\gxP\,\dlt\gxQ)\glxxLR \eQi
}
\nonumber
\\
& = &
 -\Tr (\hgxQi \dlt \hgxQ)
\lrbc{ P \glxxLR \eQi }
\;-\;
\gxP \glxxLR \lrbcs{ (\dlt \gxQi)\, \eQi }
\label{5.48yy5}
\qqq
Indeed, consider the gluing in the \lhs of this equation.
The legs of $\dlt\gxQ$ can
either be glued to the same strut of $\eQi$ or to the different struts.
In the first case, the first term of the \rhs of
\ex{5.48yy5} is reproduced. In the second case, the other two legs of
the two struts of $\eQi$ which are glued to the strut $\dlt\gxQ$ will be
glued to the legs of $P$ and in view of \ex{5.48yy4} this reproduces exactly
the second term of the \rhs of \ex{5.48yy5}.              \qed

Finally we want to prove the invariance of the integral\rx{5.45} under
diffeomorphisms if $\gxG(x)$ is an \SPN\ integration measure. We have to limit
ourselves to such diffeomorphisms $F$ that $\sFst G$ is again \SPN.
Recall that
we call an element of a graph algebra \emph{\nrw} if it can be presented as a
linear combination of connected graphs. It is easy to see that \nrw\
diffeomorphisms form a group which we denote as $\DiffxNn$.

\begin{lemma}
\label{lcon}
The groups $\DiffxN$ and $\DiffxNn$ are connected.
\end{lemma}
\proof
It is easy to connect a diffeomorphism $F$ to its linear part: one can use a
path
\qq
F(t) = F\lin + t(F-F\lin),\qquad 0\leq t\leq 1.
\label{5.52*a}
\qqq
Since $(F(t))\lin=F\lin$, then all maps $\uy=F(\ux;t)$ are non-degenerate and
thus belong to $\DiffxN$.
If $F$ is \nrw, then all diffeomorphisms $F(t)$ are also \nrw. $F\lin$, in its
turn, is connected to the identity diffeomorphism $\Id$, because the group $GL$
is connected.\qed


\begin{theorem}[diffeomorphism invariance of an integral of a measure]
\label{t5.10}
If $\gxG(\uy)$ is an \SPN\ integration measure and $\uy=F(\ux)$ is \nrw, then
$\sFst G(\ux)$ is also \SPN\ and
%
%
\qq
\FGi d\ux\; \sFst \gxG(\ux) = \FGi d\uy\; \gxG(\uy).
\label{5.52}
\qqq
\end{theorem}
\proof
Consider \ex{5.36} for $\sFst G(\ux)$. It is easy to check that $G(F(\ux))$ is
\SPN. Also it follows form \eex{5.29},\rx{5.43} and\rx{5.18} that if
$\uy=F(\ux)$
is \nrw, then the graphs of the determinant $\det F\p(\ux)$ of \ex{5.36} do not
contain connected tree subgraphs, so $\det F\p(\ux)$ is $\ux$-substantial and
$\sFst G(\ux)$ is \SPN.

Since, according to Lemma\rw{lcon}, the group $\DiffxNn$ is connected, then
there exists a 1-parametric family of \nrw\ diffeomorphisms $F(t)$, such that
%
$\atvb{F(t)}{t=0}=\Id$ and $\atvb{F(t)}{t=1}=F(\ux)$. Then \ex{5.52}
follows from the following calculation
\qq
\dlt\lrbc{ \FGi \sFst(t) \gxG(\ux)} =
\FGi \dlt\sFst(t) \gxG(\ux) =
\FGi \Ldv{ \widevec{ (\dlt F(t))\circ F^{-1}(t)} } \lrbcs{\sFst(t) \gxG(\ux)}
= 0,
\qqq
which is based on \eex{5.48},\rx{5.8x} and\rx{5.47}, the latter being
applicable since the vector field $\widevec{ (\dlt F(t))\circ F^{-1}(t)}$ is
\nrw\ and therefore $\ux$-substantial.\qed

%
\begin{corollary}
\label{c5.2}
If $\uy=F(\ux)$ is a \nrw\
diffeomorphism and $\msj(\ux)$ is an $\ux$-substantial and
$F$-invariant integration
measure, then for any \SPN\ function $G(\uy)$
\qq
\FGi d\ux\; G(F(\ux))\,\msj(\ux) = \FGi d\uy\; G(\uy)\,\msj(\uy).
\label{5.53}
\qqq
\end{corollary}
\proof
Apply \ex{5.52} to the \SPN\ integration measure $G(\uy)\msj(\uy)$.\qed

Obviously, the definition\rx{5.45} commutes with the first injection $\ftr$
of\rx{4.15*1}. However, its commutativity with $\mFiq$ is less obvious, since
as we have already discussed, the notion of a strut exponent is different in
$\cBCX$ and $\cDCXp$. In order to prove it, we have to establish that
\ex{5.45} sometimes holds even in the case when $\gxQ(\ux)$ is not purely
strut.


First, we introduce a few general definitions. A
function $\gxQ(\ux)$ is called \emph{quadratic} in $\ux$, if every graph of
$\gxQ(\ux)$ has exactly two $\ux$ legs ($\gxQ(\ux)$ may also depend on other
coordinates and parameters). Suppose that a matrix field $\xxM$ does not depend
on $\ux$, that is, its graphs do not have $\ux$ legs. Then we can associate to
$\xxM$ a quadratic (in $\ux$) function $\Mqf$ by replacing the $\dlux$ and
$d\ux$ labels in the graphs of $\xxM$ with corresponding labels $\ux$.
Conversely, if $\gxQ(\ux)$ is a quadratic
function of $\ux$, then we can associate to it a unique matrix field $\hgxQ$
which satisfies two properties: $\qfv{\hgxQ} = \gxQ(\ux)$
and $\hgxQ^\rT = \hgxQ$
(that is, $\hgxQ$ is symmetric). Obviously, our definition is consistent with
\ex{5.44}.

Let $\gxQ(\ux)$ be a quadratic function of $\ux$.
As usual, $\gxQ\str(\ux)$ denotes
its strut part.
\begin{theorem}
\label{tfunc}
Suppose that $\gxQ\str(\ux)$ is non-degenerate and that $\gxP(\ux)$ is
$\ux$-substantial. Then $\eQx \gxP(\ux)$ is \SPN\ and \ex{5.45} holds, where
$\gxQi(\ux)$ is defined as $\gxQi(\ux)= \qfv{\hgxQi}$.
\end{theorem}
\proof
According to the definition of formal gaussian integral, in our case
\qq
\FGi
d\ux\;
\eQx \gxP(\ux)
& = &
\FGi d\ux\;
e^{\hlf \gxQ\str}\;\lrbc{e^{\hlf(\gxQ-\gxQ\str)} \gxP}
\nonumber
\\
& = &
\lrbcs{ \det\hgxQ\str }^{-1/2}
\lrbs{
\lrbc{e^{\hlf(\gxQ-\gxQ\str)} \gxP} \glxxLR \e^{-\hlf\,\gxQ\str}
}.
\label{5.54}
\qqq
We introduce the functions
\qq
&\gxQ(x;t) = \gxQ(x) + t\lrbcs{\gxQ\str(x) - \gxQ(x)},
\qquad
\gxP(\ux;t) = e^{ {t\over 2} \big(\gxQ(\ux) - \gxQ\str(\ux)\big) } \gxP(\ux).
\label{5.56}
\qqq
depending on a $\IQ$-parameter $t$ and consider an expression
\qq
A(t) =
\lrbcs{ \det\hgxQ(t) }^{-1/2}
\lrbc{
\gxP(\ux;t)\glxxLR e^{-\hlf\,\gxQ(\ux;t)}
}.
\label{5.57}
\qqq
In view of \ex{5.54}, $A(1)$ equals the \lhs of \ex{5.45}, whereas $A(0)$
obviously equals the \rhs of that equation. Let us consider $\dlt A(t)$. It is
easy to see that the proof of \ex{5.48yy2} also works in the case when
$\gxQ(\ux)$ is
a quadratic function of $\ux$ which is not necessarily strut. Since
in our particular case
\qq
\gxP(\ux;t)\,\dlt\gxQ(\ux;t) + 2 \dlt\gxP(\ux;t) = 0,
\label{5.58}
\qqq
then, according to \ex{5.48yy2}, $\dlt A(t)=0$. Hence $A(0)=A(1)$ and this
proves \ex{5.45} in case when $\gxQ(\ux)$ is a quadratic function of $\ux$.\qed

Now we are ready to prove that a gaussian integral commutes with the injection
$\mFiq$.
\begin{theorem}
Let $\gxG(\ux;\ua)$ be a tensor field or an integration measure which depends
on Cartan parameters $\ua$. Suppose that $\ua\in\Xp$ and denote
$\gxGp(\ux) = \mFiq( \gxG(\ux;\ua) )$. If $\gxGp$ is \SP, then $\gxG$ is \SP.
If in addition $\gxG$ is non-degenerate, then $\gxGp$ is also non-degenerate
and
\qq
\FGi d\ux\;\gxGp(\ux) = \mFiq\lrbc{ \FGi d\ux\;\gxG(\ux;\ua) }.
\label{5.58*1}
\qqq
\end{theorem}
\proof
$\gxGp$ being \SP\ and belonging to the image of $\mFiq$
means that $\gxGp=e^{\hlf\gxQp(\ux)}\gxPp(\ux)$, where
\qq
\gxQp =
\sijoxN \lij(\ua) \tstv{x_i}{x_j},\qquad
\lij(\ua)\in\IQ[[\ua]],
\label{5.58*2}
\qqq
and $\gxPp\in\cDCX\subset\cQDCX$ is $\ux$-substantial.  Then
$\gxG(\ux;\ua) = \eQx \gxP(\ux;\ua)$, where
\qq
\gxQ(\ux) = \sijoxN \lij(\ua)\Big|_{\ua=0} \tstv{x_i}{x_j},
\label{5.58*3}
\qqq
while
\qq
\gxP(\ux;\ua) =
\exp \mFiqi\lrbc{ \hlf\sijoxN (\lij(\ua) - \lij(\ua)\Big|_{\ua=0})
\tstv{x_i}{x_j}
}
\,
\mFiqi(\gxPp(\ux))
\label{5.58*4}
\qqq
Since the exponent in \ex{5.58*4} does not contain $\cDCX$ struts and since
$\mFiqi$ in contrast to $\mFiq$ does not create new struts, then
$\gxP(\ux;\ua)$ is $\ux$-substantial and hence $\gxG(\ux;\ua)$ is \SP.

If $\gxG$ is \SPN, then $\det \mtr{\lij(\ua)}\Big|_{\ua=0}\neq 0$, hence
$\det\mtr{\lij(\ua)}\not\equiv 0$ and $\gxGp$ is also \SPN.

It remains to prove \ex{5.58*1}. Consider its \rhs:
\qq
\FGi d\ux\;\gxG(\ux;\ua) =
\FGi d\ux\; \exp\lrbc{\hlf \mFiqi(\gxQp(\ux))}
\mFiqi(\gxPp(\ux)).
\label{5.58*5}
\qqq
Since $\mFiqi(\gxQp(\ux))$ is a quadratic function of $\ux$ and its strut part
is assumed to be non-degenerate, then Theorem\rw{tfunc} says that
we can calculate the integral in the \rhs of \ex{5.58*5} with the help of
\ex{5.45} in which now $\gxQ = \mFiqi(\gxQp)$ and $\gxP=\mFiqi(\gxPp)$. The
result obviously coincides with the \lhs of \ex{5.58*1} (we can use
Theorem\rw{tdet} in order to establish the equality of `1-loop determinant
factors').\qed

\begin{remark}
\label{rq2}
\rm
Equation\rx{5.58*2} indicates that
\qq
\det\hgxQ = \det \mtr{\lij(\ua)}\Big|_{\ua=0}.
\label{5.58*6}
\qqq
Therefore it may happen that $\det\hgxQ=0$, whereas $\det\hgxQ\p\not\equiv0$,
which means that $\gxG(\ux;\ua)$ is degenerate, while $\gxGp(\ux)$ is
non-degenerate. This is similar to what we described in Remark\rw{rq1}: a
definition of a gaussian integral is finer in $\cQDCX$ than in $\cBCX$. This is
important for us, since it will turn out that the \urcc invariant is defined
through a gaussian integral which is defined only in $\cQDCX$.
\end{remark}

\nsection{A formal integration over a coadjoint orbit}
\label{xs5}
\subsection{Invariant integration measure on a coadjoint orbit}
\label{5xs.1}
In subsection\rw{2xs.3}
we used $\xG/\xT$ as a model of a coadjoint orbit, and
exponential map turned the root space
$\mfr\subset\mfg$ into a coordinate chart in the vicinity of $e\in\xG$. Now
we are going to do the same at the level of graph algebras. Since we
consider just a single coadjoint orbit (a case of multiple orbits would be
an obvious generalization), we introduce a single root coordinate color
$\crr$. Our immediate task is to define a diffeomorphism $\FLgy$ ($y$ being
a root parameter), which is an analog of diffeomorphism $\FLy$ defined by the
commutative diagram\rx{1.37}, and then find an $\FLgy$-invariant integration
measure $\msgr$.

In order to distinguish between root, Cartan and total labels and coordinates, we
assume the following convention: root coordiantes are denoted as $x,y,\ldots$,
Cartan coordinates are denoted as $a,b,\ldots$ and total coordinates are
denoted as $\bfx, \bfy,\ldots$.  For any map
$\bfy = F(\cdot)$ we define two related maps $y=\Fr(\cdot)$ and
$a=\Fc(\cdot)$ in which the total $\bfy$ leg of $F$ is replaced by a root
$y$ leg or by a Cartan $a$ leg. Also, if we replace a total label $\bfx$
in the argument of any object $A(\bfx)$ by either $y$ or $a$, then this
means that the $\bfx$ legs of $A$ are declared root or Cartan.

We have to express \ex{1.38} in purely graphical terms.
First, we have to learn how to graphically multiply the exponentials with
non-commutative exponents. This has been already described in\cx{A2},
where a BCH `forest' has been introduced. If $X$ and $Y$
are two formal
non-commutative variables, then according to BCH,
\qq
& e^X e^Y = e^Z
\nonumber
\\
&Z = X + Y + \hlf\,[X,Y] + {1\over 12}\, [X,[X,Y]] - {1\over 12}\,[Y,[X,Y]] +
\mbox{higher commutators}.
\label{6.1}
\qqq
Following\cx{A2}, we turn every term in the expression for $Z$ into a tree
graph with total edges, placing the color $\bfz$ at the `root' and
colors $\bfx$ and $\bfy$ at the ends of the `branches' of these trees. The sum
of all these trees defines a BCH map which we denote as $\bfz =
\Ev{\bfx}{\bfy}$:
\qq
\BCH =
\tstv{\bfz}{\bfx} + \tstv{\bfz}{\bfy}
+ \mbox{trees with multiple edges}.
\label{6.1*3}
\qqq
%
We also define a multiple BCH map $\bfy = \Esv{\bfx_1\ldtc \bfx_n}$ as a
logarithm of the product of $n$ exponentials. The multiple map can be
presented as a composition of binary BCH maps in any order compatible with
the associativity. The BCH maps satisfy some obvious properties coming
from the basic properties of exponentials. Among these properties are
\qq
\Ev{t_1 \bfx}{t_2 \bfx} = \Ev{(t_1+t_2) \bfx}{0} =
\tstv{\cdot}{
(t_1+t_2)\bfx},
\label{6.1*}
\qqq
which follows from the commutativity of legs of the same color,
\qq
\Ev{a}{b} = \Ev{a+b}{0} =
\tstv{\cdot}{(a+b)},
\label{6.1*1}
\qqq
which follows from $\CC$ relations expressing the commutativity
of Cartan edges (Lemma\rw{lcc}), and
\qq
&\Esv{\bfy,\bfx,-\bfy} = \mAdev{\bfy}(\bfx),
\label{6.1*2}
\\
&\mAdesv{\Ev{\bfx}{\bfy}} = \mAdev{\bfx}\circ\mAdev{\bfy}.
\label{6.1*4}
\qqq
($\mAdev{\cdot}$ map being defined in Example\rw{e5.1}), which follows
from the IHX relation and the combinatorial identity
$e^x e^y e^{-x} = e^{e^{\ad_x} y}$. Also, due to the presence of struts
in\rx{6.1*3}, the map $\BCH$ can not mix root and Cartan input in order to
produce a purely root output:
\begin{lemma}
\label{l6.1}
Consider two maps: $y=\Hr(\bfx)$ with a root output and $a=\Hc(\bfx)$ with a
Cartan output. If a Cartan part $\Fc$ of
a map $\bfy=F(\bfx)$ defined as
\qq F(\bfx) =
\Esv{\Hr(\bfx),\,\Hc(\bfx)},
\label{6.13}
\qqq
is zero, then
\qq
F(\bfx) = \Hr(\bfx)\qquad\mbox{and}\qquad \Hc(\bfx)=0.
\label{6.14}
\qqq
\end{lemma}
\proof
Consider a graph of $\Hc$ with the smallest number of edges. Then due to
the struts in\rx{6.1*3} it will produce the graph with the smallest number
of edges among the graphs in \rhs of \ex{6.13} with Cartan output. Hence,
this graph can not be cancelled to conform with the fact that the \lhs of
\ex{6.13} has only root output.\qed

Next, we have to learn how to split a single exponential into a product of two
exponentials: one of a root element and the other of a Cartan element. Consider
a label space with two coordinates $(x,a)$.
The map $(y,b) = \Sv{x}{a}$ that we want is the inverse of the map
$(x,a) = \tEv{y}{b}$ defined by the formula
\qq
\tEv{y}{b} = (\Erv{y}{b},\Ecv{y}{b}).
\label{6.2}
\qqq
The formula for the $\BCH$ map\rx{6.1*3} demonstrates
that the strut part of $\tBCH$ is the identity map $\Id$, so according to
Theorem\rw{t5.1}, $\tBCH$ is invertible and its inverse can be constructed
perturbatively as described in its proof. Now we define the map
$z = \FLgy(x)$ as
\qq
\FLgy(x) = \Srv{\Erv{y}{x}}{\Ecv{y}{x}},
\label{6.3}
\qqq
where $y=\Srv{x}{a}$ is a part of the map
$(y,b) = \Sv{x}{a}$ which has the color $y$ at its outgoing leg.
The
strut part of $z=\FLgy (x)$ is
\qq
(\FLgy)\str = \rstv{z}{x} + \rstv{z}{y}.
\label{6.3*}
\qqq
%
We will also need a similar map
\qq
\tFLgy(x) = \Scv{\Erv{y}{x}}{\Ecv{y}{x}},
\label{6.7}
\qqq
so that according to the definition of $\SPL$,
\qq
\Ev{y}{x} = \Ev{\FLgy(x)}{\tFLgy(x)}.
\label{6.8}
\qqq

\begin{lemma}
\label{l6.x2}
The diffeomorphism $z=\FLgy(x)$ is \nrw.
\end{lemma}
\proof
The BCH formula shows that $\bfz = \Ev{\bfx}{\bfy}$ is \nrw, therefore
$(x,a) = \tEv{y}{b}$ is \nrw\ and so is its inverse $(y,b) = \Sv{x}{a}$. Hence
the composition\rx{6.3} is also \nrw.\qed

Let us establish a few properties of $\FLgy$. First of all, it is easy to
verify that
\qq
\FLv{t_2 y}(t_1 y) = (t_1+t_2) y,\quad
\tFLv{t_2 y}(t_1 y) = 0.
\label{6.7*}
\qqq
Next, we establish the graph analog of \ex{1.41}, but in order to
formulate it, we need a slightly modified version of the map
$\bfz=\mAdev{\bfy}(\bfx)$. Namely, when the total parameter $\bfy$ is
replaced by a Cartan parameter $a$, and the argument is a root coordiante
$x$, then the $\CCt$ relation says that $\Big(\mAdev{a}(x)\Big)\xdc=0$
and $\bfz=\mAdev{a}(x)$ is reduced to $z=\mAdev{a}{x}$. It is easy to see that
the root-root maps $\mAdev{ta}$ still form a 1-parametric group, which is
generated by a (purely root) vector field $\mad_a$.

\begin{lemma}
\label{l6.2}
The diffeomorphisms $\FLgy$ satisfy the following composition formula
\qq
\FLv{y_2}\circ\FLv{y_1} =
\FLv{\FLv{y_2}(y_1)}\circ\mAdesv{\tFLv{y_2}(y_1)}.
\label{6.15}
\qqq
\end{lemma}
\proof
consider a map from $(x,y_1,y_2)$ to $\bfz$
\qq
\bfz = \Esvb{\,y_2,\,y_1,\,x,\,-\tFLv{y_1}(x) -
\tFLv{\tFLv{y_2}}\circ\FLv{y_1}(x)\,}.
\label{6.9}
\qqq
We will transform the \rhs by applying \ex{6.8} to pairs of adjacent root
coordinates in two different orders. First, we begin by applying \ex{6.8}
to $(y_1,x)$
\qq
\bfz & = &\Esvb{\,y_2,\,\FLv{y_1}(x),\,\tFLv{y_1}(x),\,
-\tFLv{y_1}(x) -
\tFLv{y_2}\circ \tFLv{y_1}(x)\, }
\nonumber\\
& = &
\Esvb{\,y_2,\,\FLv{y_1}(x),\,-\tFLv{y_2}\circ\FLv{y_1}(x)\, }
\nonumber\\
& = & \Esvb{\,\FLv{y_2}\circ\FLv{y_1}(x),\,
\tFLv{y_2}\circ\FLv{y_1}(x),\,
-\tFLv{y_2}\circ\FLv{y_1}(x)\,}
\nonumber\\
& = & \FLv{y_2}\circ\FLv{y_1}(x).
\label{6.10}
\qqq
On the other hand, we may first apply \ex{6.8} to the pair $(y_2,y_1)$
\qq
\bfz & = & \Esvb{\,\FLv{y_2}(y_1),\,\tFLv{y_2}(y_1),\,x,
\,-\tFLv{y_1}(x) - \tFLv{y_2}\circ\FLv{y_1}(x)\, }
\nonumber\\
& = &
\Esvb{\,\FLv{y_2}(y_1),\,\tFLv{y_2}(y_1),\,x,\,-\tFLv{y_2}(y_1),
\,\tFLv{y_2}(y_1)
- \tFLv{y_1}(x)-\tFLv{y_2}\circ\FLv{y_1}(x)\, }
\nonumber\\
& = &
\Esvb{\,\FLv{y_2}(y_1),\mAdesv{\tFLv{y_2}(y_1)}(x),\,
\tFLv{y_2}(y_1)
- \tFLv{y_1}(x)-\tFLv{y_2}\circ\FLv{y_1}(x)\, }
\nonumber\\
& = &
\BCH\Big(\FLv{\FLv{y_2}(y_1)}\circ\mAdesv{\tFLv{y_2}(y_1)}(x),\,
\tFLv{\FLv{y_2}(y_1)}\circ\mAdesv{\tFLv{y_2}(y_1)}(x)+
\tFLv{y_2}(y_1)
\nonumber\\
&&\hspace{3in}
- \tFLv{y_1}(x)-\tFLv{y_2}\circ\FLv{y_1}(x)\,\Big) .
\label{6.11}
\qqq
Comparing \eex{6.10} and\rx{6.11} we find that
\qq
\FLv{y_2}\circ\FLv{y_1}(x) & = &
\BCH\Big(\FLv{\FLv{y_2}(y_1)}\circ\mAdesv{\tFLv{y_2}(y_1)}(x),\,
\tFLv{\FLv{y_2}(y_1)}\circ\mAdesv{\tFLv{y_2}(y_1)}(x)
\nonumber\\
&&\hspace{2in}
+\tFLv{y_2}(y_1)- \tFLv{y_1}(x)-\tFLv{y_2}\circ\FLv{y_1}(x)\,\Big) .
\label{6.12}
\qqq
According to the first of \eex{6.14} of Lemma\rw{l6.1}, \ex{6.12} implies
the composition law\rx{6.15}. \qed
\begin{corollary}
\label{c6.1}
The diffeomorphisms $\FLv{ty}$ form a 1-parametric group.
\end{corollary}
\proof. This corollary follows from a simple calculation based on
\eex{1.65} and\rx{6.7*}
\qq
\FLv{t_2y}\circ\FLv{t_1y} = \FLv{\FLv{t_2y}(t_1y)}\circ
\mAdev{\tFLv{t_2y}(t_1y)} = \FLv{(t_1+t_2)y}\circ\mAdev{0}
= \FLv{(t_1+t_2)y}.
\label{6.12*}
\qqq
\qed

Let $\xiFy$ denote the vector field which generates $\FLv{ty}$. It follows
from \ex{6.3*} that  its strut part is
\qq
\lrbcs{\xiFy(x)}\str = \rstv{\del_x}{y}.
\label{6.12*x}
\qqq

Following \ex{1.39} we choose the measure for a root space to be the
inverse determinant of the derivative of $\FLgx$
%
\qq
\msgx = \lrbcs{ \det\FLgx\p (0) }^{-1}.
\label{6.4}
\qqq
Obviously, $\msjv{0} = \bfun$ and it follows from \ex{5.36} that
\qq
\msgx = \FLv{-x}\drst \msjv{0}
\label{6.12*1}
\qqq
(note that \ex{5.36} defines the action of a diffeomorphism as a
pull-back).

\begin{theorem}
\label{t6.1}
The measure\rx{6.4} is invariant under the substitution
$y=\mAda(x)$.
\end{theorem}
\proof
Since the substitutions $\mAdta$ form a 1-parametric group generated by
$\mad_a$, then it is sufficient to check that
\qq
\Ldmada \msj = \ndmada \msj + \msj\; \dvxmada = 0.
\label{6.6}
\qqq
In fact, both terms in the middle expression are equal to 0.
Indeed $\ndmada\msj=0$ because of \ex{5.8**2} of Lemma\rw{l5.4}, while
%
$\dvxmada=0$, because the symmetry of the graph $\rcrc$ sends an odd power of
$a$ to zero. \qed
\begin{theorem}
\label{t6.2}
The measure\rx{6.4} is invariant under the action of $\FLgy$.
\end{theorem}
\proof
This invariance is a simple corollary of Lemma\rw{l6.2} and
Theorem\rw{t6.1}. Indeed
\qq
\FLgy\drst \msjv{x} & = & \FLgy\drst\lrbcs{ \FLv{-x}\drst \;\msjv{0} }
=
\lrbcs{ \FLv{-x}\circ\FLgy }\drst\;\msjv{0}
=
\lrbcs{ (\FLv{-y}\circ\FLv{x})^{-1} }\drst\;\msjv{0}
\nonumber
\\
& = &
\lrbs{ \lrbcs{
\FLv{ \FLv{-y}(x) } \circ \mAdev{ \tFLv{-y}(x) }
}^{-1} }\drst\;\msjv{0}
\nonumber
\\
& = &
\lrbcs{ \mAdev{-\tFLv{-y}(x)} \circ \FLv{-\FLv{-y}(x)} }\drst\;\msjv{0}
\nonumber\\
& = &
\FLv{-\FLv{-y}(x)}\drst \lrbcs{
\mAdev{-\tFLv{-y}(x)}\drst\;\msjv{0} }
\nonumber\\
& = & \FLv{-\FLv{-y}(x)}\drst\;\msjv{0}
= \msjv {\FLv{-y}(x)} = \msjv{ \FLgy^{-1}(x) }.
\label{6.16}
\qqq
The resulting equation $\FLgy\drst\msjv{x}=\msjv{ \FLgy^{-1}(x) }$
means that the measure $\msgx$ is indeed invariant under the action of
$\FLgy$.\qed

Let us describe the structure of $\msgx$ in more details.
\begin{theorem}
\label{t6.5}
The measure $\msgx$ is an exponential $\msgx = e^{\omsgx}$, where $\omsgx$ is
a sum of wheel graphs,
and $\msgx$ is an even function of $x$: $\msjv{-x} = \msgx$.
\end{theorem}
\proof
For the purpose of calculating $\msgx$ it is better to use the inverse
function derivative formula\rx{5.13**} in order to rewrite \ex{6.4} as
\qq
\msgx = \det\FLgmx\p (x).
\label{6.16*1}
\qqq
The reason is that $\FLgmx\p$ is much easier to calculate. Indeed, $\FL$ is
defined by \ex{6.3} as a composition of two maps,
while $\Ev{-x}{x}=0$ and the
derivative $\del_{z}\Srv{z}{a}$ at $z=0$, $a=0$ is the identity matrix.
Therefore,
\qq
\FLgmx\p(x) = \atv{\lrbcs{\del_x \Erv{y}{x} } }{y=-x}.
\label{6.16*2}
\qqq
The formula for the derivative
$\atv{\lrbcs{\del_X \Ev{Y}{X} } }{Y=-X}$
of $\Ev{Y}{X}$, all of whose edges and indices are total,
follows from the combinatorial relation which is similar to that of \ex{1.44}:
\qq
e^{-X} e^{X+\Delta X} - 1\approx {1 - e^{-\ad_X}\over \ad_X}\,\Delta X
= \lrbc{ 1 + \snoi {(-\ad_X)^n\over (n+1)!} }\Delta X.
\label{6.16*3}
\qqq
Therefore,
\qq
\atv{\lrbcs{\del_X \Ev{Y}{X} } }{Y=-X} =
I + M_X, \qquad M_X=\snoi{(-\mmad_X)^n\over (n+1)!}.
\label{6.16*4}
\qqq
where $\mmad_X$ is a matrix field with total indices.
In order to express the derivative\\ $\atv{\lrbcs{\del_x \Ev{y}{x} } }{y=-x}$,
we have to replace the total legs $dX$ and $\del_X$ of $M_X$ (and $I$) with
root legs $dx$ and $\del_x$, and replace the total $X$ coordinate label with
root label $x$. We denote the corresponding operator as $\Mrx$. Then we find
from the definition\rx{5.22} that
%
%
\qq
\msgx & = & \det\FLgmx\p (x) = \det \atv{\lrbcs{\del_x \Erv{y}{x} } }{y=-x}
= \det(I+\Mrx) = \exp\lrbcs{ \Tr\lgt(I+\Mrx) }
\nonumber\\
& = &  \exp\lrbc{ -\snoi {(-1)^{n}\over n}\,\Tr(\Mrx^{\times n}) }.
\label{6.16*5}
\qqq
Since each trace $\Tr(\Mrx^{\times n})$ is obviously represented by a wheel
graph, then $\omega(x) = \log \msgx$ is indeed a sum of wheels. Note that these
wheels are different from those related to the wheeling map and Duflo
isomorphism, because $n$ internal edges of $\Tr(\Mrx^{\times n})$ are root
rather than total due to the fact that we passed from $M_X$ to $\Mrx$.

In order to prove that $\msgx$ is even, we observe that
$(\mmad_X)^\rT = - \mmad_X = \mmad_{-X}$, hence $M_X^\rT=M_{-X}$ and
$\Mrx^\rT=-\Mrmx$. Then, in view of \eex{5.15x},
\qq
\Tr(\Mrx^{\times n}) = \Tr(\Mrx^{\times n})^\rT =
\Tr(\Mrmx^{\times n})
\label{6.16*6}
\qqq
and hence $\msgx=\msjv{-x}$. \qed

Comparing \eex{6.16*5} and\rx{1.49} we see that
\qq
\msravr = \mmrcr\;\dgva\ \Tg \msgvx.
\label{6.16*7}
\qqq

For $\xN$ root coordinates $\ux$ we denote
\qq
\umsgx = \pjoxN \msjv{x_j}.
\label{6.24*5}
\qqq
Theorem\rw{t6.5} has an obvious
\begin{corollary}
\label{c6.x1}
The integration measure $\umsgx$ is $\ux$-substantial.
\end{corollary}

\subsection{Properties of the coadjoint orbit integral}
\label{5xs.2}

So far, we have described a graph analog of a coadjoint orbit defined as
$\xG/\xT$. That quotient space is related to an actual coadjoint orbit
$\corba$ by the isomorphism\rx{1.33}. Its graph analog is a
map
$\bfx = \DRax$ (where $x$ is
a root argument and $a$ is a Cartan parameter) defined as
\qq
\DRax = \mAdev{x}(a).
\label{6.17}
\qqq
\begin{lemma}
\label{l6.3}
The map $\DRa$ is \nrw\ and it satisfies the following properties
\qq
&\DRa\circ\mAdev{b} = \mAdev{b}\circ\DRa,
\label{6.18}
\\
&\DRa\circ\FLgy = \mAdev{y}\circ\DRa.
\label{6.19}
\qqq
\end{lemma}
\proof
The \nrw ness of $\bfz=\DRax$ is obvious.
To prove \ex{6.18}, we calculate its \lhs starting with the
definition\rx{6.17} and then applying \ex{5.8*3} and the relation
\qq
\mAdev{b}(a) = a,
\label{6.20}
\qqq
which follows from $\CC$ relations. Thus
\qq
\DRa\circ\mAdev{b}(x) & = & \mAdesv{\mAdev{b}(x)}(a)   =
\mAdev{b}\circ\mAdev{x}\circ\mAdev{-b}(a)
\nonumber\\
& = &
\mAdev{b}\circ\mAdev{x}(a) = \mAdev{b}\circ\DRa.
\label{6.21}
\qqq
In order to prove \ex{6.19}, we start with its \rhs
\qq
\mAdev{y}\circ\DRa (x) & = & \mAdev{y}\circ \mAdev{x}(a)
= \mAdesv{\Ev{y}{x}}(a)
\nonumber\\
& = & \mAdesv{\Ev{\FLgy(x)}{\tFLgy(x)}} (a) =
\mAdesv{\FLgy(x)}\circ\mAdesv{\tFLgy(x)} (a)
\nonumber\\
& = & \mAdesv{\FLgy(x)}(a) = \DRa\circ\FLgy.
\label{6.22}
\qqq
\qed

The substitution\rx{6.17} converts a function $G(\ubfx)$ of total coordinates
$\ubfx$ into a function $\Gruax(\ux;\ua)$ of root coordinates $\ux$ and Cartan
coordinates $\ua$ defined as
\qq
\Gruax(\ux;\ua) = G(\DRav{1}(x_1)\ldtc\DRav{\xN}(x_\xN)).
\label{6.24*3}
\qqq
Since we will be interested in the integrals of functions $\Gruax$, then
following Remark\rw{rq2}, we go further and apply the map $\mFiq$ to $\Gruax$
assuming that labels $\ua$ belong to $\Xp$:
\qq
\Grua(\ux) = \mFiq( \Gruax(\ux;\ua)).
\label{6.24*3x}
\qqq
In other words, $\ua$ are no longer parameters in $\Grua(\ux)$, rather their
legs are hidden in symmetric algebras of cohomologies of graphs.

\begin{theorem}
\label{t6.3}
Let $G(\bfz)$ be a function such that
$\Grua(x)=\mFiq(G(\DRax))$ is \SPN\ as a function of $x$.
Then
\qq
\FGi dx\;  G\lrbcs{ \mAdev{b}\circ\DRax }\msgx
& = & \FGi dx\; G(x)\,\msgx,
\label{6.23}
\\
\FGi dx\; G\lrbcs{ \mAdev{y}\circ\DRax}\msgx
& = & \FGi dx\; G(x)\,\msgx.
\label{6.24}
\qqq
\end{theorem}
\proof
We will prove the first formula, since both proofs are essentially the same.
According to \ex{6.18},
\qq
\FGi dx\; G\lrbcs{ \mAdev{b}\circ\DRax} \msgx =
\FGi dx\; G\lrbcs{ \DRa\circ\mAdev{b}(x) } \msgx.
\label{6.24*}
\qqq
Now \ex{6.23} follows from Corollary\rw{c5.2}, where
$G(y)$ is $G( \DRa(y) )$
(this corollary is applicable, since
$y=\mAdev{b}(x)$ is \nrw),
and from Theorem\rw{t6.1} which states the
$\mAdev{b}$-invariance of the measure $\msgx$.\qed

%

In what follows we will be interested in functions $G(\ubfx)$ such that
$\Grua(\ux;\ua)$ considered as functions of $\ux$ are \SPN, so that we can
define an integral $I(\ua) = \FGi d\ux\;\Grua(\ux;\ua)\,\umsgx $.

\begin{example}
\rm
\label{e6.1}
For $\xN=1$ let $\gxG(X;b) = \exp(\cstv{X}{b})$, where $b$ is a Cartan
parameter. We assume that $a,b\in\Xp$, so
$\Grua(x;a,b) = \exp(ab\blt\, +\, \hlf ab\rstv{x}{x}\, +\, \cO(x^3))$.
Obviously, $\Grua(x;a,b)$ is \SPN\ in $x$ and $Q(x) = - ab\rstv{x}{x}$.

\begin{conjecture}[Duistermaat-Heckmann for graphs]
\label{cdh}
The stationary phase
integral $\FGi dx\; \Grua(x;a,b)\,\msgx$ is `1-loop exact', that is
\qq
\FGi dx\; \exp\lrbcs{ \cstv{\DRax}{b} }\,\msgx = \exp(ab\;\blt)\,
\exp\lrbcs{-\hlf  \log(ab) \rcrc
}.
\label{6.24x*1}
\qqq
\end{conjecture}
In other words, this conjecture says that for this particular integral, the term
$$\gxP(x)
\glxxLRo
\eQixo$$ of \ex{5.45} is equal to $\gxP(0)$.
We hope to prove this conjecture in\cx{Ron} by modifying the standard
Duistermaat-Heckman proof for graph algebra case with the help of grassmanian
variables.
\end{example}

We are particularly interested in integrals which come from
$\mAdev{y}$-invariant functions $\gxG(\ubfx)$:
\qq
G(\mAdev{y}(\bfx_1)\ldtc\mAdev{y}(\bfx_\xN)) = G(\bfx_1\ldtc\bfx_\xN).
\label{6.24*2}
\qqq
Suppose that $\gxG(\ubfx)$ is $\mAdev{y}$-invariant and that
$\Grua(\uxz)$ is \SP. Denote the (twice) corresponding strut exponent
of $\Grua(\uxz)$ as
\qq
Q(\uxz) = \sijoxN \lij(\ua) \rstv{x_i}{x_j}.
\label{6.24*2x1}
\qqq
We denote
\qq
Q_k(\uxrkz) = \atv{Q(\uxz)}{x_k=0}.
\label{6.24*2x2}
\qqq
The matrix of the `quadratic form' $Q_k$ is $\mtrk{\lij(\ua)}$ which is
the $(k,k)$-minor of the full matrix $\mtr{\lij(\ua)}$.
\begin{lemma}
\label{lzro}
If $\gxG(\ubfx)$ is $\mAdev{y}$-invariant, then
\qq
\det\mtr{\lijua}=0
\label{dzr}
\qqq
and
the determinant
$\det\mtrk{\lij(\ua)}$ does not depend on $k$.
\end{lemma}
\proof
The claim of this lemma follows easily from the following relation
\qq
\sjoxN \lij(\ua) = 0,\qquad 1\leq i\leq \xN.
\label{6.24*2x5}
\qqq
In order to prove it, we observe that
since $\gxG(\ubfx)$ is $\mAdev{y}$-invariant, then, in view of \ex{6.19},
$\Grua(\uxz)$ is $\FLgy$-invariant, which means that
\qq
\Ldv{\bxiFy}\Grua(\uxz) = 0,
\label{6.24*2x3}
\qqq
where $\bxiFy(\ux) = \sjoxN\xiFy(x_j)$. Therefore
\qq
\Ldv{(\bxiFy)\str}Q(\uxz) = 0,
\label{6.24*2x4}
\qqq
where $(\bxiFy)\str$ is determined by \ex{6.12*x}. Equation\rx{6.24*2x4}
implies \ex{6.24*2x5}.\qed

Since $\det\mtrk{\lij(\ua)}$ does not depend on $k$, then neither does
$\det\hgxQ_k$, $\hgxQ$ being defined by \ex{5.44}. Thus we denote these
determinants as
\qq
\detp \mtr{\lij(\ua)} = \det \mtrk{\lij(\ua)},\qquad
\detp \hgxQ = \det \hgxQ_k.
\label{6.24*2x6}
\qqq

The following theorem will imply the independence of the universal \urcc
invariant of the choice of the broken link component.
\begin{theorem}
\label{c6.3}
Let $G(\ubfx)$ be an $\mAdev{y}$-invariant function of $\xN$ total coordinates
$\ubfx$
\qq
G(\mAdev{y}(\bfx_1)\ldtc\mAdev{y}(\bfx_\xN)) = G(\bfx_1\ldtc\bfx_\xN).
\label{6.24*2x}
\qqq
Suppose that for some $k$ ($1\leq k\leq\xN$) the function
$\atv{\Grua(\uxz)}{x_k=0}$
is non-degenerate \SP. Then the same is true for
any other $j$, $1\leq j\leq \xN$ and
\qq
\FGi d\ux\rv{j} \atv{\lrbcs{\Grua(\uxz)\,\umsgx}}{x_j=0} =
\FGi \duxrvk \atv{\lrbcs{\Grua(\uxz)\,\umsgx}}{x_k=0}.
\label{6.24*4}
\qqq
\end{theorem}
\proof
In order to simplify our notations, let us assume that $j=1$ and $k=2$.
Then assuming that the individual integral over $x_1$ is non-degenerate,
we can rewrite the \rhs of \ex{6.24} as
\qq
\lefteqn{
\FGi \dxtxN\; \jxtxN \lrbc{
\FGi dx_1\; \msjv{x_1}\; G\lrbcs{\DRav{1}(x_1),a_2,\DRav{3}(x_3)
\ldtc \DRav{\xN}(x_\xN) }
}
}
\nonumber\\
& = &
\FGi \dxtxN \; \jxtxN \Bigg(
\FGi dx_1\;\msjv{x_1}\;
\label{6.24*6}
\\
&&\hspace{-2pt}\hspace{-0.5cm}\times
G\lrbcs{\mAdev{-x_1}\circ\DRav{1}(x_1), \mAdev{-x_1}(a_1),
\mAdev{-x_1}\circ \DRav{3}(x_3)\ldtc \mAdev{-x_1}\circ\DRav{\xN}(x_\xN)
}\Bigg)
\nonumber\\
& = & \hspace{-2pt}
\FGi \dxtxN\; \jxtxN \lrbc{
\FGi \hspace{-0.7pt}
dx_2\; \msjv{x_2}\; G\lrbcs{a_1,\DRav{2}(-x_2),\DRav{3}(x_3)
\ldtc \DRav{\xN}(x_\xN) }
}
\nonumber\\
& = & \hspace{-2pt}
\FGi \dxtxN\; \jxtxN \lrbc{
\FGi \hspace{-0.7pt}
dx_2\; \msjv{x_2}\; G\lrbcs{a_1,\DRav{2}(x_2),\DRav{3}(x_3)
\ldtc \DRav{\xN}(x_\xN) }
}.
\nonumber
\qqq
Here the first equality is due to the $\mAdev{y}$-invariance of
$G(\ubfx)$, the second equality is due to \ex{6.24} of Theorem\rw{t6.3}
and the definition of $\DRa$ and the third equality is due to
parity-invariance of the gaussian integral and integration measure
$\msgx$.\qed
\begin{remark}
\label{r6.2}
\rm
Our assumption that the integral over $x_1$ is non-degenerate is not
restrictive. If the integral is degenerate, then we can add an
$\Ad$-invaraint `regularizing' term (\eg $\e\tstv{X_1}{X_2}$)
to the strut exponent so that the
integral becomes non-degenerate. After proving \ex{6.24*4} in the
regularized form we remove the regularization by setting $\e=0$.
\end{remark}

\begin{lemma}
\label{c6.2}
Let $\zA(\bfx)$ be a tensor field such that
$\zA(\DRax)$ is \SPN. Consider another
tensor field
\qq
\zchA(\bfx;\bfy) = \ndv{\mad_\bfy(\bfx)}\zA(\bfx),
\label{6.24*1}
\qqq
which
depends on $\bfx$ and on a total parameter $\bfy$. Then
\qq
\FGi dx\; \zchA(\DRax;\bfy)\,\msgx = 0.
\label{6.25}
\qqq
\end{lemma}
\proof
Since a total leg is a sum of a root leg and a Cartan leg, then it is
sufficient to prove \ex{6.25} for two cases: when $\bfy$ is replaced either by
a Cartan parameter $b$ or by a root parameter $y$. Since both proofs are
essentially the same, we will prove the $b$ case.
It is easy to see that
\qq
\zchA(\bfx;b) = \dlt\!\atv{ \zA(\mAdev{tb} (\bfx))}{t=0}.
\label{6.26}
\qqq
Then according to \ex{6.23},
\qq
\FGi dx\; \zchA(\DRax;b)\,\msgx & = &
\atv{{\dlt\FGi dx\; \zA\lrbcs{ \mAdev{tb}\circ\DRax }\msgx} }{t=0}
\nonumber\\
& = &
\atv{{\dlt\FGi dx\; \zA(\DRax)\,\msgx}}{t=0} = 0.
\label{6.27}
\qqq
\qed

For a total valued vector field $\xXi(\bfx)$ denote
\qq
\dvadv{\bfx} \xXi(\bfx) = \gllv{\del_\bfx}{\bfy}{1}
\lrbcs{ \ndv{\mad_\bfy(\bfx)}\xXi(\bfx) }
\label{6.28}
\qqq
(the name $\dvadv{}$ is a combination of $\dvrgv{}$ and $\ad$).
\begin{figure}[htb]
\begin{center}
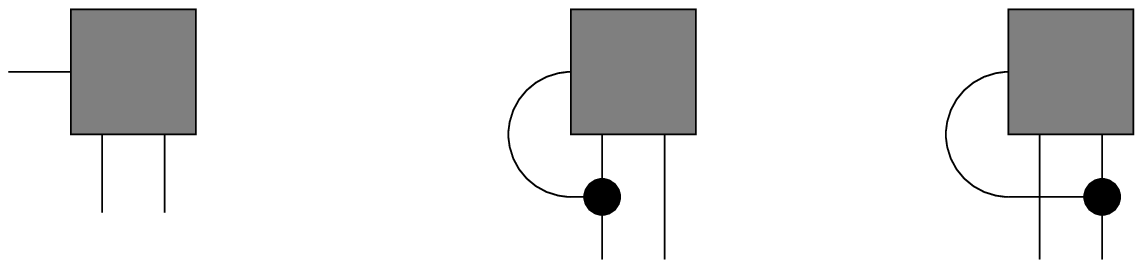
\end{center}
\caption{The definition of $\dvadv{\bfx} \xXi(\bfx)$}
\label{f5.2}
\end{figure}
\begin{remark}
\rm
\label{r6.1}
It is easy to see that the definition\rx{6.28} can be modified:
\qq
\dvadv{\bfx} \xXi(\bfx) = \dvrgv{\bfx} \lrbcs{\mmad_{\bfx}\times\xXi(\bfx)}.
\label{6.32}
\qqq
(note that $\dvrgv{\bfx}\mad_{\bfy}(\bfx)=0$, since
the resulting graph is zero because of the AS relation).
\end{remark}

The next theorem is essential for the proof of the topological invariance of
the universal \urcc invariant.
\begin{theorem}
\label{t6.4}
Let $\xXi(\bfx)$ be a total valued vector field such that $\xXi(\DRax)$ is \SPN.
Then
\qq
\FGi dx\; \dvadv{\bfx}\xXi(\DRax)\,\msgx = 0.
\label{6.29}
\qqq

\end{theorem}
\proof
We modify the \lhs of \ex{6.29} by using the definition\rx{6.28}, commuting the
gluing $\gllv{\del_\bfx}{\bfy}{1}$ with the formal gaussian integral and then
applying\rx{6.25} with $A=\xXi$
\qq
\FGi dx\; \dvadv{\bfx}\xXi(\DRax)\,\msgx & = &
\FGi dx\; \gllv{\del_\bfx}{\bfy}{1}
\lrbcs{ \ndv{\mad_\bfy(\bfx)}\xXi(\bfx) }\,\msgx
\nonumber\\
& = & \gllv{\del_\bfx}{\bfy}{1} \lrbc{
\FGi dx\; \ndv{\mad_\bfy(\bfx)}\xXi(\bfx)\,\msgx }
\nonumber\\
& = & 0.
\label{6.29*}
\qqq
\qed

\nsection{A universal \urcc invariant}
\label{xs6}

\subsection{Definition and invariance}
\label{6xs.1}

Let $\cLM$ be an $\xN$-component oriented marked Morse
link in a rational homology sphere $M$. We denoted by $\ZKBL$
its Kontsevich integral in $\cBL$. By replacing the $\cL$ labels
$1\ldtc \xN$
at the legs of
$\ZKBL$ with the corresponding variables
$\ubfx=(X_1\ldtc X_{\xN})$,
we turn Kontsevich integral
into a function $\ZKBLx$.
In subsection\rw{3xs.1*} we described an injection $\ftr$ of the algebra
$\cBL$ into $\cBCL$, which transforms a graph of $\cBL$ into the same graph
in $\cBCL$, all of whose edges are total.
Let us denote as $\ZKBCLx$ the
Kontsevich integral
of $\cLM$ in $\cBCL$:
\qq
\ZKBCLx = \ftr( \ZKBLx ).
\label{6.32xx}
\qqq
Obviously, $\ZKBCLx$ is also a function of $\xN$ total variables $\ubfx$.


\begin{lemma}
\label{l6.x1}
Kontsevich integral $\ZKBCLx$ is both $\mAdev{y}$- and $\mAdev{b}$-invariant as
a function of coordinates $\ubfx$.
\end{lemma}
\proof
Let us prove the invariance under $\mAdev{y}$, the proof for $\mAdev{b}$ is
similar. Since $\mAdev{ty}$ forms a 1-parametric group, then it is enough to
establish that
\qq
\Ldv{\mady(\ubfx)}\ZKBCLx=\nabla_{\mady(\ubfx)}\ZKBCLx=0,
\label{6.32*}
\qqq
where $\mady(\ubfx) = \sjoxN \mady(\bfx_j)$, so that $\nabla{\mady(\ubfx)}$ acts
on all coordinates $\ubfx$. Equation\rx{6.32*} follows easily from
Lemma\rw{l5.4} and from the definition\rx{6.32xx}, which establishes that all
edges of $\ZKBCLx$ are total.\qed

Our next step is to perform a substitution\rx{6.24*3} and then
convert the legs carrying Cartan coordinates $\ua$ into the elements of graph
cohomology algebras through the application of the map $\mFiq$ as in
\ex{6.24*3x}. We
define $\cDCL$ as an algebra $\cDCXp$, where $\Xp=\{\ua\}$ is a set of
$L$ Cartan labels and $X\setminus\Xp=\{\ux\}$ is a set of $L$ root labels. We
also define $\cQDCL$ as the corresponding algebra $\cQDCXp$. Now we convert
Kontsevich integral into an element of $\cDCL$ by the formula
%
%
\qq
\ZKDCrLx = \mFiq\lrbcs{\ZKBCLv{\DRav{1}(x_1)\ldtc\DRav{\xN}(x_\xN)}}.
\label{6.33}
\qqq
%


\begin{lemma}
\label{l6.4}
The function $\ZKDCrLx$ is \SP.
\end{lemma}
\proof
In view of \ex{1.8}, this function is an exponential
\qq
\ZKDCrLx = \exp\lrbcs{\WKDCrLx},
\label{6.34*1}
\qqq
where
\qq
\WKDCrLx = \mFiq\lrbc{ \WKBCLv{\DRav{1}(x_1)\ldtc\DRav{\xN}(x_\xN)} }
\label{6.34*2}
\qqq
and $\WKBCLx = \ftr(\WKLX)$.
Since $\WKL$ and $\bfx = \DRax$ are both \nrw,
then so is $\WKDCrLx$.
In particular, this means that the dependence of $\ZKDCrLx$ on struts is
purely exponential.\qed

Following \ex{5.42}, let $\gxQLxa$ denote (twice) the strut part of $\WKDCrLx$
\qq
\gxQLxa = 2(\WKDCrLx)\str
\label{6.34x1}
\qqq
We present $\gxQLxa$ as
\qq
&\gxQLxa = \sijoxN \lijua\; \rstv{x_i}{x_j},
\label{6.34x2}
\\
&
\lijua \in \IQuua\subset\QIQuua,\qquad l_{ji}(\ua) = l_{ij}(-\ua)
\nonumber
\qqq
(the origin of the coefficients $\lijua$ from the haircomb graphs of $\WKL$
will be explained in details in subsection\rw{6xs.1*}).
$\ZKDCrLxk$ is also \SP, and its quadratic exponent part
is given by
\qq
\xQLxak = \atv{\gxQLxa}{x_k=0}.
\label{6.35x1}
\qqq
%


Since Lemma\rw{l6.x1} claims that $\ZKDCrLx$ is $\mAdev{y}$-invariant, then
we are going to study the integral
\qq
\ZKDCurLa =
\FGi d\ux\rv{k} \atv{\lrbcs{\ZKDCrLx\,\umsgx}}{x_k=0}, \qquad
1\leq k\leq \xN
\label{6.35}
\qqq
(\cf \ex{6.24*}), which is an invariant of the oriented marked Morse link
$\cL\subset M$.

\begin{theorem}
For $1\leq k\leq \xN$,
the function $\ZKDCrLxk$ is \SPN\ (that is, $\mtrk{\lijua}$ is non-degenerate)
iff
\qq
\AFLut\not\equiv 0.
\label{6.34x3}
\qqq
Moreover,
\qq
{\dtplijua
\over \AFLeua}
= \pjba
(1 + \cO(\ua)).
\label{6.34x4}
\qqq
\end{theorem}
\proof
It is sufficient to prove \ex{6.34x4}.
Our proof is based on the results of\cx{Ro2}. It is easy to see that
\qq
\Tg \hQkua = \gQLuak,
\label{6.34x5}
\qqq
the operator $\gQLuak$ being defined by \ex{1.60}. Therefore
\qq
\Tg \det \hQkua = \det\gQLuak,
\label{6.34x6}
\qqq
or, more explicitly, according to \ex{4.23*1},
\qq
\lrbc{\prvldp \det\mtrk{\lij( \scp{\vua}{\vl})} }^2= \det\gQLuak.
\label{6.34x6*}
\qqq
If $\mfg=su(2)$, then $\mfh$ is 1-dimensional and there is only one positive
root $\vl$, so relation\rx{6.34x6*} becomes
\qq
\det \Qstkss = \lrbcs{\det \mtrk{\lij(\ua)} }^2,
\label{6.34x7}
\qqq
where we use an $su(2)$ $\mfhs$ coordinate $a=\scp{\va}{\vl}$.
On the other hand, it was established in\cx{Ro2} that
\qq
{\lrbcs{\det\Qstkss}^{1/2} \over  \AFLv{e^{\ua}} } =
\lrbc{\pjoxN a_j} (1 + \cO(\ua)).
\label{6.34x8}
\qqq
Equation\rx{6.34x4} follows from \eex{6.34x7} and\rx{6.34x8}.\qed

\begin{remark}
\label{rq3}
\rm
Equation\rx{6.34x4} indicates that $\dtplijua\Big|_{\ua=0}=0$, which means that
although $\ZKDCrLxk$ may be non-degenerate, its predecessor
$\ZKBCLv{\DRav{1}(x_1)\ldtc\DRav{\xN}(x_\xN)}$ is always degenerate and
therefore the integral\rx{6.35} would not be well-defined in $\cBCL$ (\cf
Remark\rw{rq2}).
\end{remark}


\begin{theorem}
\label{t6.6}
The integral\rx{6.35}
does not depend on the choice of $k$ ($1\leq j\leq \xN$) in
\ex{6.35}.
\end{theorem}
\proof
This theorem is a particular case of Theorem\rw{c6.3}.\qed
\begin{theorem}
\label{t6.7}
The integral\rx{6.35}
is an invariant of an oriented link $\cLSt$, that is, it does not
depend on a presentation of $\cL$ as a marked Morse link.
\end{theorem}

Before we prove this theorem, we have to describe how Kontsevich integral
$\ZKBLx$ depends on the choice of marked points on the components of $\cL$.
Let us introduce a cyclic smooth parameter $t$ on $\cL_1$. Then $\ZKBLxt$
depends explicitly on a position $t$ of the marked point on $\cL_1$.
\begin{lemma}
\label{l6.5}
The derivative $\del_t\ZKBLxt$ has a form
\qq
\del_t\ZKBLxt = \dvadv{X_1}\lrbcs{\xi(\ubfx)\,\ZKBLxt},
\label{6.36}
\qqq
where $\xi(\ubfx)$ is a vector field with $\del_{X_1}$ component.
The vector field $\xi(\ubfx)$ is \nrw\
and therefore
has a polynomial (actually, at most linear) dependence on the struts.
\end{lemma}
\pr{Theorem}{t6.7}
This theorem follows immediately from \ex{6.36} and
Theorem\rw{t6.4}.\qed

\pr{Lemma}{l6.5}
Consider graph algebras $\cAxy$ and $\cAz$.
There is an $\cA$-multiplication map
$\xmap{\mxyz}{\cAxy}{\cAz}$. With a slight abuse of notations we will use the
same notation for an $\cA$-multiplication map $\xmap{\mxyz}{\cBxy}{\cBz}$
induced on $\cB$ algebras by the
isomorphisms $\xchi$ (see the diagram\rx{1.17}). We will also define
$\mcxyz = \mxyz - \myxz$.
Now consider a vector field
$\xXi(x)$. Recall that its graph has a few $X$ legs and exactly one $\del_X$
leg. The relation
\qq
\xmltv{[\del_X,X]}{X}\xXi(X) = \dvadv{X} \xXi(X)
\label{6.37}
\qqq
follows easily from a combinatorial identity
\qq
Y\,X^n - X^n Y = \sum_{m=0}^{n-1} X^m\, [Y,X]\, X^{n-m-1}.
\label{6.38}
\qqq

Consider two marked points $t$ and $t+\DLt$ on a link component $\cL_1$. They
split it into two pieces (one big and one small). Thus we get a modified
Morse `link'
$\cLp$ in which a closed component $\cL$ is replaced by two segments. Consider a
Kontsevich integral $\ZKALp\in \cALp$ and its PBW-symmetrized version
$\ZKBLpe\in\cBLp$. Since
\qq
\ZKALpv{t} = \xmltv{\DLxo,X_1}{X_1} \ZKALp, \quad
\ZKALpv{\tDLt} = \xmltv{X_1,\DLxo}{X_1} \ZKALp,
\label{6.39}
\qqq
then
\qq
\Delta\ZKBLxt = \ZKBv{\cL;\ubfx;\tDLt} - \ZKBv{\cL;\ubfx;t}
= - \xmltv{[\DLxo,X_1]}{X_1} \ZKBLpe.
\label{6.40}
\qqq
As usual, $\ZKBLpe = \exp\lrbcs{\WKLpe}$, where $\WKLpe$ is \nrw.
We are interested only in the terms which are at most linear
in $\DLt$. It is easy to see from the definition of Kontsevich integral
that every $\DLx$ leg in the graphs of $\WKLpe$ carries with it at least a
factor of $\DLt$. Therefore $\atv{\WKLpe}{\DLt=0} = \WKLte$ and the graphs whose
coefficients are linear in $\DLt$ must have exactly one $\DLxo$ leg. Thus if we
take a sum of all graphs of $\WKLpe$ whose coefficients are proportional to
$\DLt$ are replace there $\DLxo$ with $\dlxo$, then the result can be declared
a vector field $\xi(\ubfx)\,\DLt$. Now
\qq
\Delta\ZKBLxt \approx - \xmltv{[\dlxo,X_1]}{X_1}\lrbcs{
\xi(\ubfx)\,\ZKBLxt } \DLt
\label{6.41}
\qqq
and \ex{6.36} follows in view of \ex{6.37} in which $\xXi(X)$ is replaced by
$\xi(\ubfx)\,\ZKBLxt$.\qed

Theorems\rw{t6.6} and\rw{t6.7} allow us to formulate the following
\begin{definition}
For a link $\cL$ in a rational homology sphere $M$ such that
$\AFLut\not\equiv 0$, the unwheeled graph \urcc\ invariant
$\ZKDCurLa$ is defined by the integral\rx{6.35}.
\end{definition}

If a link has a single component ($\xN=1$), then there is no integration
in\rx{6.35}, hence $\ZKDCurLa$ is obtained from $\ZKDCLx$ by changing its
legs from total to Cartan. In other words, for a knot $\cK$,
\qq
\ZKDCurKa = \mFiq\circ\fct(\ZKBK),
\label{6.41y}
\qqq
or, less formally, $\ZKDCurKa = \ZKBKa$, since according to our conventions,
replacing the total labels $\bfx$ of $\ZKBKx$ by Cartan labels $a$ makes its
total legs Cartan.
Theorem\rw{t4.1}, $\fct$ is an injection in a knot case, so no information
about a knot is lost in passing from the traditional Kontsevich
integral $\ZKBKx$ to its Cartan-restricted version $\ZKDCurKa$.

The dependence of $\ZKDCurLa$ on the orientation of $\cL$ is quite transparent.
For an oriented link $\cL$ let $\cLroj$ denote the same link in which we
switched the orientation of the $j$-th component $\cL_j$ and let $\bcL$
denote the link $\cL$ in which we switched the orientation of all components.

\begin{theorem}
\label{tori}
For an oriented link $\cL$ with $\ndgap$
\qq
\ZKDCurv{\cLroj} = \xPsj\lrbcs{\ZKDCurLa},\qquad
\ZKDCurv{\bcL} = \ZKDCurLa,
\label{6.41y1}
\qqq
where the parity operators $\xPsj$ are defined by \eex{4.7*} and\rx{4.13**3x}.
\end{theorem}
\proof
It is well-known (and easy to see from its definition) how Kontsevich integral
changes when orientation of link components is switched:
\qq
\ZKBv{\cLroj;\ubfx} = \ZKBv{\cL;X_1\ldtc -X_j\ldtc X_{|\Xp|} },\qquad
\ZKBv{\bcL;\ubfx} = \ZKBLx.
\label{6.41y2}
\qqq
Therefore, since $\DRv{-a}(x) = -\DRax$, then
\qq
\ZKDCrv{\cLroj;\ux} = \xPsj\lrbcs{ \ZKDCrLx },\qquad
\ZKDCrv{\bcL;\ux} = \ZKDCrLx,
\label{6.41y3}
\qqq
and \eex{6.41y1} follow by applying the integral of \ex{6.35} to these
equations.\qed

In order to make connection with $\mfg$-based \urcc invariant of
subsection\rw{2xs.4},
we have to construct a graph algebra \urcc invariant from the wheeled Kontsevich
integral\rx{1.20*}. We follow the same steps that led us to \ex{6.35} and
introduce
\qq
\ZKODCurLa =
\FGi d\ux\rv{k} \atv{\lrbcs{\ZKODCrLx\,\umsgx}}{x_k=0}, \qquad
1\leq j\leq \xN,
\label{6.42}
\qqq
where again
\qq
\ZKODCrLx = \ZKOBCLv{\DRav{1}(x_1)\ldtc\DRav{\xN}(x_\xN)}.
\label{6.43}
\qqq
\begin{theorem}
\label{t6.8}
$\ZKODCrLx$ depends neither on the choice of $k$ ($1\leq k\leq \xN$)
nor on the
presentation of $\cL$ as a marked Morse link.
\end{theorem}
\proof
The proof of the first claim is exactly the same as the proof of
Theorem\rw{t6.6}. In order to
prove the second claim, note that according to
Lemma\rw{l5.4}, the wheeling operator $\WhhL$ (which is a particular case of
gluing $\mad$-invariant graphs)
commutes with
ad-divergence $\dvadv{X_1}$. Therefore \ex{6.36}
implies a similar relation for
$\ZKOBLxt$:
\qq
\del_t\ZKOBLxt = \WhhL\; \dvadv{x_1}\lrbcs{\xi(x)\,\ZKBLxt} =
\dvadv{x_1}\WhhL\lrbcs{\xi(x)\,\ZKBLxt},
\label{6.44}
\qqq
The rest of the proof is similar to that of Theorem\rw{t6.7}.\qed

\begin{definition}
For a link $\cL$ in a rational homology sphere $M$ such that
$\AFLut\not\equiv 0$, the wheeled graph \urcc\ invariant is defined by
the integral\rx{6.42}.
\end{definition}

\begin{remark}
\rm
Since graphs of wheeling operators have even numbers of legs, then wheeled
Kontsevich integral still satisfies \eex{6.41y2}. Hence the wheeled graph \urcc
invariant satisfies \eex{6.41y1} of Theorem\rw{tori}.
\end{remark}

The wheeled graph \urcc invariant determines the Lie algebra based \urcc
invariant\rx{1.51}. In order to see this, we introduce a weight system
$\Tgth$ which acts in the same way as $\Tgh$, except that it multiplies the
graphs by $\hb^{\dgth{D}}$ rather than by $\hb^{\dgh{D}}$. Since the exponent
of \ex{1.28} comes from the \rhs of \ex{1.25*} (which is the exponent of
\ex{1.26}) through the substitution $\uvbet=\hb^{-1}\,\uvb$, and since
according to \ex{1.11x}, $\dgth{D} = \dgh{D}-\#(\mbox{1-vertices})$, then
\qq
\smnd \LLbmn\,\hb^{n-1} = \Tgth\, \WKOL.
\label{6.44*}
\qqq
where $\uvb$ are the elements of $\mfgs$ which are placed at the legs of the
graphs of $\WKOL$.
Since, according to Theorem\rx{t6.5}, $\msgx$ is an exponential of the sum of
wheel graphs, then
\qq
\dgth{\msgx} = 0.
\label{6.44*1}
\qqq
Therefore we can replace $\Tg$ by $\Tgth$ in \ex{6.16*7}:
\qq
\msravr = \mmrcr\;\dgva\ \Tg \msgvx.
\label{6.44*2}
\qqq
A combination of \eex{6.44*} and\rx{6.44*2} indicates that \ex{1.51} can be
rewritten as
\qq
\IrfvaNL & = &
\hbmnDp\,
\mohbpLg\, \dguvarho
\nonumber
\\
&&\qquad
\ilvrkz \duxvrk \Tgth
\atv{\lrbcs{\ZKODCrLx\,\umsgx}}{\vx_k=0}.
\label{6.44*3}
\qqq
Applying \ex{5.45*x3} to this equation we come to the following
\begin{theorem}
\label{t6.9}
The graph algebra \urcc invariant $\ZKODCurLa$ is `universal' in the sense that
it determines Lie algebra based \urcc invariants through the application of the
$\Tgth$ weight system
\qq
\IrfvaNL =
\hbmnDp\;\;\dguvarho\;\;\Tgth\,\lrbcs{ \ZKODCurLa}.
\label{6.45}
\qqq
\end{theorem}

Although Lie algebra based \urcc invariants are completely determined by the
universal invariant, they are still interesting in their own right: at
least, in case of $\mfg=su(2)$ they can be derived from $R$-matrix type
calculation rather than from Kontsevich integral, and they exhibit
interesting integrality properties\cx{Ro1}.

\subsection{Haircomb graphs, wheels and
Kontsevich integral formula for the Alexander polynomial}
\label{6xs.1*}

The Alexander polynomial of a link is hidden inside its \urcc invariant. In
order to see this, first, we have to express the Alexander polynomial in terms
of Kontsevich integral of the link. Let us concentrate on the strut and circle
parts of the logarithm of Kontsevich integral $\WKLX$ after the
substitution\rx{6.34*2}:
\qq
\WKDCrLx = \hlf\gxQLxa + \hxFzLua\rcrc
+ \mbox{other graphs},
\label{6.45*x}
\qqq
where $\gxQLxa$ is the strut expression\rx{6.34x2} and
$\xFzLua\in\evn{\IQuua}$ is a new formal power series. Since series $\lijua$
and $\xFzLua$ will play an important role in our calculations, we take a short
digression in order to explain explicitly how they originate from the graphs of
$\WKL$.

We begin with the strut coefficients $\lijua$.
The strut exponent $\gxQLxa$ comes from the haircomb graphs
of $\WKLX$. A \emph{haircomb} graph is a graph which consists of a single chain
connecting two 1-valent vertices and legs attached to that chain (\cf
\fg{f1.1}). The chain is called a \emph{spine} of the graph and the 1-valent
vertices are called \emph{endpoints} of the haircomb. The same graph may have
up to four different spines depending on the choice of end-points.

A haircomb graph of $\WKLX$ may contribute to $Q(\uxz)$ in two different ways.
The first way is to make a substitution $\bfx_i=\hlf\,\tmadv{x_i}^2(a_i)$ at
one of the legs and $\bfx_i=a_i$ at all other legs. Only the struts of
$\WKLX$ may
contribute in this way (struts are also haircombs). If the strut part of
$\WKLX$ is
\qq
W\istr(\cL;\ubfx) = \hlf \sijoxN \xxlij\;\; \tstv{\bfx_i}{\bfx_j},
\label{6.29zz2}
\qqq
then its first-way contribution to $Q(\uxz)$ is
\qq
\sijoxN \xxlij\; a_i a_j\;\;\rstv{x_i}{x_i}.
\label{6.29zz3}
\qqq
The second way, in which a haircomb graph of $\WKLX$
may contribute to $Q(\uxz)$,
is to perform the substitutions $\bfx_i=\tmadv{x_i}(a_i)$ at its endpoints and
the substitutions $\bfx_i=a_i$ at all other 1-vertices.
Let us quantify this contribution.
For a set of
non-negative integers $\bfn$ let $\bfD_{i,j;\bfn}$ be a set of haircomb graphs
with selected spines such that their endpoints have colors $\bfx_i,\bfx_j$ and
they have $n_k$ legs of color $\bfx_k$ attached to the spine for all $k$,
$1\leq k\leq \xN$. Let $c(\xD)$ be the coefficients of these haircomb graphs in
an expression of $\WKLX$ as an element of $\tcBX$ (we fix the signs of
$c(\xD)$ by assuming that legs are attached to the right of the spine, if it is
oriented from $\bfx_j$ to $\bfx_i$). Now we define the power series
\qq
\lpijua = -\sum_{\bfn\geq 0} \lrbc{\sum_{\xD\in\bfD_{i,j;\bfn}} c(\xD)}
a_1^{n_1}\cdots a_{\xN}^{n_{\xN}}.
\label{6.29zz4}
\qqq
Note that $l\p_{ji}(\ua) = l\p_{ij}(-\ua)$. Now the second contribution of
haircomb graphs to $Q(\uxz)$ is
\qq
\sijoxN \lpijua\;\; \rstv{x_i}{x_j}.
\label{6.29zz5}
\qqq
Since a combination of contributions\rx{6.29zz3} and \rx{6.29zz5} is equal to
the \rhs of \ex{6.24*2x1}, then we conclude that
\qq
\lij(\ua) = a_i \lrbc{a_j\,
\lpijua + \delta_{ij}\sum_{k=1}^{\xN}a_k\, \ell_{ik} }.
\label{6.29zz6}
\qqq
\begin{remark}
\rm
The contribution of diagonal struts $\ell_{ii}\tstv{\bfx_i}{\bfx_j}$ of $\WKLX$
to both expressions\rx{6.29zz3} and\rx{6.29zz5} are equal up to a sign and thus
cancel each other in the final expression\rx{6.29zz6}. Therefore, diagonal
coefficients $\ell_{ii}$ do not participate in the final
expression\rx{6.29zz6}.
\end{remark}

The circle graph contribution $\hxFzLua\rcrc$ to $\WKDCrLx$ comes from the
wheel graphs of $\WKLX$, if we substitute $\bfx_i=a_i$ at all their legs. Thus
let $\bfD_\bfn$ be the set of wheel graphs of $\cBL$ which have $n_i$
legs of color $\bfx_i$ for all $i$, $1\leq i\leq\xN$. Then
\qq
\xFzLua = 2\sum_{\bfn\geq 0\atop n_1+\cdots n_\xN \in 2\ZZ}
\sum_{\xD\in\bfD_\bfn} c(\xD),
\label{6.29zz7}
\qqq
where $c(\xD)$ are again the coefficients at graphs $\xD$ with which they
appear in $\WKLX$.

Since the tree parts of $\WKLX$ and $\WKOLX$ are the same, then
similarly to \ex{6.45*x} we can write
\qq
\WKODCrLx =  \hlf\gxQLxa + \hxFWzLua\rcrc
+ \mbox{other graphs},
\label{6.45*a}
\qqq
The circle parts of $\WKODCrLx$ and $\WKDCrLx$ are different, because the
wheeling of $\WKLX$ produces new wheel graphs by gluing the struts of $\WKLX$
to the legs of the wheels of $\Wh$ defined by \ex{1.14*}. Therefore it follows
from \eex{1.14} and\rx{1.8*} that
\qq
\xFWzLua = \xFzLua - \log\lrbcs{\lknLa},
\label{6.45*a1}
\qqq
where we defined a function $\lknLa$ as
\qq
\lknLa = \pioxN \lrbc{ \sinh\lrbc{ \hxslija}\over \hxslija }.
\label{6.45*a2}
\qqq

Now we can prove a theorem that we neglected to formulate
in\cx{Ro2}:
\begin{theorem}
\label{t6.11}
The Alexander polynomial $\AFLut$ is determined by the
haircomb tree and `1-loop' parts
of the (logarithm of) Kontsevich integral
$\WKOL$:
\qq
\AFLeua & = &
{\dtplijua\over \pja }\;
\exp\lrbcs{- \xFWzLua }
\label{6.53x}
\\
& = & {\dtplijua \over \pja} \;\;
\lknLa
\; \exp\lrbcs{ - \xFzLua}.
\nonumber
\qqq
\end{theorem}
\proof
In view of \ex{6.45*a1} it is sufficient to prove only the first equation
of\rx{6.53x}.

Since an individual polynomial $\LLbmn$ of the \lhs of \ex{6.44*} comes from
the graphs $\xD$ of $\WKOL$ which have $\echi(D)=n-1$ and $m$ legs,
then
\ex{6.44*} implies that
\qq
\smzi \LLamo = \svldp \xFWzLav.
\label{6.45*y}
\qqq
Consider now a case of $\mfg=su(2)$.
Since $su(2)$ has only one
positive root, then (in terms of the $su(2)$ $\mfhs$ coordinate
$a=\scp{\va}{\vl}$)
\qq
\smzi \LLv{m,0}{\ua} = \xFWzLua.
\label{6.53y}
\qqq
Then \ex{6.53x} follows easily from the combination of
equations\rx{1.65y2},\rx{6.34x7} and\rx{6.45*y}.\qed

We give a more direct proof of Theorem\rw{t6.11} in Appendix\rw{A1}. It relies
neither on Cartan and root edges, nor on the properties of the Jones polynomial,
but rather uses a relation between the Alexander polynomial and the supergroup
$\Uoo$.

Note that neither $\dtplijua$, nor $\xFzLua$ is the invariant of an `unmarked'
link $\cL$. Only their combination\rx{6.53x} is independent of the `Morse
marking'.

According to \ex{6.29zz6},
the series $\lijua$ is proportional to both $a_i$
and $a_j$ and, as a result, the determinant $\dtplijua$ of \ex{6.53x} is
divisible by $\pjbai$ in $\IQuua$ (if $L>1$):
\qq
\pjbai \dtplijua \in\IQuua.
\label{6.39*}
\qqq

In view of the relation between the tree part of the Kontsevich integral and
Milnor's linking numbers of the link established by N.~Habegger and
G.~Masbaum\cx{HM}, it is easy to see that the factorization\rx{6.53x} of the
Alexander polynomial $\AFLeua$ into the determinant\rx{6.39*} and the
exponential of $- \xFWzLua$ corresponds to that of Theorem~2 of J.~Levine's
paper\cx{JL}.

\subsection{Alexander
polynomial and the structure of the universal \urcc invariant}
\label{6xs.2}

Let us consider the structure of the graph invariant $\ZKDCurLa$ defined by
\ex{6.35} in more details. Let $\cgDL\subset \gDL$ be the set of
connected graphs of
$\gDL$ which have no legs and such that their $\echi$ is strictly positive.
Also, we introduce a notation
\qq
\tAFLua = \pjba\,\AFLeua.
\label{6.45x}
\qqq
%

\begin{theorem}
\label{l6.6}
The \urcc invariant $\ZKDCurLa$ can be presented in the form
\qq
\ZKDCurLa = \exp\lrbcs{\WKDCurLa}.
\label{6.45*}
\qqq
where
\qq
\WKDCurLa & = &  \hlf\sijoxN \xxlij\, a_i a_j\, \blt
\;-\;\hlf\,\log\lrbcs{\tAFLua
}\rcrc
\nonumber
\\
&&\qquad
+ \sum_{\xD\in\cgDL
} {\xFDL\over
\preDr
\tAFLhue
}
\label{6.45*1}
\qqq
and
\qq
\xFDL\in\cHD.
\label{6.45*2}
\qqq
\end{theorem}
\proof
Consider the logarithm of the integrand of \ex{6.35}
%
\qq
\log \lrbcs{\ZKDCrLx\,\umsgx} = \WKDCrLx + \bomsgx,
\qqq
where $\bomsgx = \sjoxN \omsjv{x_j}$. Let us define  $\WKFDDCrLx\in\cDCX$
by the equation
%
\qq
\lefteqn{
\WKDCrLx + \bomsgx
}
\nonumber
\\
&&=  \hlf\sijoxN \xxlij\, a_i a_j \blt   + \hxFzLua\rcrc
+ \hlf\,\gxQLxa +  \WKFDDCrLx
\label{6.46}
\qqq
as
the `least trivial' part of this logarithm.
 Since $\ZKDCLv{\ua}$ does
not contain the graph $\lcrc$, then
$\WKFDDCrLx$ contains only the graphs $\xD$ such that either $\xD$ has at
least three legs, or $\echi(\xD)=0$ and $\xD$ has at least one leg, or
$\echi(\xD)\geq 1$.

The definition\rx{5.45} applied to the integral\rx{6.35} yields
\qq
\!\!\!
\ZKDCurLa = \exp\hlf\lrbc{\sijoxN \xxlij\, a_i a_j\, \blt
\;+\;
\xFzLua\rcrc
}
\lrbcs{ \detp\hgxQ(\cL) }^{-1/2} \,\FDLgua,
\label{6.47}
\qqq
where by definition
\qq
\lrbcs{ \detp\hgxQ(\cL) }^{-1/2} =
\exp\lrbcs{ -\hlf\,\log\dtplijua\rcrc },
\label{6.47*}
\qqq
while
\qq
\FDLgua =
\exp\lrbcs{ \WKFDDCrLx\Big|_{x_k=0} } \glxxkLR \exp\lrbcs{-\hlf\,
\xQiLxak
}.
\label{6.48}
\qqq
In view of \ex{6.53x},
\qq
\tAFLua = \dtplijua\;\lknLa\; \exp\lrbcs{- \xFzLua },
\label{6.45x1}
\qqq
so we can
rewrite \ex{6.47}
as
\qq
\ZKDCurLa = \exp\hlf\lrbc{
\sijoxN \xxlij\, a_i a_j\, \blt \; -\; \log\lrbc{\tAFLua\over\lknLa}
\rcrc
}
\FDLgua.
\label{6.48*}
\qqq

Now it remains to determine the structure of $\FDLgua$.
Both exponents in the \rhs of \ex{6.48} are narrow, therefore (see, \eg
\cx{A2}) $\FDLgua$ is also an exponential of a narrow exponent, whose terms
are constructed by gluing the struts of $\gxQLxa$ to \emph{all}
legs of the graphs of
$\WKFDDCrLx$ (possibly, to more than one
graph simultaneously, but in such a way
that a resulting graph is always connected).
\qq
\FDLgua = \exp\lrbcs{\FDLngua}.
\label{6.49}
\qqq
Therefore $\FDLngua$ is a sum of contributions of connected legless graphs and
the conditions on the graphs of $\WKFDDCrLx$ imply that $\echi$ of the graphs
of $\FDLngua$ must be not less than 1. Thus $\FDLngua$ contains only the graphs
of $\cgDL$, and it remains to check the denominators of their contributions.

Since $\ZKDCrLx,\;\umsgx\in\cDCL\subset\cQDCL$, then
$\WKFDDCrLx\in\cDCL\subset\cQDCL$ and the struts of
$\xQiLxak$
are the only
source of denominators in $\FDLngua$. These denominators are
$\dtplijua$ which
appear in the expression for the elements of the inverse matrix
$\mtrik{\lijua}$. Thus, if a root edge $\ed$ of a graph $D$ in the
expression for $\FDLngua$ is a glued strut of
$\xQiLxak$,
then it carries
the denominator $\dtplijue$, and the expression for $\FDLngua$ has a form
\qq
\FDLngua =  \sum_{\xD\in\cgDL
} {\xtFDL\over
\preDr
\dtplijue
}.
\label{6.49*}
\qqq
If we multiply the numerators and denominators of the summands in the \rhs
of this equation by $\preDr \exp\lrbcs{-\xFWzLue}$, then in view of
\ex{6.45x1} they take the form of the summands of \ex{6.45*1} if we set
\qq
\xFDL = \xtFDL \preDr\exp\lrbcs{-\xFWzLue}.
\label{6.49*1}
\qqq
\qed

Since $\WKODCLx$ has the same tree part as $\WKLx$, while their circle parts
are related by \ex{6.45*a1}, then it is easy to work out a formula for
$\ZKODCurLa$ which is similar to \eex{6.45*},\rx{6.45*1}:
%
%
\qq
\ZKODCurLa = \exp\lrbcs{\WKODCurLa}.
\label{6.50}
\qqq
where
\qq
\WKODCurLa & = &   \hlf\sijoxN \xxlij\, a_i a_j\, \blt
\;-\;\hlf\,\log\lrbcs{\tAFLua
}\rcrc
\nonumber
\\
&&\qquad
+ \sDDL
{\xFODL\over
\preDr\tAFLhue
}
\label{6.51}
\qqq
and
\qq
\xFODL\in\cHD.
\label{6.52}
\qqq

\pr{Theorem}{t1.2}
Equation\rx{1.65y3} is derived from (the exponent of) \ex{6.50} with the
help of \ex{6.45}, which relates the $\mfg$-based and universal \urcc
invariants. Indeed, it is easy to verify that
\qq
&\Tgth \exp\lrbc{\hlf\sijoxN \xxlij\, a_i a_j\, \blt} =
e^{\ohb\lkLuva},
\label{6.53}
\\
&
\Tgth \exp\lrbc{\hlf\,\log\tAFLua\rcrc} =
\FNgLua{}.
\label{6.54}
\qqq
Also
\qq
\log\FDLua =
\snoi\;\;\hb^n\!\!\!\! \sDDLn \Tg \lrbc{
{\xFODL\over\preDr\tAFLhue}}.
\label{6.55}
\qqq
The calculation of $\Tg\lrbcs{{\xFODL\over\preDr\tAFLhue}}$ is performed
according to the definition\rx{4.23}. Each assignment of roots $c$ converts
a denominator $\preDr\tAFLhue$ into\\ $\preDr\tAFLv{\scp{\uva}{c(e)}}$.
Since a graph $D\in\cgDL$ has $3\echi(\xD)$ edges, then
$\lrbs{\FNdg}^{3n}$ may serve as a common denominator for all graphs with
$\echi(\xD)=n$. Hence we can present the summand of \ex{6.55} in the form
\qq
\sDDLn \Tg \lrbc{
{\xFODL\over\preDr\tAFLhue}}
=
{\pnQnLuva \over
\lrbs{\FNdg}^{3n} },
\qquad
\pnQnLuva\in\IQuva,
\label{6.56}
\qqq
and this proves \ex{1.65y3}.\qed

\begin{remark}
\rm
The power of the Alexander polynomial in denominators of \ex{1.4} is $2n$
rather than $3n$ as in \ex{6.66}, because $su(2)$ has only two roots.
Therefore, for a graph $D\in\cgDL$ to produce a non-zero contribution to
the $su(2)$-based \urcc invariant, at least one edge at every 3-vertex
must be Cartan. As a result, out of $3\echi{D}$ edges only $2\echi{D}$
edges are root and $\lrbs{\FNdg}^{2n}$ may serve as common denominator.
\end{remark}

\subsection{The rationality conjecture}
\label{6xs.3}
The rational structure of the summands in \ex{1.4} for the $su(2)$ \urcc
invariant suggests that the wheeled universal \urcc invariant also has a
more restrictive form than\rx{6.51}. In order to describe it, we have to
define the graph analogues of the variables $\ut=q^{\ual}$ which appear as
arguments of the polynomials $P_n$ of \ex{1.4}. Therefore, we come back to
the definition of the spaces $\cQHD$ for graphs $\xD\in\gDXp$.

Let $\cgDXp\subset\gDXp$ be the set of `normal' (1,3)-trivalent graphs (we
exclude dots and circles).
This time, for a graph $\xD\in\cgDXp$ we consider a group algebra
\qq
\HCfeobD = \ZHobZt
\label{6.57}
\qqq
instead of the symmetric algebra $\SHCfobD$.
Similarly to \ex{4.5} we define
\qq
\cXepHD = \Invrs{(\HCfeobD)^{\otimes|\Xp|}}{\GD},
\label{6.59}
\qqq
where the action of the symmetry group $\GD$ on
$(\HCfeobD)^{\otimes|\Xp|}$ is defined in an obvious way (it includes a
sign coming from a possible change of cyclic order at 3-vertices).
Then we define the algebra $\cQXepHD$
as an extension of this algebra, which
permits denominators of the form
\qq
\preDr \pexXp,   \qquad
p_\ed(x_1\ldtc x_{|\Xp|})\in\IQ[x_1^{\pm 1}\ldtc x_{|\Xp|}^{\pm 1}],
\label{6.60}
\qqq
and the graph algebra 
\qq
\tcQDCeXp = \boDgDXpp \cQXepHD.
\label{6.59*}
\qqq

As usual, when $\Xp$ is well-known, we may drop it from
our notations.
Also, when working with a link $\cL$, we assume that
$\Xp=\{\ua\}$.

Next, we define the IHX ideal $\tcQDCeXpIHX$ similarly to the definition of
$\tcDCXIHX$ and $\tcQDCXIHX$ and finally we introduce the algebra
\qq
\cQDCeX = \tcQDCeX/\tcQDCeXpIHX.
\label{6.59*1}
\qqq
%
We expect the latter algebra to
appear in the theory of properly defined loop-filtered finite type invariants
of links.

The algebras $\cQDCX$ and $\cQDCeX$ are related.
For an integer number $o$ we define an injection
\qq
\mpcd{\HCfeobD}{\Fso}{\SHCfobD}
\label{6.61}
\qqq
by its action on the edges $\hed\in\HCfeobD$: $\Fso(\hed) = \exp(o\hed)$.
For a set of integers $\yuo=(o_1\ldtc o_{|\Xp|})$ we define an injection
\qq
\mpcd{ \lrbcs{\HCfeobD}^{\otimes|\Xp|} }{\Fsuo}{\SHCfobDtXp},
\qquad
\Fsuo  = \bigotimes_{j=1}^{|\Xp|} F_{o_j}
\label{6.62}
\qqq
and project it to an injection
\qq
\mpcd{\cXepHD}{\Fsuo}{\cXpHD}.
\label{6.63}
\qqq
This injection can be further extended to
\qq
\mpcd{\cQXepHD}{\Fsuo}{\cQHDXp}\qquad\mbox{and}\qquad
\mpcd{\tcQDCeXp}{\Fsuo}{\tcQDCXp}.
\label{6.63*}
\qqq
%
Obviously, $\Fsuo(\tcQDCeXpIHX)\subset\tcQDCXIHX$, so there is a map
\qq
\mpcd{\cQDCeXp}{\Fsuo}{\cQDCXp}.
\label{6.63*1}
\qqq
\begin{conjecture}
\label{cinj}
We conjecture the relation
\qq
\Fsuo(\tcQDCeX)\cap\tcQDCXIHX = \Fsuo(\tcQDCeXpIHX),
\label{6.63*2}
\qqq
which implies that the map\rx{6.63} is an injection.
\end{conjecture}

\begin{conjecture}[rationality]
\label{6c.1}
Let $\cL$ be a link with at least 2 components
in a rational homology sphere such that
$\AFLut\not\equiv 0$. Let $\yuo$ denote the orders of link components as
elements of $H_1(M;\ZZ)$. Then for every graph $D\in\cgDL$ there exists an
element $P_D(\cL)\in\cLepHD$ such that
\qq
\WKODCurLa & = &   \hlf\sijoxN \xxlij\, a_i a_j\, \blt
\;-\;\hlf\,\log\lrbcs{\tAFLua
}\rcrc
\nonumber
\\
&&\qquad
+ \Fsuo\lrbc{\sDDL
{ P_D(\cL)\over
\preDr\AFLv{\uhed}
}}.
\label{6.64}
\qqq
\end{conjecture}
\begin{remark}
\rm
We expect that the same conjecture holds also for $\WKDCurLa$.
\end{remark}

\begin{remark}
\rm
If the map\rx{6.63} is an injection, as implied by Conjecture\rw{cinj}, then
the element
\qq
\sDDL
{P_D(\cL)\over
\preDr\AFLv{\uhed}
}\in\cQDCeL
\label{6.65*}
\qqq
is an invariant of the oriented link $\cL$ (actually, a similar statement for
knots was proved in\cx{GaKr} without proving the injectivity of\rx{6.63}).
\end{remark}


\pr{Conjecture}{3c.1}
Equation\rx{6.66} is derived from \ex{6.64} in exactly the same way as
\ex{1.65y3} is derived from \ex{6.51}.\qed

\noindent
{\bf Acknowledgements.}
I am very thankful to S.~Garoufalidis, R.~Lawrence, D.~Thurston and A.~Vaintrob
for discussing this work. I am especially thankful to D.~Thurston for helping
me to improve the exposition. This work was supported by NSF Grants DMS-0196235
and DMS-0196131.

\appendix

\nsection{A superalgebra $\uoo$ weight system and the Alexander polynomial of
a link}
\label{A1}

The proof of the formula\rx{6.53x} which expresses the Alexander polynomial of
a link in terms of its Kontsevich integral, is too convoluted. It uses the
`rational' expansion of the colored Jones polynomial described in\cx{Ro1}
plus the equivalence of this expansion to the $su(2)$-based \urcc invariant
established in\cx{Ro2}. One would certainly expect that a simple
formula\rx{6.53x} has a short direct proof. Indeed, such proofs exist and they
are based on calculations with the weight system which produces the Alexander
polynomial of a link from its Kontsevich integral. One possible approach is to
describe this weight system in terms of spanning trees of 3-valent graphs. This
calculation is being carried out by G.~Masbaum and A.~Vaintrob\cx{MV}. We
suggest an alternative (and, perhaphs, simpler) approach which is based
on a relation between the Alexander polynomial and the $\Uoo$
Lie supergroup. Kauffman and Saleur\cx{KS} were the first to notice that if the
$U_q \Uoo$ supergroup $R$-matrix is used instead of the $SU_q(2)$ $R$-matrix in
the formula for the Jones polynomial, then one obtains the Alexander
polynomial. Following this idea, A.~Vaintrob\cx{VA} showed that an application
of the $\uoo$ weight system to Kontsevich integral of a link also produces
the Alexander polynomial. Our formula\rx{6.53x} may be thought of as a result
of calculating a fermionic (Berezin) integral in the Reshetikhin formula for
$\uoo$. However, since we are not aware of Kirillov's integral formula for
superalgebra characters, which would have immediately proved \ex{6.53x}, then we
will do the corresponding weight caclculations explicitly.


The superalgebra $\uoo$ has four generators: two bosonic ones $\saE, \saF$
and two fermionic ones $\saPp$, $\saPm$. The non-zero
(super-)commutation relations
are
\qq
[\saF,\saPpm] = \pm \saPpm,\qquad \{\saPp,\saPm\}=\saE,
\label{A.1}
\qqq
where $\{\saPp,\saPm\}=\saPp\saPm + \saPm\saPp$. The non-degenerate Killing
form $\scp{\cdot}{\cdot}$ is symmetric on bosonic generators, but anti-symmetric on
fermionic ones:
\qq
\scp{\saF}{\saE} = \scp{\saE}{\saF} = 1,\qquad \scp{\saPp}{\saPm} =
-\scp{\saPm}{\saPp} = 1.
\label{A.2}
\qqq
Thus the non-zero components of the inverse metric tensor
$\zhi\in \Srs^2 \uoo$ and structure constant tensor
$\zf\in\bigwedge^3\lrbcs{\uoo}^\ast$ are
\qq
&(\zhi)^{EF}=(\zhi)^{FE}=(\zhi)^{-+}=1,\qquad (\zhi)^{+-}=-1,
\label{A.3}
\\
&
\zf_{F+-} = \zf_{F-+}=\zf_{+-F}=\zf_{-+F}=1,\qquad\zf_{+F-}=\zf_{-F+}=-1.
\label{A.4}
\qqq

The definition of the maps $\Tg$ of\rx{1.11*} for a superalgebra $\mfg$
is slightly more complicated than for an ordinary Lie algebra. For a graph
$\xD\in\cAX$ choose a `fermionic' subgraph $\xDF$, which consists of mutually
nonintersecting chains. These chains may be either closed (cycles) or end at
1-vertices of $\xD$. We call such 1-vertices fermionic. Let $\neDf$ be the
number of edges in $\xDF$ (we call them fermionic) and let $\neDb$ be the
number of edges in $\xD\setminus\xDF$ (we call them bosonic).

For a pair $\xD,\xDF$
we define an element $\TgDF\in\UgX$. First of all, we choose an
orientation on all chains of $\xDF$ (ultimately $\TgDF$ will not depend
on it). As a linear space $\mfg$ splits into a direct sum of its
bosonic and fermionic subspaces $\mfg=\mfgb\oplus\mfgf$. Correspondingly, the
inverse Killing form tensor $\himfg$ also splits: $\himfg=\himfgb+\himfgf$,
where $\himfgb\in S^2\mfgb$ and $\himfgf\in\bigwedge^2\mfgf$. Thus we assign
$\himfgb$ to the bosonic edges and $\himfgf$ to the oriented fermionic edges of
$\xD$ in order to contract them with the $\mfg$ structure tensors assignes to
3-vertices, as prescribed by the index contraction map\rx{1.11*1}. This
contraction produces an element
$\xCD\lrbcs{\fmfg^{\otimes\neDt}\otimes(\himfgb)^{\otimes\neDb}
\otimes(\himfgf)^{\otimes\neDf}}\in\UgX$,
which we still have to multiply by
sign factors in order to get $\TgDF$. The first sign factor is $(-1)^{\#\xDF}$,
where $\#\xDF$ is the number of connected components of $\xDF$. When passing
through 3-vertices, fermionic chains induce cyclic order on incident edges. Let
$\nor$ be the number of 3-vertices on fermionic chains on which the induces
order is opposite to that of the original order of $\xD$. Then the second sign
factor is $(-1)^{\nor}$. In order to define the third ordering we pick a linear
order on the set of open fermionic chains. Since these chains are oriented,
this order induces a linear order on the set of fermionic 1-vertices (the final
vertex of an open chain is considered to be immediately following the initial
vertex of that chain). The same set of fermionic 1-vertices has an alternative
linear order coming from the linear ordering of the elements of the set $\xX$
and from the orientation of the segments of $\uparrow_{\xX}$. Let $\norl$ be
the number of pairs of fermionic 1-vertices, on which both orders differ. Then
the third sign factor is $(-1)^\norl$. Thus we define
\qq
&\TgDF = (-1)^{\#\xDF+\nor+\norl}\;
\xCD\lrbcs{\fmfg^{\otimes\neDt}\otimes(\himfgb)^{\otimes\neDb}
\otimes(\himfgf)^{\otimes\neDf}},
\label{A.5}
\\
&\Tg(\xD) = \sum_{\xDF\subset\xD} \TgDF,
\label{A.6}
\qqq
where the sum in the latter formula goes over all possible fermionic subgraphs
of $\xD$.

For an element $x\in\UgX$ consider its representation in a $\mfg$-module
$\TVuva = \TVmvaoX$, where $\Vmva{j}$ are irreducible $\mfg$-modules with
(shifted) highest weights $\val_j$. A supertrace over a linear superspace
$V=\Vbsn\oplus\Vfrm$ is a difference between the traces over its bosonic and
fermionic subspaces: $\STr_V = \Tr_{\Vbsn} - \Tr_{\Vfrm}$. For $1\leq
k\leq\orX$, let $\STrvuak x\in\End(\Vmva{k})$ denote the supertrace of $x$ over
all the spaces $\Vmva{j}$ except the space $\Vmva{k}$. Since $x$ is
$\mfg$-invariant, then so is $\STrvuak x$. Since we assumed that $\Vmva{k}$ is
irreducible, this means that $\STrvuak x = CI$, where $I$ is the identity
operator acting on $\Vmva{k}$, while $C$ is a constant. We will denote this
constant as $\STrsvuak x$. It can be presented as a complete trace over
$\TVuva$: if we choose $\pO\in\End(\Vmva{k})$ such that $\STrvuak \pO=1$, and
denote
$\pO_k = I^{\otimes (k-1)}\otimes
\pO\otimes I^{\otimes(\orX-k-1)}\in\End{\TVuva}$,
then
\qq
\STrsvuak x = \STrvua \pO_k\, x.
\label{A.6*1}
\qqq

Now we turn back to the case of $\mfg=\uoo$.
Consider a family of 2-dimensional $\uoo$-modules $\Vuoa$ parametrized by a
formal parameter $\gra$. Their basis is formed by two vectors: $\vgrap$ (boson)
and $\vgram$ (fermion), the action of the algebra generators being
\qq
&\saE\,\vgrapm = a\,\vgrapm,\qquad \saF\,\vgrapm = \pm {1\over 2}\,\vgrapm,
\nonumber
\\
&\saPm\,\vgrap = \vgram,\qquad \saPm\,\vgram = 0,
\nonumber
\\
&\saPp\,\vgram = a\,\vgrap,\qquad\saPp\,\vgrap =0.
\label{A.7}
\qqq

\begin{theorem}[A.~Vaintrob]
For an $\xN$-component link $\cL\in S^3$
and for a tensor product of
$\uoo$-modules $\TVua=\Vm{a_1}\otimes\cdots\otimes\Vm{a_{\xN}}$,
\qq
\AFLeua =  a_k^{-1} \STrsuak\, \Tuoo\lrbcs{\ZKAL}
\label{A.8}
\qqq
for any $k$, $1\leq k\leq \xN$.
\end{theorem}

If we define $\pO\in\End{\Vuoa}$ by its action on the basis vectors as
\qq
\pO\,\vgrap = \vgrap,\qquad \pO\,\vgram = 0,
\label{A.8*1}
\qqq
then $\STr_{\Vuoa} \pO=1$ and,
according to \ex{A.6*1}, equation\rx{A.8} is equivalent to
\qq
\AFLeua = a_k^{-1} \STrua\lrbs{\pOk\,\Tuoo\lrbcs{\ZKAL}}.
\label{A.8*2}
\qqq
\pr{Theorem}{t6.11}
We will prove the second line of \ex{6.53x} by modifying the \rhs of \ex{A.8*2}.
We will simplify
$\Tuoo\lrbcs{\ZKAL}$ and find a way to calculate the supertrace $\STrua$
in a way, which is similar to Kirillov's integral formula.
Actually, we will work thourgh the
symmetric algebra $\Suoo$, that is, we will describe
$\Tuoo\lrbcs{\ZKBL}\in\Suoo$ and then describe a trick to calculate a
composition of maps
\qq
\begin{CD}
\Suoo @>{\buoo}>>\Uuoo @>{\STrVa}>> \IC,
\end{CD}
\label{A.9}
\qqq
where $\buoo$ is the PBW isomorphism (\cf \ex{1.13}).

Let $\xD$ be a graph of $\cBL$ and $\xDF$ -- one of its possible
fermionic subgraphs. Suppose that $\xD$ has a bosonic edge which is incident to
two 3-vertices. According to \ex{A.3}, the tensor $\zhi_\bsn$ has only
two non-zero matrix elements: $EF$ and $FE$, while the only bosonic index of
the non-zero structure tensor elements\rx{A.4} is $F$. Therefore, the $E$ index
of $\zhi_9\bsn$ can not be matched by structure tensors at 3-vertices, and
$\Tuoo(\xD,\xDF)=0$. Hence, if $\Tuoo(\xD,\xDF)\neq 0$, then every bosonic
edge of $\xD$ must be a leg. This means that there are only three types of
pairs $\xD,\xDF$, $\xD$ being connected, with non-zero $\uoo$ weights: a
fermionic strut, a haircomb with fermionic spine and a wheel with fermionic
cycle and bosonic legs. Let us find the $\uoo$ weights of these graphs. We
present $\SuootN$ as a polynomial algebra of $4\xN$ variables
$\saE_i,\saF_i,\saPp_i,\saPm_i$, $1\leq i\leq\xN$. Now if $\xD$ is a bosonic
strut $\tstv{\bfx_i}{\bfx_j}$, then $\Tuoo(\xD,\xDF)=\saE_i\saF_j +
\saF_i\saE_j$. If $\xD$ is a haircomb graph of \fg{f1.1} with ferminic spine,
then
$\Tuoo(\xD,\xDF)=(\saPm_i\saPp_j-
(-1)^{n_1+\cdots+n_{\xN}}\saPp_i\saPm_j)E_1^{n_1}\cdots
E_{\xN}^{n_{\xN}}$, where $n_i$ ($1\leq i\leq\xN$) are the numbers of bosonic
legs of colors $\cL_i$ attached to the fermionic spine. If $\xD$ is a wheel
with fermionic circle, then $\Tuoo(\xD,\xDF)=-(1-(-1)^{n_1+\cdots+n_{\xN}})
E_1^{n_1}\cdots E_{\xN}^{n_{\xN}}$.

We are going to convert $\Tuoo(\xD)$
into an element of $\Uuoo$ and then evaluate its supertrace in the tensor
product of 2-dimensional representation\rx{A.7}. Since $\saE=a I$ in $V_a$,
then for the purpose of our calculations we can replace $E_i$ by $a_i$ in our
formulas fro $\Tuoo(\xD,\xDF)$. Also, we know that the coefficients at the struts
in $\WKLX$ are linking numbers $\xxlij$ of $\cL$, while the coefficients at its
haircomb graphs form generating series $\lpijua$ of \ex{6.29zz4} and the
coefficients at wheels form the generating series $\xFzLua$ (see the discussion
at the beginning of subsection\rw{6xs.1*}).  Therefore, for the purpose of
evaluating the \rhs of \ex{A.8*2}, we can replace $\Tuoo\lrbcs{\ZKBL}\in\SuootN$
by a rather simple exponential
\qq
\AFLeua = a_k^{-1} e^{-\xFzLua}\, \STrua\pOCk\,\buoo\lrbs{
\exp\lrbc{
\sijoxN\lrbcs{\xxlij a_i \saF_j  - \lpijua\, \saPm_i\saPp_j}
}}.
\label{A.10}
\qqq

Now we describe a useful trick for calculating the composition of maps\rx{A.9}.
This trick works for any Lie (super-)algebra $\mfg$ and for any $\mfg$-module
$V$. There is a canonical isomorphism $\xmap{\xgmg}{\Sg}{\Dcstg}$ between the
symmetric algebra $\Sg$ and the algebra of differential operators on $\mfg$
with constant coefficients $\Dcstg$. Namely, for $\vb\in\mfg\subset\Sg$,
$\xgm(\vb)$ is the Lie derivative along the constant vector field on $\mfg$,
which is equal to $\vb$ at every point. Then it is well-known that for any
$x\in\Sg$,
\qq
\STrV \bmg (x) =  \xgmg(x) \atv{\lrbcs{\STrV e^{\vY}}}{\vY=0},
\label{A.11}
\qqq
where in the \rhs the differential operator $\xgmg(x)$ is acting on
the function $\STrV e^{\vY}$ of $\vY\in\mfg$ and the resulting function is
restricted at $\vY=0$. We will also need a slight generalization of this
equation:
for any $\pO\in\End(V)$,
\qq
\STrV \lrbcs{\pO\,\bmg(x) }
= \xgmg(x) \atv{\lrbcs{\STrV \pO\, e^{\vY}}}{\vY=0}.
\label{A.12}
\qqq

We introduce one bosonic coordinate $\gxf$
and two fermionic coordinates $\etp,\etm$ on $\uoo$
\qq
\vY = \gxf\saF + \etp \saPp + \etm \saPm
\label{A.13}
\qqq
(we did not introduce a coordinate for $\saE$, because all operators $\saE_j$
are already replaced by the corresponding variables $a_j$). Then \eex{A.11}
and\rx{A.12} allow us to rewrite \ex{A.10} as
\qq
\lefteqn{
\AFLeua = a_k^{-1} e^{-\xFzLua}\,
}
\nonumber\\
&&\times
\atv{
\exp\lrbc{
\sijoxN\lrbcs{\xxlij\, a_i\, \dgxfj
  - \lpijua\, \detmi\detpj}
}
\STrVak \lrbcs{\pOC\, e^{\vY_k} }
\prod_{1\leq i\leq \xN\atop i\neq k} \STrVai e^{\vY_i}
}{\vbfY=0}.
%
\label{A.14}
\qqq
Now we calculate the traces by using a QFT trick. We introduce two
(fermionic creation and annihilation) operators $\sapp,\sapm$, which satisfy
the anti-commutation relations
\qq
\{\sapp,\sapp\} = \{\sapm,\sapm\}=0, \qquad \{\sapp,\sapm\}=I
\label{A.15}
\qqq
and act on $\Vuoa$
as
\qq
&
\begin{array}{cc}
\sapp\,\vgrap=0, & \sapp\,\vgram=\vgrap,\\
\sapm\,\vgrap=\vgram,&
\sapm\,\vgram=0.
\end{array}
\label{A.16}
\qqq
Comparing the matrix elements\rx{A.7} and\rx{A.16} we conclude that in $\Vuoa$
\qq
\saPm = \sapm,\qquad \saPp = a\sapp,\qquad \saF = \hlf\,I - \sapm\sapp.
\label{A.17}
\qqq
and
\qq
e^{\vY} = e^{\phi/2}\exp(-\phi\,\sapm\sapp + a\,\etp\sapp + \etm\sapm).
\label{A.18}
\qqq
Let us define new operators
\qq
\saprp = \sapp + {1\over \phi}\,\etm,\qquad
\saprm = \sapm - {a\over \phi}\,\etp,
\label{A.19}
\qqq
so that \ex{A.18} becomes
\qq
\e^{\vY} = e^{\phi/2}\, \eaphe \exp(-\phi\, \saprm\saprp).
\label{A.20}
\qqq
Since the operators $\saprp,\saprm$ satisfy the same anticommutation
relations\rx{A.15} as $\sapp,\sapm$, then they have the same matrix
elements\rx{A.16} in an appropriate basis of $\Vuoa$. Those basis vectors are
eigenvectors of $\exp(-\phi\,\saprm\saprp)$ with eigenvalues $1$ and
$e^{-\phi}$,
so
\qq
\STrVa \exp(-\phi\, \saprm\saprp) = 1 - e^{-\phi}
\label{A.21}
\qqq
and
\qq
\STrVa e^{\vY} = (e^{\phi/2} - e^{-\phi/2})\,\eaphe.
\label{A.22}
\qqq
Since for two fermionic variables $\savpp,\savpm$
\qq
\int \exp(\phi\,\savpp\savpm + a\,\etp\savpp + \etm\savpm)\;
d\savpm d\savpp =  \phi\;\eaphe,
\label{A.22*1}
\qqq
then \ex{A.22} can be rewritten as
\qq
\STrVa e^{\vY} = \shlfph
\int \exp(\phi\,\savpp\savpm + a\,\etp\savpp + \etm\savpm)\;
d\savpm d\savpp.
\label{A.22*2}
\qqq
If we interpret $\savpp,\savpm$ as coordinates on the coadjoint orbit of
$\uoo$, then this equation represents Kirillov's integral formula for the
character of this superalgebra.

Comparing the matrix elements\rx{A.8*1} and\rx{A.16}, we find that in $\Vuoa$
\qq
\pOC = \sapp\sapm  = (\saprp - {1\over \phi}\,\etm)
(\saprm + {a\over \phi}\,\etp),
\label{A.23}
\qqq
so that
\qq
\STrVa \pOC\,e^{\vY}  & = & e^{\phi/2}\,\eaphe \STrVa\lrbs{
(\saprp - {1\over \phi}\,\etm)
(\saprm + {a\over \phi}\,\etp) \exp(-\phi\, \saprm\saprp)
}
\nonumber
\\
& = &
\shlfph - \lrbc{\shlfph - e^{\phi/2}}\,\eaphe,
\label{A.24}
\qqq
%
%
Since
\qq
\int\savpp\savpm\, \exp(\phi\,\savpp\savpm + a\,\etp\savpp + \etm\savpm)\;
d\savpm d\savpp = 1,
\label{A.25}
\qqq
then we can rewrite \ex{A.24} in terms of convenient fermionic integrals
\qq
\STrVa\pOC\,e^{\vY} & = & \shlfph\,  \Bigg(
\int\savpp\savpm\, \exp(\phi\,\savpp\savpm + a\,\etp\savpp + \etm\savpm)\;
d\savpm d\savpp
\nonumber
\\
&&\qquad - \xfB(\phi)
\int\exp(\phi\,\savpp\savpm + a\,\etp\savpp + \etm\savpm)\;
d\savpm d\savpp\Bigg),
\label{A.26}
\qqq
where we defined a shortcut function
\qq
\xfB(x)
= {1\over x}\lrbc{1 - {x\over 1-e^{-x}} }.
\label{A.27}
\qqq

Let us substitute \eex{A.22*2} and\rx{A.26} into \ex{A.14}. Since for any
analytic function $f(\phi)$ and for any constant $c$ we have a Taylor series
formula
\qq
\atv{\exp(c\,\del_\phi)\,f(\phi)}{\phi=0} = f(c)
\label{A.28}
\qqq
and since
\qq
\lefteqn{
\atv{
\exp\lrbc{-\sijoxN \lpijua\, \detmi\detpj } \exp\lrbc{ \sioxN (\etp_i \savppi +
\etm_i \savpmi)}
}{\ueta=0}
}
\nonumber
\\
&&\hspace{3in} =
\exp\lrbc{ \sijoxN \lpijua\,\savppj\,\savpmi
},
\label{A.29}
\qqq
then in view of \ex{6.29zz6}
\qq
\AFLeua & = & a_k^{-1}\, e^{-\xFzLua}\,
\lknLa\,\Bigg[\int\savppk\,\savpmk\,\exp
\lrbc{
\sijoxN a_i^{-1}\,\lijua\,\savppj\,\savpmi
} d\usapm\,d\usapp
\nonumber
\\
&&\qquad
- \xfB\lrbc{\sioxN \xxl_{ik}\, a_i }
\int\exp
\lrbc{
\sijoxN a_i^{-1}\,\lijua\,\savppj\,\savpmi
} d\usapm\,d\usapp
\Bigg],
\label{A.30}
\qqq
where  $d\usapm\,d\usapp = d\savpv{-,1}\,d\savpv{+,1}\cdots
d\savpv{-,\xN}\,d\savpv{+,\xN}$. Since for any matrix $\mtr{A_{ij}}$
\qq
&
\int \exp\lrbc{\sijoxN A_{ij}\,\savppj\,\savpmi} d\usapm\,d\usapp =
\det\mtr{A_{ij}},
\label{A.31}
\\
&
\int \savppk\,\savpmk\,\exp\lrbc{\sijoxN A_{ij}\,\savppj\,\savpmi}
d\usapm\,d\usapp = \det\mtrk{A_{ij}},
\label{A.32}
\qqq
where $\mtrk{A_{ij}}$ is the $(k,k)$-minor of $\mtr{A_{ij}}$, then \ex{A.30}
can be rewritten as
\qq
\AFLeua & = & a_k^{-1}\, e^{-\xFzLua}\,
\lknLa\,
\nonumber
\\
&&\qquad\times
\lrbs{
\det\mtrk{ a_i^{-1}\lijua } - \xfB\lrbc{\sioxN \xxl_{ik}\, a_i }
\det\mtr{ a_i^{-1}\lijua }
}
\label{A.33}
\qqq
Then \ex{6.53x} follows, since
\qq
\det\mtrk{ a_i^{-1}\lijua} =
{\det\mtrk{\lijua} \over \prod_{1\leq j\leq \xN\atop
j\neq k} a_j},\qquad
\det\mtr{ a_i^{-1}\lijua} = {\det\mtr{\lijua} \over \pjoxN a_j}
\label{A.34}
\qqq
and the latter expression is zero in view of \ex{dzr}.\qed

\begin{remark}
\rm
Equation\rx{dzr} can be proved without using the machinery of Cartan and root
edges. Indeed, consider the \rhs of \ex{A.8*2} if we choose $\pO=I$. Since
$\STrVa I=\sdim \Vuoa = 0$, then
\qq
\STrua\lrbs{\Tuoo\lrbcs{\ZKAL}} = 0.
\label{A.35}
\qqq
If we calculate the \lhs of this equation in the same way that we did it
for the \rhs of \ex{A.8*2}, then we find that similarly to \ex{A.33}
\qq
\STrua\lrbs{\Tuoo\lrbcs{\ZKAL}} =e^{-\xFzLua}\,\lknLa\,\det\mtr{a_i^{-1}\lijua}.
\label{A.36}
\qqq
Now a combination of \eex{A.35} and\rx{A.36} indicates that
\qq
\det\mtr{a_i^{-1}\lijua} = 0,
\label{A.37}
\qqq
which is what we need to derive \ex{6.53x} from \ex{A.33}.
\end{remark}

\end{document}

*******************************

The spaces $\cDCX$ and $\cQDCX$ have algebra structure described explicitly
in\cx{RoC}. Let us review its definition. For two graphs $D_1,D_2\in\gDX$ let
$D=D_1\bcup D_2$ be their disconnected union.
Then $\HobD = \HobDo\oplus\HobDt$ and there is a pair of natural algebra
homomorphisms
\qq
\cHDo\otimes\cHDt\longrightarrow \cHD,\quad
\cQHDo\otimes\cQHDt\longrightarrow \cQHD,
\label{4.16}
\qqq
which define the product $x_1 x_2$ of $x_i\in\cHv{D_i}$
or $x_i\in\cQHv{D_i}$ ($i=1,2$).
It is easy to see that the isomorphism\rx{4.11} is the isomorphism of algebras.

**************************

In one of three IHX graphs of $\tcBCXIHXt$ $\Xp$ legs are
also incident to the same 3-vertex, so this graph is $\CC$ and
$\tcBCXIHXt$ relations express the fact that $\Xp$ legs are `commutative':
two graphs which differ only in the order in which their $\Xp$ legs are
attached to their root edges, are equal in the quotient space
$\tcBCX/\tcBCXIHXt$.

There is a natural map of the $\orXp$-th tensor power of $\HCfobD$
into $\cBCX$
\qq
\begin{CD}
\SHCfobDtXp @>\mF>> \cBCX.
\end{CD}
\label{4.4}
\qqq
In order to define it, we associate the elements of the label set $\Xp$ with
$\orXp$ factors of the tensor product $\SHCfobDtXp$. Then any element of
$\SHCfobDtXp$ can be presented (not uniquely) as a polynomial of labeled
oriented root edges $\edx$ of $D$, $x\in\Xp$ being
a color. $\mF$ maps a monomial
$\prod_{\ed\in\bEDr} \prod_{x\in\Xp}\lrbc{\edx}^{\mex}$
into a graph
constructed by attaching $\mex$ Cartan legs with label $x$ for every $x\in\Xp$
on the
left side of every oriented edge $\ed\in\bEDr$, where $\bEDr\subset\bED$ is the
set of root edges of $D$.
The order of legs of different labels
attached to the same edge does not matter in view of
Cartan commutativity lemma\rw{lcc}.
The same lemma for vertices, one of whose
Cartan edges is an $\Xp$ leg, together with the IHX relation, which
involves one Cartan leg, guarantees that the image of $\mF$ does not
depend on the choice of a representation of the element of $\SHCfobD$.

The map\rx{4.4} restricts to
\qq
\begin{CD}
\cXpHD @>\mF>> \cBCX.
\end{CD}
\label{4.4*}
\qqq

\begin{theorem}
\label{t4.2}
The map\rx{4.4*} extends to the isomorphism
\qq
\begin{CD}
\cDCX @>{\mF}>> \cBCX,
\end{CD}
\label{4.11}
\qqq
which preserves the grading $\dgth{}$.
\end{theorem}
\proof
The proof is completeley similar to that of Theorem~3.1 of\cx{RoC} with leg
commutativity lemma~3.3 of\cx{RoC}
replaced by Cartan leg commutativity lemma\rw{lcc}. The
preservation of grading is obvious, since attaching legs to the edges of graphs
of $\cDCX$ via map\rx{4.4} does not change their Euler characteristic.\qed


**********

Let $\xD$ be a (1,3)-valent graph of $\cBCX$ which has no $\Xp$ Cartan
legs. We denote by $\txD$ a graph constructed from $\xD$ by removing all
of its Cartan edges and `dissolving' their incident vertices. We will work
with a relative rational cohomology space $\HCfobDd$ and a tensor power
$\SHCfobDtXpd$ of its symmetric algebra.

*****************************************

%
\exp\lrbs{

********************

Let us also define a similar series $\xFWzLua\in\evn{\IQuua}$ by the formula
\qq
\mFi\circ\fct\lrbcs{\WKOL} = \hxFWzLua\rcrc
+ \mbox{other graphs}.
\label{6.45*a}
\qqq
The difference between the series $\hxFzLua$ and $\hxFWzLua$ comes from the
fact that wheeling of $\WKL$ produces new 1-loop graphs, when struts of $\WKL$
are glued to the wheels of $\Wh$ defined by \ex{1.14*}.

******************************

Let us define a series
$\xFzLua\in\evn{\IQuua}$ by the formula
\qq
\mFi\circ\fct\lrbcs{\WKL} = \hxFzLua\rcrc
+ \mbox{other graphs},
\qqq
in which $\mFi\circ\fct\lrbcs{\WKL}$ is obtained by declaring all legs of
$\WKL$ to be Cartan and converting them into symmetric algebra of graph
cohomology. Alternatively, $\hxFzLua$ is a coefficient at $\rcrc$ in
$\WKDCrLx$.

*****************************

\subsection{Strut exponent and haircomb graphs}
\label{5xs.3}

Since the quadratic exponent\rx{6.24*2x1} plays a prominent role in all our
calculations, we want to give a precise description of their origin in terms of
the graphs of the original function $\gxG(\ubfx)\in\cBCX$. In particular, we
are interested in the case when $\gxG(\ubfx)$ belongs to $\cBX$ and is mapped
into $\cBCX$ by the injection $\ftr$ of\rx{4.2}. Also, we assume that $\gxG$ is
an exponential
\qq
\gxG(\ubfx) = \exp W(\ubfx),\qquad W(\ubfx)\in\cBX,
\label{6.29zz1}
\qqq
where $W(\ubfx)$ is narrow (Kontsevich integral of a link satisfies this
property). Then the strut exponent\rx{6.24*2x1} comes from the haircomb graphs
of $W$. A \emph{haircomb} graph is a graph which consists of a single chain
connecting two 1-valent vertices and legs attached to that chain. The chain is
called a \emph{spine} of the graph and the 1-valent vertices are called
\emph{endpoints} of the haircomb. The same graph may have up to four different
spines depending on the choice of end-points.

A haircomb graph of $W$ may contribute to $Q(\uxz)$ in two different ways. The
first way is to make a substitution $\bfx_i=\hlf\,\tmadv{x_i}^2(a_i)$ at one of
the legs and $\bfx_i=a_i$ at all other legs. Only the struts of $W$ may
contribute in this way (struts are also haircombs). If the strut part of $W$ is
\qq
W\istr(\ubfx) = \hlf \sijoxX \xxlij\;\; \tstv{\bfx_i}{\bfx_j},
\label{6.29zz2}
\qqq
then its first-way contribution to $Q(\uxz)$ is
\qq
\sijoxX \xxlij\; a_i a_j\;\;\rstv{x_i}{x_i}.
\label{6.29zz3}
\qqq
The second way, in which a haircomb graph of $W$ may contribute to $Q(\uxz)$,
is to perform the substitutions $\bfx_i=\tmadv{x_i}(a_i)$ at its endpoints and
the substitutions $\bfx_i=a_i$ at all other 1-vertices.
Let us quantify this contribution.
For a set of
non-negative integers $\bfn$ let $\bfD_{i,j;\bfn}$ be a set of haircomb graphs
with selected spines such that their endpoints have colors $\bfx_i,\bfx_j$ and
they have $n_k$ legs of color $\bfx_k$ attached to the spine for all $k$,
$1\leq k\leq \orX$. Let $c(\xD)$ be the coefficients of these haircomb graphs in
an expression of $W(\ubfx)$ as an element of $\tcBX$ (we fix the signs of
$c(\xD)$ by assuming that legs are attached to the right of the spine, if it is
oriented from $\bfx_j$ to $\bfx_i$). Now we define the power series
\qq
\lpijua = -\sum_{\bfn\geq 0} \lrbc{\sum_{\xD\in\bfD_{i,j;\bfn}} c(\xD)}
a_1^{n_1}\cdots a_{\orX}^{n_{\orX}}.
\label{6.29zz4}
\qqq
Note that $l\p_{ji}(\ua) = l\p_{ij}(-\ua)$. Now the second contribution of
haircomb graphs to $Q(\uxz)$ is
\qq
\sijoxX \lpijua\;\; \rstv{x_i}{x_j}.
\label{6.29zz5}
\qqq
Since a combination of contributions\rx{6.29zz3} and \rx{6.29zz5} is equal to
the \rhs of \ex{6.24*2x1}, then we conclude that
\qq
\lij(\ua) = a_i a_j\lrbc{\lpijua + \delta_{ij}\sum_{k=1}^{\orX} \ell_{ik} }.
\label{6.29zz6}
\qqq
\begin{remark}
\rm
The contribution of diagonal struts $\ell_{ii}\;\tstv{\bfx_i}{\bfx_j}$ of $W$
to both expressions\rx{6.29zz3} and\rx{6.29zz5} are equal up to a sign and thus
cancel each other in the final expression\rx{6.29zz6}. Therefore, diagonal
coefficients $\ell_{ii}$ do not participate in the final
expression\rx{6.29zz6}.
\end{remark}

****************************************

\appendix

\nsection{A superalgebra $u(1|1)$ weight system and the Alexander polynomial of
a link}

The proof of the formula\rx{6.53x} which expresses the Alexander polynomial of
a link, is too convoluted. It uses the `rational' expansion of the colored
Jones polynomial established in\cx{Ro1} plus the equivalence of this expansion
to the $su(2)$-based \urcc invariant established in\cx{Ro2}. One would
certainly expect that a simple formula\rx{6.53x} has a short direct proof. We
are going to sketch a `semi-proof' of this equation based on an observation of
L.~Kauffman and H.~Saleur\cx{KS} that the Alexander polynomial is related to
Lie superalgebra $u(1|1)$ in the same way that the Jones polynomial is related
to Lie algebra $su(2)$. Kauffman and Saleur showed that if a link is
presented as a closure of a braid, then its Alexander polynomial can be
calculated by taking the trace of the product $u(1|1)$-based $R$-matrices which
represent the braid. We are going to apply the $u(1|1)$ weight system to the
Kontsevich integral of the link and show that it yields the \rhs of \ex{6.53x}.

The superalgebra $u(1|1)$ has four generators: two bosonic ones $\saE, \saF$
and two fermionic ones $\saPp$, $\saPm$. The non-trivial
(super-)commutation relations
are
\qq
[\saF,\saPpm] = \pm \saPpm,\qquad \{\saPp,\saPm\}=\saE,
\label{A.1}
\qqq
where $\{\saPp,\saPm\}=\saPp\saPm + \saPm\saPp$. The non-degenerate Killing
form $\scp{\cdot}{\cdot}$
$(\cdot,\cdot)$ is symmetric on bosonic generators, but anti-symmetric on
fermionic ones:

**************************************************************

Then we will show that $\tftri$ and $\tfcti$ can be extended to the maps

 of the maps $\tftr$, $\tfct$ and then

\noindent
\emph{Proof of injectivity of\rx{4.2}.}
Let us define a left-inverse of the injection $\tftr$. The graphs of $\tcBCX$,
each of whose edges is either total or Cartan, form a basis there.
We define a map
$\xmap{\tftri}{\tcBCX}{\tcBX}$ by its action on this basis: if a all edges of a
graph $\xD\in\tcBCX$ are total, then it is mapped to the same graph in $\tcBX$,
while if $\xD$ has at least one Cartan edge, then it is mapped to zero. Thus
\qq
\tftri\circ\tftr = I.
\label{4.2xt1}
\qqq
If we prove that $\tftri$ can be extended to a map $\xmap{\ftri}{\cBCX}{\cBX}$,
then it would follow that $\ftri\circ\ftr=I$, which would imply the injectivity
of $\ftr$. Thus it is sufficient to prove that
$\tftri(\tcBCXIHX),\tftri(\tcBCXCC)\subset\tcBXIHX$. The first inclusion is
obvious. As for the second one, we observe that $\tftri(\tcBCXCCo)=0$, because
the $\CCo$ graphs contain Cartan edges. Almost all $\CCt$ graphs also contain
Cartan edges and therefore are also mapped to 0. The only exception are the
$\CCt$ subgraphs which have no Cartan connecting legs and a single root
connecting leg. This connecting leg is a \emph{bridge}, that is, an edge, whose
removal would increase the number of connected components of the graph.
However, a $\tftri$ image of a graph which has a bridge connecting two
subgraphs, one of which has no legs, is zero.

Let us present this latter as a difference between a total
edge and a Cartan leg.

Obviously, a map
$\ximap{\ftr}{\tcBX}{\tcBCXoC}$ at the level of `unquotiened' spaces is an
injection.
%
%
The graphs of $\tcBCXoC$, each of whose edges is either total or Cartan, form a
basis there.
Let us split this basis into two subsets. The first subset contains the graphs,
all of whose edges are total. Obviously, the span of this subset is equal to
$\ftr(\tcBX)$. The second subset contains basis graphs which have at least one
Cartan edge. We denote their span as $\VC$.
%

The
whole space $\tcBCXoC$ is a direct sum of $\ftr(\tcBX)$ and $\VC$.
%
It is
easy to see that $\tcBCXIHX$ splits into a direct sum of two subspaces:
$\ftr(\tcBXIHX)\subset\ftr(\tcBX)$ and $\VCIHX\subset\VC$. Therefore, since
$\cBX = \tcBX/\tcBXIHX$, then
the following is an injection
\qq
\ximap{\ftr}{\cBX}{(\tcBCXoC)/\tcBCXIHX}.
\label{4.2*}
\qqq
If it were true
that $\tcBCXCC\subset\VC$, then the fact that the map\rx{4.2} is an injection
would follow immediately. This is almost the case, since it is easy to see that
any graph (possibly having root edges), which has at least one Cartain edge,
belongs to $\VC$. Thus there is only one type of $\CCt$ relations, which does
not belong ot $\VC$:
%
this happens when there are no Cartan connecting legs and all the edges
inside the box of \fg{f4.3} are total. However, it is not hard to see that a
subgraph, which has just one connecting
leg (of any type), while all its internal edges
are total, is equal to zero because of the IHX relations, and so it can be
neglected after we took a quotient over $\tcBCXIHX$.\qed

In order to prove the injectivity of\rx{4.2x} we need a
couple of lemmas. Recall that an edge of a graph is called a \emph{bridge}, if
it connects two otherwise disconnected components.

\begin{lemma}
\label{r-bridge}
If all legs of a graph $\xD\in\tcBCX$ are Cartan, and $\xD$ has a root bridge,
then $\xD\in\tcBCXCCt$.
\end{lemma}
\proof
If we cut the bridge edge into two legs thus splitting $\xD$ into two proper
subgraphs, then, obviously, both of these subgraphs are $\CCt$.\qed

\noindent
\emph{Proof of injectivity of\rx{4.2x}.}
We assume that the set $\xX$ contains just one element. We want to be careful
and, contrary to our tradition, use different notations for linear maps
constructed from each other by taking quotients. Thus let
\qq
\xmap{\fcto}{\tcBX}{\tcBCX}
\label{4.2w1}
\qqq
be the original map which ultimately induces\rx{4.2x}.
Consider another map
\qq
\xmap{\tfcto}{\tcBX}{\tcBCX},
\label{4.2z1}
\qqq
which maps a graph of $\tcBX$ into the same graph
in $\tcBCX$, in which all non-bridge edges are total, while all bridges
(including legs) are Cartan. Obviously, $\tfct$ is injective. Let us choose a
basis of $\tcBCX$ which is formed by graphs whose edges satisfy a \emph{basis
condition}: if an edge is not a bridge, then it is either total or Cartan; if an
edge is a bridge, then it is either root or Cartan. Let us also
split $\tcBCX$ into a direct sum of two spaces
\qq
\tcBCX = \tfcto(\tcBX) \oplus V,
\label{4.2z*1}
\qqq
%
where
$V$ is a span of basis graphs, which either have a Cartan non-bridge edge or a
root bridge edge. It is easy to see that the remaining basis vectors of
$\tcBCX$ form a basis of $\tfcto(\tcBX)$.

First, we show that
\qq
\tcBCXCCt\subset V.
\label{4.2w2}
\qqq
Indeed, $\tcBCXCCt$ can be presented as a span of
$\CCt$ graphs all of whose edges, except possibly one root edge, satisfy the
basis condition. If the exceptional root edge of such a $\CCt$ graph (which is
the root edge of \fg{f4.3}) is
a bridge, then the graph by definition belongs to
$V$. If root edge is non-bridge, then present it as a difference between the
total and Cartan edges. Both resulting graphs are in $V$, because at least one
of Cartan edges of \fg{f4.3} has to be non-bridge (otherwise, a proper subgraph
would have exactly one non-bridge leg, which is impossible).

Second, we show that the difference between $\fcto$ and $\tfcto$ belongs to
$\tcBCXCCt$:
\qq
\tfcto(D) - \fcto(D) \in \tcBCXCCt,\qquad D\in\tcBX.
\label{4.2w3}
\qqq
Indeed, the difference between $\tfcto(D)$ and $\fcto(D)$ is that the former
makes bridges Cartan, while the latter makes them total. Since a difference
between a total edge and a Cartan edge is a root edge, then
$\tfcto(D)-\fcto(D)$ can be presented as a sum of graphs, which have at least
one root bridge, while all their legs are Cartan. According to
Lemma\rw{r-bridge}, these graphs belong to $\tcBCXCCt$.

Now consider induced maps
\qq
\xmap{\fctt,\tfctt}{\tcBX}{\tcBCX/\tcBCXCCt}.
\label{4.2u1}
\qqq
In view of\rx{4.2w2}, $\tfctt$ is injective and in view of\rx{4.2w3},
\qq
\fctt=\tfctt.
\label{4.2u2}
\qqq
Thus, $\fctt$ is also injective. It remains to check the effect of
factoring over $\tcBCXCCo$ and $\tcBCXIHX$.

First, we show that $\tcBCXCCo$ splits into a direct sum of its
intersections with $\tfcto(\tcBCX)$ and with $V$:
\qq
\tcBCXCCo = (\tcBCXCCo\bcap\tfcto(\tcBCX))\oplus
(\tcBCXCCo\bcap V).
\label{4.2w4}
\qqq
and
\qq
\tcBCXCCo\bcap\tfcto(\tcBCX)\subset \tfcto(\tcBCXIHX).
\label{4.2w5}
\qqq
Indeed, $\tcBCXCCo$ is a span of basis graphs which have a 3-vertex
incident to 3 Cartan edges. Some of them are in $\tfcto(\tcBCX)$ and some
are in $V$. If a basis $\CCo$ graph is in $\tfcto(\tcBCX)$, then it means
that all of its Cartan edges are bridges. In particular, 3 Cartan edges,
which are incident to the same vertex of \fg{f4.3}, are bridges. Therefore,
the $\tfcto$ preimage of this graph in $\tcBX$
has 3 bridges incident to the same
vertex, hence, in view of Lemma\rw{3-bridge}, it belongs to $\tcBXIHX$.

In view of\rx{4.2z*1} and\rx{4.2w2}, $\tcBCX$ splits
\qq
\tcBCX/\tcBCXCCt = \tfcto(\tcBCX)\oplus V\p,\qquad
V\p = V/(\tcBCXCCt).
\label{4.2w8}
\qqq
Since $\tcBCXCCo$ splits in $\tcBCX$, then, in view of\rw{4.2w2} it splits
also in $\tcBCX/\tcBCXCCt$
\qq
\lefteqn{
\tcBCXCCo/(\tcBCXCCo\bcap\tcBCXCCt)
}
\nonumber
\\
& = &
\tcBCXCCo\bcap\tfcto(\tcBCX) \oplus
\lrbcs{
\tcBCXCCo/(\tcBCXCCo\bcap\tcBCXCCt)
}
\bcap V\p,
\label{4.2w6}
\qqq

The IHX subspace $\tcBCXIHX$ also splits in the quotient space
$\tcBCX/\tcBCXCCt$
\qq
\lefteqn{
\tcBCXIHX/(\tcBCXIHX\bcap\tcBCXCCt)
}
\nonumber
\\
& = &
(\tcBCXIHX\bcap\tfcto(\tcBCX))
\oplus
\lrbcs{
\tcBCXIHX/(\tcBCXIHX\bcap\tcBCXCCt)
}
\bcap V\p
\label{4.2w7}
\qqq
and
\qq
\tcBCXIHX\bcap\tfctt(\tcBCX)\subset\fctt(\tcBXIHX).
\label{4.2w8*}
\qqq
Indeed, if the intermediate edge in all 3 graphs of the IHX triplet of
\fg{f4.2} are not bridges, then since those edges are total, the graphs
belong simultaneously either to $\tfcto(\tcBX)$ or to $V$ (assuming that
other edges satisfy the basis condition). If an intermediate edge is a
bridge, then we can declare it Cartan, because the difference between a
total and a Cartan edge is a root edge and, according to
Lemma\rw{r-bridge}, a graph with
a root bridge and Cartan legs belongs to $\tcBCXCCt$, which we have
factored out. It is easy to verify that if intermediate bridge edges in
IHX triplets are declared Cartan, then again all 3 graphs belong either to
$\tfcto(\tcBX)$ or to $V$ simultaneously. Relation\rx{4.2w8*} follows by
factoring an easily verifiable relation
\qq
\fcto(\tcBX)\bcap\tcBCXIHX = \fcto(\tcBXIHX)
\label{4.2w9}
\qqq
combined with \ex{4.2u2}.

Since both $\tcBCXCCo$ and $\tcBCXIHX$ split with respect to the
splitting\rx{4.2w8}, since their intersections with $\tfctt(\tcBX)$ are
subspaces of $\tfctt(\tcBXIHX)$ and in view of \ex{4.2u2} we conclude that
the map\rx{4.2x}, constructed by factoring the domain of the map\rx{4.2u1}
over $\tcBXIHX$ and its range over $\tcBCXCCo+\tcBCXIHX$, is injective.\qed

***************************************************************

 An individual polynomial $\LLbmn$ comes from the graphs
$\xD$ of $\WKOL$ which have $\echi(D)=n-1$ and $m$ legs. In particular,
this means that if we define
$\xFzLua\in\evn{\IQuua}$
by the formula
\qq
\mFi\circ\fct\lrbcs{\WKL} = \hxFzLua\rcrc
+ \mbox{other graphs},
\label{6.45*x}
\qqq
in which $\mFi\circ\fct\lrbcs{\WKL}$ is obtained by declaring all legs of
$\WKL$ to be Cartan and converting them into symmetric algebra of graph
cohomology,
then
\qq
\smzi \LLamo = \svldp \xFzLav,
\label{6.45*y}
\qqq
(indeed, since wheeling does not affect the coefficient at the
graph $\rcrc$, then
$\WKL$ can be replaced by $\WKOL$
in \ex{6.45*x}).

At this point we can prove a theorem that we neglected to formulate in\cx{Ro2}:
\begin{theorem}
\label{t6.11}
The Alexander polynomial $\AFLut$ is determined by the
haircomb tree and `1-loop' parts
of the (logarithm of) Kontsevich integral
$\WKOL$:
\qq
\AFLeua & = & \pjbai \dtplijua\; \exp\lrbcs{- \xFzLua }.
\label{6.53x}
\qqq
\end{theorem}
\proof
Consider \ex{6.45*y} in case of $\mfg=su(2)$. Since $su(2)$ has only one
positive root, then (in terms of the $su(2)$ $\mfhs$ coordinate
$a=\scp{\va}{\vl}$)
\qq
\smzi \LLv{m,0}{\ua} = \xFzLua.
\label{6.53y}
\qqq
Then \ex{6.53x} follows easily from the combination of
equations\rx{1.65y2},\rx{6.34x7} and\rx{6.45*y}.\qed

Note that neither $\dtplijua$, nor $\xFzLua$ is the invariant of an `unmarked'
link $\cL$. Only their combination\rx{6.53x} is independent of the `Morse
marking'.

*******************************************

\begin{remark}
\label{rq3}
\rm
We
will be interested in the integrals of functions $\frua$, so following
Remark\rw{rq2} we will always assume that the integral is calculated in
$\cQDCXp$ (\eg, if the integrand comes from $\cBCX$, then it is mapped by
$\mFiq$ into $\cQDCXp$), while variables $\ua$ belong to $\Xp$. We may keep
parameters $\ua$ in the integrand
as a reminder, but they are actually not parameters, but elements of $\Xp$,
their legs being hidden in the symmetric algebras of
cohomologies of graphs.
\end{remark}

***********************************

We also denote
\qq
\umsgx = \pjoxN \msjv{x_j}.
\label{6.24*5}
\qqq
Theorem\rw{t6.5} has an obvious
\begin{corollary}
\label{c6.x1}
The integration measure $\umsgx$ is $\ux$-substantial.
\end{corollary}

Consider graph calculus in $\cBX$ (that is, all edges are total).
If a function $f(x)$ depends only on a single coordinate $x$
(and does not depend on
any other parameters), then
\qq
\ndv{\mad_y} f(x) = 0.
\label{5.8**2}
\qqq

****************************************************

Suppose that a
map $\bfy=F(\bfx)$ is such that $\Fc=0$. If
\qq F(\bfx) =
\Esv{\Hr(\bfx),\,\Hc(\bfx)}
\label{6.13}
\qqq
for two other maps: with
a root output $y=\Hr(\bfx)$ and with a Cartan output $a=\Hc(\bfx)$, then
\qq
F(\bfx) = \Hr(\bfx)\qquad\mbox{and}\qquad \Hc(\bfx)=0.
\label{6.14}
\qqq

*************************

If for an element $x\in\cQHD$ and for a root edge $\ed$ of $D$ the polynomial
\qq
\pexXp\Big|_{\hedx=0\;\;{\rm for\;\; all}\;\; x\in\Xp} \neq 0
\label{4.13*}
\qqq
in all denominators\rx{4.11*} of $x$, then $x$
is non-singular at $\ed$.

**************************

It is easy to see that the weight system\rx{4.3*} commutes with a linear map
\qq
\xmap{\ftr}{\cBX}{\cBCX} , \qquad \xX=\xXt,
\label{4.2}
\qqq
which maps each graph of $\cBX$ to the same graph in $\cBCX$, all of its edges
being total.

\subsection{Injections}
\label{3xs.1*}

Let us consider two linear maps
\qq
&\xmap{\ftr}{\cBX}{\cBCX} , \qquad \xX=\xXt,
\label{4.2}
\\
&\xmap{\fct}{\cBX}{\cBCX}, \qquad \xX=\xXc,
\label{4.2x}
\qqq
which map each graph $\xD$ of $\cBX$ into the same graph of $\cBCX$ such that
all edges of $\ftr(xD)$ are total, while all internal edges of $\fct(\xD)$ are
total and all legs are Cartan. It is easy to see that the weight
system\rx{4.3*} commutes with $\ftr$.

In addition to\rx{4.2}, there is another interesting map
\qq
\xmap{\fct}{\cBX}{\cBCX}, \qquad \xX=\xXc,
\label{4.2x}
\qqq
which maps each graph of $\cBX$ to the same graph of $\cBCX$, all of its
internal edges being total, while all of its legs are
declared to be Cartan.

***********************************

There is an obvious map
\qq
\xmap{\ftr}{\cBX}{\cBCX} , \qquad \xX=\xXt,
\label{4.2}
\qqq
which maps each graph of $\cBX$ to the same graph in $\cBCX$, all of its edges
being total.
There exists a map
%
\qq
\xmap{\Tg}{\cBCX}{\oSgXT}
\label{4.3*}
\qqq
similar to that of\rx{1.11*} which
commutes with\rx{4.2}. We just have to modify\rx{1.11*2}.
Projectors $\prjr$ and $\prjh$ project naturally $\stg$ onto $\str$ and
$\sth$. Thus, we construct an element $\hiD\in\stgx$ by taking a tensor
product of $\himfr = \prjr(\himfg)$ for every root edge and $\himfh =
\prjh(\himfg)$ for every Cartan edge and then we modify the
definition\rx{1.11*2}
\qq
\fmfg^{\otimes\neDt}\otimes\hiD
\mathop{{\longmapsto}}^{\xCD} \tTg(D)\in \tmfc
\subset\tmfgx,
\label{4.3}
\qqq
where $\nehDo$ and $\nerDo$ are the numbers of Cartan and root legs of $D$.
The symmetrization over 1-valent
vertices which have the same label projects $\tTg(D)$ to
$\Tg(D)\in\oSgXT$.

For special circle
graphs we define $\Tg(\rcrc) = \dim\mfr$ and $\Tg(\lcrc)=\dim\mfh$.
%


***********************************

It is easy to see that
$\tcBCXCCt\subset\tcBCX$ does not intersect with $\tfct(\tcBCX)$. Indeed, if a
graph $\xD$ is in the image of $\tfct$, then all its Cartan edges are bridge.
Suppose that $\xD$ contains a proper $\CCt$ subgraph. Since its Cartan
edges are bridges, then so must be the root edge, but all bridges must be
Cartan. This is a contradiction. Since $\tcBCXCCt$ does not intersect with
$\tfct(\tcBCX)$, then
\qq
\xmap{\tfct}{\tcBX}{\tcBCX/\tcBCXCCt}
\label{4.2z2}
\qqq
is injective. This map coincides with the map
\qq
\xmap{\fct}{\tcBX}{\tcBCX/\tcBCXCCt}.
\label{4.2z3}
\qqq
Indeed, the only difference between these maps is whether internal bridges are
declared Cartan or total. A difference between a total edge and Cartan edge is
a root edge. However, if all legs of a graph in $\tcBCX$ are Cartan, then a
graph with a root edge is $\CCt$ (breaking this edge into two legs splits the
graph into two proper $\CCt$ subgraphs) and the equality between the
maps\rx{4.2z2} and\rx{4.2z3} follows. Now we can use these maps
interchangeably.

Obviously, the map\rx{4.2z3} maps $\tcBXIHX\subset\tcBX$ in
$\tcBCXIHX\subset\tcBC$. Therefore the maps\rx{4.2z2} and\rx{4.2z3} descend to
a pair of equal injections
\qq
\xmap{\fct,\tfct}{\cBX}{\tcBCX/(\tcBCXCCt\oplus\tcBCXIHX)}.
\label{4.2z4}
\qqq

Finally, the subspace $\tcBCXCCo/(\tcBCXCCo\bcap(\tcBCXCCt\oplus\tcBCXIHX))$
does not intersect with the image of $\tfct$ of\rx{4.2z4}. Indeed, if a graph
has a 3-vertex which is incident to three Cartan

This map maps $\tcBXIHX$ into
$\tcBCXIHX\subset (\tcBCXIHX\bcap\tcBCXCCt)$.
Indeed, the combination of graphs \fg{f4.2} in $\tcBX$ maps into a similar
combination

Now we prove that the map\rx{4.2x} is an injection, if the set $\xX$
contains just one element.

An edge of a graph is called a \emph{bridge}, if its removal splits a connected
component to which it belonged into two disconnected pieces. We consider the
legs as bridges. Consider a map
$\xmap{\tfct}{\tcBX}{\tcBCX}$ which maps a graph of $\tcBX$ into the same graph
in $\tcBCX$ in which all non-bridge edges are total, while all bridges
(including legs) are Cartan. Obviously, $\tfct$ is injection.

*************************

where
\qq
V\p = V/(\tcBCXCCt),
\label{4.2w7}
\qqq
so that in view of\rx{4.2z*1} and\rx{4.2w2},
\qq
\tcBCX/\tcBCXCCt = \tfcto(\tcBCX)\oplus V\p.
\label{4.2w8}
\qqq

*******************************

We can consider a case of a $t$ dependent preexponential factor and a
$t$-dependent quadratic form separately.

If $\gxQ(x)$ does not depend on $t$, while $\gxP(x;t)$ does depend on it,
then \ex{5.48} follows easily from the definition\rx{5.45}. Now suppose
that $\gxQ(x;t)$ is $t$-dependent, while $\gxP(x;t)$ is $t$-independent. In
this case \ex{5.48} requires an inclusion of the 1-loop determinant
in the definition\rx{5.45}.

Let us calculate explicitly both sides of \ex{5.48}. First,  the \rhs
\qq
\lefteqn{
\dlt\lrbs{ \lrbcs{ \det\hgxQ }^{-1/2} \lrbc{ P \glxxLR \eQi } }
}
\label{5.49}
\\
&&
\hspace{-0.5in}=
-\hlf \lrbcs{ \det\hgxQ }^{-1/2} \Bigg[ \Tr (\hgxQi \dlt \hgxQ)
\lrbc{ P \glxxLR \eQi }
+
\lrbs{
\gxP \glxxLR \lrbcs{ (\dlt \gxQi)\, \eQi } }\Bigg],
\nonumber
\qqq
where $\dlt\gxQi$ can be determined with the help of the formula
\qq
\widehat{\dlt \gxQi} = - \hgxQi\times\dlt\hgxQ\times\hgxQi.
\label{5.50}
\qqq
Next, we differentiate the \lhs of \ex{5.45}
\qq
\dlt \lrbc{ \FGi d\ux \; \gxG(x;t) } & = &
\hlf \FGi d\ux\; \eQ\,\gxP\,\dlt\gxQ
\nonumber\\
& = &
\hlf\,
\lrbcs{ \det\hgxQ }^{-1/2}
\lrbc{ (P\,\dlt\gxQ) \glxxLR \eQi }.
\label{5.51}
\qqq
Consider the gluing in the last expression. The legs of $\dlt\gxQ$ can
either be glued to the same strut of $\eQi$ or to the different struts.
In the first case the first term of the \rhs of
\ex{5.49} is reproduced. In the second case, the other two legs of
the two struts of $\eQi$ which are glued to the strut $\dlt\gxQ$ will be
glued to the legs of $P$ and in view of \ex{5.50} this reproduces exactly
the second term of the \rhs of \ex{5.49}.              \qed

**************************************************

we consider a deformation of the
picture of $\Tr(\Mrx^{\times n})$. As a result of this transformation we get
the same graph, except that all $x$ legs now face inwards. In view of the AS
relation, turning the legs outwards amounts to replacing $x$ with $-x$.
Therefore
\qq
\Tr(\Mrx^{\times n}) = \Tr(\Mrmx^{\times n})
\label{6.16*6x}
\qqq
and hence $\msgx=\msjv{-x}$. \qed

*************************

By using the results of\cx{Ro2} we will show in Section\rw{s6} that
%
\qq
{\det \gQLuak \over
\FNdg
} = 1 + \cO(\uva) \in\IQuva,
\label{1.66}
\qqq
where
%
Relation\rx{1.66} means that we can replace the denominators $\det\gQLuak$
by the functions $\FNdg$ which have an advantage of being
invariants of $\cL$ (that is, they do not depend on presenting $\cL$ as
a dotted Morse link) at a cost of modifying the formal series $\pQnLuva$.

Finally, we will prove the following

************************

A relation\rx{6.45} between $\mfg$-based and universal \urcc invariants
dictates a relation between the parts of \eex{1.62} and\rx{6.51} which appear
with similar powers of $\hb$. In view of \ex{6.34x6*}, we find that

\qq
\LLrarkzo = \svldp \xFzLav
\label{6.53xxx}
\qqq
and
\qq
\log\FDLua =
\snoi\;\;\hb^n\!\!\!\! \sDDLn \Tg \lrbc{
{\xtFODL\over
\preDr \dtplijue
}
}.
\label{6.54x}
\qqq

A combination of Lemma\rw{l6.6} with $su(2)$ results of\cx{Ro2} leads to the
following two theorems.

\begin{lemma}
All graphs $D$ of $\FDLngua$ have no legs and $\echi(\xD)\geq 1$.
\end{lemma}
The absence of legs follows from the definition of gluing $\glxxkLRb$, which
must glue all available legs. The fact that $\echi{\xD}\geq 1$ follows
from the fact that every graph of $\WKFDDCrLx$ either has at least 3 legs
or its $\echi$ is non-negative.\qed

*******************

\begin{theorem}
\label{t1.2}
If $\FNgLua{}{\not\equiv} 0$, then
The $\mfg$-based \urcc invariant
defined by \eex{1.62} and\rx{1.64} is well-defined and
does
not depend on a presentation of $\cL$ as a dotted Morse link and on a
choice of $\xk$ ($1\leq \xk\leq \xN$) thus being a topological invariant
of $\cL$. It can be presented in a form
\qq
\IrhLuva  =
e^{\ohb\lkLuva}\, \hbpLg\, \dguva\,
\lrbcs{\FNgLua{}}^{-1}  \,
\exp\lrbc{\snoi {\pnQnLuva \over\FNgLua{3n} }\,\hb^n
},
\label{1.68}
\qqq
where
\qq
\pnQnLuva\in\IQuva
\label{1.69}
\qqq
are link invariants.
\end{theorem}

%
\lrbcs{ \detp\hgxQ(\cL;\ua) }

%

**********************************************

If both $\Dr$
and $\Dc$ contain proper $\CC$ subgraphs,
then we define $\cHD$ to be 0-dimensional. If only $\Dr$ contains proper
$\CC$ subgraphs,
then we define $\cHD$ to be the same as $\cHDc$.
If $\Dr$ does not contain proper $\CC$ subgraphs, then we define
$\cHD$ to be isomorphic to $\cHDr$ and we define the inclusion
of $\cHD$ into $\cHDr\oplus\cHDc$ as
\qq
\cHD\hookrightarrow\cHDr\oplus\cHDc,
\quad
x\longmapsto x + \fed(x),
\label{4.9}
\qqq
where $\fed$ is the map\rx{4.8}.

*******************************************************

Note that the intermediate
edge there is total.
The subspace $\tcBCXCC$ is a span of two
`Cartan-commutator' (or just CC) relations depicted in \fg{f4.3}. The first
relation $\CCo$ says that Cartan edges `commute' among themselves.
We call a subgraph $D\p$ of $D$
\emph{proper} if all vertices of $D\p$ are either 1-valent or 3-valent.
The
second relation $\CCt$
says that if a graph has a proper non-strut subgraph,
all of whose legs are Cartan except
for one leg which is root, then the original graph is zero.

*****************************************************

So far, in defining the universal \urcc invariant  we have tacitly assumed that
the strut part of the integrands $\ZKDCrLx$ and $\ZKODCrLx$ is non-degenerate.
Let us look at it more closely. The strut graphs of $\WKDCrLx$ come from
the haircomb graphs of $\WKL$ (see \fg{f1.1}). Since wheeling $\WhhL$ only
increases the Euler characteristic of graphs, then the strut part of
$\WKODCrLx$ is exactly the same as that of $\WKDCrLx$. Let us denote it as
\qq
\sijoxN  \xlij(\ua)\;\;x_i\rstrt x_j.
\label{6.46}
\qqq
A contribution of a haircomb graph which is not a strut, to the series
$\xlij(\ua)$ is at least of order $a^3$, so the dominant quadratic part of
$\xlij(\ua)$ comes exclusively from the struts of $\WKL$. The strut part
of $\WKL$ is
\qq
\sijoxN \xxlij\;\; x_i\tstrt x_j,
\label{6.47}
\qqq
where $\xxlij$ are linking numbers of $\cL$. Thus
\qq
\xlij(\ua) = \xlijt(\ua) + \cO(a^3),\qquad
\xlijt(\ua) =
\left\{
\displaystyle
\begin{array}{cl}
\displaystyle
- \xxlij\, a_i a_j & \mbox{if $i\neq j$}
\vspace{3mm}
\\
\displaystyle
\sum\limits_{1\leq k\leq \xN\atop k\neq i} \xxlij\, a_i a_k &
\mbox{if $i=j$.}
\end{array}
\right.
\label{6.48}
\qqq
Note that although the matrix $\mtr{\xlij(\ua)}$ depends on the marking of
$\cL$, the leading quadratic part $\mtr{\xlijt(\ua)}$ is a topological
invariant.

For a matrix $\mtr{A_{ij}}$ we denote by $\mtrrv{A_{ij}}{k}$ its minor
constructed by removing $k$-th row and $k$-th column. The quadratic form
governing the integral\rx{6.35} is $\mtrrv{\xlij(\ua)}{k}$. According to
\ex{6.48},
\qq
\det\mtrrv{\xlij(\ua)}{k} = \det\mtrrv{\xlijt(\ua)}{k} +
\cO(a^{2\xN-1}),
\label{6.49}
\qqq
while the first term in the \rhs is of order $\cO(a^{2\xN-2})$. Therefore,
if $\mtrrv{\xlijt(\ua)}{k}$ is non-degenerate, then the integrals\rx{6.35}
and\rx{6.42} are well-defined. The requirement
\qq
\det\mtrrv{\xlijt(\ua)}{k} \not\equiv 0
\label{6.50}
\qqq
is a simple topological condition on $\cL$. Let a \emph{linking graph} of
$\cL$ be a graph whose vertices are link components and two vertices are
adjacent iff the corresponding link components have a non-zero linking
number. A link is called \emph{algebraically connected} if its linking
graph is connected. It is easy to check (see \eg\cx{Ro2}) that a
condition\rx{6.50} is satisfied for all algebraically connected links.

Still, we can do much better if we use the full results of\cx{Ro2} where
we established an equivalence between the $su(2)$-based \urcc invariant
and a similar invariant derived from the $su(2)$ $R$-matrix formula for
the colored Jones polynomial. Namely, we showed in\cx{Ro2} that
\qq
{\det\mtrrv{\xlij(\ua)}{k} \over \AFLv{e^{\ua}}}
= 1 + \cO(a).
\label{6.51}
\qqq
Therefore we come to
\begin{theorem}
\label{t6.10}
The integrals\rx{6.35} and\rx{6.42} are well-defined iff
$\AFLut\not\equiv 0$.
\end{theorem}

The results of\cx{Ro2} also allow us to provide an explicit expression
for the Alexander polynomial $\AFLut$ in terms of Kontsevich integral
$\ZKBL$. Namely, let us take the 1-loop part of $\WKL$ and declare all
their legs to be Cartan. Only the wheel graphs survive this procedure and
the result can be presented as a formal power series
$\WK_0(\ua)\;\; \bigcirc$. Then
\qq
\AFLv{e^{\ua}} = \det\mtrrv{\xlij(\ua)}{k} \;\exp\lrbcs{ -\WK_0(\ua) }.
\label{6.52}
\qqq

Finally, a combination of \eex{6.45} and\rx{1.4} leads us to
\begin{conjecture}[rationality of Kontsevich integral of links]
Let $\cL$ be a link with a non-degenerate Alexander polynomial. Then there
exist

\end{conjecture}

******************************************************************

\qq
\lefteqn{
\dlt\lrbs{ \lrbcs{ \det\hgxQ }^{-1/2} \lrbc{ P \glxxLR \eQi } }
}
\nonumber
\\
&&\hspace{-0.3in} =
- \hlf \lrbcs{ \det\hgxQ }^{-1/2}
\Tr (\hgxQi \dlt \hgxQ)
 - \hlf \lrbcs{ \det\hgxQ }^{-1/2}
\lrbs{
\gxP \glxxLR \lrbcs{ (\dlt \gxQi)\, \eQi } },
\label{5.49x}
\qqq

Theorems\rw{6.6}--\ref{6.8} establish $\ZKDCrLx$ and $\ZKODCrLx$ as topological
invariants of links in $S^3$. The next theorem is a result of an easy
comparison of calculations of Section{s3} and this Section.
\begin{theorem}
\label{t6.9x}
$\ZKODCrLx$ determines the $\mfg$-based \urcc invariant for any simple Lie
algebra $\mfg$
\qq
\IrfvaNL = \Tg\, \ZKODCurLa
\label{6.45x}
\qqq
\end{theorem}

, where
\qq
\dgh{\xD} = \# (\mbox{edges}) - \# (\mbox{3-vertices})
= \hlf\; \#(\mbox{all vertices})
= \#(\mbox{1-vertices}) + \echi(\xD),
\label{1.11}
\qqq
where $\echi(D) = \#(\mbox{edges}) - \#(\mbox{all vertices})$ is the
(\emph{opposite})
Euler characteristic of $D$.

From the point of view of a graph algebra, the vector field $\xXi(\bfy)$
can also be interpreted as a function $\Gfz(\bfy)$ which depends on a
parameter $\bfz$. Then \ex{6.25} says that
\qq
\FGi dx\; \tGfz(\DRax;\bfu)\,\msgx = 0,
\label{6.30}
\qqq
where $\tGfz(\bfy;\bfu)$ is defined according to \ex{6.24*1}. The graphs of
$\tGfz(\bfy;\bfu)$ have a single $\bfz$ leg coming from $\Gfz(\bfy)$ and a
single $\bfu$ leg coming from $\mad_{\bfu}$. It is easy to see that $\tG(\bfy)$
of \ex{6.28} can be constructed from $\tGfz$ by gluing the $\bfz$ and $\bfu$
legs together:
\qq
\tG(\bfy) = \gllv{\bfz}{\bfu}{1}\lrbcs{ \tGfz(\bfy;\bfu) }.
\label{6.31}
\qqq
Since gluing $\gllv{\bfz}{\bfu}{1}$ commutes with the stationary phase
integration of \ex{6.30}, then we get \ex{6.29} by applying
$\gllv{\bfz}{\bfu}{1}$ to both sides of \ex{6.30}.\qed

Let $\xXi(\bfx)$ be a total valued vector field. Consider a vector field
$\mmad_{\bfy}\times\xXi(\bfx)$ and denote by $\dvadv{\bfx}\xXi(\bfx)$ its
divergence evaluated at $\bfy=\bfx$
(that is, after calculating the divergence, we replace the color
$\bfy$ with $\bfx$; the name $\dvadv{}$ is a combination of $\dvrgv{}$
and $\ad$)
\qq
\dvadv{\bfy} \xXi(\bfy) =
\atv{
\dvrgv{\bfx} \lrbcs{\mmad_{\bfy}\times\xXi(\bfx)}
}{\bfy=\bfx}.
\label{6.28x}
\qqq

It is easy to see that $\dvrgv{\bfx}\mad_{\bfy}(\bfx)=0$ (the resulting graphs
resulting graph is zero because of the AS relation), and as a result, the
definition\rx{6.28} could be modified:

Suppose that in the vector field $\mad_{\bfy}{\bfu}$, $\bfy$ is interpreted as
a coordinate and $\bfu$ is interpreted as a parameter, so that the `vector' leg
carries the color $\del_{\bfy}$. Then
it is not hard to see that $\dvrgv{\bfy}\mad_{\bfy}(\bfu)=0$ (the
resulting graph is zero because of the AS relation), and as a result, the
definition\rx{6.28} could be modified:
\qq
\dvadv{\bfy} \xXi(\bfy) = \dvrgv{\bfy} \lrbcs{\mmad_{\bfy}\times\xXi(\bfy)}.
\label{6.32}
\qqq
\end{remark}

We introduce the following conventions for the grpahs. If we place a linear
combination of colors at a leg, then this denotes a linear combination of
graphs with individual colors and corresponding coefficients. Also, if we place
a root (or Cartan) color at a total leg, then this means that this leg is
declared root (or Cartan). In order to distinguish between root, Cartan and
total colors we will denote root colors as $x,y,\ldots$, Cartan colors as
$a,b,\ldots$ and total colors as $\bfx, \bfy,\ldots$. For any map
$\bfy = F(\cdot)$ we define two related maps $y=\Fr(\cdot)$ and
$a=\Fc(\cdot)$ in which the total $\bfy$ leg of $F$ is replaced by a root
$y$ leg or by a Cartan $a$ leg. Also, if we replace a total color $\bfx$
in the argument of any object $A(\bfx)$ by either $y$ or $a$, then this
means that the $\bfx$ legs of $A$ are declared root or Cartan.

We are going to integrate over the variables $\ux$. By definition, a gaussian
(or stationary phase, or just \SP) integrand is an integration measure of the
form
\qq
\gxG(x) = \eQx \gxP(x),
\label{5.42z}
\qqq
where
\qq
\gxQ(x) = \sijoxN \lji\; x_j\rstrt x_i,
\label{5.43z}
\qqq
and the \SP\ preexponential factor $\gxP(x)$ has a polynomial dependence on the
struts $\ux\rstrt\ux$. It is easy to check that the expression\rx{5.42} can be
an integration measure, that is, if a substitution $x=F(y)$ is performed on it
according to the rule\rx{5.36}, then it still retains the form\rx{5.42}.

If the operator
\qq
\hgxQ = \sijoxN \lji\; \dlxj \rstrt dx_i
\label{5.44z}
\qqq
is non-degenerate,

\begin{theorem}
For a matrix field $M$, $\det( \Tg M) = \Tg (\det M)$.
\end{theorem}
\proof
Let us first proove this for a strut matrix field\rx{5.17}. Let us assume for
simplicity that all coordinates are root (a general case is easy to consider,
since struts do not mix root and cartan coordinates, so that a general strut
matrix field has a block-diagonal form\rx{5.17*}). A weight system $\Tg$
applied to a
root strut matrix field
\qq
M = \sijoxN
\lji\;
\rstv{\dlxj}{dx_i}
\label{5.35*}
\qqq
takes value in $\QShoXp\otimes\End\lrbcs{\bigoplus_{j=1}^{\xN}\mfr_j}$, where
$\mfr_j$ are $\xN$ copies of the space $\mfr$. In fact, the weight system of a
matrix
field\rx{5.35*} is diagonal with respect to the root spaces $\mfrvl$, so
\qq
\Tg M = \svlDg \TgMvl,\qquad
\TgMvl = \mtr{\lij(\vl)} \in \QShoXp\otimes
\End\lrbc{
\bigoplus_{j=1}^{\xN} \mfrvlj }.
\label{5.35*1}
\qqq
%
As a result,
\qq
\det(\Tg M) = \pvlDg \det\TgMvl = \pvlDg \det\mtr{\lij(\vl)}.
\label{5.35*2}
\qqq
A combination of the definition of determinant\rx{5.18} and the
formula\rx{4.23*} for the application of $\Tg$ to $\cQHv{\rcrc}$ yield the same
expression in view of \ex{5.19} applied to $\mtr{\lij(\vl)}$. Thus we proved
the theorem for strut matrix fields.

A general definition of a determinant is based on formulas\rx{5.22}
and\rx{5.29}. The application of weight system $\Tg$ converts \ex{5.22} into
\ex{5.19} and converts \ex{5.29} into a multiplicativity of the ordinary
determinant. Therefore $\Tg$ `commutes' with these defining equations, and
hence it commutes with graph determinant.\qed

Note that
\qq
\LLbetmz = \LLbetmo = 0,\qquad\mbox{if $m<2$}.
\label{1.27*}
\qqq
Indeed, the polynomials $\LLbetmz$ come from the tree graphs ($\echi(D)=-1$) and
tree graphs have at least two 1-valent vertices, so $\LLbetmz=0$ if $m< 2$.
The polynomials $\LLbetmo$ come from the graphs with $\echi(D)=0$. Kontsevich
integral $\WKL$ contains no graphs without 1-valent vertices. Such graphs
appear after the wheeling, but it is not hard to see that for them
$\echi(D)\geq 2$, so none of them contributes to $\LLbetmo$. Also, any graph
with a single 1-vertex is zero in $\cB$, hence\rx{1.27}.

\qq
\begin{array}{ccccc}
&& \cAL &&
\vspace{0.1in}
\\
& \mapne{\xchi} && \mapse{\Tg} &
\vspace{0.1in}
\\
\cBL & \mapright{\Tg} & \SgL & \mapright{\bmg} & \UgL
\vspace{0.1in}
\\
\mapdown{\WhhL}& &\mapdown{\dflL} & & \mapdown{\TrVua}
\vspace{0.1in}
\\
\cBL & \mapright{\Tg} & \SgL &
\mapright{\iOuval}
& \IC
\end{array}
\label{1.17}
\qqq


\begin{remark}
\rm
\label{r6.1}
Suppose that in the vector field $\mad_{\bfy}{\bfu}$, $\bfy$ is interpreted as
a coordinate and $\bfu$ is interpreted as a parameter, so that the `vector' leg
carries the color $\del_{\bfy}$. Then
it is not hard to see that $\dvrgv{\bfy}\mad_{\bfy}(\bfu)=0$ (the
resulting graph is zero because of the AS relation), and as a result, the
definition\rx{6.28} could be modified:
\qq
\tG(\bfy) = \dvrgv{\bfy} \lrbcs{\mad_{\bfy}(G(\bfy))}.
\label{6.32}
\qqq
\end{remark}

Let $\bfz = G(\bfy)$ be a substitution such that $\bfz = G\circ\DRax$ is a
\SP\ substitution. Consider a total vector field
$\mad_{\bfw}(G(\bfy))$ and denote by $\tG(\bfy)$ its divergence evaluated at
$\bfw = \bfy$ (that is, after calculating the diverence, we replace the color
$\bfw$ with $\bfy$)
\qq
\tG(\bfy) =
\atv{
\dvrgv{\bfy} \lrbcs{\mad_{\bfw}(G(\bfy))}
}{\bfw=\bfy}
\label{6.28}
\qqq
Then
\qq
\FGi dx\; \tG(\DRax)\,\msgx = 0.
\label{6.29}
\qqq

Let $\cLSt$ be an $\xN$-component link. We denoted by $\ZKBL$ and $\ZKOBL$
its Kontsevich integral and wheeled Kontsevich integral in $\cBL$. As
explained in Section\rw{s4}, the algebra $\cBL$ can be injected into
$\cBCL$, a graph of $\cBL$ becoming the same graph in $\cBCL$, all of
whose edges are total. Let us denote as $\ZKBCLx$ and $\ZKOBCLx$ the
corresponding Kontsevich integrals. According to the definitions of
Section\rw{s5}, they are functions of $\xN$ total colors $\ubfx$.

The algebra $\cBCL$ is isomorphic to the algebra $\cDCL$ in which Cartan
legs of all $\xN$ colors are converted into cohomology of the
graphs. We denote the corresponding Kontsevich integrals as $\ZKDCLx$ and
$\ZKODCLx$. We will use the same notation for Kontsevich integrals in a
bigger algebra $\cQDCL$. Note that since all edges of the graphs of
$\ZKBCLx$ and $\ZKOBCLx$ are total, then so far there are no Cartan legs,
and the corresponding elements of algebrals $\cHD$ in the expressions for
Kontsevich integrals are of zero degree.


$(ST^2S)_{0 \lambda}$

where 0 is the identity representation, and $\lambda$
is an integrable representation of level k.

Equivalently, there should be a nice way to sum up

\qq
\sum_{\lambda}  \chi_{\lambda}(h_\mu)  exp[4 \pi i c_2(\lambda)/(k+h)]
\chi_{\lambda}(h_\mu')
\qqq
where  $\chi_\lambda$  is the character in the representation
$\lambda$   of the finite dimensional group,
$h_{\mu} = \exp[ 2\pi i (\mu+\rho)/(k+h)],$
and the sum is over the integrable representations of level $k$.


\begin{theorem}
\label{t6.4}
Let $\bfz = G(\bfy)$ be a substitution such that $\bfz = G\circ\DRax$ is a
\SP\ substitution. Consider a total vector field
$\mad_{\bfy}(G(\bfy))$ and denote by $\tG(\bfy)$ its divergence:
\qq
\tG(\bfy) = \dvrgv{\bfy} \lrbcs{\mad_{\bfy}(G(\bfy))}
\label{6.25x}
\qqq
Then
\qq
\FGi dx\; \tG(\DRax)\,\msgx = 0.
\label{6.26x}
\qqq
\end{theorem}
\proof

We follow the combinatorial proof of Theorem\rw{t1.1} and establish an
analog of \ex{1.41}. Define
\qq
\tFLgy(x) = \Scv{\Erv{y}{x}}{\Ecv{y}{x}}
\label{6.7}
\qqq
(\cf \ex{6.3}). Then, according to the definition of $\SPL$,
\qq
\Ev{y}{x} = \Ev{\FLgy(x)}{\tFLgy(x)}.
\label{6.8}
\qqq
Now consider a map from $(x,y_1,y_2)$ to $\bfz$
\qq
\bfz = \Esvb{\,y_2,\,y_1,\,x,\,-\tFLv{y_1}(x) -
\tFLv{\tFLv{y_2}}\circ\FLv{y_1}(x)\,}.
\label{6.9}
\qqq
We will transform the \rhs by applying \ex{6.8} to pairs of adjacent root
coordinates in two different orders. First, we begin by applying \ex{6.8}
to $(y_1,x)$
\qq
\bfz & = &\Esvb{\,y_2,\,\FLv{y_1}(x),\,\tFLv{y_1}(x),\,
-\tFLv{y_1}(x) -
\tFLv{y_2}\circ \tFLv{y_1}(x)\, }
\nonumber\\
& = &
\Esvb{\,y_2,\,\FLv{y_1}(x),\,-\tFLv{y_2}\circ\FLv{y_1}(x)\, }
\nonumber\\
& = & \Esvb{\,\FLv{y_2}\circ\FLv{y_1}(x),\,
\tFLv{y_2}\circ\FLv{y_1}(x),\,
-\tFLv{y_2}\circ\FLv{y_1}(x)\,}
\nonumber\\
& = & \FLv{y_2}\circ\FLv{y_1}(x).
\label{6.10}
\qqq
On the other hand, we may first apply \ex{6.8} to the pair $(y_2,y_1)$
\qq
\bfz & = & \Esvb{\,\FLv{y_2}(y_1),\,\tFLv{y_2}(y_1),\,x,
\,-\tFLv{y_1}(x) - \tFLv{y_2}\circ\FLv{y_1}(x)\, }
\nonumber\\
& = &
\Esvb{\,\FLv{y_2}(y_1),\,\tFLv{y_2}(y_1),\,x,\,-\tFLv{y_2}(y_1),
\,\tFLv{y_2}(y_1)
- \tFLv{y_1}(x)-\tFLv{y_2}\circ\FLv{y_1}(x)\, }
\nonumber\\
& = &
\Esvb{\,\FLv{y_2}(y_1),\mAdesv{\tFLv{y_2}(y_1)}(x),\,
\tFLv{y_2}(y_1)
- \tFLv{y_1}(x)-\tFLv{y_2}\circ\FLv{y_1}(x)\, }
\nonumber\\
& = &
\BCH\Big(\FLv{\FLv{y_2}(y_1)}\circ\mAdesv{\tFLv{y_2}(y_1)}(x),\,
\tFLv{\FLv{y_2}(y_1)}\circ\mAdesv{\tFLv{y_2}(y_1)}(x)+
\tFLv{y_2}(y_1)
\nonumber\\
&&\hspace{3in}
- \tFLv{y_1}(x)-\tFLv{y_2}\circ\FLv{y_1}(x)\,\Big) .
\label{6.11}
\qqq
Thus we find that
\qq
\FLv{y_2}\circ\FLv{y_1}(x) & = &
\BCH\Big(\FLv{\FLv{y_2}(y_1)}\circ\mAdesv{\tFLv{y_2}(y_1)}(x),\,
\tFLv{\FLv{y_2}(y_1)}\circ\mAdesv{\tFLv{y_2}(y_1)}(x)
\nonumber\\
&&\hspace{2in}
+\tFLv{y_2}(y_1)- \tFLv{y_1}(x)-\tFLv{y_2}\circ\FLv{y_1}(x)\,\Big) .
\label{6.12}
\qqq

According to this lemma, \ex{6.12} implies the equality between maps
\qq
\FLv{y_2}\circ\FLv{y_1} =
\FLv{\FLv{y_2}(y_1)}\circ\mAdesv{\tFLv{y_2}(y_1)}.
\label{6.15}
\qqq
which map $x$ to $z$ and depend on two root parameters $y_1$ and $y_2$.

Consider the action of the $\Ad$ map defined in Example\rw{e5.1} on a root
coordinate: $\bfy = \mAda(x)$. Since parameter $a$ is Cartan, then in view
of condition $\CCo$,
$(\mAda(x))_\mfh=0$, so $\mAda(x)$
is reduced to $y = \mAda(x)$ and can be presented as $e^a\;y\rstrt x$,
while the maps $\mAdta$ form a 1-parametric
subgroup generated by the (purely root) vector field $\mad_a$.

Now we need a lemma which claims that
due to the presence of struts
in\rx{6.1*3}, the map $\BCH$ can not mix root and Cartan input in order to
produce a purely root output. \begin{lemma} \label{l6.1} Suppose that a
map $\bfy=F(\bfx)$ is such that $F_\mfh=0$. If
\qq F(\bfx) =
\Esv{\Hr(\bfx),\,\Hc(\bfx)}
\label{6.13}
\qqq
for two other maps: with
a root output $y=\Hr(\bfx)$ and with a Cartan output $a=\Hc(\bfx)$, then
\qq
F(\bfx) = \Hr(\bfx)\qquad\mbox{and}\qquad \Hc(\bfx)=0.
\label{6.14}
\qqq
\end{lemma}
\proof
Consider a graph of $\Hc$ with the smallest number of edges. Then due to
the struts in\rx{6.1*3} it will produce the graph with the smallest number
of edges among the graphs in \rhs of \ex{6.13} with Cartan output. Hence,
this graph can not be cancelled to conform with the fact that the \lhs of
\ex{6.13} has only root output.\qed

Let us define yet another map $z = \tFLga(x)$, where $a$ is a Cartan parameter,
\qq
\tFLga = e^a\;\;z\rstrt x.
\label{6.5}
\qqq
and we identified the color $a$ with the variable assiciated with the edge as
an element of $\HobD$, $D=z \rstrt x$.
\begin{theorem}
\label{t6.1}
The measure\rx{6.4} is invariant under the substitution\rx{6.5}.
\end{theorem}
\proof
It is easy to see that the substitutions $\tFLgta$, $t$ being a \Qpr, form a
1-parametric group, generated by the vector field
$\xia = a\;\; \dlz\rstrt dx$. Therefore, it is sufficient to
check that
\qq
\Ldxia \msj = \ndxia \msj + \msj\, \dvxxia = 0.
\label{6.6}
\qqq
In fact, both terms in the middle expression are equal to 0. Indeed,
$\ndxia$ acts by attaching an $a$ Cartan leg to an $x$ leg in all possible
ways. Since all the legs of the graphs of $j$ are colored by $x$, then
this amounts to a `link relation' which still works on graphs with Cartan
and root edges if the leg which is being attached is Cartan. Also,
$\dvxxia=0$, because the symmetry of the graph sends an odd coefficient
$a$ to zero. \qed

In other
words, we have to find a map from a pair of coordinates $(x,a)$, $x$ being
root and $b$ being Cartan into a similar pair
\qq
(y,b) = \Sv{x}{a}
\label{6.2}
\qqq
satisfies the property
\qq
\tEv{ \Srv{x}{a} }{ \Scv{x}{a} } = \Id,
\label{6.3}
\qqq
where $\Sr$ and $\Sc$ denote $y$ and $b$ parts of the whole map\rx{6.2}
and $(x,a) = \tEv{y}{b}$ is the same map as $\BCH$ except that it splits
$\Ev{y}{b}$ into a root and a Cartan part.

Relation\rx{6.3} means that the map $\SPL$ is the inverse of a map
$(x,a) = \tEv{y}{b}$ which is the same as $\BCH$ except that it splits
$\Ev{y}{b}$ into a root and Cartan parts.

We leave it for the reader to check that vector fields form the Lie algebra of
$\DiffxN$ and that the Lie bracket of vector fields is given by the familiar
Lie derivative formula
\qq
[\xi_1,\xi_2] = \ndv{\xi_1} \xi_2 - \ndv{\xi_2} \xi_1.
\label{5.8}
\qqq

collection of \emph{matrix
elements} $\yxMij(\ux)$ ($\oijxNi$), $\yxMij(\ux)$ being a graph all of whose
legs are coordinates $\ux$ and parameters except for two legs, one carrying the
color $dx_j$ and the other carrying the color $\dlyi$. We will use an
abbreviated notation $\Mij$ for $\xxMij$. A derivative of a map $\uy=\uF(\ux)$
is a matrix field $\yxFpij(\ux)$, where a matrix element $\yxMij$ is represented
by a sum of graphs constructed in all possible ways from the graph of $\uF_i$
by replacing the color $y_i$ by $\dlyi$ and one of the colors $x_j$ by $dx_j$.
Now, the action of a diffeomorphism $\uy=\uF(\ux)$ on a vector field $\xi(\ux)$
consists of an index contraction which comes in the form of leg gluing
\qq
\yxFpij(\ux) \gllsv{dx_j}{\dlxj}{1} \xi(\ux)
\qqq

A \emph{function} $F$ is a graph whose legs are colored by the variables
$x,y,\ldots$ A \emph{1-form} $\omega$ is a graph all of whose legs are
colored by variables except for a single leg, which is colored by a
differential $dx$, $dy$,\ldots A \emph{partial differential} $\del_x F$ of
a function $F$ is a gradient field, which is a sum over all ways of
replacing a color $x$ at one leg of $F$ by a differential $dx$. A
\emph{differential} of a function $dF$ is a sum of its partial
differentials over all variables.

In fact, we define a partial differential $\del_x A$ of any object $A$
presented by a graph: $\del_x A$ is a sum over the ways of changing a
color at one leg of $A$ from $x$ to $dx$.

A \emph{vector field} $\xi$ is a graph all of whose legs are colored by
the variables except for one leg which is colored by a derivative
$\del_x$, $\del_y$,\ldots We define a pairing $\prcox$ between a
1-form $\omega$ and a vector field $\xi$ in an obvious way: if $\omega$
and $\xi$ are each represented by a single graph, then $\prcox$ is a graph
constructed by gluing the differential leg of $\omega$ with the derivative
leg of $\xi$ if the variables match, and $\prcox=0$ otherwise. If $\omega$
and $\xi$ are sums of graphs, then the pairing is defined by its
bilinearity.

If $D$ does not contain $\CC$ subgraphs, then we consider a space
$\SHCfobD$.

It is a
(symmetric) algebra and since $\SHCfobDtXp$ is a tensor product of algebras,
then it also has an algebra structure. We denote by $\QeSHCfobDtX$ an extended
algebra which includes fractions, whose denominators are of the form
\qq
\preDr \pexXp,
\label{4.11*}
\qqq
where $\pexXp\in\SHCfobDtXp$
are polynomials of $(\hedx)_{x\in\Xp}$ (that is, each polynomial depends on a
particular root edge of all possible colors). Note that all individual
polynomials $\pexXp$ are non-zero in $\SHCfobDtXp$.

Suppose that two graphs have one
1-valent vertex marked. Then gluing them means joining and `dissolving' these
vertices. We will also call this procedure `a gluing of legs'.
Let $D_1$ and $D_2$ be two graphs of $\cBX$ such that $D_1$ has at least $m$
vertices of color $x_1$ and $D_2$ has at least $m$ vertices of color $x_2$.
Then the gluing $D_1 \gllv{x_1}{x_2}{m} D_2$ is the sum of the graphs
constructed by choosing $m$ 1-vertices $x_1$ of $D_1$ and $x_2$ of $D_2$ and
gluing them together pairwise (the sum goes over all the ways to choose the
legs and glue them pairwise).

Let us describe the gluing of root edges in the space $\cQDCX$. Suppose that the
graphs $D_1,D_2\in\gDX$ have at least $m$ root legs of colors $x_1$ and $x_2$.
Let $\Dpj$ $j=1,\ldots$
be the graphs constructed by gluing the legs in all
possible ways. A map
\qq
\begin{CD}
\cQHDo\otimes\cQHDt @>\gllv{x_1}{x_2}{m}>> \bigoplus_{j}\cQHv{\Dpj}.
\label{4.17x}
\end{CD}
\qqq
is a sum of individual maps
\qq
\begin{CD}
\cQHDo\otimes\cQHDt @>\gllv{x_1}{x_2}{m|j}>> \cQHv{\Dpj}.
\label{4.18x}
\end{CD}
\qqq
These maps are defined in the following way. If a graph $\Dpj$ has a
$\CCo$ or $\CCt$ subgraph, then the space $\cQHv{\Dpj}$ is
zero-dimensional and we define the image of\rx{4.18} as zero. If $\Dpj$ does
not contain CC subgraphs, then we consider
a natural map
\qq
\Hhobv{\Dpj}\longrightarrow \Hhobv{D_1}\oplus\Hhobv{D_2}
\label{4.19x}
\qqq
which simply cuts the cycles at gluing points. Then we extend the dual map
\qq
\Hobv{D_1}\oplus\Hobv{D_2}\longrightarrow\Hobv{\Dpj}
\label{4.20x}
\qqq
to the algebra homomorphism\rx{4.18}.

Now we define $\tcQDCXIHX$ by copying the definition of $\tcQDCXIHX$. For $D$
being a graph with a 4-valent vertex, we consider all elements $x\in\cQHD$
such that $f_i(x)\in\cQHDev{i}$, $i=1,2,3$, where $\ed$ are the `extra' edges in
the graphs $D_i$ which come from the resolution of the 4-valent vertex of $D$.

$\hed\in\HCobD$ (in view of Remark\rw{r4.1*1} this means that if $\ed$ is
declared Cartan, then $D$ contains a $\CCt$ subgraph) or if
$x$ can be presented as a sum of fractions whose denominators do not
belong to the kernel of the unsymmetrized version
\qq
\begin{CD}
\SHCfobDtXp   @>\fed>> \SHCfobDptXp
\end{CD}
\label{4.13*1x}
\qqq
of the map $\fed$ of\rx{4.8}.

, note that any element of $\SHCfobDtX$ can be presented
(not uniquely) as a polynomial of oriented root edges $e_j$ of $D$. $\mF$ maps
a monomial $\prod e_j^{m_j}$ into a graph constructed from $D$ by attaching
$m_j$ Cartan legs to every oriented edge $e_j$ on its left side. The CC
relation for vertices one of whose Cartan edges is a leg, guarantees that the
image of $\mF$ does not depend on the choice of a representation of the element
of $\SHCfobD$.

The space $\cBCX$ is produced by factoring $\tcBCX$ by three relations.
The first relation is the familiar IHX relation of \fg{f4.2}. Note that the
intermediate edge is total. We called the other two relations
`Cartan-commutator' or simply CC.
Relation
CC1 says that Cartan edges `commute' among themselves. Relation CC2 says that
if a graph has a subgraph, all of whose connecting edges are Cartan except
for one edge which is root, then the original graph is zero.

\qq
\cBC = \tcBC/(\IHX+\CC).
\label{4.1}
\qqq

Let $D$ be a (1,3)-valent graph of $\cBCX$ all of whose edges are root.
Thinking of $D$ as a cell complex, we consider a rational
relative cohomology space $\HobD$. The oriented egdes of $D$ represent integral
elements of $\HobD$ (a pairing between an edge and a cycle is a coefficient at
that egde in a representation of the cycle as a linear combination of egdes).
There is a natural map of a symmetric algebra $\SHobD$ into $\cBCX$
\qq
\begin{CD}
\SHobD @>\mF>> \cBCX
\end{CD}
\label{4.4}
\qqq
which is defined by its action on a product $\prod e_j^{m_j}$ of oriented
edges $e_j$ of $D$. For every edge $e_j$ we attach $m_j$ Cartan legs on its
left side. This way a product of oriented edges of $D$ becomes a graph of
$\cBCX$.

Let $\GD$ be the symmetry group of $D$:
$\GD$ maps edges to edges, 3-vertices to 3-vertices and 1-vertices to
1-vertices of the same color. The elements of $\GD$ may change orientations
at 3-valent vertices. The natural action of $\GD$ on $\HobD$ can be extended to
a group-like action on the symmetric algebra $\SHobD$. We modify this action by
multiplying the action of $g\in\GD$ by a sign factor $(-1)^{|g|}$, where $|g|$
is the number of 3-vertices of $D$ whose cyclic order was changed by $g$.
We need the $\GD$-invariant part of $\SHobD$
\qq
\cHD = \Invrs{\SHobD}{\GD}.
\label{4.5}
\qqq

Now suppose that $D$ is a (1,3)-valent graph of $\cBCX$ such that it has no
Cartan legs.
Define
$\dslD$ as a graph obtained from $D$ by `dissolving' its internal Cartan edges:
we remove each Cartan edge of $D$ and join together the root edges at its
incident 3-valent vertices so that these vertices disappear. The group $\GD$
still acts on $\dslD$, $\Hobv{\dslD}$ and $\SHobv{\dslD}$, so
we define
\qq
\cHD = \Invrs{\SHobv{\dslD}}{\GD}
\label{4.6}
\qqq

Note that the polynomials $\LLbmz$ come from the tree graphs ($\echi(D)=-1$).
Tree graphs have at least two 1-valent vertices, so $\LLbmz=0$ if $m< 2$.
The polynomials $\LLbmo$ come from the graphs with $\echi(D)=0$. Kontsevich
integral $\WKBL$ contains no graphs without 1-valent vertices. Such graphs
appear after the wheeling, but it is not hard to see that for them
$\echi(D)\geq 2$, so none of them contributes to $\LLbmo$. Since any graph with
a single 1-vertex is zero in $\cB$, then $\LLmbo=0$ if $m<2$.

\begin{lemma}
If $D$ is a tree graph and
$\vbom\in\mfh$, then $\TgDbm=0$ unless $D$ consists of a single egde and two
1-valent vertices (such graphs are called `struts').
\end{lemma}
\proof
If $D$ is not a strut, then it contains a 3-valent vertex connected by edges to
two 1-valent vertices. This means that the expression of $\TgDbm$ contains a
commutator of the elements $\vbj$ attached to those 1-valent vertices. Then
$\TgDbm=0$, because all elements of $\mfh$ commute.\qed

The quantum invariants of knots, links and 3-manifolds, such as the
Jones polynominal and the Witten-Reshetikhin-Turaev invariant, were
discovered about 10 years a However, their interpretation in terms
of classical 3-dimensional topology still remains a mystery.

Let us compare the skein relation definition of the Jones polynomial
to that of a much older Alexander-Conway polynomial. The
single-variable Alexander-Conway polynomial $\APLt\in\ZZti$ is a
unique invariant of links in $S^3$ which satisfies the following two
properties. First, the normalization condition:
\qq
\APtbas{\unknot} = 1.
\label{1.1}
\qqq
Second, if $\cLp$, $\cLm$ and $\cLz$ are three links whose regular
projection on a plane is the same except at one spot (see \fg{f1}),
\begin{figure}[hbt]
\leavevmode \centerline{\epsfbox{pct1.eps}} \caption{The links
$\cLp$, $\cLm$ and $\cLz$}
\label{f1}
\end{figure}
then
\qq
\APtbas{\cLp} - \APtbas{\cLm} = (t^{1/2}-t^{-1/2})\,\APtbas{\cLz}.
\label{1.2}
\qqq
\def\pio{ \pi_1 }
This definition is purely combinatorial and it is a bit unnatural
from the 3-dimensional point of view, since it requires a projection
of a link. However, there exist alternative definitions of the
Alexander-Conway polynomial of a knot $\cK$ which are purely topological.
One derives $\APKt$ from the structure of the knot group
$\pi_1(S^3\setminus \cK)$, and the variable $t$ represents the action of the
homology $\pio/[\pio,\pio]$ onto the quotient
$[\pio,\pio]/[[\pio,\pio],[\pio,\pio]]$, where $\pio$ is the group of the knot
($\pio = \pio(S^3\setminus \cK)$). The other definition
relates the Alexander polynomial to the
Reidemeister torsion of a local system in the knot
complement, the variable $t$ being the twist acquired by that system
along the meridian of $\cK$. From both definitions of $\APKt$ it is
clear that $t$ is intimately related to the meridian of $\cK$.

The Jones polynomial of links $\JtL\in\ZZqhi$ can also be defined by
skein relations. It is the unique invariant which satisfies the
following two properties: the normalization condition
\qq
\Jt{\unknot} = q^{1/2} + q^{-1/2}
\label{1.3}
\qqq
and the skein relation
\qq
q\,\Jt{\cLp} - q^{-1}\,\Jt{\cLm} = (q^{1/2} - q^{-1/2})\,\Jt{\cLz},
\label{1.4}
\qqq
where the links $\cLp$, $\cLm$ and $\cLz$ are the same as those in
\ex{1.2}. Despite an obvious similarity between \eex{1.1},\rx{1.2}
and \eex{1.3},\rx{1.4}, there exists no intrepretation of $\JtL$ in
terms of the ``classical'' objects of 3-dimensional topology, such as
the fundamental group of the knot complement. In particular, there is
no indication that the variable $q$ has any connection to the
meridian of $\cK$.

A new hope for a topological interpretation of $\JtL$ emerged when
J.~Birman, X.-S.~Lin and D.~Bar-Natan discovered that both the
Alexander-Conway and Jones polynomials are packed with Vassiliev
invariants. Consider the expansions
\qq
&
\APKt = \snzi a_n(\cK)\,(t-1)^n,
\label{1.5}\\
&
\JtK = \snzi b_n(\cK)\,(q-1)^n.
\label{1.6}
\qqq
It is not hard to see from the skein relations\rx{1.2} and\rx{1.4}
that the coefficients $\a_n(\cK)$ and $\b_n(\cK)$ are Vassiliev
invariants of degree $n$. However, Vassiliev invariants by definition
are related to the topology of ``the space of all maps
$S^1\longrightarrow S^3$, rather than to the topology of knots
themselves. The latter relation is still missing, although some bits
of it are known, such as the relation between the tree Vassiliev
invariants and Milnor's linking numbers (see \cx{HM} and references therein).

By looking at \ex{1.5} we may say that the Alexander-Conway
polynomial presents a way of assembling some Vassiliev invariants of
knots into a polynomial which has a clear interpretation in terms of
the classical 3-dimensional topology. At the same time, the Jones
polynomial assembles some other Vassiliev invariants into another
polynomial whose topological origin is rather obscure. Therefore one
may wonder if there is a way of reassembling all Vassiliev invariants
into the polynomials which would be similar to the Alexander-Conway
polynomial rather than to the Jones polynomial in terms of their
topological interpretation.

In Sections\rw{gs} and\rw{s3} we present
an algorithm of assembling Vassiliev invariants coming from the Kontsevich
integral of
a knot into a sequence of functions of a variable $t$.
In Section\rw{s4}
we conjecture
that these functions are rational: their denominators are powers of
the Alexander-Conway polynomial of $t$ while their numerators are new
polynomial invariants of knots. Since these new polynomials depend on
the same variable $t$, we expect them to have a topological
interpretation in which, similarly to the case of the
Alexander-Conway polynomial, $t$ will also be related to the meridian
of a knot.

Since the first version of this paper was written and reported, Andrew Kricker
has proved the rationality conjecture in his paper\cx{Kr}.

Kontsevich integral is related to the colored Jones (HOMFLY) polynomial of the
knot through the application of a Lie algebra weight system. In Section\rw{s5}
we explain how to apply this weight system to the `repackaged' Vassiliev
invariants. Then we show how the rational structure of Kontsevich integral
appears as a rational the Jones polynomial. In Section\rw{s6} we use these
results to extract the first non-trivial knot polynomial related to the
$\theta$-graph from the expansion of the $SU(3)$ colored Jones polynomial. In
the Appendix we present a table of these `2-loop' polynomials for knots with up
to 7 crossings.

\nsection{Graph spaces}
\label{gs}

We are going to define an algebra $\cD$ based on 3-valent graphs, but first let
us recall the definition of the algebra $\cB$ of (1,3)-valent graphs related
to Vassiliev invariants of a knot.
Each 3-valent vertex of a graph
is endowed with a cyclic ordering of 3 egdes attached to it. When we draw a
picture of a graph, we assume that this ordering is counterclockwise.

A graph $D$ has 2 degrees. They are defined as
\qq
\dego (D) & = & \#\mbox{1-vertices},
\label{1.7*}
\\
\degt (D) & = & \#\mbox{chords} - \#\mbox{3-vertices} =
\chi(D) + \dego (D),
\label{1.7}
\qqq
where $\chi(D)$ is the Euler characteristic of $D$
(more precisely, $\chi(D)$ denotes the Euler characteristic with the
\emph{opposite} sign).

Let $\tcB_{m,n}$ be a
formal vector space (over $\IC$) whose basis elements are in a
one-to-one correspondence with (1,3)-valent graphs of degrees $m$ and $n$
\qq
\tcB_{m,n} = \span{D\,|\,\dego(D)=m,\,\degt(D)=n}.
\qqq
Together all such spaces form a graded space $\tcB$
\qq
\tcB = \bopnzi \tcB_n, \qquad\mbox{where}\;\;
\tcB_n =
\bopmzi\tcB_{m,n}.
\label{1.8}
\qqq
The space $\tcB$ has two important subspaces:
$\tcBas$ and $\tcBihx$.
$\tcBas$ is spanned by the sums $D_1+D_2$, where $D_1$ and $D_2$ are the same
graphs except that they have different cyclic orders
at one 3-valent vertex:
%
\qq
\tcBas = \span{D_1+D_2\quad\mbox{for all pairs $D_1,D_2$}}.
\label{1.9}
\qqq
In order to define $\tcBihx$, consider a space $\tcBs$ whose basis
vectors are graphs with 1-valent and 3-valent vertices and exactly
one 4-valent vertex. We define a linear map
$\dihx:\,\tcBs\longrightarrow \tcB$ by its action on the individual
graphs of $D\in\tcBs$
\qq
\dihx:\; D\mapsto D_1-D_2+D_3,
\label{1.10}
\qqq
where all four graphs $D,D_1,D_2,D_3$ are the same except at one
spot, where they differ according to \fg{f2}.
\begin{figure}[hbt]
\centerline{\epsfbox{pct2.eps}} \caption{The graphs $D$, $D_1$,
$D_2$ and $D_3$}
\label{f2}
\end{figure}
Then we define the second subspace $\tcBihx\subset\tcB$ as the
image of $\dihx$.

Now we introduce a space
\qq
\cB = \tcB / (\tcBas + \tcBihx).
\label{1.11}
\qqq
%
Since the graphs $D_1,D_2$ of\rx{1.9} and the graphs $D_1,D_2,D_3$
of\rx{1.10} have the
same degrees\rx{1.7*},\rxw{1.7} among themselves,
then both subspaces $\tcBas$ and $\tcBihx$ respect the
gradings\rx{1.8} and as a result the space $\cB$ is also graded
\qq
\cB = \bopnzi \cB_n, \qquad\mbox{where}\;\;
\cB_n =
\bopmzi\cB_{m,n}.
\label{1.12}
\qqq
It is well-known that the dual space $\cBs$ is isomorphic to the
space of all Vassiliev invariants of knots, and the grading
$\cBs=\bopnzi \cBs_n$ corresponds to the grading of Vassiliev
invariants.

The space $\cB$ can be endowed with a commutative
algebra structure. We define the product
of two graphs in $\tcB$ as their disjoint union. It is easy to see that this
product respects the gradings\rx{1.7*},\rx{1.7} and that the subspace
$\tcBas + \tcBihx$ is the ideal in algebra $\tcB$. Therefore, the quotient
space $\cB$ is also an algebra.

We are going to introduce another algebra $\cD$ which is isomorphic to
$\cB$. This construction has been known to some people\cx{pr}. It
appeared as an attempt to better understand the structure of $\cB$
and, in particular, to evaluate the dimension of the spaces
$\cB_n$. I am especially indebted to A.~Vaintrob for illuminating
discussions on the structure of $\cD$. I introduce
the algebra $\cD$ in order to formulate a conjecture about the structure
of Kontsevich integral, which was motivated by the study
of the Melvin-Morton expansion of the colored Jones polynomial as it
comes of $R$-matrix expression and which has now been proved by
A.~Kricker\cx{Kr}.

We begin by defining a bigger space $\tcD$.
Let $D$ be a graph with
3-valent vertices and no 1-valent vertices. We think of this graph as
a $CW$-complex and consider a space of its rational cohomologies
$\Hor{D}$. Let $\SgrD$ be the group of symmetry of a graph $D$ (it maps
3-vertices to 3-vertices and edges to edges) and let $\SgroD\subset\SgrD$ be
its subgroup which preserves the cyclic order of the edges at the vertices.
$\SgrD$ acts naturally on $\Hor{D}$ and this group action can be extended to
the symmetric algebra $\SHor{D}$. We denote by $\Hso{D}$ the $\SgroD$-invariant
part of the latter space:
%
\qq
\Hso{D} = \bopmzi \Hsmo{D},\qquad
\Hsio{m}{D} = \Invgdo{\SmHor{D}},
\label{1.13}
\qqq
while $\PgoD$ is the corresponding projector
\qq
\PgoD:\,\SHor{D} \longrightarrow \Hso{D},\qquad
\PgoD(x) = {1\over \rnk{\SgroD}} \sum_{g\in \SgroD} g(x),
\label{1.14}
\qqq
where $\rnk{\SgroD}$ denotes the number of elements in $\SgroD$.
Now we define a linear space $\tcDo$ as
\qq
\tcDo = \bopmnzi \tcDio{m,n},\qquad\mbox{where}\;\;
\tcDio{m,n} = \bigoplus_{D:\,\chi(D)=n} \Hsmo{D}.
\label{1.15}
\qqq

The space $\tcDo$ has an associative, commutative
algebra structure. First, note that for a disjoint
union $D_1\cup D_2$ of two graphs $D_1,D_2$
\qq
\Hor{D_1\cup D_2} = \Hor{D_1} \oplus \Hor{D_2}
\qqq
and therefore
\qq
\SHor{D_1\cup D_2} = \SHor{D_1} \otimes \SHor{D_2}
\qqq
as algebras. The latter equation allows us to define a product of two elements
$x_i\in\Hor{D_i}$, $i=1,2$ as a projection of their tensor product in
$\SHor{D_1\cup D_2}$
\qq
x_1 x_2 = \Pgio{D_1\cup D_2} (x_1\otimes x_2) \in \Hso{D_1\cup D_2}.
\label{1.15*}
\qqq
If the graphs $D_1$, $D_2$ do not have isomorphic connected components, then
$\Sgroi{D_1\cup D_2} = \Sgroi{D_1}\times \Sgroi{D_2}$
and the projector in \ex{1.15*} may
be omitted: $x_1 x_2 = x_1\otimes x_2$. The commutativity of the
product\rx{1.15*} is obvious. Associativity follows from a relation
\qq
(x_1 x_2) x_3 = x_1 (x_2 x_3) = \Pgio{D_1\cup D_2 \cup D_3}
(x_1\otimes x_2\otimes x_3)
\qqq
Finally, since the product\rx{1.15*}
respects both gradings\rx{1.15}, then $\tcDo$ is a graded algebra.

Next, we define the subspace $\tcDASo\subset\tcDo$
which comes from the change of
cyclic order at 3-valent vertices.
The definition of the symmetric algebra $S^*\Hor{D}$ is independent of this
cyclic order. Therefore if we take a graph $D_1$ and
change the cyclic order at one of its vertices, thus producing a new graph
$D_2$, then there is a natural isomorphism of cohomologies
\qq
\hfAS:\; \Hor{D_1}\longrightarrow \Hor{D_2},
\label{1.16}
\qqq
because $D_2$ was constructed in such a way that there is a natural isomorphism
between $D_1$ and $D_2$ as $CW$-complexes (generally, there could be more than
one isomorphism due to the symmetry group $G_{D_1}$). The isomorphism\rx{1.16}
can be extended to an isomorphism of symmetric algebras
\qq
\hfAS:\; \SHor{D_1}\longrightarrow \SHor{D_2},
\label{1.17}
\qqq
%
Let $\tVAS$ be the graph of this map
\qq
\tVAS = \{ (x,y) |\; y = \hfAS(x) \}\subset \SHor{D_1}\oplus
\SHor{D_2}.
\label{1.18}
\qqq
We denote by $\VAS$ its projection onto $\Hso{D_1}\oplus \Hso{D_2}$
\qq
\VAS = \Pgio{D_1}\Pgio{D_2}(\tVAS) \subset \Hso{D_1}\oplus \Hso{D_2}.
\label{1.19}
\qqq
We define the subspace $\tcDASo\subset\tcDo$ as the sum of all the spaces $\VAS$
for all 3-valent diagrams $D_1$ and all choices of vertices of $D_1$ where we
change the orientation. It is easy to check that $\tcDASo$ is an ideal: for any
element $x\in \SHor{D_1}$, and for any element $y\in\Hso{D_3}$
\qq
(\Pgio{D_1}(x) + \Pgio{D_2} \hfAS(x))\, y & = &
\Pgio{D_1\cup D_3}(\Pgio{D_1}(x)\otimes y) +
\Pgio{D_2\cup D_3}(\Pgio{D_2}\hfAS(x)\otimes y)
\nonumber\\
& = &
\Pgio{D_1\cup D_3}(x\otimes y) + \Pgio{D_2\cup D_3} \hfAS(x\otimes y),
\qqq
because obviously
\qq
\Pgio{D_i\cup D_3} (\Pgio{D_i}\otimes I) = \Pgio{D_i\cup D_3},\qquad i=1,2.
\qqq
and
\qq
\hfAS(x\otimes y) = \hfAS(x) \otimes y,
\qqq
where in the \lhs $\hfAS$ comes from the change of cyclic order at a vertex in
the whole graph $D_1\cup D_3$.

Finally, we define a subspace $\tcDIHXo\subset\tcDo$. Let $D$ be a graph with
3-valent vertices and exactly one 4-valent vertex, and with fixed cyclic
order at every vertex.
By adding an extra edge to $D$, we ``resolve'' the 4-valent vertex in
3 different ways, thus converting $D$ into one of the 3-valent graphs
$D_1,D_2,D_3$ of \fg{f2}. A removal of this extra edge generates 3
natural maps of rational homologies
\qq
\hf_i:\; \oHor{D_i} \longrightarrow \oHor{D},\qquad i=1,2,3.
\label{1.20}
\qqq
We extend the dual maps
$\hfs_i:\; \Hor{D}\longrightarrow\Hor{D_i}$ as
algebra homomorphisms
\qq
\hfs_i:\; \SHor{D}\longrightarrow\SHor{D_i},\qquad i=1,2,3.
\label{1.21}
\qqq
We define the map
$\hdihx:\; \SHor{D} \longrightarrow \bigoplus_{i=1}^3\Hso{D_i}$
by the formula (\cf \ex{1.10})
\qq
\hdihx = - \Pgio{D_1} \hfs_1 + \Pgio{D_2} \hfs_2 - \Pgio{D_3} \hfs_3.
\label{1.22}
\qqq
The subspace $\tcDIHXo$ is the sum of the images of all the maps
$\hdihx$ for all the graphs $D$. It is easy to check that similarly to
$\tcDASo$, $\tcDIHXo$ is also an ideal in $\tcDo$.

Now we define the quotient space
\qq
\cD = \tcDo / (\tcDASo + \tcDIHXo).
\label{1.23}
\qqq
Since the graphs $D_1$, $D_2$ of \rxw{1.16} and $D$, $D_1$, $D_2$, $D_3$
of \rxw{1.20} have the same Euler characteristic among themselves and since the
maps\rx{1.17} and\rx{1.21} preserve the grading of symmetric algebras, then
$\cD$ is a graded algebra:
%
\qq
\cD = \bopmnzi \cD_{m,n},
\label{1.24}
\qqq
where the spaces $\cD_{m,n}$ are the quotients of
the spaces $\tcD_{m,n}$.
\qq
\cD_{m,n} = \tcD_{m,n} / \cD_{m,n}\cap (\tcDAS + \tcDIHX).
\label{1.24*}
\qqq
%

This description of the algebra $\cD$
makes it easy to establish its
isomorphism with the algebra $\cB$, but there exists a slightly different
description of $\cD$ which suits better for the formulation of our conjecture
about the structure of Kontsevich integral. Recall that $\SgrD$ denotes the
full symmetry group of a 3-valent graph $D$ (including the maps which do not
preserve the cyclic order at the vertices). As we have mentioned, $\SgrD$ acts
naturally on $\SHor{D}$. We modify this action by multiplying the action of an
element $g\in\SHor{D}$ by $(-1)^{|g|}$, where $|g|$ denotes the number of
vertices of $D$ whose cyclic order is changed by $g$. Now instead of\rx{1.13}
we define
\qq
\Hsx{D} = \bopmzi \Hsmx{D},\qquad
\Hsix{m}{D} = \Invgdx{\SmHor{D}},
\label{1.13*1}
\qqq
while $\PgxD$ is the corresponding projector
\qq
\PgxD:\,\SHor{D} \longrightarrow \Hsx{D},\qquad
\PgxD(x) = {1\over \rnk{\SgrD}} \sum_{g\in \SgrD} g(x),
\label{1.14*1}
\qqq
Let $\bD$ be a set of all 3-valent graphs with a particular cyclic order of
edges at vertices chosen for every graph (so that each isomorphism class of
3-valent graphs is represented in $\bD$ exactly once). Define
\qq
\tcD = \bopmnzi \tcDix{m,n},\qquad\mbox{where}\;\;
\tcDix{m,n} = \bigoplus_{D:\,\chi(D)=n} \Hsmx{D}
\label{1.14*2}
\qqq
(\cf \ex{1.15}).
%
%
If we choose a different set $\bD\p$, then there is a natural isomorphism
between $\tcDbd$ and $\tcD_{\bD\p}$. Namely, if $D_1\in \bD$ and
$D_2\in \bD\p$ represent the same 3-valent graph (but possibly with different
cyclic orders), then we identify the spaces $\Hsx{D_1}$ and $\Hsx{D_2}$ by an
identity map with an extra sign factor $(-1)^{\dsd}$,
where $\dsd$ is the number of
vertices in the graphs $D_1$, $D_2$ which have different cyclic orders. In the
future we will sometimes denote $\tcDbd$ simply as $\tcD$, assuming that the
choice of cyclic order for every 3-valent graph was somehow fixed.

\begin{lemma}
There is a natural isomorphism $\tcDbd \cong \tcDo/\tcDASo$.
\end{lemma}
\proof
Since $\SgroD\subset \SgrD$, then $\Hsx{D}\subset\Hso{D}$. As a result,
$\tcDbd$ may be considered a subspace of $\tcDo$ and thus we have a map
$f\,:\;\tcDbd\longrightarrow \tcDo/\tcDASo$. On the other hand, one can
construct a natural map $g\,:\;\tcDo\longrightarrow \tcDbd$ in the following
way: if a 3-valent graph $D_1$ is isomorphic to a graph $D_2\in\bD$, then $g$
maps $\Hso{D_1}$ to $\Hsx{D_2}\subset\Hso{D_1}$ as
$(-1)^{\dsd}\Pgix{D_1}$. Obviously, $\tcDASo\subset \ker\,g$, so we have a map
$h\,:\;\tcDo/\tcDASo \longrightarrow \tcDbd$. We leave it to the reader to
check that $f$ and $h$ constitute an isomorphism.\qed

After constructing an isomorphism
$h\,:\;\tcDo/\tcDASo\longrightarrow \tcDbd$ we define the space $\tcDIHXx$
simply as the image of
$\tcDIHXo/ (\tcDIHXo\cap\tcDASo)$. Thus we proved the following
\begin{theorem}
\label{t2.1}
There is a natural isomorphism $\cD \cong \tcDbd/\tcDIHXx$.
\end{theorem}
The grading subspaces
$\cD_{m,n}$ turn out to be the quotients $\tcD_{m,n}/(\tcD_{m,n}\cup\tcDIHXx)$.

The advantage of this description of $\cD$ is that it allows us to
work with rather natural spaces $\Invgdx{\SmHor{D}}$ instead of bigger
and less symmetric spaces $\Invgdo{\SmHor{D}}$.

\nsection{Isomorphism between $\cB$ and $\cD$}
\label{s3}

\begin{theorem}
\label{t3.1}
There exists a canonical isomorphism of algebras
\qq
\hxA:\;\cB\longrightarrow \cD,
\label{3.1}
\qqq
which respects the grading
\qq
\hxA:\;\cB_{m,n}\longrightarrow \cD_{m,n-m}.
\label{3.2}
\qqq
\end{theorem}
\begin{corollary}
If $m>n$, then $\cB_{m,n} = \empt$.
\end{corollary}

Before we prove this theorem, we have to establish some facts
concerning the structure of the space $\cB$. We call an edge of a
(1,3)-valent graph \emph{a leg} if this edge is connected to a
1-valent vertex. All other edges are called \emph{internal}.

\begin{lemma}
\label{l3.1}
If two legs of a (1,3)-valent graph $D$ are attached to the same 3-valent
vertex, then $D\in\tcBas$.
\end{lemma}

\proof
Suppose that a (1,3)-valent graph $D$ contains such a 3-valent
vertex. Since the 1-valent vertices of our graphs are not ordered in
any way, then changing the cyclic order at that 3-valent vertex does
not change the graph. Therefore $2D\in\tcBas$ and this proves the
lemma.\qed

Let us call a (1,3)-valent graph \emph{restricted} if each of its
3-valent vertices contains at most one leg. Let $\tcBr$ be a formal
space whose basis vectors are restricted graphs. We introduce
familiar subspaces. The subspaces $\tcBasrii\subset \tcBr$, $i=0,1$
are spanned by the sums of restricted diagrams $D_1,D_2$ which differ
in the ordering at a 3-valent vertex which is attached to $i$ legs.
The subspaces $\tcBihxii\subset \tcBr$, $i=0,1$ are spanned by the
images of the map\rx{1.10} acting on the (3,4)-valent diagrams
whose
single
4-valent vertex contains $i$ legs. Then Lemma\rw{l3.1} has a simple
corollary:
\qq
\cB_{m,n} = \tcBrz_{m,n} / (\tcBasrz + \tcBihxz),\qquad
\mbox{where}\;\; \tcBrz_{m,n} = \tcBr_{m,n}/ (\tcBasro + \tcBihxo) .
\label{3.3}
\qqq
Indeed, this relation follows from the fact that if the 4-valent
vertex of a (3,4)-valent graph $D$ has at least two legs, then the
intersection of the image of the corresponding operator\rx{1.10} with
the space $\tcBr$ is trivial.
Also, it is easy to see that $\tcBasro$ and $\tcBihxo$ are ideals in $\tcBr$,
so the quotient
\qq
\tcBrz = \bopmnzi \tcBrz_{m,n} = \tcBr / (\tcBasro + \tcBihxo)
\qqq
has a graded algebra structure.

Now we begin to construct the isomorphism. Let $D$ be a 3-valent
graph with $N$ edges and cyclic order at vertices.
Thinking of $D$ as a $CW$-complex, let $\Co$ be the space of 1-chains. In other
words, $\Co$ is an $N$-dimensional vector space spanned by the oriented edges
of $D$, if we assume that an edge with the opposite orientation is equal to the
opposite of the edge as an element of $\Co$. Thus, if we pick an orientation on
the edges of $D$, then $\Co$ has a natural basis $e_j$, $1\leq j\leq N$ of the
edges of $D$. We will also need the dual space $\Cos$ with the dual basis
$f_j$, $1\leq j\leq N$. The symmetry group of the graph $G_D$ acts on both
spaces $\Co$ and $\Cos$.


Next, consider a vector space whose basis is formed by
$m$-legged (1,3)-valent restricted graphs such that if we remove
their legs, then we get the 3-valent graph
$D$. We denote the quotient of this space by
its intersection with $\tcBasro$ as $\tcBmd$. We also have to
consider a bigger space. Suppose that we index the edges of $D$ and
then attach $m$ legs to its edges in order to produce restricted graphs. These
(1,3)-valent
graphs still carry the indexing of the edges of $D$. If we factor
this space by its intersection with the obvious analog of $\tcBasro$,
then we get the space $\cchBmdp$. The symmetry group $\SgroD$ of the
graph $D$ acts on $\cchBmdp$ by mapping the edges of $D$ together with their
legs, while preserving the cyclic order at the vertices.
The invariant subspace of this action
is canonically isomorphic to $\tcBmd$:
\qq
\tcBmd = \Invgdo{\cchBmdp}.
\label{3.4*}
\qqq

Let us introduce a
multi-index notation
\qq
\um = (m_1,\ldots,m_N), \qquad |\um| = \sjoN m_j.
\label{3.4}
\qqq
For $N$ non-negative numbers $\um$ and for a choice of orientation of the
edges of $D$
construct a diagram $D_\um$ in
the following way: for every $j$, $1\leq j\leq N$ attach $m_j$ legs
to $D$ on the left side of the edge $e_j$ (the notion of the left
side is well-defined since $e_j$ is oriented). It is easy to see that
all graphs $D_\um$, $|\um|=m$ form a basis of the space $\cchBmdp$,
because after we took the quotient over the analog of the space
$\tcBasro$, we can flip the legs of the graphs of $\cchBmdp$ to a
particular side of each edge of $D$ (at the cost of changing the
signs of the corresponding vectors of $\cchBmdp$).

There is a natural isomorphism
$\xA:\;\cchBmdp \longrightarrow \SmCos$ which acts on the basis
vectors as
\qq
\hxA:\;D_\um\mapsto\pjoN f_j^{m_j}.
\label{3.5}
\qqq

Suppose that the 3-valent graph $D$ has $\Nz$ vertices $v_j$,
$1\leq j\leq \Nz$. Consider the $\Nz$-dimensional space $\Cz$ of 0-chains whose
basis vectors are in a one-to-one correspondence with these vertices.
Then there is a natural boundary map $\xdel:\;\Co\longrightarrow\Cz$.
Let $\chCos$ be the space of 1-cocycles, it is
the subspace of $\Cos$ whose elements annihilate the
kernel of $\xdel$. Apparently,
\qq
\Hor{D} = \Cos/\chCos.
\label{3.6}
\qqq

Let $\cchBihxrodm = \cchBmdp\cap\cchBihxro$, where the space
$\cchBihxro$
is the analog of the space $\tcBihxro$ for the graphs which come from
3-valent graphs with indexed edges.

\begin{lemma}
\label{l3.2}
The map $\hxA$ establishes an isomorphism between the spaces
$\cchBihxrodo$ and $\chCos$.
\end{lemma}

\proof
For $1\leq j\leq \Nz$, denote as $V_j$ the image in $\cchBodp$ of
the operator\rx{1.10} associated with the vertex $v_j$ of $D$ (that
is, one of the two 3-valent vertices in each of the graphs of \fg{f2}
is $v_j$, while the other vertex is attached to a leg). Then the
space $\cchBihxrodo$ is spanned by all the spaces $V_j$.

For $1\leq j\leq \Nz$ and for $x\in\Co$ let $\xdel_j(x)$ be the
coefficient in front of $v_j\in\Cz$ in the expansion of $\xdel(x)$
with respect to the basis $v$.
Then $\ker\xdel = \bcap_{j=1}^{\Nz} \ker\xdel_j$
and, as a result, the
space $\chCos$ is spanned by the spaces
$V\p_j\subset\Cos$ which annihilate the spaces
$\ker\xdel_j\subset\Co$. It is very easy to see that for every
$j$, $\hxA$ establishes an isomorphism between the corresponding
spaces $V_j$ and $V\p_j$. This proves the lemma.\qed

\begin{lemma}
\label{l3.3}
$\hxA$ establishes the isomorphism between the spaces
$\cchBmdp/\cchBihxrodm$ and
\\
$\SmHord$.
\end{lemma}

To prove this lemma we need a simple fact from linear algebra.
\begin{lemma}
\label{l3.4}
Let $V$ be a finite dimensional vector space and $W$ be its subspace.
Denote by $\Ps$ a symmetrizing projector
$\Ps:\;V^{\otimes m}\longrightarrow S^mV$. Then
\qq
S^m V / \Ps(S^{m-1}V\otimes W) = S^m(V/W).
\label{3.7}
\qqq
\end{lemma}
\proof
We leave the proof to the reader.

\pr{Lemma}{l3.3}
It is easy to see that $\hxA$ maps the space $\cchBihxrodm$ onto
$\Ps(S^{m-1}\Cos\otimes \chCos)$. Then the claim of the lemma follows
from \eex{3.6} and\rx{3.7} if we set $V=\Cos$ and $W=\chCos$ in the
latter equation.\qed

Consider a space $\tcBihxrodm = \tcBmd \cap \tcBihxro$.

\begin{lemma}
\label{l3.5}
There is a natural isomorphism between the quotient spaces
\qq
\tcBmd / \tcBihxrodm = \Invgdo{\cchBmdp /\cchBihxrodm}.
\label{3.8}
\qqq
\end{lemma}

In order to prove this isomorphism we need another linear algebra
lemma.

\begin{lemma}
\label{l3.6}
Let $V$ be a finite-dimensional representation of a finite group $G$.
Let $W\subset V$ be a subspace, which is invariant under the action
of $G$. Then there is a natural isomorphism
\qq
\Invg{V} / \Invg{W}  = \Invg{V/W}.
\label{3.9}
\qqq
\end{lemma}
\proof
For example, one could use the fact that a finite-dimensional
representation of $G$ is a sum of irreducible representations. We
leave the details to the reader.\qed

\pr{lemma}{l3.5}
The cyclic order preserving
symmetry group $\SgroD$ of the 3-valent graph $D$ acts on the space
$\cchBmdp$. Obviously, the symmetrization over this action projects
$\cchBmdp$ onto $\tcBmd$. Thus
\qq
\tcBmd = \Invgdo{\cchBmdp}.
\label{3.10}
\qqq
At the same time, the subspace $\cchBihxrodm$ is invariant under the
action of $\SgroD$ and
\qq
\tcBihxrodm = \Invgdo{\cchBihxrodm}.
\label{3.11}
\qqq
Then \ex{3.8} follows from \ex{3.9} in view of the relations\rx{3.10}
and\rx{3.11}.\qed

Let us introduce a notation $\cBmd = \tcBmd / \tcBihxrodm$.

\begin{corollary}
\label{c3.2}
The map $\hxA$ establishes the isomorphism between the spaces $\cBmd$
and $\Hsio{m}{D}$ (see \ex{1.13} for the definition of the latter
space).
\end{corollary}

\proof This isomorphism follows from the combination of
Lemmas\rw{l3.3} and\rw{l3.5}.\qed

We leave it for the reader to check that the isomorphism $\hxA$ intertwines the
maps
\qq
\cB_{m_1}(D_1)\otimes\cB_{m_2}(D_2) &\longrightarrow &
\cB_{m_1 + m_2} (D_1\cup D_2)
\nonumber\\
\Hsio{m_1}{D_1}\otimes\Hsio{m_2}{D_2} &\longrightarrow &
\Hsio{m_1 + m_2}{D_1\cup D_2}
\qqq
which come from the multiplications in the algebras $\tcB$ and $\tcD$ as
defined in Section\rw{gs}.


\pr{Theorem}{t3.1}
According the definition\rx{3.3} of the space $\tcBrz_{m,n}$,
\qq
\tcBrz_{m,n+m} = \bigoplus_{D:\;\chi(D)=n} \cBmd,
\label{3.12}
\qqq
while by its definition
\qq
\tcDo_{m,n} = \bigoplus_{D:\;\chi(D)=n} \Hsio{m}{D}.
\label{3.13}
\qqq
It is easy to see that $\hxA$ establishes the isomorphisms
\qq
\hxA:\; \tcBasrz\cap\tcBrz_{m,n+m}\longrightarrow
\tcDASo\cap\tcDo_{m,n},\qquad
\tcBihxz\cap\tcBrz_{m,n+m}\longrightarrow\tcDIHXo\cap\tcDo_{m,n}.
\label{3.14}
\qqq
Then \ex{3.2} follows from \eex{3.3} and\rx{1.24*} together with the
isomorphism of Corollary\rw{c3.2}.\qed

\nsection{Rationality Conjecture}
\label{s4}

Recall that Kontsevich integral of a knot $\cK\in S^3$ is a sequence
of vectors $\IB_{m,n}(\cK)\in \cB_{m,n}$, $m\geq 0$, $n\geq m$
depending on the topological class of $\cK$. The space $\cB_{0,0}$ is
1-dimensional, its basis vector is the empty graph, so it can be
naturally identified with $\IC$. It is known that
$\IB_{0,0}(\cK) = 1$.

We combine the vectors $\IB_n(\cK)$ into a formal
power series of a formal variable $\hb$
\qq
\IB(\cK;\hb) = 1 + \smnoge \IB_{m,n}(\cK)\,\hb^n\in\cB.
\label{4.1}
\qqq
Prior to formulating a conjecture about the structure of
$\IB(\cK;\hb)$
we have to apply to it some transformations. First, we apply
the wheeling map $\hOm:\;\cB\longrightarrow\cB$, described in\cx{Wh},
in order to produce
\qq
\IO(\cK;\hb) = \hOm(\, \IB(\cK;\hb)) =
1 + \smnog \IO_{m,n}(\cK)\,\hb^{m+n} \in \cD,\qquad
\IO_{m,n}\in\cB_{m,n}.
\label{4.1*}
\qqq
Then we apply the isomorphism $\hxA$, which maps Kontsevich integral
from $\cB$ to $\cD$. More precisely, we choose a set
$\bD$ of 3-valent graphs $D$ such
that each type of a graph (without distinguishing them by cyclic
order at vertices) is represented there exactly once, and then we map $\cB$ to
$\cD_{\bD}$ as described at the end of Section\rw{gs}.
Thus we get
\qq
\ID(\cK;\hb) = \hxA(\, \IB(\cK;\hb)) =
1 + \smnog \ID_{m,n}(\cK)\,\hb^{m+n} \in \cD,
\qquad \ID_{m,n}(\cK)\in \cD_{m,n}.
\label{4.2}
\qqq
%


By using the algebra structure of $\cD$ and manipulating the formal power
series in $\hb$ we can define the logarithm of Kontsevich integral
\qq
&
\Il(\cK;\hb) = \log \ID(\cK;\hb) =
\smnog \Il_{m,n}(\cK)\,\hb^{m+n} \in \cD,
\qquad \Il_{m,n}(\cK)\in \cD_{m,n}.
\label{4.3*}
\qqq
through the formula
\qq
\log(1+x) = \sum_{n=1}^{\infty} (-1)^{n-1}{x^n\over n}
\qqq
The advantage of the logarithm $\Il(\cK;\hb)$ is that it can be expressed
exclusively in terms of \emph{connected} 3-valent graphs.

Kontsevich integral $\Il(\cK;\hb)$ belongs to the quotient space\rx{1.23}. Let
$\tIl(\cK;\hb)$ be a representative of $\Il(\cK;\hb)$ in the space
$\tcD$ (Of course,
it is defined only up to an element
of $\tcDIHX$).
We present
$\tIl(\cK;\hb)$ as
\qq
\tIl(\cK;\hb) = \sDc \smzi x_m(\cK,D)\, \hb^{\chi(D)+m},
\label{4.4}
\qqq
where $\bDc\subset \bD$ is a subset of connected 3-valent graphs and
$x_m(\cK,D)\in \Hsm{D}$.

Now we are almost ready to formulate our conjecture. Let $V$ be a vector space.
For $x\in V$ we define $e^x\in \Ss V$ by the power series $e^x =\snzi x^n/n!$.
If $\Lambda$ is a lattice in $V$, then we extend this exponential map to
an injection of a group algebra $\Exp:\;\IQ[\Lambda]\rightarrow \Ss V$. For a
graph $D$, $\Hoz{D}$ forms a lattice in $\Hor{D}$.
We denote
\qq
\HsQ{D} = \Exp
\Qalghz
\subset\Hs{D}.
\qqq
In other words, $\HsQ{D}$ is $\SgrD$-invariant part of the rational span of the
exponents of the elements of $\Hoz{D}$ and $\Exp$ establishes its isomorphism
with $\Qalghz$.

Now recall that if $D$ has $N$ edges, then $e_j$ ($1\leq j\leq N$)
denote the oriented edges forming a basis in the space of 1-chains $\Co$, while
$f_j$, $1\leq j\leq N$ form the dual basis in the dual space $\Cos$.
In view of \ex{3.6} we can think of $f_j$ as elements of $\Hor{D}$.

\begin{lemma}
The product of the Alexander-Conway polynomial of $e^{f_j}$ is an element of
the algebra
$\HsQ{D}$:
\qq
\pjoN \APbas{\cK}{\exp(f_j)}\in\HsQ{D}.
\label{4.4*}
\qqq
and its inverse is a well-defined element of $\Hs{D}$
\end{lemma}
\proof
To prove relation\rx{4.4*}, we have to show that its \lhs is $G_D$-invariant.
The elements of the group $G_D$ not only permute $f_j$, $1\leq j\leq N$, but
they may also reverse the orientation of some edges of $D$ and thus change the
signs of corresponding $f_j$. However, the relation
\qq
\APbas{\cK}{1/t} = \APbas{\cK}{t},
\label{4.5*}
\qqq
guarantees that this change of sign does not affect the expression\rx{4.4*}
and hence it is $G_D$-invariant.
At the same time, the Alexander-Conway polynomial satisfies the
property $\APbas{\cK}{1}=1$ which guarantees that the inverse of\rx{4.4*} can
be inverted within $\Hs{D}$.\qed


Let us introduce a notation
\qq
\Il(\cK,D) = \smzi x_m(\cK,D) \in \Invgd{\SHor{D}}.
\label{4.5*1}
\qqq
The only 3-valent graph $D$ with $\chi(D)=0$ is a
 circle.
The value
of $\Il(\cK,\crcl)$ has been established by D.~Bar-Natan and
S.~Garoufalidis in\cx{BG}
\qq
\Il(\cK,\crcl)
= \hlfv
\lrbs{
\log \lrbc{\sinh(f/2) \over (f/2) } -
\log \APbas{\cK}{\exp(f)}
},
\label{4.6}
\qqq
where $f$ represents the integral generator of $\Hor{\crcl}$. Our conjecture
deals with the value of $\Il(\cK,D)$ for graphs with $\chi(D)\geq 1$. Recall
that such graphs have exactly $N=3\chi(D)$ edges.

\begin{conjecture}
\label{c4.1}
The representative $\tIl(\cK;\hb)\in\tcD$ of Kontsevich integral
$\Il(\cK;\hb)\in\cD$ can be chosen in such a way that for any
$D\in\bD$, $\chi(D)\geq 1$ there exists an element
$y(\cK,D)\in\HsQ{D}$
such that
\qq
\Il(\cK,D)
= {y(\cK,D) \over \pjochD \APbas{\cK}{\exp(f_j)} }.
\label{4.5}
\qqq
\end{conjecture}
\begin{remark}
\rm
Andrew Kricker has proved this conjecture in his paper\cx{Kr}.
\end{remark}

\begin{remark}
\rm
D.~Thurston presented arguments which show that if Conjecture\rw{c4.1} is true
as it is formulated, then it should also be true if one defines $\ID(\cK,\hb)$ directly
as an image of $\IB(\cK,\hb)$ under the isomorphism $\hxA$ without appying the wheeling
map $\hOm$ of \ex{4.1*}.
\end{remark}

\begin{remark}
\rm
It is convenient to introduce some other notations in relation to \ex{4.5}. Let
$p(\cK,D)\in \Qalghz$ be such that $\Expalg{p(\cK,D)} = y(\cK,D)$. Also, if
we index the edges of $D$ in such a way that $\fbs$ form a
basis of
$\Hoz{D}$ and $\Hor{D}$, then we can write $\Il(\cK,D)$ and
$y(\cK,D)$ more explicitly as
\qq
\Il(\cK,D) & = & \Il(\cK,D;\fbs),
\label{4.5*2}\\
p(\cK,D) & = & p(\cK,D;\fbs),
\label{4.5*2*1}\\
y(\cK,D) & = & p(\cK,D;\efbs),
\label{4.5*2*2}
\qqq
where
\qq
&\Il(\cK,D;\xbs)\in\IQ[[\xbs]],
\label{4.5*3*1}\\
& p(\cK,D;\tbs)\in
\IQtbs.
\label{4.5*3}
\qqq

\end{remark}

\nsection{Rational structure of the Jones polynomial}
\label{s5}

There is a well-known relation between the Kontsevich integral and the colored
Jones polynomial of a knot, so the rationality conjecture\rw{c4.1} should
manifests itself in the structure of the latter object. In fact, this
manifestation observed in\cx{Ro}, served for us as evidence which led to the
rationality conjecture. Another advantage in establishing a relation between
\ex{4.5} and the rational expansion of the Jones polynomial\cx{Ro} is that
at present it is much easier to calculate the colored Jones polynomial than
Kontsevich integral. Therefore, working out the rational expansion of\cx{Ro} is
a practical way of finding the polynomials $y(\cK,D)$ of \ex{4.5}.

Let us recall the exact relation between the Kontsevich integral and a colored
Jones (or, more generally, HOMFLY) polynomial based on a simple Lie algebra
$\gg$. We equip $\gg$ with the ad-invariant
%
scalar product
normalized in such a way that long roots have length $\sqrt{2}$ (this
scalar product allows us to identify the dual space $\gg^*$ with
$\gg$ itself). Let
$\val\in\gh$ be the hightest weight of a representation of $\gg$,
shifted by $\vr$ (which is half the sum of positive roots of $\gg$).
Reshetikhin and Turaev associate to this data a polynomial
$\JvaK\in\ZZ[q^{\pm 1/2}]$. If we substitute
\qq
q = e^{\hb},
\label{5.1*}
\qqq
then we
can expand $\JvaK$ in power series of $\hb$
\qq
\JvaK = \snzi \pnvaK\,\hb^n,
\label{5.1}
\qqq
whose coefficients $\pnvaK$ are polynomials of $\val$.
The same series\rx{5.1} can be deduced from the value of Kontsevich
integral.

The data $\gg,\val$ defines an element in the dual space
$\cBs$, which is called \emph{the weight system}. We will define it
in such a way that it will be
suitable for application to
$\IO(\cK;\hb)$.
The first steps in the definition of the weight systems are fairly
standard. Let $\vx_a$, $1\leq a\leq \dim \gg$ be a basis of $\gg$.
Define the structure constants $f_{abc}$ by the relation
\qq
[\vx_a,\vx_b] = \sum_{c=1}^{\dim \gg} \fabco\, \vx_c.
\label{5.2}
\qqq
We can raise and lower the indices of $\fabco$ with the help of the
metric tensor
\qq
h_{ab} = \vx_a\cdot \vx_b
\label{5.3}
\qqq
and its inverse $h^{ab}$.

Let $D$ be a (1,3)-valent graph, $\dego(D)=m$, $\degt(D)=n+m$.
Suppose that if we strip off its legs, then we get a 3-valent graph
$\Dz$. Let us orient the edges of $\Dz$ and assign
orientation to the edges of $D$ in such a way that it is compatible
with the orientation of $\Dz$ and legs are oriented in the direction
from 1-valent vertex to 3-valent vertex. Next, we assign the tensors
$f$ to 3-valent vertices, assigning their indices to attached edges
according to the cyclic ordering. We use the upper indices for the
incoming edges and lower indices for the outgoing edges. Finally, we
take the product of all tensors $f$ assigned to 3-valent vertices,
contract each pair of indices of $f$'s along each internal edge,
while contracting each index assigned to a leg with $\a_a$
($\val = \sum_{a=1}^{\dim\gg} \a_a\,\vx_a$).
Thus we get a Weyl group invariant
homogeneous polynomial $\pDval$ of $\val$ of degree $m$. It is easy
to see that it does not depend on the choice of orientation of the
edges of $\Dz$. For a fixed weight $\val$, $\yp$ assignes a number to
each (1,3)-valent graph, so $\yp\in\tcB^*$. In fact, due to the
anti-symmetry of $f$ and to the Jacobi identity, satisfied by
the commutator\rx{5.2}, $\yp$ annihilates the subspaces $\tcBas$ and
$\tcBihx$ and therefore it can be projected to $\cBs$.

The usual way to proceed further is to convert $\pD$ as a Weyl group invariant
polynomial on $\gh$ into an element of $\Inv{\gg}{S^m\gg}$, then use a PBW map
to convert it into an element of $\Inv{\gg}{U\gg}$ and calculate the trace of
that element in a $\gg$ module with the highest weight $\val - \vr$, thus
obtaining another polynomial $\pval$ of $\val\in\gh$ which is the standard
weight of the graph $D$ coming from $\gh$, or thinking of it as a function on
all graphs $D$, $\yps$ is a weight system on $\cB$. Then the relation between
the expansion\rx{5.1} and Kontsevich integral is
\qq
\JvaK = \dval \lrbc{
1 + \smnog \yps(\,\IB_{m,n}(\cK),\val)\,\hb^{m+n} },
\label{5.4*}
\qqq
where $\dval$ is the dimension of the representation of $\gg$ with
the shifted highest weight $\val$. However, as explained in\cx{Wh}, the
wheeling map allows one to get the expansion\rx{5.1} straight from the weight
$\pDval$ without going through PBW map and calculating the trace:
\qq
\JvaK = \dval \lrbc{
1 + \smnog \yp(\,\IO_{m,n}(\cK),\val)\,\hb^{m+n} },
\label{5.4}
\qqq
This is the formula that we will work with, because the weight function
$\pDval$ is easy to transfer from $\cB$ to $\cD$.
The inverse of the dual isomorphism map $\hxA^*$ maps the weight
system $\yp\in \cBs$ into an element of $\cDs$, which we will call
$\ypD$. In order to see how $\ypD$ acts on $\cD$ we come back to the
calculation of $\pDval$ and modify it.

Suppose that $\gg$ has $2k$ roots $\l_1,\ldots,\l_{2k}$. Let us index them in
such a way that $\l_1,\ldots,\l_{\kg}$ are positive roots and
$\l_1,\ldots\l_{\rg}$ are simple roots, $\rg$ being the rank of $\gg$.

For a root $\l$ of $\gg$ let $\Prl$
denote the operator projecting $\gg$ onto the root space
$V_\l\subset \gg$. We also introduce an operator $\Ph$, projecting
$\gg$ onto $\gh$. Let us assign a root of $\gg$
or the Cartan subalgebra to each internal edge of $D$. Let $\tbfS$ be
a set of all such assignments. For an assignment $c\in\tbfS$ we
modify the contraction of indices of tensors $f$ in the following
way: if an internal edge carries an index $a$ at the beginning and
index $b$ at the end, then instead of contracting them (that is,
instead of setting $a=b$ and taking a sum over their values) we bring
in an extra factor $P_b^a$, where $P$ is the projector corresponding
to the subspace assigned to that edge by $c$, and then contract the
pairs of indices $a$ and $b$ independently. In other words, we
project Lie algebras $\gg$ flowing along the internal edges of $D$
onto root spaces and Cartan subalgebras. Let us denote the resulting
number as $\pcDval$. Since the sum of projectors $\Ph$ and $\Prl$ for
all roots $\l$ of $\gg$ is equal to the identity operator, then
\qq
\pDval = \sum_{c\in\tbfS} \pcDval.
\label{5.5}
\qqq

The sum in the \rhs of this equation can be simplified.
Since $\val\in \gh$, then
\qq
[\val,\vy] = (\val\cdot\l)\,\vy\qquad\mbox{if $\vy\in V_\l$},
\qquad [\val,\vy] = 0 \qquad\mbox{if $\vy\in\gh$}.
\label{5.6}
\qqq
Therefore, $\pcDval=0$ unless the following two conditions are
met. First, $c$ must assign the same projector to internal edges of
$D$ which correspond to the same edge of $\Dz$. Second, there is a
\emph{compatibility requirement} at every 3-valent vertex: Cartan
subalgebra can be assigned to at most one of its edges and the sum of
the roots on incoming edges is equal to the sum of the roots on
outgoing edges.
Thus we can replace
the set $\tbfS$ in \ex{5.5} with the set $\bfS$ of `compatible' assignments
whose elements assign
subspaces to the edges of $\Dz$ in such a way that the compatibility
condition is satisfied at all of its vertices.

Equations\rx{5.6} also indicate that the effect of leg contractions
is easy to take into account in the calculation of $\pcDval$. If a leg is
attached to at least one edge, to which a Cartan subalgebra is
assigned, then $\pcDval=0$. Otherwise, if $m_j$ legs are attached on
the left side of an oriented edge $e_j$ of $\Dz$ to which a root
$\l$ is assigned, then they contribute a factor of
$(\val\cdot\l)^{m_j}$. Let $\lcj$ denote the root of $\gg$ assigned
by $c\in\bfS$ to the edge $e_j$ of $\Dz$. If $c$ assigns $\gh$ to
$e_j$, then we set $\lcj=0$. With these notations we see that
\qq
\pcDval = \pcDz \pjoN (\val\cdot\lcj)^{m_j},
\label{5.7}
\qqq
where $\pcDz = \pcDzval$ (we had to introduce this new notation
because the graph $\Dz$ has no legs and as a result $\pcDzval$ does
not depend on $\val$). Note that in \ex{5.7} we adopted a convention
that $0^0=1$.

The isomorphism\rx{3.5} completes the translation of $\pcDval$ into
the language of 3-valent graphs. For an assignment $c\in\bfS$
consider a linear combination of edges
\qq
\ecval = \sjoN (\val\cdot\lcj)\,e_j\in \Co.
\label{5.8}
\qqq
According to the compatibility condition satisfied by $c$,
$\ecval\in\ker\xdel=\oHor{D}$. Therefore, we can evaluate an element
$x\in \SsHor{D}$ on $\ecval$ and get a number (or a formal series)
$x(\ecval)$. Equations\rx{3.5},\rx{5.7} and\rx{5.8} indicate that for
an
element $x\in\cchBmdpz/\tcBihxri{1}(D_0,m)$,
\qq
\pcxval = \pcDz\,(\hxA\,x)(\ecval).
\label{5.9}
\qqq
Then, according to \ex{5.5}, after taking a sum over the assignments
of $\bfS$, we come to the following relation: for any $x\in\tcBmdz$,
\qq
\pxval = \pDxval,
\label{5.10}
\qqq
where
\qq
\pDyval=\sum_{c\in\bfS}\pcDz\, y(\ecval),
\qquad y \in \Hsio{m}{D_0}.
\label{5.11}
\qqq
Thus \ex{5.11} defines the element $\ypD\in\cDs$ corresponding to
$\yp\in\cBs$.

Applying \ex{5.10} to \ex{5.4}, we find that
\qq
\JvaK = \dval \lrbc{
1 + \smnog \ypD(\,\ID_{m,n}(\cK),\val)\,\hb^{m+n} }.
\label{5.12}
\qqq
It is easy to see that the weight system $\ypD$ behaves nicely under
the multiplication of elements of $\cD$:
$\ypD(xy,\val) = \ypD(x,\val)\,\ypD(y,\val)$ for any $x,y\in\cD$.
Therefore the analog of \ex{5.12} holds for the modified
integral\rx{4.3*}
\qq
\lgdval
= \smnog \ypD(\,\Il_{m,n}(\cK),\val)\,\hb^{m+n},
\label{5.13}
\qqq
and for its representative\rx{4.4} in the space $\tcD$
\qq
\lgdval = \sD \smzi \ypD(\,x_m(\cK,D),\val)\, \hb^{\chi(D)+m}.
\label{5.14}
\qqq
%
By using the formula\rx{5.11} for the
weight system, we can rewrite \ex{5.14} as
\qq
\lgdval = \sD \sum_{c\in\bfS}
\smzi \pcD\,x_m(\cK,D)(e_{c,\val})\, \hb^{\chi(D)+m},
\label{5.15}
\qqq
where $x_m(\cK,D)(e_{c,\val})$ denotes the evaluation of the element
$x_m(\cK,D)\in \Hsm{D}$ on $e_{c,\val}\in \oHor{D}$. According to
\ex{5.8}, $e_{c,\val}$ is a linear function of $\val$, while
$x_m(\cK,D)(e_{c,\val})$ is the homogeneous polynomial of
$e_{c,\val}$
of degree $m$. Therefore, \ex{5.15} can be further modified as
\qq
\lgdval & = & \sD \hb^{\chi(D)}\sass
\pcD\smzi x_m(\cK,D)(e_{c,\hb\val})
\nonumber\\
& = &
\sD \hb^{\chi(D)}\sum_{c\in\bfS}
\pcD\Il(\cK,D) (e_{c,\hb\val})
\label{5.16}
\qqq
the last line coming from \ex{4.5*1}.
Since by the definition of the
dual basis $f_j(e_i) = \delta_{ij}$, then according to \ex{5.8},
$f_j(e_{c,\hb\val}) = \hb\,(\val\cdot\lcj)$ and as a result, in view
of\rx{5.1*},
\qq
\Il(\cK,D) (e_{c,\hb\val}) =
\IlKi{D}{\aockh}
\label{5*.1}
\qqq
(see \ex{4.5*2}). Equation\rx{4.6} allows us to write the contribution of the
`1-loop' graph ($\chi(D)=0$) explicitly. Assignments
$c$ simply put different roots
on the circle,
$\pccir=1$ and
\qq
\lefteqn{
\sass \pccir \IlKi{\crcl}{\hb(\val\cdot\lci{1})}
}
\label{5*.2}
\\
&&
\hspace{1in}
=
\sum_{j=1}^k \log \lrbc{ q^{(\val\cdot\l_j)/2} -
q^{-(\val\cdot\l_j)/2}
\over \hb \,(\val\cdot\l_j)}
- \sum_{j=1}^k \log \APbas{\cK}{q^{\val\cdot\l_j}}.
\nonumber
\qqq
Thus if we exponentiate both sides of \ex{5.16} and use the
formulas\rx{5*.1},\rx{5*.2} and the dimension formula
\qq
d_{\val} = \pjok { \val\cdot \l_j \over \vr \cdot \l_j},
\qqq
then we find that
\qq
\JvaK  =
{\dqval \over \Dg(\cK;\qaok)}\;\;\Cqg \;
\exp\lrbc{\snoi \JlKn\,\hb^n }
\label{5*.3}
\qqq
where
\qq
&\JlKn = \sDln \;\;\sass\pcD\,\IlKi{D}{\aockh},
\label{5*.3*1}
\qqq
while
\qq
&\dqval = \pjok {q^{(\val\cdot\l_j/2)} - q^{-(\val\cdot\l_j)/2}
\over q^{(\vr\cdot\l_j/2)} - q^{-(\vr\cdot\l_j)/2} }
\label{5.20*2}
\qqq
is called the quantum dimension of the $\gg$-module with highest weight
$\val-\vr$ and
\qq
&\Cqg = \sum_{j=1}^k \log \lrbc{ q^{(\vr\cdot\l_j)/2} -
q^{-(\vr\cdot\l_j)/2}
\over \hb \,(\val\cdot\l_j)} = 1 + \cO(\hb^2),
\label{5.20*4}\\
&\Dg(\cK;\qaok) = \pjok \APbas{\cK}{q^{\val\cdot\l_k}}.
\label{5.20*5}
\qqq

Now let us apply Conjecture\rw{c4.1} to the \rhs of \ex{5*.3*1}. According to
\eex{4.5} and\rx{4.5*2},
\qq
\IlKi{D}{\aockh}
=
{
p(\cK,D;q^{\val\cdot\lci{1}},\ldots,q^{\val\cdot\lci{\chi(D)+1}})
\over \pjochD \APbas{\cK}{q^{\val\cdot\lcj}} }
\qqq
Therefore if we bring all terms in the sums of \ex{5*.3*1} to the common
denominator $$\Dg^{3n}(\cK;\qaok),$$ then we find that $\JlKn$ has a rational
form
\qq
\JlKn = {\plgn(\cK;\qaok)\over \Dg^{3n}(\cK;\qaok)},\qquad
\plgn(\cK;\tok)\in\IQtorg.
\label{5*.4}
\qqq
Then substituting this formula to \ex{5*.3}, exponentiating the formal power
series and expanding $\Cqg$ in powers of $\hb$ we come to the following
\begin{crl}
For a knot $\cK$ and a simple algebra $\gg$ there exist the
polynomials
\qq
&p_n(\cK;t_1,\ldots,t_{\rg}) \in \IQtorg,
&\qquad n\geq 0,
\label{5.19*}
\qqq
such that
\qq
\JvaK & = &
{\dqval \over \Dg(\cK;\qaok)}
\lrbc{1+
\snoi   {p_n(\cK;\qaok) \over \Dg^{3n}(\cK;\qaok) }\,\hb^n
}.
\label{5.20*}
\qqq
\end{crl}

We can check this prediction for the case of $\gg=su(2)$. In fact, in this case
the power of $\Dg$ in denominators\rx{5.20*} can be reduced. Indeed, the algebra
$su(2)$ has only one positive root. As a result, the elements of
$\bfS$ assign the subspaces of $su(2)$ to the edges of a graph $D$
in such a way that for any three edges attached to the same vertex,
two are assigned a root space and the third is assigned the Cartan
subalgebra.
Therefore, of $3\chi(D)$ edges that a graph $D$ has, $\chi(D)$ edges always
carry a Cartan subalgebra and only $2\chi(D)$ edges carry the root spaces.
%
Therefore, in case of $su(2)$ \ex{5.20*} is reduced to
\qq
\JaK = {[\a]\over \APKqa}
\lrbc{ 1 +
\snoi
{p_n(\cK;\qa) \over \APKqan }\,\hb^n
},
\label{5.22}
\qqq
where $\a$ is the dimension of the $su(2)$ module attached to the knot $\cK$
and
\qq
[\a] = {q^{\a/2} - q^{-\a/2}\over q^{1/2} - q^{-1/2} }
\label{5.25}
\qqq
is its quantum dimension.


Equation\rx{5.22} can be verified directly. We proved in\cx{Ro} that
for a knot $\cK$ in $S^3$ there exist the polynomials
\qq
\PnKt \in \ZZ[t^{\pm 1}],\;\;n\geq 1,
\label{5.23}
\qqq
such that the expansion\rx{5.1} can be rewritten as
\qq
\JaK = {[\a]\over \APKqa}
\lrbc{1+
\snoi {\PnKqa \over \APKqan}\,h^n,
},
\label{5.24}
\qqq
where
\qq
h = q-1 = e^\hb - 1.
\label{5.25*}
\qqq
It is easy to see that \ex{5.22} follows easily from \ex{5.24}.

\nsection{2-loop invariant and the $SU(3)$ colored Jones polynomial}
\label{s6}

Let us describe more precisely the implications of Conjecture\rw{c4.1} for the
value of Kontsevich integral at the level of `2-loop' graphs, \ie the graphs
with $\chi(D)=1$. There are only 2 such connected graphs in $\bDc$:
the theta-graph $D_1$
and the dumbbell $D_2$ of \fg{f3}.
\begin{figure}[hbt]
\leavevmode \centerline{
\hspace*{1.1in}
\epsfbox{pct5.eps}} \caption{The 2-loop
graphs $D_1$, $D_2$ and $D_3$}
\label{f3}
\end{figure}
Therefore, we can present the 2-loop part of the Kontsevich integral\rx{4.3*}
as
\qq
\smzi\Il_{m,1}(\cK) = \Il(\cK,D_1;\fDot{1}) +
\Il(\cK,D_2;\fDot{2})
\label{6.1}
\qqq
(\cf \eex{4.4},\rx{4.5*1} and\rx{4.5*2}), where we used a notation $\fD{i}{j}$
instead of simply $f_i$ in order to distinguish the dual edges coming from
different graphs. The formal power series in the \rhs of \ex{6.1} are not
themselves the invariants of $\cK$. They become the invariants only after the
factorization over the subspace $\tcDIHXx$ (see Theorem\rw{t2.1} and preceding
discussion). Let us describe the IHX indeterminacy in these power series more
precisely. The graph $D_3$ of Fig.\rw{f3} is the only connected 2-loop
graph with a 4-valent vertex. Applying the operator $\hdihx$ of\rx{1.22} to an
element $z(\fDot{3})\in\SHor{D_3}$ we get
\qq
\lefteqn{
{2\over 3}[z(\fDot{1}) + z(\fD{2}{1},-\fD{1}{1}-\fD{2}{1})
+ z(-\fD{1}{1}-\fD{2}{1},\fD{1}{1}) ]
-z(\fDot{2})
}&&
\nonumber\\
&&\hspace{4in}\in\bigoplus_{i=1}^2 \,\Invgdix{i}{\SHor{D_i}}.
\label{6.2}
\qqq
%
In this formula we assumed for simplicity of notation that
$z(x_1,x_2)\in\IQ[[x_1,x_2]]$ already has the symmetries
\qq
z(x_1,x_2) = z(x_2,x_1) = -z(-x_1,x_2),
\label{6.3}
\qqq
which makes the additional symmetrization of the expression\rx{6.2}
unnecessary. Expression\rx{6.2} indicates that by using the IHX freedom we can
bring the expression\rx{6.1} to the form
\qq
&\smzi\Il_{m,1}(\cK) = \zthK{\fDot{1}}\in\Invgdix{D_1}{\SHor{D_1}},
\label{6.4}
\qqq
where
\qq
\lefteqn{
\zthK{x_1,x_2} = \Il(\cK,D_1;x_1,x_2)
}
\label{6.5}\\
&& +
{2\over 3}[
\Il(\cK,D_2;x_1,x_2) + \Il(\cK,D_2;x_2,-x_1-x_2) +
\Il(\cK,D_2;-x_1-x_2,x_1)],
\nonumber
\qqq
thus eliminating the graph $D_2$ from Kontsevich integral. At the same time,
expression\rx{6.2} shows that $\zthK{x_1,x_2}$ of \ex{6.4} is the
$IHX$-invariant combination and therefore it is the only 2-loop invariant of
$\cK$.

The rationality conjecture implies that $\zthK{x_1,x_2}$ also has a
rational structure. Indeed, according to the conjecture, one can use the IHX
freedom in order to bring the terms in the \rhs of \ex{6.1} to the following
form:
\qq
\Il(\cK,D_1;\xot) & = &
{p(\cK,D_1;\exot) \over \APKex{1}\APKex{2}\APKe{x_1+x_2}}
\nonumber\\
\Il(\cK,D_2;\xot) & = & {p(\cK,D_2;\exot) \over \APKex{1}\APKex{2}}.
\label{6.6}
\qqq
Then according to \ex{6.5}, $\zthKx$ has a form
\qq
&\zthKx = {\pthKex\over \APKex{1}\APKex{2}\APKe{x_1+x_2}},
\label{6.7}
\qqq
where
the polynomial $\pthKt\in \IQtot$
is an invariant of $\cK$. Both this polynomial and a rational function
\qq
\zthsKt = {\pthKt\over \APKti{1}\APKti{2}\APK{t_1 t_2} }
\label{6.7*}
\qqq
have the symmetries
%
\qq
f(\tot) = f(t_2,t_1) = f( (t_1 t_2)^{-1},t_2) = f(t_1^{-1},t_2^{-1})
\label{6.8}
\qqq
implied by
the symmetry group $\Sgri{D_1}$. Finally, we rewrite \ex{6.4} with
the help of \ex{6.7}
\qq
\smzi\Il_{m,1}(\cK) =
\zthK{\fDot{1}} =
{\pthKef \over \APK{\efD{1}}\APK{\efD{2}}\APK{e^{-f_{1,D_1}-f_{2,D_1}}}  }.
\label{6.8*}
\qqq

It is easy to see from its definition that Kontsevich integral\rx{4.1} does not
contain (1,3)-valent graphs without legs. The wheeling map $\hOm$ produces such
graphs, however their Euler characteristic is at least 2. Therefore,
$\Il_{0,1}(\cK)=0$ in\rx{4.3*} and this means that
\qq
\zthK{0,0} = \zthsK{1,1} = 0.
\label{6.8*1}
\qqq

The polynomial $\pthKt$ can be extracted from the colored $SU(3)$ Jones
polynomial as described in Section\rw{s5}. In\cx{Rop} we will prove
a slightly strengthened version of \ex{5.20*}
for the groups $SU(n)$:
\qq
&
\JvaK  =
{\dqval\over \Dg(\cK;\qaok)}
\lrbc{1 +
\snoi   {P_n(\cK;\qaok) \over \Dg^{3n}(\cK;\qaok) }\,h^n
},
\label{6.9*}
\qqq
\qq
&
P_n(\cK;\tok) \in \ZZtok
\nonumber
\qqq
(note that here we used an expansion parameter $h=e^{\hb}-1$ instead of $\hb$
and as a result obtained the polynomials with \emph{integer} coefficients). For
the case of $SU(3)$ this formula implies that
\qq
&
\JvaK = {\dqval \over \Dg(\cK;\qaot)}\lrbc{1 + \hb\,
\FoK{\qaot} + \cO(\hb^2)},
\label{6.10*}
\qqq
where
\qq
& \FoK{\tot} =
{P_1(\cK;\tot)\over
\lrbs{
\APKti{1}\APKti{2}\APK{t_1 t_2}
}^3
}.
\label{6.11*}
\qqq
As we explained in Section\rw{s5}, a similar formula\rx{5*.3} can be obtained
by applying the $su(3)$ weight system to the logarithm of the Kontsevich
integral of $\cK$. Comparing \eex{6.10*} and\rx{5*.3}
and taking into account that $\Cqg = 1 + \cO(\hb^2)$,
we see that
\qq
\FoK{\qaot} = \JlKo.
\label{6.12*}
\qqq
Equations\rx{5*.3*1},\rx{6.7*} and\rx{6.8*} show that
\qq
\JlKo = \sass\pcDo \zthsK{\qaotc}
\label{6.13*}
\qqq
and the sum in this formula goes over the compatible assignments of root spaces
and Cartan subalgebra to the egdes of the
$\theta$-shaped graph $D_1$. There are two types of
such assignments. The first
one assigns two opposite roots to two edges and Cartan subalgebra to the third
edge, so
$\pcDo=2$. There are 3 choices of pairs of roots, and within each choice
there are 6 distinct assignments which give the same contributions due to the
symmetries\rx{6.8}. Therefore, the total contribution of the first assignment
to the \rhs of \ex{6.13*}
is
\qq
12
\Big(
\zthsK{\qaic{1},1} + \zthsK{\qaic{2},1} +
\zthK{\qalot,1 }
\Big).
\label{6.14*}
\qqq
Assignments of the second type put 3 different roots on the edges of $D_1$, so
$\pcDo=1$ There are 2 choices of compatible triplets of roots, and there are 6
ways to assign each triplet to the egdes of $D_1$. Thus we have 12 assignments
of the second type, and each gives the same contribution due to the
symmetries\rx{6.8}. Therefore, the total contribution of the second assignment
to the \rhs of \ex{6.13*} is
\qq
12\,\zthsK{\qaot}.
\label{6.15*}
\qqq
Thus putting the sum of\rx{6.14*} and\rx{6.15*} in the \rhs of \ex{6.13*} we
find from \ex{6.12*} that
\qq
\FoK{\tot} =
12\,\Big(
\zthsK{t_1,1} + \zthsK{t_2,1} + \zthsK{(t_1 t_2)^{-1},1} + \zthsKt
\Big).
\label{6.16*}
\qqq
It is easy to solve this equation for $\zthsKt$. By setting $t_2=1$ and using
\ex{6.8*1} and the symmetries\rx{6.8} we get
\qq
\FoK{t_1,1} = 36 \zthsK{t_1,1},
\label{6.17*}
\qqq
hence
\qq
\zthsKt = {1\over 36}\;
\Big(
3\FoK{\tot} - \FoK{t_1,1} - \FoK{t_2,1} - \FoK{(t_1 t_2)^{-1},1}
\Big).
\label{6.18*}
\qqq
In\cx{Rop} we will present a relatively efficient way of calculating
$\FoK{\tot}$. We have already written a Maple V program\cx{RoM}
which implements this
algorithm. For a knot presented as a cyclic closure of a braid, this program
calculates $\APKt$, $P_1(\cK;\tot)$ of \ex{6.11*} and then it finds $\pthKt$
through \ex{6.18*}.

\nsection{Discussion}
\label{s7}

Since the first version of this paper was written, A.~Kricker\cx{Kr} has proved
Conjecture\rw{c4.1}. In fact, he proved it for a more general case of knots in
integer homology spheres, where an analog of Kontsevich integral for knots is
defined with the help of the LMO invariant\cx{LMO} or its \AA rhus
version\cx{Aa}. This knot invariant lies in the same space $\cB$, so the
previous discussion equally applies in that case. The analog of the colored
Jones polynomial is the so-called \emph{trivial connection contribution to the
colored Jones polynomial} defined for $SU(2)$ in\cx{RoII} for knots in rational
homology spheres. It also has a rational structure\rx{5.24}.

Naturally, one wants to extend the rationality conjecture to the most
general case of links in rational homology spheres. Unfortunately, Kricker's
proof works only for integer homology spheres, so it can not be generalized
easily to homologically non-trivial knots in rational homology spheres, for
which the analog of the rationality conjecture can be easily formulated in
accordance with the $SU(2)$ results of\cx{RoII} -- one just has to use
fractional exponents $e^{f_i/h(\cK)}$ in \ex{4.5*2*2}, where $h(\cK)$ is the
order of the homology element represented by the knot $\cK$.

Generalizing the rationality conjecture to links is not a straightforward
exercise, because the arguments of Section\rw{s3} hinge upon Lemma\rw{l3.1}.
For this lemma to work, the legs of a (1,3)-valent graph have to be
interchangeable (or, in other words, `commutative'). In case of links however,
legs are attached to different components, and as a result, one may have a
non-zero graph in $\cB$ which has two legs attached to the same 3-valent
vertex, if they come from different link components.

S.~Garoufalidis and A.~Kricker\cx{GK} have circumvented this difficulty in
the case of boundary links and proved an analog of the rationality property.
However, the rationality property of the $SU(2)$ colored Jones polynomial of
links described in\cx{RoII} suggests a different approach. Namely, there is a
graph space map which sends Kontsevich integral of a link into a close relative
of the space $\cD$. This map is similar to the \AA rhus map\cx{Aa}. It
implements diagrammatically the stationary phase integration performed
in\cx{RoII} and ultimately it `makes' all legs commutative. Similarly to the
\AA rhus map, one would have to prove that the image of the map is a link
invariant, and this is work in progress.

Despite the fact that polynomials\rx{4.5*3} share the variables $\tbs$ with the
Alexander polynomial, their topological interpretation remains unclear. First
of all, because of the $IHX$ indeterminacy, the rational expressions\rx{4.5}
are not knot invariants. Only their linear combinations which are insensitive
to the $IHX$ transformations are true invariants of knots. We explained this
point in details in Section\rw{s6} for 2-loop graphs. In that case we presented
an explicit linear combination\rx{6.5} which is invariant and which yields a
2-loop invariant polynomial $\pthKt$. So just as a beginning, it would be
interesting to establish its topological interpretation.

In the framework of the quantum Chern-Simons field theory and in the framework
of the theory of finite type (Vassiliev) invariants, $\zthKx$ and the polynomial
$\pthKt$ are analogs of the Casson-Walker invariant of rational homology
spheres, so one might try to related $\pthKt$ to the moduli space of flat
connections in the knot complement for an appropriate Lie group. At a simpler
(`1-loop') level the Alexander polynomial $\APKt$ is the analog of the order of
integer homology $\HoZ{M}$ of a rational homology sphere. The order of
$\HoZ{M}$ is equal to the number of flat $U(1)$ connections on $M$. At the same
time, at least for fibered knots, $\APKt$ is related to the monodromy map
acting on the moduli space $\mduos$ of flat $U(1)$ connections on the Seifert
surface $\Sigma$ of $\cK$. Namely, the monodromy map
$f:\;\Sigma\longrightarrow \Sigma$ which defines the structure of a fiber
bundle for the knot complement $S^3\setminus\cK$, acts also on $\mduos$ and
\qq
\APKt = \sum_{n=0}^{2g(\Sigma)}(-1)^n t^{n-g(\Sigma)}\Tr_{H^n(\mduos)} \fs,
\label{7.1}
\qqq
where $\fs$ denotes the action of $f$ on $H^n(\mduos)$. Since the
Casson-Walker invariant `counts' the number of flat $SU(2)$ connections on a
rational homology sphere, then one might expect that $\pthKt$ can be expressed
somehow similarly to \ex{7.1} through the action of the monodromy $f$ on moduli
spaces of flat connections of other Lie groups. Unfortunately, no
such interpretation exists at present.

\noindent
{\bf Acknowledgements.}
I am very thankful to D.~Thurston and A.~Vaintrob for discussing this
work. I am especially indebted to A.~Vaintrob for numerous
discussions of the properties of the space $\cD$ and to D.~Thurston for substantive
discussions of the conjecture and for explaining the effects of
unwheeling procedure at the diagrammatic level.
This work was supported by NSF Grants
DMS-0196235 and DMS-0196131.

\appendix
\nsection{The 2-loop polynomial $\pthKt$ for knots with up to 8 crossing}

\begin{table}
\begin{center}
\begin{tabular}{||c|l|l||} \hline\hline
Knot & $\APKt$ & $12\pthKt$ \\
\hline\hline
$3_1$ & $t-1$ & $-t_1^2 t_2 + t_1^2$
\\ \hline
$4_1$ & $t^2 - 3t + 5$ & $0$
\\ \hline
$5_1$ & $t^2 - t + 1$ &
$2t_1^4 t_2^2 - 2t_1^4 t_2 + 2 t_1^4 - t_1^2 t_2 + t_1^2$
\\ \hline
$5_2$ & $2t - 3$ & $-13 t_1^2 t_2 + 9 t_1^2 + 6 t_1 - 12$
\\ \hline
$6_1$ & $-2t+5$ & $3t_1^2 t_2 - t_1^2 - 6t_1 + 24$
\\ \hline
$6_2$ & $-t^2 + 3t - 3$ & $-3t_1^4 t_2^2 +
3t_1^4 t_2 - t_1^4 - 6 t_1^3 - 11 t_1^2 t_2 + 15t_1^2$
\\ \hline
$6_3$ & $t^2 - 3t + 5$ & $0$
\\ \hline
$7_1$ & $ t^3 - t^2 + t - 1$ &
$-3t_1^6 t_2^3 + 3t_1^6 t_2^2 - 3t_1^6 t_2 + 2t_1^4 t_2^2 + 3t_1^6
 - 2 t_1^4 t_2 + 2 t_1^4 - t_1^2 t_2 + t_1^2$
\\ \hline
$7_2$ & $3t - 5$ & $-58t_1^2 t_2+ 36 t_1^2 + 36 t_1 - 84$
\\ \hline
$7_3$ & $2 t^2 - 3 t + 3$ & $-25 t_1^4 t_2^2 + 25 t_1^4 t_2 - 17 t_1^4 + 7
t_1^2 t_2 - 12 t_1^3 + t_1^2 - 6 t_1 + 12$
\\ \hline
$7_4$ & $4t-7$ & $ 136 t_1^2 t_2 - 80 t_1^2 - 96 t_1 + 240$
\\ \hline
$7_5$ & $2t^2 - 4t + 5$ &
$41 t_1^4 t_2^2 - 33 t_1^4 t_2 - 16 t_1^3 t_2 + 17 t_1^4 + 12 t_1^2 t_2
+32 t_1^3 + 4 t_1^2 - 14 t_1 + 36$
\\ \hline
$7_6$ & $-t^2 + 5t - 7$ &
$-7 t_1^4 t_2^2 + 5 t_1^4 t_2 + 10 t_1^3 t_2 - t_1^4 - 20 t_1^3
- 98 t_1^2 t_2  + 64 t_1^2 + 50 t_1 - 108$
\\ \hline
$7_7$ & $t^2 - 5t + 9$ & $ - 5 t_1^2 t_2 + t_1^2 + 12 t_1 -48 $
\\ \hline
$8_1$ & $-3t + 7$ &
$23 t_1^2 t_2 - 9 t_1^2 - 36 t_1 + 132$
\\ \hline
$8_2$ & $-t^3 + 3t^2 - 3t$ &
$6t_1^6 t_2^3 - 6 t_1^6 t_2^2 + 6 t_1^6 t_2 + 20 t_1^4 t_2^2 - 2 t_1^6
-20 t_1^4 t_2 - 12 t_1^5 + 30 t_1^4$
\\
& $+ 3$& $ -11 t_1^2 t_2 - 6 t_1^3 + 15 t_1^2$
\\ \hline
$8_3$ & $ -4t + 9$ & $ 0 $
\\ \hline
$8_4$ & $-t^2 + 3t - 3$ &
$-3t_1^4 t_2^2 + 3 t_1^4 t_2 - t_1^4 - 11 t_1^2 t_2 - 6 t_1^3 + 15 t_1^2$
\\ \hline
$8_5$ & $-t^3 + 3t^2 - 4t$ &
$-10 t_1^6 t_2^3 + 8 t_1^6 t_2^2 +6 t_1^5 t_2^2 - 6 t_1^6 t_2 - 29 t_1^4 t_2^2
- 6 t_1^5 t_2 + 2 t_1^6 + 12 t_1^4 t_2$
\\
& $+5$ &
$ + 12 t_1^5 + 13 t_1^3 t_2 - 15 t_1^4 + 15 t_1^2 t_2 + 6 t_1^3
 - 43 t_1^2 + 16 t_1 - 48$
\\ \hline
$8_6$ & $-2t^2 + 6t - 7$ &
$-31 t_1^4 t_2^2 + 27 t_1^4 t_2 + 12 t_1^3 t_2 - 9 t_1^4 - 111 t_1^2 t_2
- 54 t_1^3 + 111 t_1^2 + 18 t_1$
\\
& & $-12$
\\ \hline
\end{tabular}
\end{center}
\caption{The Alexander polynomial $\APKt$ and the
2-loop polynomial $12\pthKt$ presented by monomials in funamental
domains}
\label{t6.1}
\end{table}

\begin{table}
\begin{center}
\begin{tabular}{||c|l|l||}
\hline\hline
Knot & $\APKt$ & $12\pthKt$
\\ \hline\hline
$8_7$ & $ t^3 - 3t^2 + 5t - 5$ &
$5t_1^6 t_2^3 - 5 t_1^6 t_2^2 + 3 t_1^6 t_2 - t_1^4 t_2^2 + 6 t_1^5 t_2 - t_1^6
- 7 t_1^4 t_2 - 6 t_1^5$
\\
&  & $ + 4 t_1^3 t_2 - 3 t_1^4 + 19 t_1^2 t_2 + 16 t_1^3 - 31 t_1^2$
\\ \hline
$8_8$ & $2t^2 - 6t + 9$ &
$-5t_1^4 t_2^2 + 3 t_1^4 t_2 + 6 t_1^3 t_2 - t_1^4 - 5t_1^2 t_2 - 6 t_1^3
- 9 t_1^2 + 18 t_1 - 60$
\\ \hline
$8_9$ & $-t^3 + 3t^2 - 5t+7$ &
$0$
\\ \hline
$8_{10}$ & $t^3 - 3t^2 + 6t-7$ &
$7t_1^6 t_2^3 - 6t_1^6 t_2^2 - 3t_1^5 t_2^2 + 3t_1^6 t_2 - 2t_1^4 t_2^2
+ 9 t_1^5 t_2 - t_1^6 + 2 t_1^4 t_2$
\\
&  &
$ - 6 t_1^5- 5 t_1^3 t_2 - 14 t_1^4
+ 48 t_1^2 t_2 + 20 t_1^3 - 40 t_1^2 - 18 t_1 + 36$
\\ \hline
$8_{11}$ & $-2t^2 + 7t - 9$ &
$ - 39 t_1^4 t_2^2 + 31 t_1^4 t_2 + 28 t_1^3 t_2 - 9 t_1^4 - 206 t_1^2 t_2
- 76 t_1^3 + 160 t_1^2$
\\
& & $ + 74 t_1 - 132$
\\ \hline
$8_{12}$ & $t^2 - 7t + 13$ & $0$
\\ \hline
$8_{13}$ & $2t^2 - 7t + 11$ &
$-5t_1^4 t_2^2 + 3 t_1^4 t_2 + 6 t_1^3 t_2 - t_1^4 - 7 t_1^2 t_2 - 6 t_1^3
-9t_1^2 + 24t_1 - 84$
\\ \hline
$8_{14}$ & $-2t^2 + 8t - 11$ &
$-47 t_1^4 t_2^2 + 35 t_1^4 t_2 + 48 t_1^3 t_2 - 9t_1^4 - 356 t_1^2 t_2
- 102 t_1^3 + 236 t_1^2$
\\
& & $ + 168 t_1-336$
\\ \hline
$8_{15}$ & $3t^2 - 8t+11$ &
$203 t_1^4 t_2 - 148 t_1^4 t_2 - 145 t_1^3 t_2 + 57 t_1^4 + 375 t_1^2 t_2
+ 240 t_1^3 $
\\
& & $- 111 t_1^2 - 304 t_1 + 756$
\\ \hline
$8_{16}$ & $t^3 - 4 t^2 + 8t - 9$ &
$-9t_1^6 t_2^3 + 8 t_1^6 t_2^2 + 4 t_1^5 t_2^2 - 4t_1^6 t_2 - 26 t_1^4 t_2^2
- 16 t_1^5 t_2 + t_1^6$
\\
& & $ + 30 t_1^4 t_2 + 12 t_1^5 + 7 t_1^3 t_2 + 6 t_1^4
- 90 t_1^2 t_2 - 66 t_1^3 + 106 t_1^2 + 6 t_1$
\\
& & $ + 12$
\\ \hline
$8_{17}$ & $-t^3 + 4t^2 - 8t + 11$ & $0$
\\ \hline
$8_{18}$ & $-t^3 + 5t^2 - 10t+13$ & $0$
\\ \hline \hline
\end{tabular}
\end{center}
\caption{Table 1 continued}
\label{t6.1*}
\end{table}

\begin{table}
\begin{center}
\begin{tabular}{||c|l|l||}
\hline\hline
Knot & $\APKt$ & $12\pthKt$
\\ \hline\hline
$3_1$ & $u-1$ &
$u_1^2 - 3u_2 - 2u_1-6$
\\ \hline
$4_1$ & $u^2 - 3u + 3$ & $0$
\\ \hline
$5_1$ & $u^2 - u - 1$ &
$2u_1^4 - 10 u_2 u_1^2 - 4u_1^3 + 10 u_2^2 + 10 u_2 u_1 - 23 u_1^2 + 53u_2
+ 26u_1$
\\
&& $ + 66$
\\ \hline
$5_2$ & $2u-3$ &
$9u_1^2 - 31 u_2 - 12 u_1 - 66$
\\ \hline
$6_1$ & $-2u+5$ &
$-u_1^2 + 5u_2 - 4u_1 + 30$
\\ \hline
$6_2$ & $-u^2 + 3u -1$ &
$-u_1^4 + 7u_2 u_1^2 - 11 u_2^2 - 5 u_2 u_1 + 31 u_1^2 - 73 u_2 - 34 u_1 - 114$
\\ \hline
$6_3$ & $u^2 - 3u + 3$ & $0$
\\ \hline
$7_1$ & $u^3 - u^2 - 2u + 1$ &
$3u_1^6 - 21 u_2 u_1^4 - 6 u_1^5 + 42 u_2^2 u_1^2 + 21 u_1 u_1^3 - 52 u_1^4
- 21 u_2^3$
\\
&& $ + 215 u_2 u_1^2 + 62 u_1^3 - 152 u_2^2 - 16 u_2 u_1 + 268 u_1^2
- 358 u_2$
\\
&& $ - 64 u_1 - 276$
\\ \hline
$7_2$ & $3u-5$ & $ 36 u_1^2 - 130 u_2 - 36 u_1 - 300$
\\ \hline
$7_3$ & $2u^2 - 3u - 1$ &
$- 17u_1^4 + 93 u_2 u_1^2 + 38 u_1^3 - 109 u_2^2 - 121 u_2 u_1 + 221 u_1^2
- 559 u_2$
\\
&& $ - 314 u_1 - 702$
\\ \hline
$7_4$ & $4u-7$ & $-80u_1^2 + 296 u_2 + 64 u_1 + 720$
\\ \hline
$7_5$ & $2u^2 - 4u + 1$ &
$17 u_1^4 - 101 u_2 u_1^2 - 50 u_1^3 + 141 u_2^2 + 165 u_2 u_1 - 200 u_1^2
+ 624 u_2$
\\
&& $ + 392 u_1 + 672$
\\ \hline
$7_6$ & $-u^2 + 5u - 5$ &
$-u_1^4 + 9 u_2 u_1^2 - 6u_1^3 - 19u_2^2 + 17 u_2 u_1 + 64 u_1^2 - 194 u_2
- 16 u_1$
\\
&& $ - 324$
\\ \hline
$7_7$ & $u^2 - 5u + 7$ &
$u_1^2 - 7u_2 + 10u_1 - 54$
\\ \hline
$8_1$ & $-3u+7$ &
$-ou_1^2 + 41u_2 - 18 u_1 + 186$
\\ \hline
$8_2$ & $-u^3 + 3u^2 - 3$ &
$-2u_1^6 + 18 u_2 u_1^4 - 48 u_2^2 u_1^2 - 6 u_2 u_1^3 + 66 u_1^4 + 34 u_2^3
+ 12 u_2^2 u_1$
\\
&& $- 314 u_2 u_1^2 - 54 u_1^3 + 324 u_2^2 + 114 u_2 u_1 - 463 u_1^2
+ 977 u_2$
\\
&& $ + 248 u_1 + 894$
\\ \hline\hline
\end{tabular}
\end{center}
\caption{The Alexander polynomial $\APKt$ and the
2-loop polynomial $12\pthKt$ expressed in terms of symmetric polynomials $u$
and $u_1$, $u_2$}
\label{t6.2}
\end{table}

\begin{table}
\begin{center}
\begin{tabular}{||c|l|l||}
\hline\hline
Knot & $\APKt$ & $12\pthKt$
\\ \hline\hline
$8_3$ & $-4u+9$ & $0$
\\ \hline
$8_4$ & $-u^2+3u-1$ &
$-u_1^4 + 7u_2u_1^2 - 11 u_2^2 - 5u_2 u_1 + 31 u_1^2 - 73 u_2
- 34 u_1 -114$
\\ \hline
$8_5$ & $-u^3 + 3u^2 - u -1$ &
$2u_1^6 - 18 u_2 u_1^4 - 4u_1^5 + 50 u_2^2 u_1^2 + 32 u_2 u_1^3 -59 u_1^4 -
42 u_2^3$
\\
&& $ - 52 u_2^2 u_1 + 288 u_2  u_1^2 + 132 u_1^3 - 359 u_2^2 - 341 u_2 u_1
+ 351 u_1^2$
\\
&& $ -983u_2 - 528 u_1 -834$
\\ \hline
$8_6$ & $-2u^2 + 6u -3$ &
$-9u_1^4 + 63 u_2 u_1^2 + 8 u_1^3 - 103 u_2^2 - 57 u_2 u_1 + 231 u_1^2
- 631 u_2$
\\
&& $ - 288 u_1 - 822$
\\ \hline
$8_7$ & $u^3 - 3u^2 + 2u +1$ &
$-u_1^6 + u_2 u_1^4 + 4 u_1^5 - 26 u_2^2 u_1^2 - 23 u_2 u_1^3 + 11 u_1^4
+ 23 u_2^3$
\\
&& $ + 34 u_2^2 u_1 - 94 u_2 u_1^2 - 34 u_1^3 + 157 u_2^2 + 133 u_2 u_1
- 92 u_1^2$
\\
&& $ + 346 u_2 + 112 u_1 + 276$
\\ \hline
$8_8$ & $2u^2 - 6u + 5$ &
$-u_1^4 + 7 u_2 u_1^2 + 4 u_1^3 - 13 u_2^2 - 11 u_2 u_1 - 5u_1^2 - 15 u_2
+ 8 u_1$
\\
&& $ + 6$
\\ \hline
$8_9$ & $-u^3 + 3u^2 - 2u + 1$ & $0$
\\ \hline
$8_{10}$ & $u^3 - 3u^2 + 3u - 1$ &
$-u_1^6 + 9 u_2 u_1^4 + 6 u_1^5 - 27 u_2^2 u_1^2 - 36 u_2 u_1^3 + 4 u_1^4
+ 27 u_2^3$
\\
&& $ + 54 u_2^2 u_1 - 62 u_2 u_1^2 - 58 u_1^3 + 152 u_2^2 + 198 u_2 u_1
- 9 u_1^2$
\\
&& $ + 259 u_2 + 156 u_1 + 138$
\\ \hline
$8_{11}$ & $-2u^2 + 7u - 5$ &
$-9u_1^4 + 67 u_2 u_1^2 + 2u_1^3 - 119 u_2^2 - 35 u_2 u_1 + 256 u_1^2
- 730 u_2$
\\
&& $ - 260 u_1 - 972$
\\ \hline
$8_{12}$ & $u^2 - 7u + 11$ & $0$
\\ \hline
$8_{13}$ & $2u^2 - 7u + 7$ &
$-u_1^4 + 7 u_2 u_1^2 + 4u_1^3 - 13 u_2^2 - 11 u_2 u_1 - 5u_1^2 - 17 u_2
+ 14 u_1$
\\
&& $ - 18$
\\ \hline
$8_{14}$ & $-2u^2 + 8u - 7$ &
$-9u_1^4 + 71u_2 u_1^2 - 8u_1^3 - 135 u_2^2 + 3u_2 u_1 + 300 u_1^2 - 916 u_2$
\\
&& $- 204 u_1 - 1272$
\\ \hline
$8_{15}$ & $3u^2 - 8u + 5$ &
$57u_1^4 -376 u_2 u_1^2 - 166 u_1^3 + 613u_2^2 + 569 u_2 u_1 - 687 u_1^2$
\\
&& $+ 2543u_2 +1258 u_1 + 2574$
\\ \hline\hline
\end{tabular}
\end{center}
\caption{Table 3 continued}
\label{t6.2*}
\end{table}

\begin{table}
\begin{center}
\begin{tabular}{||c|l|l||}
\hline\hline
Knot & $\APKt$ & $12\pthKt$
\\ \hline\hline
$8_{16}$ & $u^3 - 4u^2 + 5u - 1$ &
$u_1^6 - 10 u_2 u_1^4 - 4u_1^5 + 33 u_2^2 u_1^2 + 25 u_2 u_1^3 - 12 u_1^4
-35 u_2^3$
\\
&& $ - 43 u_2^2 u_1 + 128 u_2 u_1^2 + 22 u_1^3 - 244 u_2^2 - 158 u_2 u_1
+ 189 u_1^2$
\\
&& $ - 607 u_2 - 182 u_1 - 558$
\\ \hline
$8_{17}$ & $-u^3 + 4u^2 - 5u + 3$ & $0$
\\ \hline
$8_{18}$ & $-u^3 + 5u^2 - 7u + 3$ & $0$
\\ \hline\hline
\end{tabular}
\end{center}
\caption{Table 3 continued}
\label{t6.2*1}
\end{table}


Here are the results of calculating the polynomials $\pthKt$ for the first few
knots (with up to 8 crossings). We present these results in two different ways.
First, as we know, $\pth(\cK)\in\Qalghzo$ and relations\rx{6.8} come from the
symmetry $\Sgri{D_1}$. More explicitly, $\Hoz{D_1}$ looks like $su(3)$ root
lattice with elements $f_1$ and $f_2$ (and variables $t_1$, $t_2$)
corresponding to the simple roots (see~\fg{fx}).
\begin{figure}
\vspace{0.2in}
\input tstpct5.pst
\caption{Fundamental domain of the symmetry\rx{6.8}}
\label{fx}
\end{figure}
The
symmetry group $\Sgri{D_1}$ is the symmetry of this lattice (which preserves
the origin). So instead of writing the whole polynomial $\pthKt$ we may list
just the monomials belonging to a fundamental domain of $\Sgri{D_1}$. From our
$su(3)$ lattice description it is easy to see that we may choose a fundamental
domain to include the monomials
\qq
t_1^{m_1} t_2^{m_2},\qquad m_1,m_2\geq 0,\; m_1 \geq 2 m_2.
\label{6.20*}
\qqq
Then the other monomials will be determined by the symmetries\rx{6.8}.
Similarly, in view of the symmetry\rx{4.5*} it is enough to list only the
monomials of $\APKt$ with non-negative powers of $t$. Thus in \tb{t6.1} we
present the `fundamental domain' parts of the Alexander polynomial $\APKt$ and
(scaled) 2-loop polynomial $12\pthKt$.

An alternative way of describing $\pthKt$ comes from the observation that the ring
of Laurent polynomials with the symmetries\rx{6.8} can be written as
$\IQ[u_1,u_2]$, where
\qq
u_1(t_1,t_2) &=&
t_1 + t_1^{-1} + t_2 + t_2^{-1} + t_1 t_2 + t_1^{-1} t_2^{-1},
\nonumber\\
u_2(t_1,t_2) & = &
t_1^2 t_2 + t_1^{-2} t_2^{-1} + t_1 t_2^{2} + t_1^{-1} t_2^{-2}
+ t_1 t_2^{-1} + t_1^{-1} t_2.
\label{6.21*}
\qqq
So in \tb{t6.2} we present the expressions for the Alexander polynomial $\APKt$
in terms of $u=t+t^{-1}$ and for the (scaled) 2-loop polynomial $12\pthKt$ in
terms of $u_1$ and $u_2$.

\begin{remark}
\rm
If $\cK\p$ is the mirror image of $\cK$, then
$\pth(\cK\p;\tot) = -\pthKt$, hence $\pthKt=0$ for amphicheiral knots.
\end{remark}
\begin{remark}
\rm
As we see, experimental evidence suggests that
\qq
12\pthKt\in\ZZtot.
\label{6.19*}
\qqq
\end{remark}
\begin{remark}
\rm
The degree of the Alexander polynomial is bounded by the genus of the knot
$g(\cK)$
\qq
\deg \APKt \leq g(\cK),
\label{6.22*}
\qqq
In view of the symmetries\rx{6.8} (which come from $\Sgri{D_1}$), the reasonable
measure of the degree of $\pthKt$ is the $t_1$ degree of its fundamental domain
part. Let us denote it simply as $\deg\pthKt$. Then \tb{t6.1} suggests a
similar inequality
\qq
\deg\pthKt\leq 2g(\cK).
\label{6.23*}
\qqq

\end{remark}

\end{document}

\qq
\APbas{\cK}{\exp(f_j)}(e_{c,\hb\val}) & = &
\APbas{\cK}{q^{\val\cdot\lcj}},
\label{5.17}
\\
 \pcD\,y(\cK,D)(e_{c,\hb\val}) & = &
\pcD p(\cK,D;q^{\val\cdot\lci{1}},\ldots,q^{\val\cdot\lci{\chi(D)+1}})
\nonumber\\
& = & p_c(\cK,D;q^{\val\cdot\l_1},\ldots,q^{\val\cdot\l_{\rg}})
\label{5.18}
\qqq
(\cf \ex{4.5*2}),
where $\l_1,\ldots,\l_{\rg}$ are all simple roots of $\gg$ and
\qq
p_c(\cK,D;t_1,\ldots,t_{\rg}) \in \IQ[t_1,\ldots,t_{\rg}]
\label{5.19}
\qqq
is a polynomial depending on the weight assignment $c$.

Applying the weight system to the contribution of the circle graph,
described by \ex{4.6}, is rather straightforward. We just have to
take a sum over assignments of the roots of $\gg$ to the circle. Note
that in view of \ex{4.5*}, an assignment of a positive root $\l$
yields the same number as the assignment of its opposite.

Thus we proved that Conjecture\rw{c4.1} leads to the following formula
\qq
\lgdval & = & \sum_{j=1}^k \log \lrbc{ q^{(\val\cdot\l_j)/2} -
q^{-(\val\cdot\l_j)/2}
\over \hb \,(\val\cdot\l_j)}
- \sum_{j=1}^k \log \APbas{\cK}{q^{\val\cdot\l_j}}
\label{5.20}
\\
&&\qquad
+
\sDl \hb^{\chi(D)}\sum_{c\in\bfS}
{p_c(\cK,D;q^{\val\cdot\l_1},\ldots,q^{\val\cdot\l_{\rg}})
\over \pjoN \APbas{\cK}{q^{\val\cdot\lcj}} }.
\nonumber
\qqq
By using the dimension formula
\qq
d_{\val} = \pjok { \val\cdot \l_j \over \vr \cdot \l_j},
\qqq
we can rewrite \ex{5.20} as
\qq
\lgdval & = &\log\lrbc{\dqval/ \dval} + \log\Cqg
- \log \Dg(\cK;\qaok)
\label{5.20*3}
\\
&&\qquad
+
\sDl \hb^{\chi(D)}\sum_{c\in\bfS}
{p_c(\cK,D;q^{\val\cdot\l_1},\ldots,q^{\val\cdot\l_{\rg}})
\over \pjoN \APbas{\cK}{q^{\val\cdot\lcj}} }.
\nonumber
\qqq
where
\qq
\dqval = \pjok {q^{(\val\cdot\l_j/2)} - q^{-(\val\cdot\l_j)/2}
\over q^{(\vr\cdot\l_j/2)} - q^{-(\vr\cdot\l_j)/2} }
\label{5.20*2}
\qqq
is called the quantum dimension of the $\gg$-module with highest weight
$\val-\vr$ and
\qq
&\Cqg = \sum_{j=1}^k \log \lrbc{ q^{(\vr\cdot\l_j)/2} -
q^{-(\vr\cdot\l_j)/2}
\over \hb \,(\val\cdot\l_j)} = 1 + \cO(\hb^2),
\label{5.20*4}\\
&\Dg(\cK;\qaok) = \pjok \APbas{\cK}{q^{\val\cdot\l_k}}.
\label{5.20*5}
\qqq

The formula\rx{5.20*3}
can be further simplified in order to adapt it to the actual expansion as it
comes from the colored Jones polynomial. Namely, we take into account an obvious
relation between the Euler characteristic and the
number of edges in a 3-valent graph: $N = 3 \chi(D)$, $\chi(D)\geq 1$.
Therefore, after bringing the terms in the \rhs of \ex{5.20} coming from the
graphs with the same value of $\chi(D)$ to a common denominator
$\Dg^{3\chi(D)}(\tok)$,
exponentiating both sides and expanding $\Cqg$ in powers of $\hb$
we come to the following:

\qq
&\hspace{-0.5in}
\lgdval =   \log({\dqval/\dval}) -
\log\Dg(\cK;\qaot)
+
F_1(\cK;\qaot)
\,\hb
+ \cO(\hb^2),
\label{6.10}\\
& F_1(\cK;\tot) =
{P_1(\cK;\tot)\over
\lrbs{
\APKti{1}\APKti{2}\APK{t_1 t_2}
}^3
}
\label{6.11}
\qqq
As we explained in Section\rw{s5}, the same formula\rx{6.10} can be obtained by
applying the $su(3)$ weight system to the logarithm of Kontsevich
integral\rx{4.3*} (see \ex{5.20*3}). Comparing \eex{5.20*3} and\rx{6.10} and
taking into account that $\log \Cqg = \cO(\hb^2)$, we find that
\qq
F_1(\cK;\qaot) =
\sDlo \sum_{c\in\bfS}
{p_c(\cK,D;\qaot)
\over \pjoN \APbas{\cK}{q^{\val\cdot\lcj}} }.
\label{6.12}
\qqq
In fact, as we know, there are only two connected 3-valent graphs. Moreover,
since we choose the value\rx{6.8*} for the representative of the 2-loop part of
Kontsevich integral in the space $\tcD$, then \ex{6.12} is reduced to
\qq
F_1(\cK;\qaot) =
\sum_{c\in\bfS}
{p_c(\cK,D_1;\qaot)
\over\pjoth \APbas{\cK}{q^{\val\cdot\lcj}} }.
\label{6.13}
\qqq
In this formula $\bfS$ is the set of compatible assignments of root
spaces and Cartan subalgebra of $su(3)$ to the edges of the graph $D_1$ (see
the definition after \ex{5.6}), while the polynomial $p_c(\cK,D_1;\tot)$ is
determined by the polynomial $\pthKt$ through the formula\rx{5.18} (we should
just substitute $\pthKt$ for
$p(\cK,D;q^{\val\cdot\lci{1}},\ldots,q^{\val\cdot\lci{\chi(D)+1}})$).

There are two types of compatible assignments for the graph $D_1$. The first

\qq
\pthKt = \pthK{t_2,t_1} = \pthK{(t_1 t_2)^{-1},t_2} =
\pthK{t_1^{-1},t_2^{-1}}
\label{6.8}
\qqq

********************************************

Since we are ultimately interested in $\zthsKt$, we want to use \ex{6.16*} in
order to express it in terms of $\FoK{\tot}$. Set

Thus adding up the contributions\rx{6.14*} and\rx{6.15*} we transform
\ex{6.13*} into
%
%

18 such assignments, and each of them gives the same
contribution to the \rhs of \ex{6.13} due to the symmetries\rx{6.8}. Since
$\pcDo=$, then all these assignments contribute

For $SU(3)$ this formula is
\qq
\JvaK & = & \lrbc{
\pjoth {q^{(\val\cdot\l_j/2)} - q^{-(\val\cdot\l_j)/2}
\over q^{1/2} - q^{-1/2}}
}
\snzi   {p_n(\cK,D;\qaok) \over \Dg^{3n+1}(\qaok) }\,\hb^n.
\label{6.9}
\qqq

\qq
\JvaK & = & \lrbc{
\pjok {q^{(\val\cdot\l_j/2)} - q^{-(\val\cdot\l_j)/2}
\over q^{(\vr\cdot\l_j/2)} - q^{-(\vr\cdot\l_j)/2} }
}
\snzi   {p_n(\cK;\qaok) \over \Dg^{3n+1}(\cK;\qaok) }\,\hb^n.
\nonumber
\qqq
\qq
&
\JvaK  =  \lrbc{
\pjok {q^{(\val\cdot\l_j/2)} - q^{-(\val\cdot\l_j)/2}
\over q^{(\val\cdot\l_j/2)} - q^{-(\val\cdot\l_j)/2}
}
}
\snzi   {P_n(\cK;\qaok) \over \Dg^{3n+1}(\cK;\qaok) }\,h^n,
\label{5.20*1}\\
&
P_n(\cK;\tok) \in \ZZ[\tok]
\nonumber
\qqq
and exponentiating both sides of \ex{5.20} we can rewrite it as
\qq
\lefteqn{
\JvaK
=
{\dqval \over \Dg(\cK;\qaok)}\;\Cqg
}
&&\qquad
\hspace{1in}\times
\exp\lrbc{
\sDl \hb^{\chi(D)}\sum_{c\in\bfS}
{p_c(\cK,D;q^{\val\cdot\l_1},\ldots,q^{\val\cdot\l_{\rg}})
\over \pjoN \APbas{\cK}{q^{\val\cdot\lcj}} }
}
\qqq

\qq
&\JvaK =   {\dqval\over \Dg(\cK;\qaot)}
\lrbc{
1 +
F_1(\cK;\qaot)
\,\hb
+ \cO(\hb^2)
}.
\\
& F_1(\tot) =
{P_1(\cK;\tot)\over
\lrbs{
\APKti{1}\APKti{2}\APK{t_1 t_2}
}^3
}
\qqq

Therefore, according to Conjecture\rw{c4.1},
\qq
\smzi\Il_{m,1}(\cK) =
{y(\cK,D_1) \over \pjov{3} \APbas{\cK}{\exp(f_j)} }
+
{y(\cK,D_2) \over \pjov{2} \APbas{\cK}{\exp(f_j)} }
\qqq

\begin{table}
\begin{center}
\begin{tabular}{|c|l|l|} \hline
Knot & $12\pthKt$ in fundamental domain & $12\pthKt$ in $\IQ[u_1,u_2]$ \\
\hline\hline
$3_1$ & $-t_1^2 t_2 + t_1^2$ & $u_1^2 - 3 u_1 u_2 - 2 u_1 - 6$
\\ \hline
$4_1$ & $0$ & $0$
\\ \hline
$5_1$ & $2t_1^4 t_2^2 - 2t_1^4 t_2 + 2 t_1^4 - t_1^2 t_2 + t_1^2$ &
$2 u_1^4 - 10 u_1^3 u_2 + 10 u_1^2 u_2^2 - 4 u_1^3 + 10 u_1^2 u_3$
\\
 & & $ - 23 u_1^2 + 53 u_1 u_3 + 26 u_1 + 66$
\\ \hline
$5_2$ & $-13 t_1^2 t_2 + 9 t_1^2 + 6 t_1 - 12$ &
$9u_1^2 - 31 u_1 u_2 - 12 u_1 - 66$
\\ \hline
$6_1$ & $3t_1^2 t_2 - t_1^2 - 6t_1 + 24$ & $-u_1^2 + 5u_1 u_2 - 4 u_1 + 30$
\\ \hline
$6_2$ & $-3t_1^4 t_2^2 + 3t_1^4 t_2 - t_1^4 - 6 t_1^3 - 11 t_1^2 t_2 + 15t_1^2$
& $-u_1^4 + 7u_1^3 u_2 - 11u_1^2 u_2^2 - 5 u_1^2 u_2 + 31u_1^2$
\\
& & $-73u_1 u_3 -34 u_1 - 114$
\\ \hline
$6_3$ & $0$ & $0$
\\ \hline
$7_1$ & $-3t_1^6 t_2^3 + 3t_1^6 t_2^2 - 3t_1^6 t_2 + 2t_1^4 t_2^2
+ 3t_1^6$
& $3 u_1^6 - 21 u_1^5 u_2 + 42 u_1^4 u_2^2 - 21 u_1^3 u_2^3 - 6 u_1^5$
\\
& $ - 2 t_1^4 t_2 + 2 t_1^4 - t_1^2 t_2 + t_1^2$ &
$+
21 u_1^4 u_2 - 52 u_1^4 + 215 u_1^3 u_2 - 152 u_1^2 u_2^2 $
\\
& & $+ 62 u_1^3 -
26 u_1^2 u_2 + 268 u_1^2 - 358 u_1 u_3$
\\
& & $ - 64 u_1 - 276$
\\ \hline
$7_2$ & $-58t_1^2 t_2+ 36 t_1^2 + 36 t_1 - 84$ &
$36 u_1^2 - 130 u_1 u_2 - 36 u_1 - 300$
\\ \hline
$7_3$ & $-25 t_1^4 t_2^2 + 25 t_1^4 t_2 - 17 t_1^4 + 7 t_1^2 t_2 - 12 t_1^3$
& $ - 17 u_1^4 + 93 u_1^2 u_2 - 109 u_1^2 u_2^2 + 38 u_1^3$
\\
&
$+ t_1^2 - 6 t_1 + 12 $
&
$ - 121 u_1^2 u_2 + 221 u_1^2 - 559 u_1 u_2 - 314 u_1$
\\
& & $ - 702 $
\\ \hline
$ 7_4$ &
$ 136 t_1^2 t_2 - 80 t_1^2 - 96 t_1 + 240$ &
$ - 80 u_1^2 + 296 u_1 u_2 + 64 u_1 + 720$
\\ \hline
$7_5$ &
$41 t_1^4 t_2^2 - 33 t_1^4 t_2 - 16 t_1^3 t_2 + 17 t_1^4 + 12 t_1^2 t_2$
& $ 17 u_1^4 - 101 u_1^3 u_2 + 141 u_1^2 u_2^2 - 50 u_1^3 $
\\
& $+32 t_1^3 + 4 t_1^2 - 14 t_1 + 36$ &
$+ 165 u_1^2 u_2 - 200 u_1^2 + 624 u_1 u_2 + 392 u_1$
\\
& & $ + 672$
\\ \hline
$7_6$ & $-7 t_1^4 t_2^2 + 5 t_1^4 t_2 + 10 t_1^3 t_2 - t_1^4 - 20 t_1^3$ &
$- u_1^4 + 9 u_1^3 u_2 - 19 u_1^2 u_2^2 - 6 u_1^3$
\\
& $- 98 t_1^2 t_2  + 64 t_1^2 + 50 t_1 - 108$ &
$ + 17 u_1^2 u_2+ 64 u_1^2 - 194 u_1 u_2 - 16 u_1$
\\
& & $-324$
\\ \hline
$7_7$ & $ - 5 t_1^2 t_2 + t_1^2 + 12 t_1 -48 $ &
$u_1^2 - 7 u_1 u_2 + 10 u_1 - 54$
\\ \hline
\end{tabular}
\end{center}
\caption{The 2-loop polynomial $12\pthKt$ presented by monomials in funamental
domain and as a polynomial of $u_1$ and $u_2$}
\label{t6.1}
\end{table}

\begin{table}
\begin{center}
\begin{tabular}{|c|l|l|} \hline
Knot & $\APKt$ & $12\pthKt$ \\
\hline\hline
$3_1$ & $t-1$ & $-t_1^2 t_2 + t_1^2$
\\ \hline
$4_1$ & $t^2 - 3t + 5$ & $0$
\\ \hline
$5_1$ & $t^2 - t + 1$ &
$2t_1^4 t_2^2 - 2t_1^4 t_2 + 2 t_1^4 - t_1^2 t_2 + t_1^2$
\\ \hline
$5_2$ & $2t - 3$ & $-13 t_1^2 t_2 + 9 t_1^2 + 6 t_1 - 12$
\\ \hline
$6_1$ & $-2t+5$ & $3t_1^2 t_2 - t_1^2 - 6t_1 + 24$
\\ \hline
$6_2$ & $-t^2 + 3t - 3$ & $-3t_1^4 t_2^2 +
3t_1^4 t_2 - t_1^4 - 6 t_1^3 - 11 t_1^2 t_2 + 15t_1^2$
\\ \hline
$6_3$ & $t^2 - 3t + 5$ & $0$
\\ \hline
$7_1$ & $ t^3 - t^2 + t - 1$ &
$-3t_1^6 t_2^3 + 3t_1^6 t_2^2 - 3t_1^6 t_2 + 2t_1^4 t_2^2 + 3t_1^6
 - 2 t_1^4 t_2 + 2 t_1^4 - t_1^2 t_2 + t_1^2$
\\ \hline
$7_2$ & $3t - 5$ & $-58t_1^2 t_2+ 36 t_1^2 + 36 t_1 - 84$
\\ \hline
$7_3$ & $2 t^2 - 3 t + 3$ & $-25 t_1^4 t_2^2 + 25 t_1^4 t_2 - 17 t_1^4 + 7
t_1^2 t_2 - 12 t_1^3 + t_1^2 - 6 t_1 + 12$
\\ \hline
$7_4$ & $4t-7$ & $ 136 t_1^2 t_2 - 80 t_1^2 - 96 t_1 + 240$
\\ \hline
$7_5$ & $2t^2 - 4t + 5$ &
$41 t_1^4 t_2^2 - 33 t_1^4 t_2 - 16 t_1^3 t_2 + 17 t_1^4 + 12 t_1^2 t_2
+32 t_1^3 + 4 t_1^2 - 14 t_1 + 36$
\\ \hline
$7_6$ & $-t^2 + 5t - 7$ &
$-7 t_1^4 t_2^2 + 5 t_1^4 t_2 + 10 t_1^3 t_2 - t_1^4 - 20 t_1^3
- 98 t_1^2 t_2  + 64 t_1^2 + 50 t_1 - 108$
\\ \hline
$7_7$ & $t^2 - 5t + 9$ & $ - 5 t_1^2 t_2 + t_1^2 + 12 t_1 -48 $
\\ \hline
\end{tabular}
\end{center}
\caption{The Alexander polynomial $\APKt$ and the
2-loop polynomial $12\pthKt$ presented by monomials in funamental
domains}
\label{t6.1}
\end{table}

\begin{table}
\begin{center}
\begin{tabular}{|c|l|l|} \hline
Knot & $\APKt$ & $12\pthKt$
\\ \hline
$3_1$ & $u-1$ &
$u_1^2 - 3u_2 u_2 - 2u_1 - 6$
\\ \hline
$4_1$ & $u^2 - 3u + 3$ & $0$
\\ \hline
$5_1$ & $u^2 - u - 1$ &
$2u_1^4 - 10 u_1^3 u_2 + 10 u_1^2 u_2^2 - 4 u_13 + 10 u_1^2 u_2 - 23 u_1^2
+ 53 u_1 u_2 $
\\
& & $+ 26 u_1 +66$
\\ \hline
$5_2$ & $2u-3$ &
$ 9u_1^2 - 31u_1 u_2 - 12 u_1 - 66$
\\ \hline
$6_1$ & $-2u+5$ &
$-u_1^2 + 5u_1 u_2 - 4 u_1 + 30$
\\ \hline
$6_3$ & $u^2 - 3u + 3$ & $0$
\\ \hline
$7_1$ & $u^3 - u^2 - 2u + 1$ &
$3u_1^6 - 21 u_1^5 u_2 + 42 u_1^4 u_2^2 - 21 u_1^3 u_2^3
- 6 u_1^5 + 21 u_1^4 u_2 - 52 u_1^4$
\\
&& $ + 215 u_1^3 u_2 - 152 u_1^2 u_2^2 +
62 u_1^3 - 26 u_1^2 u_2 + 268 u_1^2 - 358 u_1 u_2$
\\
&& $ - 64 u_1 - 276$
\\ \hline
$7_2$ & $3u-5$ & $36u_1^2 - 130u_1 u_2 - 36 u_1 - 300$
\\ \hline
$7_3$ & $2u^2 - 3u - 1$ &
$-17u_1^4 + 93u_1^3 u_2 - 109 u_1^2 u_2^2 + 38u_1^3 - 121 u_1^2 u_2
+ 221 u_1^2$
\\
&& $ - 559 u_1 u_2 - 314 u_1 - 702$
\\ \hline
$7_4$ & $4u-7$ & $-80u_1^2 + 296 u_1 u_2 + 64 u_1 + 720$
\\ \hline
$7_5$ & $2u^2 - 4u + 1$ &
$17 u_1^4 - 101 u_1^3 u_2 + 141 u_1^2 u_2^2 - 50 u_1^3 + 165 u_1^2 u_2
- 200 u_1^2$
\\
&& $ + 624 u_1 u_2 + 392 u_1 + 672$
\\ \hline
$7_6$ & $-u^2 + 5u - 5$ &
$-u_1^4 + 9 u_1^3 u_2 - 19 u_1^2 u_2^2 - 6 u_1^3 + 17 u_1^2 u_2 + 64 u_1^2
- 194 u_1 u_2$
\\
&& $ - 16 u_1 - 324$
\\ \hline
$7_7$ & $u^2 - 5u + 7$ &
$u_1^2 - 7 u_1 u_2 + 10 u_1 - 54$
\\ \hline
\end{tabular}
\end{center}
\caption{The Alexander polynomial $\APKt$ and the
2-loop polynomial $12\pthKt$ expressed in terms of symmetric polynomials $u$
and $u_1$, $u_2$}
\label{t6.2}
\end{table}

For example, let $G_0$ be the fundamental group of a knot complement
$G_0 = \pi_1(S^3\setminus \cK)$. Consider the commutators
$G_1=[G_0,G_0]$, $G_2 = [G_1,G_1]$ and abelian quotients
$G_0\p = G_0/G_1$, $G_1\p = G_1/G_2$. Obviously,
$G_0\p = H_1(S^3\setminus \cK,\ZZ) = \ZZ$.
Denote by $\bt$ the generator of $G_0\p$, it represents the meridian
of $\cK$. The group $G_0\p$ acts on
$G_1\p$ by conjugation: for $\bx\in G_0\p$, $\by\in G_1\p$,
$\bx:\;\by\mapsto \bx^{-1}\by\bx$. Now the Alexander-Conway
polynomial of $\cK$ is defined (up to a factor $\pm t^n$, $n\in\ZZ$) as
the simplest polynomial such that $\APbas{\cK}{\bt}$ maps $G_1\p$
into 0. Another definition relates the Alexander polynomial to the
Reidemeister torsion of a local system in the knot
complement, the variable $t$ being the twist acquired by that system
along the meridian of $\cK$. From both definitions of $\APKt$ it is
clear that $t$ is intimately related to the meridian of $\cK$.

Here are the results of calculating the polynomials $\pthKt$ for the first few
knots (with up to 8 crossings). We present these results in two different ways.
First, as we know, $\pth(\cK)\in\Qalghzo$ and relations\rx{6.8} come from the
symmetry $\Sgri{D_1}$. More explicitly, $\Hoz{D_1}$ looks like $su(3)$ root
lattice with elements $f_1$ and $f_2$ (and variables $t_1$, $t_2$)
corresponding to the simple roots (see~\fg{fx}).
\begin{figure}
\vspace{0.2in}
\input tstpct5.pst
\caption{Fundamental domain of the symmetry\rx{6.8}}
\label{fx}
\end{figure}
The
symmetry group $\Sgri{D_1}$ is the symmetry of this lattice (which preserves
the origin). So instead of writing the whole polynomial $\pthKt$ we may list
just the monomials belonging to a fundamental domain of $\Sgri{D_1}$. From our
$su(3)$ lattice description it is easy to see that we may choose a fundamental
domain to include the monomials
\qq
t_1^{m_1} t_2^{m_2},\qquad m_1,m_2\geq 0,\; m_1 \geq 2 m_2.
\label{6.20*}
\qqq
Then the other monomials will be determined by the symmetries\rx{6.8}.
Similarly, in view of the symmetry\rx{4.5*} it is enough to list only the
monomials of $\APKt$ with non-negative powers of $t$. Thus in \tb{t6.1} we
present the `fundamental domain' parts of the Alexander polynomial $\APKt$ and
(scaled) 2-loop polynomial $12\pthKt$.

An alternative way of describing $\pthKt$ comes from the observation that the ring
of Laurent polynomials with the symmetries\rx{6.8} can be written as
$\IQ[u_1,u_2]$, where
\qq
u_1(t_1,t_2) &=&
t_1 + t_1^{-1} + t_2 + t_2^{-1} + t_1 t_2 + t_1^{-1} t_2^{-1},
\nonumber\\
u_2(t_1,t_2) & = &
t_1^2 t_2 + t_1^{-2} t_2^{-1} + t_1 t_2^{2} + t_1^{-1} t_2^{-2}
+ t_1 t_2^{-1} + t_1^{-1} t_2.
\label{6.21*}
\qqq
So in \tb{t6.2} we present the expressions for the Alexander polynomial $\APKt$
in terms of $u=t+t^{-1}$ and for the (scaled) 2-loop polynomial $12\pthKt$ in
terms of $u_1$ and $u_2$.

\begin{remark}
\rm
If $\cK\p$ is the mirror image of $\cK$, then
$\pth(\cK\p;\tot) = -\pthKt$, hence $\pthKt=0$ for amphicheiral knots.
\end{remark}
\begin{remark}
\rm
As we see, experimental evidence suggests that
\qq
12\pthKt\in\ZZtot.
\label{6.19*}
\qqq
\end{remark}
\begin{remark}
\rm
The degree of the Alexander polynomial is bounded by the genus of the knot
$g(\cK)$
\qq
\deg \APKt \leq g(\cK),
\label{6.22*}
\qqq
In view of the symmetries\rx{6.8} (which come from $\Sgri{D_1}$), the reasonable
measure of the degree of $\pthKt$ is the $t_1$ degree of its fundamental domain
part. Let us denote it simply as $\deg\pthKt$. Then \tb{t6.1} suggests a
similar inequality
\qq
\deg\pthKt\leq 2g(\cK).
\label{6.23*}
\qqq

\end{remark}